\documentclass[11pt]{amsart}
\usepackage{amssymb}
\usepackage{amsmath}
\usepackage{mathrsfs}
\usepackage{fancyhdr}
\usepackage[british]{babel}
\usepackage{geometry}
\usepackage{enumitem}
\usepackage{algpseudocode}
\usepackage{dsfont}
\usepackage{centernot}
\usepackage{xstring}
\usepackage{colortbl}
\usepackage[section]{algorithm}
\usepackage{graphicx}
\usepackage[font={footnotesize}]{caption}
\usepackage[usenames,dvipsnames,table]{xcolor}
\usepackage{tikz}
\usepackage[h]{esvect}
\usepackage[
		bookmarksopen=true,
		bookmarksopenlevel=1,
		colorlinks=true,
		linkcolor=darkblue,
        linktoc=page,
		citecolor=darkblue,
]{hyperref}
\usepackage{bbm}
\usepackage{subcaption}

\title[Decompositions of quasirandom hypergraphs]{Decompositions of quasirandom hypergraphs into hypergraphs of bounded degree}
\date{\today}
\author[Stefan~Ehard]{Stefan Ehard}
\address[Stefan~Ehard]{Institut f\"ur Optimierung und Operations Research, Universit\"at Ulm,
Germany
}
\email{stefan.ehard@uni-ulm.de}

\author[Felix~Joos]{Felix Joos}
\address[Felix~Joos]{Institut f\"ur Informatik, Universit\"at Heidelberg, Germany }
\email{joos@informatik.uni-heidelberg.de}

\thanks{The research leading to these results was partially supported by the Deutsche Forschungsgemeinschaft (DFG, German Research Foundation) -- 428212407 (F.~Joos).
}

\newtheorem{theorem}[algorithm]{Theorem}

\newtheorem{lemma}[algorithm]{Lemma}
\newtheorem{cor}[algorithm]{Corollary}
\newtheorem{fact}[algorithm]{Fact}
\newtheorem{defin}[algorithm]{Definition}

\newtheorem{problem}[algorithm]{Problem}
\newtheorem{conj}[algorithm]{Conjecture}

\newtheorem{question}[algorithm]{Question}

\newtheoremstyle{partstyle}{10pt}{5pt}{\em}{0pt}{\em}{.}{5pt}{}
\theoremstyle{partstyle}
\newtheorem{proofpart}{Part}

\newcounter{stepenv}
\newenvironment{stepenv}[1][]{\refstepcounter{stepenv}}{}

\newcounter{step}[stepenv]
\newenvironment{step}[1][]{\refstepcounter{step}\par\medskip\noindent%
        \textit{\underline{Step~{\thestep}.#1} } \itshape\rmfamily}{\medskip}
				
\newcounter{substep}[step]
\renewcommand{\thesubstep}{\thestep.\arabic{substep}}
\newenvironment{substep}[1][]{\refstepcounter{substep}\par\medskip\noindent%
        \emph{\underline{Step~{\thesubstep}.#1} } \itshape\rmfamily}{\medskip}

\newcounter{subsubstep}[substep]
\renewcommand{\thesubsubstep}{\thesubstep.\arabic{subsubstep}}

\newcounter{claim}[stepenv]
\newenvironment{claim}[1][]{\refstepcounter{claim}\par\medskip\noindent%
        \textit{Claim~\theclaim. #1} \itshape\rmfamily}{\medskip}

\numberwithin{equation}{section}

\definecolor{darkblue}{rgb}{0,0,0.5}

\def\noproof{{\unskip\nobreak\hfill\penalty50\hskip2em\hbox{}\nobreak\hfill%
       $\square$\parfillskip=0pt\finalhyphendemerits=0\par}\goodbreak}
\def\endproof{\noproof\bigskip}

\def\noclaimproof{{\unskip\nobreak\hfill\penalty50\hskip2em\hbox{}\nobreak\hfill%
       $-$\parfillskip=0pt\finalhyphendemerits=0\par}\goodbreak}

\def\endclaimproof{\noclaimproof\medskip}

\newdimen\margin
\def\textno#1&#2\par{
   \margin=\hsize
   \advance\margin by -4\parindent
          \setbox1=\hbox{\sl#1}
   \ifdim\wd1 < \margin
      $$\box1\eqno#2$$
   \else
      \bigbreak
      \hbox to \hsize{\indent$\vcenter{\advance\hsize by -3\parindent
      \it\noindent#1}\hfil#2$}
      \bigbreak
   \fi}

\def\lateproof#1{\removelastskip\penalty55\medskip\noindent\begin{stepenv}\end{stepenv}{\bf Proof of #1. }} 

\DeclareMathOperator{\dg}{deg}

\DeclareMathOperator{\plog}{polylog}

\def\claimproof{\removelastskip\penalty55\medskip\noindent{\em Proof of claim: }}

\makeatletter
\DeclareFontEncoding{LS1}{}{}
\DeclareFontSubstitution{LS1}{stix}{m}{n}
\DeclareMathAlphabet{\mathscrs}{LS1}{stixscr}{m}{n}
\makeatother

\makeatletter
\g@addto@macro \normalsize {%
 \setlength\abovedisplayskip{4pt plus 2pt minus 2pt}%
 \setlength\belowdisplayskip{4pt plus 2pt minus 2pt}%
}
\makeatother

\begin{document}

\newcommand{\new}[1]{\textcolor{red}{#1}}
\def\COMMENT#1{}
\def\TASK#1{}

\newcommand{\todo}[1]{\begin{center}\textbf{to do:} #1 \end{center}}

\def\eps{{\varepsilon}}
\def\heps{{\hat{\varepsilon}}}

\newcommand{\ex}{\mathbb{E}}
\newcommand{\pr}{\mathbb{P}}
\newcommand{\cB}{\mathcal{B}}
\newcommand{\cA}{\mathcal{A}}
\newcommand{\cE}{\mathcal{E}}
\newcommand{\cS}{\mathcal{S}}
\newcommand{\cF}{\mathcal{F}}
\newcommand{\cG}{\mathcal{G}}
\newcommand{\bL}{\mathbb{L}}
\newcommand{\bF}{\mathbb{F}}
\newcommand{\bZ}{\mathbb{Z}}
\newcommand{\cH}{\mathcal{H}}
\newcommand{\cC}{\mathcal{C}}
\newcommand{\cM}{\mathcal{M}}
\newcommand{\bN}{\mathbb{N}}
\newcommand{\bR}{\mathbb{R}}
\newcommand{\bS}{\mathbb{S}}
\def\O{\mathcal{O}}
\newcommand{\cP}{\mathcal{P}}
\newcommand{\cQ}{\mathcal{Q}}
\newcommand{\cR}{\mathcal{R}}
\newcommand{\cJ}{\mathcal{J}}
\newcommand{\cL}{\mathcal{L}}
\newcommand{\cK}{\mathcal{K}}
\newcommand{\cD}{\mathcal{D}}
\newcommand{\cI}{\mathcal{I}}
\newcommand{\cN}{\mathcal{N}}
\newcommand{\cV}{\mathcal{V}}
\newcommand{\cT}{\mathcal{T}}
\newcommand{\cU}{\mathcal{U}}
\newcommand{\cX}{\mathcal{X}}
\newcommand{\cZ}{\mathcal{Z}}
\newcommand{\cW}{\mathcal{W}}
\newcommand{\1}{{\bf 1}_{n\not\equiv \delta}}
\newcommand{\eul}{{\rm e}}
\newcommand{\Erd}{Erd\H{o}s}
\newcommand{\cupdot}{\mathbin{\mathaccent\cdot\cup}}
\newcommand{\whp}{whp }
\newcommand{\bX}{\mathcal{X}}
\newcommand{\bV}{\mathcal{V}}
\newcommand{\hbX}{\widehat{\mathcal{X}}}
\newcommand{\hbV}{\widehat{\mathcal{V}}}
\newcommand{\hpsi}{\widehat{\psi}}
\newcommand{\hA}{\widehat{A}}
\newcommand{\hcA}{\widehat{\cA}}
\newcommand{\hsA}{\widehat{\sA}}
\newcommand{\tA}{\widetilde{A}}
\newcommand{\tB}{\widetilde{B}}
\newcommand{\tE}{\widetilde{E}}
\newcommand{\tH}{\widetilde{H}}
\newcommand{\tG}{\widetilde{G}}
\newcommand{\tN}{\widetilde{N}}
\newcommand{\hB}{\widehat{B}}
\newcommand{\hX}{\widehat{X}}
\newcommand{\hV}{\widehat{V}}
\newcommand{\tX}{\widetilde{X}}
\newcommand{\tV}{\widetilde{V}}
\newcommand{\cbI}{\overline{\mathcal{I}^\alpha}}
\newcommand{\hAj}{\widehat{A}^\sigma_j}
\newcommand{\hVM}{V^M}
\newcommand{\bc}{\mathscrs{c}}
\newcommand{\bd}{\mathbf{d}}
\newcommand{\bp}{\mathbf{p}}
\newcommand{\bpp}{\mathbf{p}^{2nd}}
\newcommand{\bpA}{\mathbf{p}^{A}}
\newcommand{\bppA}{\mathbf{p}^{A,2nd}}
\newcommand{\bpB}{\mathbf{p}^{B} }
\newcommand{\bppB}{\mathbf{p}^{B,2nd}}
\newcommand{\bpZ}{\mathbf{p}^{Z} }
\newcommand{\bppZ}{\mathbf{p}^{Z,2nd}}
\newcommand{\bx}{\mathbf{x}}
\newcommand{\by}{\mathscrs{y}}
\newcommand{\be}{\mathbf{e}}
\newcommand{\sX}{\mathscr{X}}
\newcommand{\shX}{\widehat{\mathscr{X}}}
\newcommand{\sA}{\mathcal{A}} 
\newcommand{\sB}{\mathscr{B}}
\newcommand{\sP}{\mathscr{P}}
\newcommand{\stA}{{\mathscr{\widetilde{A}}}}
\newcommand{\bs}{\boldsymbol}

\newcommand{\eH}{\mathscrs{e}}
\newcommand{\ea}{\mathscrs{e}^{\ast}}
\newcommand{\fH}{\mathscrs{f}}
\newcommand{\vG}{\mathscrs{v}}
\newcommand{\gG}{\mathscrs{g}}
\newcommand{\hG}{\mathscrs{h}}
\newcommand{\rR}{\mathscrs{r}}
\newcommand{\aA}{\mathscrs{a}}
\newcommand{\bA}{\mathscrs{b}}
\newcommand{\ab}{{\mathscrs{a}\mathscrs{b}}}
\newcommand{\mm}{\mathscrs{m}}
\newcommand{\xX}{\mathscrs{x}}
\newcommand{\pP}{\mathscrs{p}}
\newcommand{\zZ}{\mathscrs{z}}

\newcommand{\VR}{{V(R)}}
\newcommand{\VRsm}{{V(R)\sm\{0\}}}
\newcommand{\JX}{J_{X}}
\newcommand{\JV}{J_{V}}
\newcommand{\JXV}{J_{XV}}
\newcommand{\last}{J^{last}}
\newcommand{\rc}{[r_\circ]}
\newcommand{\fst}{$1^{st}$}
\newcommand{\scd}{$2^{nd}$}

\newcommand{\supp}{{\rm supp}}

\newcommand{\bars}{\bar{\sigma}}

\newcommand{\cl}{{=}}

\newcommand{\doublesquig}{%
  \mathrel{%
    \vcenter{\offinterlineskip
      \ialign{##\cr$\rightsquigarrow$\cr\noalign{\kern-1.5pt}$\rightsquigarrow$\cr}%
    }%
  }%
}

\newcommand{\defn}{\emph}

\newcommand\restrict[1]{\raisebox{-.5ex}{$|$}_{#1}}

\newcommand{\prob}[1]{\mathrm{\mathbb{P}}\left[#1\right]}
\newcommand{\expn}[1]{\mathrm{\mathbb{E}}\left[#1\right]}
\def\gnp{G_{n,p}}
\def\G{\mathcal{G}}
\def\lflr{\left\lfloor}
\def\rflr{\right\rfloor}
\def\lcl{\left\lceil}
\def\rcl{\right\rceil}

\newcommand{\qbinom}[2]{\binom{#1}{#2}_{\!q}}
\newcommand{\binomdim}[2]{\binom{#1}{#2}_{\!\dim}}

\newcommand{\grass}{\mathrm{Gr}}

\newcommand{\brackets}[1]{\left(#1\right)}
\def\sm{\setminus}
\newcommand{\Set}[1]{\{#1\}}
\newcommand{\set}[2]{\{#1\,:\;#2\}}
\newcommand{\krq}[2]{K^{(#1)}_{#2}}
\newcommand{\ind}[1]{$\mathbf{S}(#1)$}
\newcommand{\indC}[1]{$\mathbf{C1}(#1)$}
\newcommand{\indCC}[1]{$\mathbf{C2}(#1)$}
\newcommand{\indcov}[1]{$(\#)_{#1}$}
\def\In{\subseteq}
\newcommand{\IND}{\mathbbm{1}}
\newcommand{\norm}[1]{\|#1\|}
\newcommand{\normv}[1]{\|#1\|_v}
\newcommand{\normV}[2]{\|#1\|_{#2}}

\begin{abstract} 
\noindent
We prove that any quasirandom uniform hypergraph $H$ can be approximately decomposed into any collection of bounded degree hypergraphs with almost as many edges.
In fact, our results also apply to multipartite hypergraphs and even to the sparse setting when the density of~$H$ quickly tends to $0$ in terms of the number of vertices of $H$.
Our results answer and address questions of Kim, K\"uhn, Osthus and Tyomkyn;
and Glock, K\"uhn and Osthus as well as Keevash.

The provided approximate decompositions exhibit strong quasirandom properties which is very useful for forthcoming applications.
Our results also imply approximate solutions to natural hypergraph versions of long-standing graph decomposition problems, as well as several decomposition results for (quasi)random simplicial complexes into various more elementary simplicial complexes such as triangulations of spheres and other manifolds.
\end{abstract}

\maketitle

\vspace{-.3cm}
\section{Introduction}
Decompositions of large mathematical objects into smaller pieces have been an integral part in nearly every field of mathematics.
We also find numerous instances of this theme in \mbox{combinatorics} with many connections to other related fields.
Potentially the most prominent recent examples are the achievements surrounding the resolution of the Existence conjecture by Keevash~\cite{keevash:14} and its generalizations by Glock, K\"uhn, Lo and Osthus as well as Keevash himself~\cite{GKLO:ta,keevash:18b}.
These results include decompositions of quasirandom uniform (multipartite) hypergraphs (the so-called host hypergraphs) into cliques, and more generally, into arbitrary but fixed hypergraphs subject to ``obvious'' divisibility conditions and the host hypergraphs being sufficiently large.

The analogous question for graphs has been famously solved by Wilson~\cite{wilson:72a,wilson:72b,wilson:75}.
His results have initiated a very vibrant research area with connections to Latin squares and arrays, geometry, design theory and combinatorial probability theory.

The following three conjectures {have been} a driving force regarding graph decompositions and concern families of sparse graphs whose number of vertices grows with the number of vertices of the host graph; sometimes the sparse graphs even have the same number of vertices as the host graph, a case which turns out to be particularly challenging.
Before we state the conjectures, let us introduce some notation.
We say a family $\cH$ of hypergraphs \emph{packs} into a hypergraph $G$ if the members of~$\cH$ can be found edge-disjointly in~$G$.
The family $\cH$ \emph{decomposes} $G$ if {additionally} $\sum_{H\in \cH}e(H)=e(G)$, {where $e(J)$ denotes the number of edges of a hypergraph $J$}.
We also say there is a \emph{decomposition} of~$G$ into copies of a hypergraph~$H$
if the edge set of $G$ can be partitioned into copies of~$H$.
{Given these definitions, we can state the aforementioned conjectures.}
\begin{enumerate}[label=(\Roman*)]
	\item\label{conj:R} (Ringel, 1963) {\it For all $n$ and all trees $T$ on $n+1$ vertices, there is a decomposition of $K_{2n+1}$ into $2n+1$ copies of $T$.}
	\item\label{conj:GL} (Tree packing conjecture -- Gy\'arf\'as and Lehel, 1976) 
	{\it For all $n$ and all sequences of trees $T_1,\ldots,T_n$ where $T_i$ has $i$ vertices, there is a decomposition of $K_n$ into $T_1,\ldots,T_n$.}
	\item\label{conj:OW} (Oberwolfach problem -- Ringel, 1967) 
	{\it For all odd $n$ and all $2$-regular graphs $F$ on $n$ vertices,
	there is a decomposition of $K_{n}$ into copies of $F$.}
\end{enumerate}
Until recently (considering that the conjectures were posed in the 1960s and 1970s), 
no substantial progress had been made, but {lately} a number of striking and exciting advances have been achieved mostly in conjunction with very elaborated analyses of random processes.
This includes the resolution of Conjectures~\ref{conj:R} and~\ref{conj:GL} for bounded degree trees~\cite{JKKO:19}, the complete resolutions of Conjecture~\ref{conj:R} in~\cite{KS:20b,MPS:20}  and Conjecture~\ref{conj:OW} in~\cite{GJKKO:18,KS:20a}, and the resolution of Conjecture~\ref{conj:GL} for families of trees with many leaves and a weak restriction concerning the maximum degree~\cite{ABCT:19}.
All these results hold for sufficiently large $n$.
Preceding these results, there is a collection of \mbox{approximate} decomposition results, that is, where a few edges of the host graph are not covered, under various and quite general conditions, see~\cite{ABHP:17,BHPT:16,FLM:17,KKOT:19,MRS:16}.
The importance of these results should not be underestimated as many of those play a key role in the actual decomposition results.
{In particular,
Kim, K\"uhn, Osthus and Tyomkyn~\cite{KKOT:19} established a powerful result for approximate decompositions of quasirandom (multipartite) graphs into spanning bounded degree graphs by proving a blow-up lemma for approximate decompositions.
This greatly extends the famous blow-up lemma due to Koml\'os, S\'ark\"ozy and Szemer\'edi~\cite{KSS:97}.
There are two natural ways for further extensions of the main result of~\cite{KKOT:19}, which we discuss in the following.
}

\medskip

\noindent
\textbf{From graphs to hypergraphs:} Although there are numerous (approximate) decomposition results for spanning structures in graphs, the situation for hypergraphs is notably different.
There are only a few {approximate decomposition} results concerning various types of cycles~\cite{BF:12,FK:12,FKL:12, JK:21} as well as Keevash's results concerning factors, that is, the vertex-disjoint collection of a hypergraph of fixed size, which follow from results in~\cite{keevash:18b}.
{
Even for embedding just a single spanning structure, there are several obstacles for generalizations to hypergraphs; for instance,
it took serious effort to lift the usual blow-up lemma for graphs to hypergraphs, a feat eventually achieved by Keevash~\cite{keevash:11}.
In view of the blow-up lemma for approximate decompositions of graphs, Kim, K\"uhn, Osthus and Tyomkyn~\cite{KKOT:19} as well as Keevash~\cite{keevash:18c} 
implicitly posed the following question.}
\begin{question}[Keevash; Kim, K\"uhn, Osthus and Tyomkyn]\label{questionHyper}
Does every quasirandom (multipartite) hypergraph admit an approximate decomposition into any collection of spanning bounded~degree hypergraphs (with slightly fewer edges)?
\end{question}

\noindent
\textbf{From the dense to the sparse regime:}
{There are various attempts to transfer nowadays classical results on dense graphs to the sparse regime.
Early instances include works of Thomason~\cite{thomason:87b,thomason:87} as well as Chung and Graham~\cite{CG:02} on sparse quasirandom graphs (see also references therein), and more recently, Allen, B\"ottcher, H\`an, Kohayakawa and Person~\cite{ABHKP:16} proved sparse versions of the usual blow-up lemma.
In light of this and their result in~\cite{KKOT:19}, Kim, K\"uhn, Osthus and Tyomkyn remarked that it would be very interesting to solve the following question.}

\begin{question}[Kim, K\"uhn, Osthus and Tyomkyn]\label{question}
Does every quasirandom graph $G$ admit an approximate decomposition into any collection of spanning bounded degree graphs (with slightly fewer edges) even when the density of $G$ tends to $0$ as the order of $G$ tends to infinity?
\end{question}

{Observe that the sparse regime exhibits notable further challenges; 
for instance, it is not even known whether every quasirandom sparse hypergraph contains every bounded degree spanning hypergraph as a subgraph -- a question that was also remarked upon by Allen, B\"ottcher, H\`an, Kohayakawa and Person in~\cite[Section~7.5]{ABHKP:16}. See also~\cite{HHM:20} for sparse embedding results of linear hypergraphs.}

\subsection{The main result -- a simplified version}

{
In this paper we simultaneously make decisive progress in both aspects outlined above by providing a versatile tool for approximately decomposing quasirandom hypergraphs into families of spanning bounded degree hypergraphs.
That is, we prove that any (sparse) (multipartite) quasirandom hypergraph admits an approximate decomposition into any family of spanning bounded degree hypergraphs with slightly fewer edges,
thereby answering (in a strong form) both Questions~\ref{questionHyper} and~\ref{question} in the affirmative.
}

{
Our main result also has some direct applications; for instance, it yields an asymptotic solution to a hypergraph Oberwolfach problem asked by Glock, K\"uhn and Osthus~\cite{GKO:ta}, see Section~\ref{sec:hypergraphs}.
}

{In order to state} a simplified version of our main result we introduce some terminology.
A hypergraph~$G$ is called \emph{$k$-uniform} or a \emph{$k$-graph} if all edges have size $k$.
The neighbourhood $N_G(S)$ of a $(k-1)$-set $S$ of vertices is the set of vertices that form an edge together with $S$.
Suppose $\eps>0$, $t\in \bN$, $d\in (0,1]$ and $G$ has $n$ vertices.
We say $G$ is \emph{$(\eps,t,d)$-typical} if 
$|\bigcap_{S\in\cS}N_G(S)|=(1\pm\eps)d^{|\cS|}n$
for all (non-empty) sets $\cS$ of $(k-1)$-sets of $V(G)$ with $|\cS|\leq t$.
We refer to the (vertex) degree of a vertex $v$ as the number of edges containing~$v$ and let $\Delta(G)$ be the maximum degree in $G$.

The following result is a direct consequence of our main result.
{In view of Question~\ref{question}, note that we allow the density to decay at a rate that is polynomial in terms of the order.}

\begin{theorem}\label{thm:simple}
For all $\alpha\in(0,1]$,
there exist $n_0,t\in \bN$ and $\eps>0$ such that the following holds for all $n\geq n_0$.
Suppose $G$ is an $(\eps,t,d)$-typical $k$-graph on $n$ vertices with $k\leq\alpha^{-1}$ as well as $d\geq n^{-\eps}$, and 
{$\cH$} is a family of $k$-graphs {$H$} on $n$ vertices with $\Delta(H)\leq\alpha^{-1}$ {for all $H\in\cH$}
and $\sum_{H\in\cH}e(H)\leq (1-\alpha)e(G)$.
Then {$\cH$} packs into~$G$.
\end{theorem}

Observe that the binomial random $k$-graph is with high probability $(\eps,t,d)$-typical whenever $k,\eps,t$ are fixed and $d\geq n^{-\eps}$,
and consequently, with high probability these $k$-graphs can be approximately decomposed into any list of bounded degree hypergraphs with almost as many edges.

\subsection{Applications}
Apart from being interesting on its own, as outlined above, we believe that our main results will be in addition useful for forthcoming applications.
We illustrate and discuss some immediate applications in the concluding Section~\ref{sec:conclusion}, but of course, more complex applications are beyond the scope of this article.
A general strategy for using our results to obtain an actual decomposition could be as follows: in a first step, use our results to embed most parts of the hypergraphs while retaining well-chosen pieces for a second step (these parts could be leaves as in~\cite{ABCT:19}, short cycles all of the same length as in~\cite{GJKKO:18} or many short paths as in~\cite{KS:20a}).
We can ensure that the `free places' for the pieces in the second step of the procedure are well-spread among the hypergraphs, which is a necessary property for such an approach to work (for details, see our notion of \emph{testers} in Section~\ref{sec:main results}).
See~\cite{ABCT:19} for an example where a similar result as Theorem~\ref{thm:simple} has been applied to find an actual decomposition (in the setting of graphs).

In addition, Theorem~\ref{thm:simple} can be applied to several natural questions on hypergraph decompositions including a question of Glock, K\"uhn and Osthus,
which we consider in Section~\ref{sec:hypergraphs}.
Closely linked to hypergraphs are simplicial complexes which are equivalent to downward closed hypergraphs $H$; that is, whenever $\eH$ is an edge of~$H$,
then $H$ contains also all subsets of $\eH$ as edges.
Gowers was one of the first who suggested the investigation of topological analogues of 1-dimensional graph structures in higher dimensions as we may consider a Hamilton cycle in a graph as a spanning 1-dimensional simplicial complex that is homeomorphic to $\mathbb{S}^1$
and hence a Hamilton cycle in higher dimensions may be viewed as spanning $k$-dimensional simplicial complex that is homeomorphic to~$\mathbb{S}^k$.
In particular, Linial has considered further questions of this type under the term `high-dimensional \mbox{combinatorics}' and has achieved several new insights.
We discuss implications of Theorem~\ref{thm:simple} to these types of questions in Section~\ref{sec:simplicial complexes}.

\subsection{Multipartite graphs and the main result}\label{sec:main results}
Let us now turn to the statement of our main result.
In the past it turned out that (approximate) decomposition results of multipartite (hyper)graphs are highly desirable for several applications (to name only two examples, see Kim, K\"uhn, Osthus and Tyomkyn~\cite{KKOT:19} and Keevash~\cite{keevash:18b}, which are for instance applied in~\cite{GJKKO:18,JKKO:19,KS:20a,KS:20b}).
In view of this, we provide all our tools also for the multipartite setting.
Thus, for our notion of quasirandomness, we lift the hypergraph blow-up lemma due to Keevash~\cite{keevash:11} {for embedding a single hypergraph} 
to an approximate decomposition result; yet, note that our notion of quasirandomness is stronger than the usual notion that comes with the regularity lemma.

Whether we have a strong control over the actual (approximate) decomposition, locally and globally, is another decisive factor in the strength of the tool for further applications.
Therefore, we make a considerable effort to implement two types of versatile test functions
with respect to which the decomposition behaves random-like.

{We use standard terminology and refer to Section~\ref{sec:notation} for more notational details.}
We say that a family/multiset of $k$-graphs $\cH=\{H_1,\ldots, H_s\}$ \defn{packs} into a $k$-graph~$G$
if there is a function $\phi: \bigcup_{H\in \cH}V(H)\to V(G)$ such that $\phi|_{V(H)}$ is injective and $\phi$ injectively maps edges onto edges. 
In such a case, we call $\phi$ a \emph{packing of~$\cH$ into~$G$}.
Our general aim is to pack a collection $\cH$ of multipartite $k$-graphs into a quasirandom host $k$-graph $G$ having the same multipartite structure, which is captured by a so-called `reduced graph' $R$.
For graphs and where $\cH$ consists of a single graph only,
this result has been proven first by Koml\'os, S\'ark\"ozy and Szemer\'edi~\cite{KSS:97}
and is famously known as the blow-up lemma.
In view of this, 
we say
$(H,G,R,\cX,\cV)$ is a \defn{blow-up instance} of size $(n,k,r)$ if\looseness=-2
\begin{itemize}
	\item $H,G,R$ are $k$-graphs where $V(R)=[r]$;
	\item $\cX=(X_i)_{i\in[r]}$ is a vertex partition of $H$ such that $|\eH\cap X_i|\leq 1$ for all $\eH\in E(H), i\in[r]$;
	\item $\cV=(V_i)_{i\in[r]}$ is a vertex partition of $G$ such that $|V_i|=|X_i|=(1\pm1/2)n$ for all $i\in [r]$; 
	\item $H[X_{i_1},\ldots,X_{i_k}]$ is empty whenever $\{i_1,\ldots,i_k \}\notin E(R)$.
\end{itemize}
We also refer to $(\cH,G,R,\cX,\cV)$ as a \defn{blow-up instance}
if $\cH$ is a collection of $k$-graphs and $\cX$ is a collection of vertex partitions $(X_i^H)_{i\in [r],H\in \cH}$ 
so that $(H,G,R,(X_i^H)_{i\in[r]},\cV)$ is a blow-up instance for every $H\in \cH$.

Given a blow-up instance $\sB=(\cH,G,R,\cX,\cV)$ of size $(n,k,r)$, we generalize the notion of typicality to the multipartite setting given by the reduced graph $R$. 
For finite and disjoint sets $A_1,\ldots,A_\ell$, we write $\bigsqcup_{i\in[\ell]}A_i:=\{\{a_1,\ldots,a_\ell\}\colon a_i\in A_i \text{ for all $i\in[\ell]$} \}$.
We say $G$ is \defn{$(\eps,t,d)$-typical with respect to~$R$} if for all $i\in[r]$ and all sets $\cS\subseteq
\bigcup_{\rR\in E(R)\colon i\in\rR} \bigsqcup_{j\in \rR\sm\{i\}}V_{j}$
with $|\cS|\leq t$, we have $|V_i\cap \bigcap_{S\in\cS}N_G(S)|=(1\pm\eps)d^{|\cS|}|V_i|$.
We say the blow-up instance $\sB$ is \defn{$(\eps,t,d)$-typical} if $G$ is $(\eps,t,d)$-typical with respect to $R$ and $|V_i|=|X_i^H|=(1\pm\eps)n$ for all $H\in\cH, i\in[r]$.
We say $\sB$ is \defn{$\Delta$-bounded} if $\Delta(R),\Delta(H)\leq \Delta$ for each $H\in \cH$.

For the readers' convenience, we first state a simplified version of our main result for multipartite graphs.

\begin{theorem}\label{thm:multi}
For all $\alpha\in(0,1]$,
there exist $n_0,t\in \bN$ and $\eps>0$ such that the following holds for all $n\geq n_0$.
	Suppose $(\cH,G,R,\cX,\cV)$ is an  $(\eps,t,d)$-typical and $\alpha^{-1}$-bounded blow-up instance of size $(n,k,r)$ with $k\leq\alpha^{-1}$, $r\leq n^{\log n}$, $d\geq n^{-\eps}$, 
	$|\cH|\leq n^{k+1}$,\COMMENT{This is not restrictive and $|\cH|\leq n^{k+1}$ is chosen arbitrarily (with some room to establish Theorem~\ref{thm:normalwithweights}).} and $\sum_{H\in \cH} e_H(X_{i_1}^H,\ldots,X_{i_k}^H) \leq (1-\alpha)dn^k$ for all $\{i_1,\ldots,i_k \} \in E(R)$.
	Then there is a packing $\phi$ of $\cH$ into $G$ such that $\phi(X_i^H)= V_i$ for all $i\in [r]$ and $H\in \cH$.
\end{theorem}

In numerous applications of the original blow-up lemma for graphs, it has been essential that it provides additional features.
In view of this, we make a substantial effort to also include further tools in our results that allow us to control the structure of the packings and will be very useful for future applications.
We achieve this with two different types of {what we call} \emph{testers}.
The first tester is a so-called \defn{set tester}; with the setting as in Theorem~\ref{thm:multi},
we can fix a set $Y\In X_i^{H}$ and $W\In V_i$ for some $i\in [r]$ and $H\in \cH$.
Then we find a packing $\phi$ such that $|W\cap \phi(Y)|=|W||Y|/n \pm \alpha n$.
Moreover, we can even fix several sets~$Y_j$ as above in multiple $k$-graphs $H_j$ in $\cH$ and {the size of} their common intersection with~$W$ is as {large as} we would expect it to be in an idealized random packing.\looseness=-1

The second type of tester is a so-called \defn{vertex tester}.
In the simplest form, we fix a vertex $c\in V_i$ with $i\in [r]$ and define a weight function on $\bigcup_{H\in \cH}X_i^H$.
Then we find a packing such that the weight of the vertices embedded onto $c$ is roughly the total weight divided by $n$.
However, for many applications it is not enough to control single vertices -- this is one reason why it is difficult to apply the hypergraph blow-up lemma due to Keevash~\cite{keevash:11}.
Here, we make a considerable effort to provide a tool to deal with larger sets.
To this end, for a vertex tester, we can also fix $c_i\in V_i$ for $i\in I$ for some $(k-1)$-set $I\In \rR\in E(R)$ and define a weight function on the $(k-1)$-sets that could be potentially embedded onto $\{c_i\}_{i\in I}$.
Then our main result yields an embedding where the weight of the actually embedded $(k-1)$-sets onto $\{c_i\}_{i\in I}$ is the appropriate proportion of the total weight assigned.

Let us now formally define these two types of testers.
Suppose $\sB$ is a blow-up instance as above.
We say $(W,Y_1,\ldots,Y_m)$ is an \defn{$\ell$-set tester} for $\sB$ 
if $m\leq \ell$ and there exist $i\in [r]$ and distinct $H_1,\ldots,H_m\in \cH$ such that $W\In V_i$ and $Y_j\In X_i^{H_j}$ for all $j\in [m]$.
We say $(\omega,\bc)$ is an \defn{$\ell$-vertex tester} for $\sB$ with centres $\bc=\{c_i \}_{i\in I}$ in $I\subseteq [r]$, if 
\begin{itemize}
\item $|I|\leq k-1$\COMMENT{\label{comment size I vertex tester}We want to put weights on $(k-1)$-sets, so usually, $|I|=k-1$. In fact, we also use the property that $|I|\leq k-1$ explicitly after display~\eqref{eq:relation weight w_0 w_ver}.} and $I\subseteq\rR$ for some $\rR\in E(R)$, and $c_i\in V_i$ for each $i\in I$, and
\item $\omega$ is a weight function on the $|I|$-tuples $\sX_{\sqcup {I}}:=\bigcup_{H\in \cH}(\bigsqcup_{i\in I} X_i^H)$ with $\omega\colon\sX_{\sqcup I} \to [0,\ell]$ and {whenever $|I|\geq 2$, we have  $\supp(\omega)=\omega^{-1}((0,\ell])\subseteq \{\xX\in\sX_{\sqcup I} \colon \xX=\eH\cap \sX_{\sqcup {I}} \text{ for some $\eH\in\cH$} \}$.}\COMMENT{This last condition on $\omega^{-1}((0,\ell])$ means that we only allow weight on tuples that really lie within an edge of $\cH$.}
\end{itemize}
For an $\ell$-tuple function $\omega\colon\binom{X}{\ell}\to\bR_{\ge 0}$ on a finite set $X$, we define $\omega(X'):=\sum_{S\in \binom{X'}{\ell}}\omega(S)$ for any $X'\In X$.
The following theorem is our main result.

\begin{theorem}\label{thm:main_new}
Suppose the assumptions of Theorem~\ref{thm:multi} hold and
suppose $\cW_{set},\cW_{ver}$ are sets of  $\alpha^{-1}$-set testers and $\alpha^{-1}$-vertex testers of size at most $n^{2\log n}$, respectively.
Then there is a packing~$\phi$ of $\cH$ into $G$ in Theorem~\ref{thm:multi} such that 
\begin{enumerate}[label=\rm (\roman*)]
	\item\label{item:set testers} $|W\cap \bigcap_{j\in [m]}\phi(Y_j)|= |W||Y_1|\cdots |Y_m|/n^m \pm \alpha n$ for all $(W,Y_1,\ldots,Y_m)\in \cW_{set}$;
	\item\label{item: vertex tester} $\omega(\phi^{-1}(\bc))=(1\pm\alpha)\omega(\sX_{\sqcup I})/n^{|I|}\pm  n^{\alpha}$ for all $(\omega,\bc)\in \cW_{ver}$ with centres $\bc$ in~$I$.
\end{enumerate}
\end{theorem}

In the same line, we also provide these two types of testers if $G$ is a (non-multipartite) quasirandom $k$-graph and the result in fact follows from Theorem~\ref{thm:main_new}.
The definition is adapted in the obvious way.
Suppose the vertex set of $G$ is $V$ and we aim to pack $\cH$ into $G$ where the vertex set of $H\in\cH$ is denoted by $X^H$.
For set testers $(W,Y_1,\ldots,Y_m)$, we proceed as above but select $W\In V$ and $Y_j\In X^{H_j}$.
For vertex testers $(\omega,\bc)$ with $\{c_i\}_{i\in I}$ and $|I|\leq k-1$,
we also proceed as above but require that $\{c_i\}_{i\in I}$ is an $|I|$-set in $V$
and $\omega$ is a function from the union over all $H\in \cH$ of the ordered $|I|$-sets of $X^H$ into $[0,\ell]$, where {$\supp(\omega)$} contains only $|I|$-tuples that are contained in an edge of $H$ {if $|I|\geq 2$}.
With this we obtain the following result.

\begin{theorem}\label{thm:normalwithweights}
For all $\alpha\in(0,1]$,
there exist $n_0,t\in \bN$ and $\eps>0$ such that the following holds for all $n\geq n_0$.
Suppose $G$ is an $(\eps,t,d)$-typical $k$-graph on $n$ vertices with $k\leq\alpha^{-1}$, $d\geq n^{-\eps}$ and 
$\cH$ is a family of $k$-graphs on $n$ vertices with $\Delta(H)\leq \alpha^{-1}$ for all $H\in \cH$ and $|\cH|\leq n^{k}$
such that $\sum_{H\in \cH}e(H)\leq (1-\alpha)e(G)$.
Suppose $\cW_{set},\cW_{ver}$ are sets of  $\alpha^{-1}$-set testers and $\alpha^{-1}$-vertex testers of size at most $n^{\log n}$, respectively.
Then there is a packing $\phi$ of $\cH$ into $G$ such that
\begin{enumerate}[label=\rm (\roman*)]
	\item\label{item:set testers2} $|W\cap \bigcap_{\ell\in [m]}\phi(Y_\ell)|= |W||Y_1|\cdots |Y_m|/n^m \pm \alpha n$ for all $(W,Y_1,\ldots,Y_m)\in \cW_{set}$;
	\item\label{item: vertex tester2} $\omega(\phi^{-1}(\bc))=(1\pm\alpha)\omega(\bigcup_{H\in \cH}V(H))/n^{|I|}\pm n^{\alpha}$ for all $(\omega,\bc)\in \cW_{ver}$ with centres $\bc$ in~$I$.
\end{enumerate}
\end{theorem}

\section{Proof overview}

In the following we outline our highlevel approach for the proof of our main result in the multipartite setting, that is, assuming we aim to pack $k$-graphs $H\in\cH$ with vertex partition $(X_1^H,\ldots,X_r^H)$ into~$G$ with vertex partition $(V_1,\ldots,V_r)$.
Our general approach is to consider each cluster $V_i$ in turn and embed simultaneously almost all vertices of $\bigcup_{H\in \cH}X_i^H$ onto $V_i$.
Afterwards we complete the embedding with another procedure.

Let us turn to a more detailed description.
At the very beginning, we will use a simple reduction to the case where each $H[X_{i_1}^H,\ldots,X_{i_k}^H]$ induces a matching for all $\{i_1,\ldots,i_k\}\in E(R)$ and $H\in\cH$
by simply splitting the clusters into smaller clusters (see Lemma~\ref{lem:refining partitions}).
This makes the analysis simpler and cleaner, and we assume this setting from now on.

The proof of the main result consists of two central parts.
The first part is an iterative procedure that considers each cluster $V_i$ of $G$ in turn 
and provides a partial packing that embeds almost all vertices of the graphs in $\cH$ onto vertices in~$G$ (at the beginning we fix some suitable ordering as~$r$ may be even larger than $n$, in particular, $r$ is not bounded by a function in terms of $\eps$).
In the $i$-th step, we embed almost all vertices of $\bigcup_{H\in \cH}X_i^H$ onto $V_i$ by finding simultaneously for each $H\in\cH$ an almost perfect matching within a `candidacy graph' $A_i^H$, which is an auxiliary bipartite graph between $X_i^H$ and $V_i$ such that $xv\in E(A_i^H)$ only if $v$ is still a suitable image for $x$ with respect to the embedding obtained in previous steps.
Of course, we have to guarantee that this indeed yields an edge-disjoint packing; 
that is, when we map $x_1\in X_i^{H_1}$ and $x_2\in X_i^{H_2}$ on the same vertex $v\in V_i$, 
then we have to ensure that for all $\eH_j\in E(H_j)$ with $x_j\in \eH_j$, $j\in[2]$, 
the $(k-1)$-sets $\eH_1\setminus \{x_1\}$ and $\eH_2\setminus \{x_2\}$ are not embedded onto the same $(k-1)$-set (provided they are already embedded).
We achieve this by defining an auxiliary hypergraph $\cH_{aux}$ with respect to the candidacy graphs to which we apply {Theorem~\ref{thm:hypermatching} --- a result on pseudorandom hypergraph matchings due to Glock and the authors~\cite{EGJ:19a} which we state in Section~\ref{sec:hypermatching}.}
There will be a bijection between matchings in $\cH_{aux}$ and valid embeddings of $\bigcup_{H\in \cH}X_i^H$ into $V_i$.
This is one of the main ingredients in the first part of our proof.

Let us give more details for the construction of $\cH_{aux}$.
Assume we are in the $i$-th step of the partial packing procedure and already found a partial packing $\phi_\circ$ of the initial $i-1$ clusters.
We define a labelling $\psi$ on the edges of the candidacy 2-graphs such that for every edge $xv\in E(A_i^H)$ and $H\in\cH$, the labelling $\psi(xv)$ contains the set of $G$-edges that are used when we extend $\phi_\circ$ by embedding $x$ onto~$v$. 
Let us for simplicity assume that only one $G$-edge would be covered when we embed~$x$ onto~$v$, say, $\psi(xv)=\gG_{xv}\in E(G)$.
Since the packing will map multiple vertices of $\bigcup_{H\in \cH}X_i^H$ onto the same vertex $v$ in $V_i$, we consider disjoint copies $(V_i^H)_{H\in\cH}$ of $V_i$ where the vertex $v^H$ is the copy of~$v$. 
For every $xv\in E(A_i^H)$ and $H\in\cH$, we define the $3$-set $\hG_{xv}:=\{x,v^H,\psi(xv)\}$ and let $\cH_{aux}$ be the $3$-graph with vertex set $\bigcup_{H\in \cH}(X_i^H \cup V_i^H) \cup E(G)$ and edge set $\{\hG_{xv}\colon xv\in E(A_i^H) \text{ for some } H\in\cH  \} $.
It is easy to see that there is a one-to-one correspondence between matchings in $\cH_{aux}$ and valid embeddings of $\bigcup_{H\in \cH}X_i^H$ into $V_i$ such that no $G$-edge is used more than once.
We used a similar hypergraph construction in our recent proof in~\cite{EJ:19} of the `blow-up lemma for approximate decompositions'.

Of course, whether we can iteratively apply this procedure depends on the choice of the partial packing in each step. 
Hence, with the aim of avoiding a future failure of the process, 
we have to maintain several pseudorandom properties throughout the entire process.
For instance, we have to ensure that there are many candidates available in each step; in more detail, we guarantee that the updated candidacy graphs after each step remain super-regular even though they naturally become sparser after each embedding step.

Unfortunately,
it is not enough to consider only candidacy graphs between pairs of clusters.
For clusters indexed by elements in $I$ where $|I|\leq k-1$, 
we consider candidacy graphs $A_I^H$ on the clusters $\bigcup_{i\in I}(X_i^H\cup V_i)$ with edges of size $2|I|$.
An edge $\aA$ in $A_I^H$ will then indicate whether the entire set of~$|I|$ vertices in $\aA\cap \bigcup_{i\in I}X_i^H$ can (still) be mapped onto $\aA\cap \bigcup_{i\in I}V_i$.
{We further discuss the purpose of these candidacy graphs in Section~\ref{sec:candidacy}, where we define them precisely.}

\smallskip
The main source yielding a smooth trajectory of our partial packing procedure is the aforementioned Theorem~\ref{thm:hypermatching} from~\cite{EGJ:19a}, which provides a tool that gives rise to a pseudorandom matching in $\cH_{aux}$ with respect to tuple-weight functions.
One difficulty of the first stage of our proof is the careful definition of these tuple-weight functions.
For instance, we have to ensure that we can indeed iteratively apply Theorem~\ref{thm:hypermatching}.
Moreover, we have to guarantee that we can turn the partial packing into a complete one in the second part of our proof. 
To that end, we define very flexible but complex weight functions on tuples of (hyper)edges of the candidacy graphs that we call \emph{edge testers}.
Dealing with hypergraphs, and especially with hypergraphs with vanishing density, makes it significantly more complex to control the weight of these edge testers during our partial packing procedure than for a similar approach for simple graphs.

One single embedding step is performed by our so-called `Approximate Packing Lemma' (see Section~\ref{sec:approx packing}).
The process where we iteratively apply our Approximate Packing Lemma is described in Section~\ref{sec:mainpart} and will provide a partial packing that maps almost all vertices of the graphs in $\cH$ onto vertices in $G$.

\smallskip
The second part of the proof deals with embedding the remaining vertices and turning the obtained partial packing into a complete one.
Our general strategy is to apply a randomized procedure where we unembed several vertices that we already embedded in the first part and
find the desired packing by using a small edge-slice $G_B$ of $G$ put aside at the beginning (that is, we did not use the edges of $G_B$ for the partial packing in the first stage).
Of course, we again have to track which vertices are still suitable images during the completion and respect the partial packing of the first part. 
To that end, for each $H\in\cH$ and $i\in[r]$, we track a second type of candidacy graphs~$B_i^H$ between $X_i^H$ and $V_i$ with respect to $G_B$, where $xv\in E(B_i^H)$ only if we could map $x$ onto $v$ during the completion. 
In fact, we track these candidacy graphs already during the partial packing procedure and carefully control several quantities using our edge testers.
In the completion step, we can then apply a randomized matching procedure within the candidacy graphs $B_i^H$ to turn the partial packing into a complete one.

\smallskip
As in many other results that were originally proven for graphs and later lifted to $k$-graphs for $k\geq 3$,
we have to overcome numerous difficulties that are specific to hypergraphs.
In our case this includes for example the much more complex intersection structure among hyperedges, which in turn complicates the analysis of our partial packing procedure considerably.
To this end, several novel ideas are needed.

Let us highlight one obstacle. 
Suppose we are in the $i$-th step of the iteration where we aim to embed essentially all vertices $\bigcup_{H\in \cH}X_i^H$ onto $V_i$,
then all $x\in \bigcup_{H\in \cH}X_i^H$ have to be grouped according to the edge intersection pattern of all edges that contain $x$ and the edges that intersect these edges with respect to the clusters that have been considered earlier.
In this context, we will define the patterns of edges in~$H$ in Section~\ref{sec:patterns}.

As we alluded to earlier, a strong control over the actual packing is of importance for further applications when an entire decomposition is sought.
In particular, it is often not enough to control how many vertices of a certain set are mapped to a particular vertex, but how many $(k-1)$-sets are embedded to a particular $(k-1)$-set.
To illustrate this, it is significantly stronger to claim that
$\Delta_{k-1}(G-\phi(\cH))\leq \alpha n$ than $\Delta(G-\phi(\cH))\leq \alpha n^{k-1}$ (where $\Delta_{k-1}(G)$ refers to the maximum number of edges containing a particular $(k-1)$-set in $G$).
It is already complicated enough to control such quantities for the multipartite setting,
but in order to transfer this ability to general quasirandom $k$-graphs,
it is necessary to allow vanishing densities of magnitude $o(\log^{-k} n)$.
This is due to the following observation.
Suppose $\cP$ is a partition of $[n]$ and say a set $S$ is \emph{crossing} (with respect to $\cP$) if each part of $\cP$ contains at most one element of $S$.
For $k\geq 3$, we need at least $\plog n$ partitions of $[n]$, which are non-trivial but have only a constant number of parts, such that all $(k-1)$-sets of~$[n]$ are (roughly) equally often crossing for the partitions,
whereas for $k=2$ one partition suffices (because $1$-element sets are always crossing).

Unfortunately, considering sparse $k$-graphs adds another complexity level to the problem.
To be more precise, $o(n)$ and $o(dn)$ no longer mean the same where $d\geq n^{-\eps}$ refers to the density,
and thus, terms of size $o(n)$ can no longer be ignored.
Essentially at all stages of the proof a substantially more careful analysis is needed to make sure that several quantities are not only $o(n)$ but $o(d^mn)$ 
because in many natural auxiliary (hyper)graphs considered in our proof vertices have typically~$d^mn$ neighbours, where $m\in \bN$ grows as we proceed in our procedure.\vspace{-.2cm}

\section{Preliminaries}
In this section we clarify notation and collect some important tools that we will frequently use throughout the paper. \vspace{-.2cm}

\subsection{Notation}\label{sec:notation}
Let us introduce some general notation and (hyper)graph terminology. 
Let $\ell\in\mathbb{N}$. We write $[\ell]_0:=[\ell]\cup\{0\}=\{0,1,\ldots,\ell\}$ and $-[\ell]:=\{-\ell,\ldots,-1\}$, where $[0]:=\emptyset$.
We refer to a set of cardinality $\ell$ as an $\ell$-set.
For finite and disjoint sets $A_1,\ldots,A_\ell$ and $I\subseteq [\ell]$, we write $A_{\cup I}:=\bigcup_{i\in I}A_i$.
For a tuple $\mathbf{a}=(a_1,\ldots,a_\ell)\in\bR^{\ell}$ and $I\subseteq [\ell]$, we write $\mathbf{a}_I:=(a_i)_{i\in I}$ and $\norm{\mathbf{a}}:=\sum_{i\in[\ell]}|a_i|$.
For a finite set $A$,
we write~$2^A$ for the powerset of~$A$ and $\binom{A}{\ell}$ for the set of all $\ell$-subsets of~$A$. 
For a graph $G$, let $\binom{E(G)}{\ell}^\cl\subseteq \binom{E(G)}{\ell}$ be the set of all matchings of size $\ell$ in $G$.
For finite sets $A_1,\ldots,A_\ell$, we write $\bigsqcup_{i\in[\ell]}A_i:=\{\{a_1,\ldots,a_\ell\}\colon a_i\in A_i \text{ for all $i\in[\ell]$} \}$,\COMMENT{The cartesian product $A_1\times\ldots\times A_\ell$ for unordered pairs.} and conversely, whenever we write $\{a_1,\ldots,a_\ell\}\in\bigsqcup_{i\in[\ell]}A_i$, we tacitly assume that $a_i\in A_i$ for all $i\in[\ell]$.
For $I\subseteq[\ell]$, we write $A_{\sqcup I}:=\bigsqcup_{i\in I}A_i$.
Whenever we consider an index set $\{i_1,\ldots,i_\ell\}\subseteq\mathbb{Z}$, we tacitly assume that $i_1\leq i_2\leq\ldots\leq i_\ell$.
For a real-valued function $f\colon A\to \bR_{\geq 0}$, let $\supp(f):=\{a\in A\colon f(a)>0 \}$ be its support.
For a function $g\colon A\to B$, let $g(A'):=\bigcup_{a\in A'\cap A}g(a)$.
We often write $a\ll b \ll c$ in our statements meaning that $a,b,c\in(0,1]$ and there are increasing functions $f,g:(0,1]\to (0,1]$
such that whenever $a\leq f(b)$ and $b \leq g(c)$,
then the subsequent result holds.
{We also write $c=a\pm b$ if $a-b\leq c\leq a+b$.}
For $a\in (0,1]$ and $\mathbf{b}\in(0,1]^\ell$, we write $\mathbf{b}\geq a$ whenever $b_i\geq a$ for all $b_i\in\mathbf{b}$, $i\in[\ell]$.
Whenever we consider a parameter with letter $d$, usually used for the density, we tacitly assume that $d\in(0,1]$.
For the sake of a clearer presentation, we avoid roundings whenever it does not affect the argument.\looseness=-2 

For a $k$-graph $G$, let $V(G)$ ad $E(G)$ denote the vertex set and edge set, respectively.
For $X\subseteq V(G)$, let~$G[X]$ be the $k$-graph induced on $X$, and for pairwise disjoint subsets $V_1,\ldots,V_k\subseteq V(G)$, let $G[V_1,\ldots,V_k]$ be the $k$-partite subgraph of $G$ induced between $V_1,\ldots,V_k$.
Let $e(G)$ denote the number of edges in~$G$ and let $e_G(V_1,\ldots,V_k):=e(G[V_1,\ldots,V_k])$.
For set $S$ of at most $k-1$ vertices of $G$, we define the neighbourhood of $S$ in $G$ by $N_G(S):=\{\gG\sm S\colon \gG\in E(G), S\subseteq\gG \}$ and let $\dg_G(S):=|N_G(S)|$. 
Note that $N_G(S)$ is a set of $(k-|S|)$-tuples. 
For $m\in[k-1]$, let $\Delta_{m}(G):=\max_{S\in \binom{V(G)}{m}}\dg_G(S)$ denote the maximum $m$-degree of~$G$.
We usually write $\Delta(G)$ instead of $\Delta_1(G)$.
Further, we say that $u$ is a $G$-vertex if $u\in V(G)$, and $u$ is a $G$-neighbour of $v\in V(G)$ if $u$ and $v$ are contained in an edge of $G$.
If $G$ is a $2$-graph, $u$ and $v$ are vertices of $G$, and $S\subseteq V(G)$, we let $N_G[v]:=N_G(v)\cup\{v\}$, $N_G(u\land v):= N_G(u)\cap N_G(v)$,
and $N_G(S):=(\bigcup_{v\in S} N_G(v))\sm S$.
We frequently treat collections of (hyper)graphs as the (hyper)graph obtained by taking the disjoint union of all members.\looseness=-1

We say a $k$-graph $G$ on $n$ vertices is \defn{$(\eps,t,d)$-typical} if for all non-empty sets $\cS\In\binom{V(G)}{k-1}$ with $|\cS|\leq t$, we have $|\bigcap_{S\in\cS}N_G(S)|=(1\pm\eps)d^{|\cS|}n$.
We often write $N_G(\cS):=\bigcap_{S\in\cS}N_G(S)$.
Throughout the paper, we usually denote a $(k-1)$-set with the letter~$S$, and a set of $(k-1)$-sets with the letter~$\cS$.\looseness=-1

For a $k$-graph $G$, we denote by $G_\ast$ the $2$-graph with vertex set $V(G)$ and edge set $\bigcup_{\gG\in E(G)}\binom{\gG}{2}$. 
That is, $G_\ast$ arises from $G$ by replacing each hyperedge in $G$ with a clique of  size $k$.

For (hyper)graphs $G$ and $H$, we write $G-H$ to denote the (hyper)graph with vertex set $V(G)$ and edge set $E(G)\sm E(H)$.

For a graph $G$ and a finite set $\cE$, we call $\psi\colon E(G)\to2^\cE$ an \emph{edge set labelling} of~$G$. 
A~label $\alpha\in\cE$ \emph{appears} on an edge $e$ if $\alpha\in\psi(e)$. 
Let $\norm{\psi}$ be the maximum number of labels that appear on any edge of $G$.
We define the \emph{maximum degree $\Delta_\psi(G)$ of $\psi$} as the maximum number of edges of $G$ on which any fixed label appears, and the \emph{maximum codegree $\Delta^c_\psi(G)$ of $\psi$} as the maximum number of edges of $G$ on which any two fixed labels appear together.

For a graph $G$ and $m\in\bN$, let $G^m$ denote the $m$-th power of $G$, that is, the $2$-graph which is obtained from $G$ by adding all edges between vertices whose distance in $G$ is at most $m$.

\subsection{Concentration inequalities}\label{sec:prob tools}
To verify the existence of subgraphs with certain properties we frequently consider random subgraphs and use McDiarmid's inequality to verify that specific random variables are highly concentrated around their mean.

\begin{theorem}[McDiarmid's inequality, see~\cite{mcdiarmid:89}\COMMENT{Lemma~1.2}] \label{thm:McDiarmid}
Suppose $X_1,\dots,X_m$ are independent random variables and suppose $b_1,\dots,b_m\in [0,B]$.
Suppose $X$ is a real-valued random variable determined by $X_1,\dots,X_m$ such that changing the outcome of $X_i$ changes $X$ by at most $b_i$ for all $i\in [m]$.
Then, for all $t>0$, we have 
\vspace{-.1cm}
$$\prob{|X-\expn{X}|\ge t} \le 2 \exp\left({-\frac{2t^2}{B\sum_{i=1}^m b_i}}\right).$$
\end{theorem}

We also state the following convenient form of Freedman's inequality~\cite{freedman:75}.

\begin{lemma}\label{lem:convenient freedman}
Suppose $X,X_1,\ldots,X_m$ are real-valued random variables with $X=\sum_{i\in [m]} X_i$ such that $|X_i|\leq B$ and $\sum_{i\in [m]}\mathrm{\mathbb{E}}_i\left[|X_i|\right]\leq \mu$, where $\mathrm{\mathbb{E}}_i\left[X_i\right]$ denotes the expectation conditional on any given values of $X_j$ for $j<i$.
Then 
\vspace{-.1cm}
$$\prob{|X|>2\mu}\leq 2\exp\left(-\frac{\mu}{4B}\right).$$
\end{lemma}

\subsection{Graph regularity}\label{sec:graph regularity}
In this section we introduce the quasirandom notion of (sparse) $\eps$-regularity for $2$-graphs and collect some important results.
For a bipartite graph $G$ with vertex partition $(V_1,V_2)$,
we define the \defn{density} of the pair $W_1,W_2$ with $W_i\In V_i$ by $d_G(W_1,W_2):=e_G(W_1,W_2)/|W_1||W_2|$.
We say $G$ is \defn{$(\eps,d)$-regular} if $d_G(W_1,W_2)=(1\pm\eps)d$ for all $W_i\In V_i$ with $|W_i|\geq \eps |V_i|$,
and $G$ is \defn{$(\eps,d)$-super-regular} if in addition $|N_G(v)\cap V_{3-i}|= (1\pm \eps)d|V_{3-i}|$ for each $i\in [2]$ and $v\in V_i$.

The following two standard results concern the robustness of $\eps$-regular graphs.\nopagebreak
\begin{fact} \label{fact:regularity}\nopagebreak
Suppose $\eps,d\in(0,1]$ and $G$ is an $(\eps,d)$-regular bipartite graph with vertex partition~$(A,B)$ and $Y\In B$ with $|Y|\ge \eps|B|$. Then all but at most $2\eps|A|$ vertices of $A$ have $(1\pm \eps)d|Y|$ neighbours in~$Y$.
\end{fact}
\COMMENT{Proof:
\\
Let $A^{high}\subseteq A$ be the vertices in $A$ that have more than $(1+\eps)d|Y|$ neighbours in $Y$ and assume for a contradiction that, $|A^{high}|>\eps|A|$.
This implies
$$e_G(A^{high},Y)>(1+\eps)d|A^{high}||Y|$$
but
$$d_G(A^{high},Y)=\frac{e_G(A^{high},Y)}{|A^{high}||Y|}\leq (1+\eps)d.$$
This implies the contradiction.
Analogously for $A^{low}$. 
}

\begin{fact}\label{fact:regularity robust}
Suppose $1/n \ll \eps$ and $d\geq n^{-\eps}$.
Suppose $G$ is an $(\eps,d)$-super-regular bipartite graph with vertex partition $(A,B)$, where $\eps^{1/6} n\le |A|,|B|\le n$. If $H$ is a subgraph of $G$ with $\Delta(H)\le \eps dn$ and $X\In A\cup B$ with $|X|\le \eps dn$, then $G[A\sm X,B\sm X]-H$ is $(\eps^{1/3},d)$-super-regular.\COMMENT{Modified for sublinear density $d$. 
Proof: The number of edges between two vertex sets $Z_1$ and $Z_2$ is at least 
$(1-\eps)d|Z_1||Z_2|- (\eps d n + \eps d n)(|Z_1|+|Z_2|)\geq (1- \eps - 4\eps^{1/2})d|Z_1||Z_2|\geq (1-\eps^{1/3})d|Z_1||Z_2|$ whenever $|Z_i|\geq \eps^{1/2}n= \eps^{1/3} \cdot \eps^{1/6}n$.}
\end{fact}

The following result from~\cite{ABSS:20} is useful to establish $\eps$-regularity for sparse graphs.
For a graph $G$, let $C_4(G)$ be the number of 4-cycles in $G$.

\begin{theorem}[{\cite[Lemma 13]{ABSS:20}}]\label{thm:new almost quasirandom}
Suppose $1/n\ll\eps$ and $d\geq n^{-\eps}$.
Suppose $G$ is a bipartite graph with vertex partition~$(A,B)$, $|A|,|B|=(1\pm\eps)n$, density $d_G(A,B)=(1\pm\eps)d$, and $C_4(G)< (1+\eps)d^4|A|^2|B|^2/4$.
Then $G$ is $(\eps^{1/13},d)$-regular.
\end{theorem}
\COMMENT{Replaced:
\\
We will also use the next result from~\cite{DLR:95}.
(In~\cite{DLR:95} it is proved in the case when $|A|=|B|$ with $16\eps^{1/5}$ instead of $\eps^{1/6}$.
The version stated below can be easily derived from this.)
\begin{theorem}\label{thm: almost quasirandom}
Suppose $1/n\ll\eps \ll \gamma$ and $d>0$.
Suppose $G$ is a bipartite graph with vertex partition~$(A,B)$ such that $|A|=n$, $\gamma n \leq |B|\leq \gamma^{-1}n$ and at least $(1-5\eps)n^2/2$ pairs $u,v\in A$ satisfy $\dg_G(u),\dg_G(v)\geq (1- \eps)d|B|$ and $|N_G(u)\cap N_G(v)|\leq (1+ \eps)d^2|B|$.
Then $G$ is $(\eps^{1/6},d)$-regular.
\COMMENT{Modified for sublinear density $d$.} 
\end{theorem}}

Further, if the common neighbourhood of most pairs in an $(\eps,d)$-super-regular graph has the appropriate size, we can establish the following useful bounds on the number of edges between subsets of vertices that we allow to be very small.
\begin{lemma}\label{lem:edges in sparse graphs}
Suppose $1/n\ll \eps$ and $d\geq n^{-\eps}$. 
Suppose $G$ is an $(\eps,d)$-super-regular graph with bipartition $(A,B)$ and $|A|=|B|=n$, and for all but at most $n^{3/2}$ pairs $\{a,a'\}$ in $A$, we have
$|N_G(a\land a')|\leq (1+\eps)d^2n$.
If $X\subseteq A$ and $Y\subseteq B$ with $n^{3/4+3\eps}\leq |X|, |Y| \leq \eps n$, then $e_G(X,Y)\leq \eps^{1/3}dn \max\{|X|,|Y|\}$.
\end{lemma}
\COMMENT{We also had the following statement that we did not use: \\If $n^{3/4+3\eps}\leq |X|$ and $|Y|=\gamma n$ for $\eps\ll\gamma$, then $e_G(X,Y)\leq(1+\eps^{1/3})d|X||Y|$.
\\Proof: 
We first establish \ref{item:sparse1}. 
We have that
\begin{align*}
&\frac{1}{2}\sum_{b\in Y} |N_G(b)\cap X|^2 =
\sum_{\{a,a'\}\in\binom{X}{2}}|N_G(a\land a')\cap Y|
\\&\leq\sum_{\{a,a'\}\in\binom{X}{2}}
\left((1+\eps)d^2n- |N_G(a\land a')\cap B\sm Y|\right)
+n^{5/2}
\\&\leq \frac{1}{2}(1+2\eps) |X|^2d^2n - \frac{1}{2}(1-\gamma-\eps^{1/2})|X|^2d^2n
\\&\leq  
\frac{1}{2}(\gamma+\eps^{1/2})|X|^2d^2n.
\end{align*}
Combining this with
\begin{align*}
\frac{1}{2}\sum_{b\in Y} |N_G(b)\cap X|^2 \geq \frac{e_G(X,Y)^2}{2|Y|},
\end{align*}
we obtain $e_G(X,Y)^2\leq (\gamma+\eps^{1/2})d^2|X|^2|Y|n$, and thus
\begin{align*}
e_G(X,Y)
\leq (1+\eps^{1/3}) d|X||Y|,
\end{align*} 
where we used that $|Y|=\gamma n$ for $\eps\ll\gamma$.
}

\begin{proof}
We have that
\begin{align*}
\frac{1}{2}\sum_{b\in Y} |N_G(b)\cap X|^2 \leq
\sum_{\{a,a'\}\in\binom{X}{2}}|N_G(a\land a')\cap Y| + \eps n^2
\leq \frac{1}{2}(1+2\eps) |X|^2d^2n.
\end{align*}
Combining this with 
$\frac{1}{2}\sum_{b\in Y} |N_G(b)\cap X|^2 \geq \frac{e_G(X,Y)^2}{2|Y|}$
yields that
\begin{align}\label{eq:bound edges codegrees}
e_G(X,Y)^2\leq (1+2\eps) d^2|X|^2|Y|n.
\end{align}

Suppose first that $|Y|\leq |X|$, and suppose for a contradiction that $e_G(X,Y)\geq \eps^{1/3}dn|X|$.
{By inserting this into~\eqref{eq:bound edges codegrees} and solving for $|Y|$,} this implies that $|Y|\geq \eps^{2/3}n/3$,
which is a contradiction to $|Y|\leq|X|\leq \eps n$.

Next, suppose $|X|\leq |Y|$, and suppose for a contradiction that $e_G(X,Y)\geq \eps^{1/3}dn|Y|$.
{By inserting this into~\eqref{eq:bound edges codegrees} and solving for $|X|$,} this implies that 
$|X|\geq {\eps^{1/3}}(|Y|n)^{1/2}/2.$
Since $|Y|\geq|X|$, this yields that $|Y|\geq \eps^{2/3}n/4$, which is a contradiction to $|Y|\leq \eps n$.
\end{proof}

We will also need the following result that is similar to~\cite[Lemma 2]{ARR:98} and guarantees that a (sparse) $(\eps,d)$-super-regular balanced bipartite graph of order $2n$ contains a spanning $m$-regular subgraph (an $m$-factor) for $m=(1-2\eps^{1/3})dn$ provided that most pairs of vertices have the appropriate number of common neighbours.
It can be proved along the same lines as in~\cite{ARR:98} by employing Lemma~\ref{lem:edges in sparse graphs}.\looseness=-1

\begin{lemma}\label{lem:m-factor}
Suppose $1/n\ll\eps$ and $d\geq n^{-\eps}$.
Suppose $G$ is an $(\eps,d)$-super-regular bipartite graph with vertex partition $(A,B)$ where $|A|=|B|=n$, and suppose that all but at most $n^{3/2}$ pairs $\{a,a'\}$ in $A$ satisfy that $|N_G(a\land a')|\leq (1+\eps)d^2n$.
Then $G$ contains an $m$-factor for $m=(1-2\eps^{1/3})dn$.
\end{lemma}

\COMMENT{As in the dense case~\cite{ARR:98}, we will use the following characterization for the proof of Lemma~\ref{lem:m-factor}.
\begin{theorem}
\label{thm:char for m-factor}
Suppose $G$ is a bipartite graph with vertex partition $(A,B)$ where $|A|=|B|=n$.
Then $G$ has an $m$-factor if and only if for all $X\subseteq A$, $Y\subseteq B$, we have 
$$m|X|+m|Y|+e_G(A\sm X, B\sm Y)\geq mn.$$
\end{theorem}
\lateproof{Lemma~\ref{lem:m-factor}}
Let $m=(1-2\eps^{1/3})dn$ and let $X$ and $Y$ be arbitrary subsets of $A$ and $B$, respectively. 
Let $\bar{X}:=A\sm X$ and $\bar{Y}:= B\sm Y$.
We establish the inequality in Theorem~\ref{thm:char for m-factor}.
The inequality is clearly satisfied if $|X|+|Y|\geq n$, and thus we may assume that $|X|+|Y|< n$ and hence, $|\bar{X}|+|\bar{Y}|>n$, and $|Y|<|\bar{X}|$.
We show that $e_G(\bar{X},\bar{Y})$ is sufficiently large.
\\
First, suppose that the size of either $\bar{X}$ or $\bar{Y}$ is less than $n^{3/4+3\eps}$, say, $|\bar{X}|\leq n^{3/4+3\eps}$. 
Since $|Y|<|\bar{X}|$, we obtain
\begin{align}\label{eq:X tiny}
e_G(\bar{X},\bar{Y})
\geq 
(1-\eps)|\bar{X}|dn-|\bar{X}||Y|
\geq 
(1-2\eps) |\bar{X}|dn.
\end{align}
Next, suppose that the size of either $\bar{X}$ or $\bar{Y}$ is at least $n^{3/4+3\eps}$ but at most $\eps n$, say, $ n^{3/4+3\eps}\leq|\bar{X}|\leq\eps n$. 
Since $|Y|<|\bar{X}|$, we obtain with Lemma~\ref{lem:edges in sparse graphs} that
$e_G(\bar{X},Y) \leq \eps^{1/3}|\bar{X}|dn$ and thus,
\begin{align}\label{eq:X small}
e_G(\bar{X},\bar{Y})
\geq 
(1-\eps)|\bar{X}|dn - \eps^{1/3}|\bar{X}|dn
\geq 
(1-2\eps^{1/3})|\bar{X}|dn.
\end{align}
By~\eqref{eq:X tiny} and~\eqref{eq:X small}, we obtain in both cases that
\begin{align*}
m|X|+m|Y|+e_G(\bar{X},\bar{Y}) 
\geq
m|X|+m|Y|
+(1-2\eps^{1/3})|\bar{X}|dn \geq mn.
\end{align*}
Now, suppose that $|\bar{X}|,|\bar{Y}|\geq \eps n$.
We obtain
\begin{align*}
m|X|+m|Y|+e_G(\bar{X},\bar{Y}) 
&\geq
m|X|+m|Y|+
(1-\eps)d|\bar{X}|\bar{Y}|
\\
&\geq
(1-2\eps^{1/3})d
\left(n|X|+n|Y|+(n-|X|)(n-|Y|) \right)
\\&
\geq(1-2\eps^{1/3})dn^2.
\end{align*}
This completes the proof.
\endproof}

\subsection{Pseudorandom hypergraph matchings}\label{sec:hypermatching}
A key ingredient in the proof of our `Approximate Packing Lemma' in Section~\ref{sec:approx packing} is the main result from~\cite{EGJ:19a} on pseudorandom hypergraph matchings.

For this we need some more notation.
Given a finite set $X$ and an integer $\ell\in \bN$, 
an \defn{$\ell$-tuple weight function on~$X$} is a function $\omega\colon \binom{X}{\ell}\to \bR_{\ge 0}$.
For a subset $X'\In X$, we then define $\omega(X'):=\sum_{S\in \binom{X'}{\ell}}\omega(S)$. 
For $m\in[\ell]$, let $\normV{\omega}{m}:=\max_{T\in \binom{X}{m}}\sum_{S\colon T\subseteq S \in \binom{X}{\ell}}\omega(S)$.\COMMENT{Previously:\\For $m\in[\ell]_0$ and a tuple $T\in \binom{X}{m}$, define 
$\omega(T):=\sum_{S\supseteq T}\omega(S)$ and let $\normV{\omega}{m}:=\max_{T\in \binom{X}{m}}\omega(T)$.\\However, the notation for $\omega(T)$ is never used again.} 
Suppose $H$ is an $r$-uniform hypergraph and $\omega$ is an $\ell$-tuple weight function on~$E(H)$. Clearly, if $\cM$ is a matching, then a tuple of edges which do not form a matching will never contribute to~$\omega(\cM)$. We thus say that $\omega$ is \defn{clean} if $\supp(\omega)\subseteq\binom{E(H)}{\ell}^=$, that is, $\omega(\cE)=0$ whenever $\cE\in\binom{E(H)}{\ell}$ is not a matching.

\begin{theorem}[\cite{EGJ:19a}]\label{thm:hypermatching}
Suppose $1/\Delta \ll \delta,1/r,1/L$, $r\geq 2$, and let $\eps:=\delta/50L^2r^2$.
Let $H$ be an $r$-uniform hypergraph with $\Delta(H)\leq \Delta$ and $\Delta_2(H)\le \Delta^{1-\delta}$ as well as $e(H)\leq \exp(\Delta^{\eps^2})$. 
Suppose that for each $\ell\in[L]$, we are given a set $\cW_\ell$ of clean $\ell$-tuple weight functions on $E(H)$ of size at most $\exp(\Delta^{\eps^2})$ such that $\omega(E(H))\ge \normV{\omega}{m}\Delta^{m+\delta}$ for all $\omega\in \cW_\ell$ and $m\in [\ell]$.
Then there exists a matching~$\cM$ in~$H$ such that $\omega(\cM)=(1\pm \Delta^{-\eps}) \omega(E(H))/\Delta^\ell$ for all $\ell\in [L]$ and $\omega \in \cW_\ell$.
\end{theorem}

\subsection{Refining partitions}\label{sec:refining partitions}
In this section we provide a useful result to refine the vertex partition of a collection $\cH$ of $k$-graphs of bounded degree such that every $H\in\cH$ only induces a matching between any $k$-set of the refined partition.
For $2$-graphs, a similar approach was already used in~\cite{RR:99} by simply applying the classical Hajnal--Szemer\'edi Theorem.
Our result is based on a random procedure which enables us to sufficiently control the weight distribution of weight functions with respect to the refined partition.
We used such a procedure already in~\cite{EJ:19}.

\begin{lemma}\label{lem:refining partitions}
	Suppose $1/n\ll\eps\ll \beta \ll \alpha,1/k$ and $r\leq n^{\log n}$.
	Suppose $\cH$ is a collection of at most~$n^{2k}$\COMMENT{This is not restrictive and $n^{2k}\geq n^k$ is chosen arbitrarily.}  $k$-graphs, $(X_i^H)_{i\in[r]}$ is a vertex partition of each $H\in\cH$, and $\Delta(\cH)\leq \alpha^{-1}$.
	Suppose $n/2\leq |X_i^H|=|X_i^{H'}| \leq 2n$ for all $H,H'\in \cH$ and $i\in [r]$.
	Suppose $\cW$ is a set of at most $n^{5\log n}$ weight functions $\omega\colon \sX_{\sqcup {I}}\to [0,\alpha^{-1}]$ with $I\subseteq[r]$, $|I|\leq k$, {and whenever $|I|\geq 2$, we have}
	$\supp(\omega)\subseteq \{\xX\in\sX_{\sqcup I} \colon \xX\subseteq \eH \text{ for some $\eH\in E(\cH)$} \}$.\COMMENT{This implies in particular that $\omega(X_{\cup[r]}^H)\leq 2\alpha^{-2}n$ for every $H\in\cH$.} 
	Then for all $H\in \cH$ and $i\in [r]$, there exists a partition $(X_{i,j}^H)_{j\in [\beta^{-1}]}$ of $X^H_i$ such that 
	\begin{enumerate}[label=\rm (\roman*)]
		\item\label{item:split1} $X^H_{i,j}$ is an independent set in $H_\ast^2$ for all $H\in \cH$, $i\in[r]$, $j\in[\beta^{-1}]$;
		\item\label{item:split2} $|X_{i,1}^H|\leq \ldots \leq |X_{i,\beta^{-1}}^H|\leq |X_{i,1}^H|+1$ for all $H\in \cH$ , $i\in[r]$;

\item\label{item:split weights one H, |I|=1}
$\omega(X_{i,j}^H)=\beta \omega(X_i^H)\pm\beta^{3/2}n$ for all $H\in\cH$, $i\in[r]$, $j\in[\beta^{-1}]$, and $\omega\in\cW$ with $\omega\colon \sX_i\to[0,\alpha^{-1}]$;\COMMENT{For $|I|=1$, we have a different error bound as the random permutations do not have any effect, therefore `$\pm \beta^{3/2}n$'.}

\item\label{item:split weights} $\omega(\bigcup_{H\in \cH}(\bigsqcup_{\ell\in[|I|]}X_{i_\ell,j_\ell}^H))
=
(1\pm  \eps)
\beta^{|I|} \omega(\sX_{\sqcup I})$ for all $\omega\in\cW$ with $\omega\colon \sX_{\sqcup {I}}\to [0,\alpha^{-1}]$, $\omega(\sX_{\sqcup {I}})\geq n^{1+\eps}$, $I=\{i_1,\ldots,i_{|I|}\}\subseteq[r]$,\COMMENT{We assume $I$ is not a multiset, i.e. $I\in\binom{[r]}{|I|}$.} $|I|\leq k$, and $j_1,\ldots,j_{|I|}\in [\beta^{-1}]$.
\COMMENT{
This will in particular also imply the required estimate for the number of $\cH$-edges:
\\
$\sum_{H\in\cH}e_H(\bigcup_{\ell\in[k]}X_{i_\ell,j_\ell}^H)= (1\pm k\beta)\beta^k\sum_{H\in \cH}e_H(X_{i_1}^H,\ldots,X_{i_k}^H)$ for all $\{i_1,\ldots,i_k\}\in\binom{[r]}{k}$ and $j_1,\ldots,j_k\in[\beta^{-1}]$.
}
	\end{enumerate}
\end{lemma}

We give an analogous statement for the non-multipartite setting, that is, when $r=1$, where the given weight functions assign weight on vertex tuples.\COMMENT{If we do not show the details of the proof of Theorem~\ref{thm:normalwithweights}, then it might suffice to say that for the proof of Theorem~\ref{thm:normalwithweights}, a non-multipartite version for weight functions on ordered vertex tuples can be stated and proved analogously as Lemma~\ref{lem:refining partitions}.}
{For a finite set $A$ and $\ell\in\bN$, let $\binom{A}{\ell}_\prec$ be the set of all $\ell$-tuples with non-repetitive entries of $A$.}\COMMENT{Needs to be checked again. Also~\ref{item:split' weights} of Lemma~\ref{lem:refining partitions r=1}.}
\begin{lemma}\label{lem:refining partitions r=1}
Suppose $1/n\ll\eps\ll \beta \ll \alpha,1/k$.
Suppose $\cH$ is a collection of at most $n^{2k}$\COMMENT{This is not restrictive and $n^{2k}\geq n^k$ is chosen arbitrarily.}  $k$-graphs  on $n$ vertices with $\Delta(\cH)\leq \alpha^{-1}$.
Suppose $\cW$ is a set of at most $n^{5\log n}$ weight functions $\omega\colon \bigcup_{H\in \cH}\binom{V(H)}{m}_\prec\to[0,\alpha^{-1}]$ with $m\in[k-1]$, and whenever $m\geq 2$, we have $\supp(\omega)\subseteq \bigcup_{H\in \cH}\{\xX\in\binom{V(H)}{m}_\prec\colon \xX\subseteq \eH \text{ for some $\eH\in E(H)$} \}$.
Then for all $H\in\cH$, there exists a partition $(X_j^H)_{j\in [\beta^{-1}]}$ of $V(H)$ such that
	\begin{enumerate}[label=\rm (\roman*)]
		\item\label{item:split1'} $X^H_{j}$ is an independent set in $H_\ast^2$ for all $H\in \cH$, $j\in[\beta^{-1}]$;
		\item\label{item:split2'} $|X_{1}^H|\leq \ldots \leq |X_{\beta^{-1}}^H|\leq |X_{1}^H|+1$ for all $H\in \cH$;

\item\label{item:split' weights one H, |I|=1}
$\omega(X_{j}^H)=\beta \omega(V(H))\pm\beta^{3/2}n$ for all $H\in\cH$, $j\in[\beta^{-1}]$, and $\omega\in\cW$ with $\omega\colon \bigcup_{H\in \cH}V(H)\to[0,\alpha^{-1}]$;\looseness=-1

\item\label{item:split' weights} $\omega(\bigcup_{H\in \cH}(X_{j_1}^H\times\cdots\times X_{j_m}^H))
=
(1\pm \beta^{1/2})
\beta^{m} \omega(V(\cH))$ for all $\omega\in\cW$ with $\omega\colon \bigcup_{H\in \cH}\binom{V(H)}{m}_\prec\to[0,\alpha^{-1}]$, $\omega(V(\cH))\geq n^{1+\eps}$, and $\{j_1,\ldots,j_{m}\}\in\binom{[\beta^{-1}]}{m}$.
	\end{enumerate}
\end{lemma}

Lemmas~\ref{lem:refining partitions} and~\ref{lem:refining partitions r=1} are very similar to the analogous result for $2$-graphs in~\cite[Lemma~3.7]{EJ:19}.
The main novelty in our setting is that we allow for weight functions on tuples of vertices.
Nevertheless, the proofs follow exactly the same strategy as in the proof our result for $2$-graphs.
In the following, we therefore just provide a short proof sketch, say, for Lemma~\ref{lem:refining partitions}.
(A detailed proof can be found in~\cite{ehard:phd}.)\looseness=-2

We first consider every $H\in\cH$ in turn and by randomly partitioning every $X_i^H$ into $\beta^{-1}$ sets, we obtain a partition $(Y_{i,j}^H)_{j\in[\beta^{-1}]}$ that essentially satisfies~\ref{item:split1}--\ref{item:split weights one H, |I|=1} with $Y_{i,j}^H$ playing the role of~$X_{i,j}^H$.
Then we perform a vertex swapping procedure to replace edges in $H_\ast^2$ between vertices in $(Y_{i,j}^H)_{j\in[\beta^{-1}]}$ and obtain $(Z_{i,j}^H)_{j\in[\beta^{-1}]}$ which satisfies~\ref{item:split1}--\ref{item:split weights one H, |I|=1}.
In the end, for every $H\in\cH$, $i\in[r]$, we randomly permute the ordering of $(Z_{i,j}^H)_{j\in[\beta^{-1}]}$ to also ensure~\ref{item:split weights}.

\COMMENT{\lateproof{Lemma~\ref{lem:refining partitions}}
We follow the same strategy as in the proof our result in~\cite[Lemma~3.7]{EJ:19} for $2$-graphs. 
Here, the main novelty is property~\ref{item:split weights} where we also allow for weight functions on tuples of vertices.
Our approach is as follows.
We first consider every $H\in \cH$ in turn and construct a partition $(Y_{i,j}^H)_{j\in[\beta^{-1}]}$ that essentially satisfies~\ref{item:split1} and~\ref{item:split2} with $Y_{i,j}^H$ playing the role of $X_{i,j}^H$. 
Then we perform a vertex swapping procedure to resolve some conflicts in $Y_{i,j}^H$ and obtain $Z_{i,j}^H$. 
In the end, we randomly permute the ordering of $(Z_{i,j}^H)_{j\in[\beta^{-1}]}$ for each $H\in \cH,i\in [r]$ to also ensure~\ref{item:split weights}.
\\
To simplify notation, we assume from now on that $|X_i^H|$ is divisible by $\beta^{-1}$ for all $H\in\cH$, $i\in [r]$,
and at the end we explain how very minor modifications yield the general case.
Recall that $H_\ast$ is the $2$-graph on $V(H)$ that arises from $H$ by replacing each hyperedge with a clique of size $k$.
For future reference, we also recall that for every $\omega\in\cW$ with $\omega\colon \sX_{\sqcup {I}}\to [0,\alpha^{-1}]$ with $I\subseteq[r]$, $|I|\leq k$, we have by assumption that
\begin{align}\label{eq:weight only on edges}
\supp(\omega)\subseteq \{\xX\in\sX_{\sqcup I} \colon \xX\subseteq\eH  \text{ for some $\eH\in\cH$} \}.
\end{align}
Note in particular that~\eqref{eq:weight only on edges} implies that  $\omega(X_{\sqcup I}^H)\leq 2k\alpha^{-2}n$ for every $H\in\cH$.
\\
Let $H\in\cH$ be fixed.
We claim that there exist partitions $(Y_{i,j}^H)_{j\in [\beta^{-1}]}$ of $X_i^H$ for each $i\in [r]$ such that for all $i\in [r],j\in [\beta^{-1}]$
\begin{enumerate}[label=(\alph*)]
	\item\label{item:correct split1} $|Y_{i,j}^H|=\beta |X_i^H|\pm\beta^2n$;
	\item\label{item:correct split2} at most $\beta^{9/5} n$ pairs of vertices in $Y_{i,j}^H$ are adjacent in $H_\ast^2$;\COMMENT{To be precise: at most $2\alpha^{-2}k\beta^2 n$ pairs. There are at most $2n\alpha^{-2}k/2$ many $P_3$'s in $H_\ast$ that start and end in $X_i^H$, and at most $n\alpha^{-1}k$ edges in $X_i^H$. Each endpoint of such a $P_3$ or $P_2$ lies in $Y_{i,j}^H$ with probability $\beta$.}
	\item\label{item:correct split3}
	$\omega(Y_{i,j}^H)=\beta \omega(X_i^H)\pm\beta^2n$ for all $\omega\in \cW$ with $\omega\colon \sX_{i}\to [0,\alpha^{-1}]$.
\end{enumerate}
Indeed, the existence of such partitions can be seen by assigning every vertex in~$X_i^H$ uniformly at random to some $Y_{i,j}^H$ for $j\in [\beta^{-1}]$. 
Together with a union bound and Theorem~\ref{thm:McDiarmid}, we conclude that~\ref{item:correct split1}--\ref{item:correct split3} hold simultaneously with positive probability (where we employ~\eqref{eq:weight only on edges} for~\ref{item:correct split3}).
\\
Next, we slightly modify these partitions $Y_{i,j}^H$ to obtain a new collection of partitions $Z_{i,j}^H$. 
For all $i\in[r]$, $j\in[\beta^{-1}]$, let $W_{i,j}^H\In Y_{i,j}^H$ be such that $|Y_{i,j}^H \setminus W_{i,j}^H|=\beta|X_i^H|- \beta^{5/3}n$ and $W_{i,j}^H$ contains all vertices in $Y_{i,j}^H$ that contain an $H_\ast^2$-neighbour in $Y_{i,j}^H$ (the sets $W_{i,j}^H$ clearly exist by~\ref{item:correct split1} and \ref{item:correct split2}).
For all $i\in[r]$, $j\in[\beta^{-1}]$, let $\{w_1^i,\ldots,w_s^i\}=W_i^H:=\bigcup_{j\in [\beta^{-1}]}W_{i,j}^H$ and observe that $s=|X_i^H|-\sum_{j\in[\beta^{-1}]}|Y_{i,j}^H\sm W_{i,j}^H|=\beta^{2/3}n$.
\\
Now for every $i\in[r]$, arbitrarily assign labels in $[\beta^{-1}]$ to the vertices in $W_i$ such that each label is used exactly~$\beta^{5/3}n$ times.
Let  $Z_{i,j}^H(0):=Y_{i,j}^H\setminus W_{i,j}^H$ for all $i\in[r]$, $j\in[\beta^{-1}]$.
To obtain the desired partitions we perform the following swap procedure for every $i\in[r]$. 
For every $t\in [s]$ in turn we do the following.
Say $w_t^i\in W_{i,j}^H$ and $w_t^i$ received label $j'$.
We select $j''\in [\beta^{-1}]\setminus \{j,j'\}$ such that $w_t^i$ has no $H_\ast^2$-neighbour in $Z_{i,j''}^H(t-1)$ and such that $Z_{i,j''}^H(t-1)$ contains a vertex $w$ that has no $H_\ast^2$-neighbour in $Z_{i,j'}^H(t-1)$.
In such a case we say that $j''$ is \defn{selected} in step $t$.
Then we define $Z_{i,j''}^H(t):= (Z_{i,j''}^H(t-1)\cup \{w_t^i\})\setminus \{w\}$, $Z_{i,j'}^H(t):=Z_{i,j'}^H(t)\cup \{w\}$ and $Z_{i,\ell}^H(t):=Z_{i,\ell}^H(t-1)$ for all $\ell\in [\beta^{-1}]\setminus \{j',j''\}$.
Note that $\beta |X_i^H| - \beta^{5/3}n \leq |Z_{i,\ell}(t)|\leq \beta |X_i^H|$ for all $t\in [s]$, $\ell\in [\beta^{-1}]$.
Observe also that we have always at least $\beta^{-1}/2$ choices to select $j''$ in step $t$.\COMMENT{There are at most $2\alpha^{-2}k\beta n$ many $P_3$'s in $H_\ast$ that start in $Y_{i,j'}$, and at most $2\alpha^{-1}k\beta n$ many $P_2$'s in $H_\ast$ that start in $Y_{i,j'}$. 
Hence, at most $6\alpha^{-2}$ entries $j''\in[\beta^{-1}]\sm\{j,j'\}$ are such that every vertex in $Y_{i,j''}^H$ has an $H_\ast^2$-neighbour in $Y_{i,j'}$, which implies that there are at least ${(3\beta^{-1})}/{4}$ good entries $j''\in[\beta^{-1}]\sm\{j,j'\}$ such that there exists a $w\in Y_{i,j''}$ which does not have an $H_\ast^2$-neighbour in $Y_{i,j'}$. 
Since $w_t$ has at most $\alpha^{-2}k$ neighbours in $H_\ast^2$, at most $\alpha^{-2}k$ of these good entries $j''\in[\beta^{-1}]\sm\{j,j'\}$ are such that $w_t$ has an $H_\ast^2$-neighbour in $Y_{i,j''}$. Hence, there are at least $\beta^{-1}/2$ choices to select $j''$ in step $t$.}
As $s= \beta^{2/3}n$, we can ensure that each $j''\in [\beta^{-1}]$ is selected, say, at most $10\beta^{5/3}n$ times.
We write $Z_{i,j}^H:=Z_{i,j}^H(s)$ and it is plain to verify that for all $j\in [\beta^{-1}]$ we have
\begin{enumerate}[label=(\alph*$'$)]
	\item\label{item:correct splitI1} $|Z_{i,j}^H|= \beta|X_i^H|$;\COMMENT{Each label $j\in[\beta^{-1}]$ is used exactly $\beta^{5/3}n$ times. In each such step $t$ where $j$ is used, the cardinality of the current set $Z_{i,j}^H(t-1)$ is increased by exactly one and remains the same during the other steps. Since $|Z_{i,j}^H(0)|=|Y_{i,j}^H\sm W_{i,j}^H|=\beta|X_i^H|-\beta^{5/3}n$, this implies that $|Z_{i,j}^H|=\beta|X_i^H|$.}
	\item\label{item:correct splitI2} $Z_{i,j}^H$ is independent in $H_\ast^2$;\COMMENT{This follows immediately from the construction of the swaps because $W_{i}$ contains all vertices that contain a $H_\ast^2$-neighbour and each $w\in W_i$ is swapped to a set such that $w$ does not have a $H_\ast^2$-neighbour in that set.}
	\item\label{item:correct splitI3}
$\omega(Z_{i,j}^H)=\beta \omega(X_{i}^H) \pm  \beta^{3/2} n/2$ for all $\omega\in \cW$ with $\omega\colon \sX_{i}\to [0,\alpha^{-1}]$.
\end{enumerate}
As $H\in \cH$ is chosen arbitrarily, the statements~\ref{item:correct splitI1}--\ref{item:correct splitI3} hold for all $H\in \cH$.
Note that~\ref{item:correct splitI3} implies~\ref{item:split weights one H, |I|=1}.
\\
It remains to show how to find permutations $\{\pi_i^H\}_{H\in\cH,i\in[r]}$ such that~\ref{item:split weights} also holds for $\omega\in\cW$ with $\omega(\sX_{\sqcup {I}})\geq n^{1+\eps}$ where $X_{i,j}^H:=Z_{i,\pi_i^H(j)}^H$.
This can be easily achieved by considering random permutations.
To see this, fix $\omega\in\cW$ with $\omega\colon \sX_{\sqcup {I}}\to [0,\alpha^{-1}]$, $I=\{i_1,\ldots,i_{|I|}\}\subseteq[r]$, $2\leq|I|\leq k$, and $j_1,\ldots,j_{|I|}\in [\beta^{-1}]$.
Note that we would expect that
$\omega(\bigcup_{H\in \cH}(\bigsqcup_{\ell\in[|I|]}X_{i_\ell,j_\ell}^H))
=
\beta^{|I|} \omega(\sX_{\sqcup I}) \pm \beta^{4/3}n$.
Hence, by Freedman's inequality\COMMENT{Add reference to precise version of Freedman's inequality.} and a union bound we obtain with probability, say, at least~$1/2$, that
$\omega(\bigcup_{H\in \cH}(\bigsqcup_{\ell\in[|I|]}X_{i_\ell,j_\ell}^H))
=
(1\pm  \eps/2)
\beta^{|I|} \omega(\sX_{\sqcup I})$
for all $\omega\in\cW$ with $\omega\colon \sX_{\sqcup {I}}\to [0,\alpha^{-1}]$, $I=\{i_1,\ldots,i_{|I|}\}\subseteq[r]$, $|I|\leq k$, and $j_1,\ldots,j_{|I|}\in [\beta^{-1}]$.
This establishes~\ref{item:split weights}.
\\
In the beginning we made the assumption that $\beta^{-1}$ divides $|X_i^H|$.
To avoid this assumption, we simply remove a set $\widetilde{X}_i^H$ of size at most $\beta^{-1}-1$ from $X_i^H$ such that $\beta^{-1}$ divides $X_i^H\setminus \widetilde{X}_i^H$ and perform the entire procedure with $X_i^H\setminus \widetilde{X}_i^H$ instead of~$X_i^H$.
To that end, for all $H\in \cH$ and $i\in [r]$, let 
$r_i^H\in[\beta^{-1}-1]_0$ be such that $r_i^H=|X_i^H| \mod \beta^{-1}$, and let
$\tX_i^H\subseteq X_i^H$ be such that
\begin{itemize}
\item $|\tX_i^H|=r_i^H$;
\item $\omega(\bigcup_{H\in \cH}\tX_{\sqcup I}^H)
\leq \eps^2\omega(\sX_{\sqcup {I}})$ for all $\omega\in\cW$ with $\omega\colon \sX_{\sqcup I}\to [0,\alpha^{-1}]$.
\end{itemize}
The existence of such sets $\tX_i^H$ can be easily seen by taking each $\tX_i^H$ as a subset of $X_i^H$ of size exactly $r_i^H$ uniformly and independently at random for all $i\in[r]$ and $H\in\cH$.
We now perform the procedure with $X_i^H\sm\tX_i^H$ instead of $X_i^H$.
That is, for all $H\in \cH$ and $i\in [r]$, we obtain a partition $(\tX_{i,j}^H)_{j\in [\beta^{-1}]}$ of $X^H_i\sm \tX^H_i$ satisfying~\ref{item:split1}--\ref{item:split weights}, where $|\tX_{i,j}^H|=|\tX_{i,j'}^H|$ for all $i\in[r],j,j'\in[\beta^{-1}]$, and~\ref{item:split weights one H, |I|=1} and ~\ref{item:split weights} hold with error term `$\pm\beta^{3/2}n/2$' and `$(1\pm\eps/2)$', respectively.
At the very end we add the vertices in $\widetilde{X}_i^H$ to the partition $(\tX_{i,j}^H)_{j\in [\beta^{-1}]}$ while preserving~\ref{item:split1} and~\ref{item:split2}.
We may do so by performing a swap argument as before.
Together with~\eqref{eq:weight only on edges}, observe that the error bounds give us enough room to spare.
\endproof}

\section{Blow-up instances and candidacy graphs}\label{sec:notation blow-up}
In this section we introduce more notation concerning blow-up instances, which will be useful throughout our packing procedure (in Sections~\ref{sec:approx packing} and~\ref{sec:mainpart}).
Let $\sB=(\cH,G,R,\cX,\cV)$ be a blow-up instance of size $(n,k,r)$ that is fixed throughout Section~\ref{sec:notation blow-up}.
Note that the reduced graph $R$ with vertex set $[r]$ gives us a natural ordering of the clusters and we assume this ordering to be fixed (for a different ordering, just relabel the cluster indices).
For simplicity, we often write  $\sX_{i}:=\bigcup_{H\in\cH}X^H_{i}$ for $i\in[r]$ and $\sX_{\sqcup {I}}:=\bigcup_{H\in \cH}(\bigsqcup_{i\in I} X_i^H)$ for $I\subseteq[r]$.
Further, we call $I\subseteq[r]$ an \defn{index set (of~$\sB$)}, if $I\subseteq \rR$ for some $\rR\in E(R)$.\COMMENT{This implies in particular that $|I|\leq k$.}

\medskip
We will introduce some important quantities that we control during our packing procedure.
For instance, we  track for each edge $\gG\in E(G)$ (and for each subset of $\gG$), the set of $\cH$-edges that still could be mapped onto~$\gG$ given a function~$\phi$ that already maps vertices of some clusters in~$\cH$ onto vertices in~$G$ (see Definition~\ref{def:X_g,p,p,phi} in Section~\ref{sec:suitable H-edges}). 
Similarly, we track for distinct edges $\gG,\hG\in E(G)$, the set of tuples of $\cH$-edges that still could be mapped together onto $\gG$ and $\hG$, respectively, with respect to $\phi$ (see Definition~\ref{def:E_g,h,phi}).
To track these quantities, {we define \emph{edge testers} (see Definitions~\ref{def:global edge tester}--\ref{def:general edge tester}) on \emph{candidacy graphs} (see Definition~\ref{def:candidacy graph}) in Sections~\ref{sec:edge testers} and~\ref{sec:candidacy}, respectively.
The definition of these edge testers depends on how the edges in $\cH$ intersect, and to that end, we define \emph{patterns} (see Definitions~\ref{def:pattern}--\ref{def:e_H,p}) in Section~\ref{sec:patterns}. }

\subsection{Candidacy graphs}\label{sec:candidacy}
For the purpose of tracking which sets of vertices in $\cH$ are still suitable images for sets of vertices in~$G$, we will consider auxiliary \defn{candidacy graphs}. 
To that end, assume we are given {$r_\circ\leq r$} and a mapping $\phi\colon \bigcup_{H\in \cH}\hX^H_{\cup \rc}\to V_{\cup \rc}$ with $\hX_q^H\subseteq X_q^H$ for $q\in \rc$ and $\phi|_{V(H)}$ is injective for each $H\in\cH$; that is, $\phi$ already embeds some $\cH$-vertices onto $G$-vertices.
We assume $r_\circ$ and $\phi$ to be fixed throughout the entire Section~\ref{sec:notation blow-up}.
We define candidacy (hyper)graphs with respect to $\phi$ and the blow-up instance $\sB=(\cH,G,R,\cX,\cV)$ such that a (hyper)edge~$\aA$ of the candidacy graph incorporates the property that the set of $\cH$-vertices in~$\aA$ can still be mapped onto the set of $G$-vertices in $\aA$ with respect to potential $\cH$-vertices that are already mapped onto $G$-vertices by~$\phi$.

\begin{defin}[Candidacy graphs $A_I^H(\phi)$]\label{def:candidacy graph}
For all $H\in\cH$ and every index set\COMMENT{Recall, we say $I\subseteq[r]$ is an \defn{index set (for $\sB$)}, if $I\subseteq \rR$ for some $\rR\in E(R)$. Hence, $|I|\leq k$.} $I\subseteq [r]$,\COMMENT{Usually, $I\subseteq[r]\sm \rc$. However, we also allow $I\subseteq[r]$ to be able to also define candidacy graphs $B_i^H(\phi)$ for $i\in \rc$, that is, for already embedded clusters.} let $A_I^H(\phi)$ be the $2|I|$-graph with vertex set $X_{\cup I}^H\cup V_{\cup I}$ and $\bigcup_{i\in I}\{x_i,v_i\}\in E(A_I^H(\phi))$ for $\{x_i,v_i\}\in X_i^H\sqcup V_i$
if all 
$\eH=\eH_\circ\cup \eH_{m}\in E(H[\hX_{\cup (\rc\sm I)}^H,X_{\cup I_m}^H])$ 
with ${m}\in[|I|]$, $I_m\in \binom{I}{m}$, $\eH_\circ\subseteq\binom{\hX_{\cup (\rc\sm I)}^H}{k-{m}}$, and $\eH_{m}=  \{x_i\}_{i\in I_m}$ satisfy 
\begin{align}\label{eq:condition updating}
\phi(\eH_\circ)\cup \{v_i\}_{i\in I_m} \in E(G[V_{\cup (\rc\sm I)},V_{\cup I_m}]).
\end{align}
We call $A_I^H(\phi)$ the \defn{candidacy graph with respect to $\phi$ and $G$.}
\end{defin}

Let us describe Definition~\ref{def:candidacy graph} in words.
Suppose first that we are given an index set $I\subseteq[r]\sm \rc$.
Then the set~$I$ contains the indices of clusters whose vertices are not yet embedded by $\phi\colon \bigcup_{H\in \cH}\hX^H_{\cup \rc}\to V_{\cup \rc}$.
If the set of vertices $\{x_i\}_{i\in I}\in X_{\sqcup I}^H$ can still be mapped onto $\{v_i\}_{i\in I}\in V_{\sqcup I}$, then $\{v_i\}_{i\in I}$ are still suitable candidates for $\{x_i\}_{i\in I}$ and we store this information in the candidacy graph $A_I^H(\phi)$ by adding the edge $\bigcup_{i\in I}\{x_i,v_i\}\in E(A_I^H(\phi))$.
Let us spell out what it means that $\{x_i\}_{i\in I}$ can still be mapped onto $\{v_i\}_{i\in I}$. It means that
all $H$-edges~$\eH$  
\begin{itemize}
\item that intersect $\{x_i\}_{i\in I}$, say in a set of $m$ vertices $\eH_m$ which lie in the clusters with indices $\{i_1,\ldots,i_m\}$,
\item and whose other $k-m$ vertices $\eH_\circ=\eH\sm\eH_{m}$ are embedded by~$\phi$,
\end{itemize}
satisfy~\eqref{eq:condition updating}, that is
\begin{itemize}
\item the embedding $\phi_\circ(\eH_\circ)$ together with the vertices of $\{v_i\}_{i\in I}$ in the clusters with $\{i_1,\ldots,i_m\}$, that is $\phi(\eH_\circ)\cup\{v_{i_1},\ldots,v_{i_{m}}\}$, forms an edge in $G$.
\end{itemize}
\usetikzlibrary{decorations.pathreplacing,shapes.misc,positioning}

We note that we allow $I\subseteq[r]$ in Definition~\ref{def:candidacy graph} (and not only $I\subseteq[r]\sm \rc$) because will use this for $I\subseteq [r]$, {$|I|=1$,} to track a second type of candidacy graphs between already embedded clusters.
At the end of our procedure we will use this second type of candidacy graphs to turn an approximate packing into a complete one.

Let us continue with another comment.
Clearly, for $I=\{i\}$, the candidacy graph $A_i^H(\phi)$ is a bipartite 2-graph.\footnote{For the sake of readability, we write $A_i^H$ instead of $A_{\{i\}}^H$.}
It is worth pointing out a crucial difference between $\bigcup_{i\in I}A_i^H(\phi)$ and $A_I^H(\phi)$: 
If $\{x_iv_i\}_{i\in I}\in \bigsqcup_{i\in I}E(A_i^H(\phi))$, then each vertex $x_i$ on its own can still be mapped onto~$v_i$, whereas if $\bigcup_{i\in I}\{x_i,v_i\}\in E(A_I^H(\phi))$, then the entire set $\{x_i\}_{i\in I}$ can still be mapped onto~$\{v_i\}_{i\in I}$.

Let $\cA_I(\phi):=\bigcup_{H\in \cH}A_I^H(\phi)$ and let $\cA(\phi)$ be the collection of all $\cA_I(\phi)$ for all index sets $I\subseteq [r]\sm \rc$.
We also refer to subgraphs of $A_I^H(\phi)$ as candidacy graphs.

\medskip
To suitably control the candidacy graphs during our approximate packing procedure, it will be important that the neighbourhood $N_{A_i^H}(x)$ in the candidacy graph $A_i^H$ of an $H$-vertex $x$  is the intersection of neighbourhoods of $(k-1)$-sets in $G$.
To that end, for $\eps>0$, $q\in\bN$, $H\in\cH$, $i\in[r]$ and a candidacy graph $A_i^H\subseteq A_i^H(\phi)$, we say
\begin{align}\label{eq:eps,t well-intersecting}
\begin{minipage}[c]{0.85\textwidth}\em
$A_i^H$ is \defn{$(\eps,q)$-well-intersecting with respect to $G$}, if for every $x\in X_i^H$, we can find $\cS_x\subseteq\binom{V(G)}{k-1}$ with $|\cS_x|\leq q$ such that
$N_{A_i^H}(x)=V_i\cap N_G(\cS_x)$,
and every $x\in X_i^H$ is contained in at most $n^{1/4+\eps}$ pairs $\{x,x'\}\in\binom{X_i^H}{2}$ such that $\cS_x\cap \cS_{x'}\neq\emptyset$.
\end{minipage}\ignorespacesafterend
\end{align}
\COMMENT{
The last condition in~\eqref{eq:eps,t well-intersecting} implies that all but at most $n^{5/4+\eps}$ pairs have the appropriate codegree.
}
We note that during our packing procedure the sets $\cS_x$ will be uniquely determined and thus,~\eqref{eq:eps,t well-intersecting} is indeed well-defined.

\subsection{Patterns}\label{sec:patterns}
The behaviour of several parameters in our packing procedure depends on the intersection pattern of the edges in $\cH$, that is, how edges in $\cH$ intersect and overlap.
To this end, we associate two vectors in~$\bN_0^r$ with certain sets of vertices in $H\in\cH$ that we call \defn{\fst-pattern} and \defn{\scd-pattern}.
{Even though we need the precise definitions of these patterns at certain points throughout the paper, it mostly suffices to remember that every set of vertices and every edge in~$\cH$ has a unique \fst- and \scd-pattern. This allows us to track certain quantities with respect to their patterns. 
We proceed to the precise definitions of \fst- and \scd-patterns.}

In order to conveniently define these vectors for a given set of vertices in $H\in\cH$, we consider supergraphs of $H$ (namely, for $Z=B$ in Definition~\ref{def:pattern} below).
We do this because we define candidacy graphs in Section~\ref{sec:candidacy}, and we will in fact consider two collections $\cA$ and~$\cB$ of candidacy graphs (as mentioned in the proof overview), and thus, we also have to distinguish two types of \fst-patterns and \scd-patterns for both $Z\in\{A,B\}$. 
To that end, it is more convenient to imagine that the clusters associated with the candidacy graphs in $\cB$ are copies of the original cluster.
For all $H\in\cH$, $J\subseteq[r]$, and $j\in J$, let $X_j^{H,B}$ be a disjoint copy of $X_j^H$, and let $\pi$ be the bijection that maps a vertex in $X_j^H$ to its copy in $X_j^{H,B}$. 
Let $H_J$ be the supergraph of $H$ with vertex set $V(H_J):=V(H)\cup X_{\cup J}^{H,B}$ and edge set $E(H_J):=E(H)\cup\{\pi(x)\cup(\eH\sm\{x\})\colon \eH\in E(H), x\in \eH\cap X_{\cup J}^H \}$.

We now define \fst-patterns and \scd-patterns and give an illustration in Figure~\ref{fig:pattern}.
{Note that $H_\emptyset-H$ is the empty graph.}

\begin{defin}[Patterns]\label{def:pattern}
For all $Z\in\{A,B\}$, index sets $I\subseteq[r]$, $J\subseteq I$, and all $\xX=\{x_i\}_{i\in I}\in X_{\sqcup I}^H$ for some $H\in\cH$,
let $\xX':=\{x_i\}_{i\in I\sm J}\cup \{\pi(x_j)\}_{j\in J}$, $H_A:=H$ and $H_B:=H_J-H$.
We define the \defn{\fst-pattern} $\bpZ(\xX,J)\in\bN_0^{r}$ and the \defn{\scd-pattern} $\bppZ(\xX,J)\in\bN_0^{r}$ as $r$-tuples where their $\ell$-th entry $\bpZ(\xX,J)_\ell$ and $\bppZ(\xX,J)_\ell$ for $\ell\in[r]$ is given by
\begin{align}
\label{eq:def pattern}
\bpZ(\xX,J)_\ell&:=
\big|
\big\{
\fH\in E\left(H_Z\right)\colon
(\fH\cap X^H_{\ell})\sm\{x_i\}_{i\in I\sm J}\neq\emptyset,
\fH\sm X_{\cup[\ell]}^H\subseteq \xX', |\fH\sm X_{\cup[\ell]}^H|\geq 2
\big\}
\big|;
\\
\label{eq:def prepattern}
\bppZ(\xX,J)_\ell&:=
\big|
\big\{
\fH\in E\left(H_Z\right)\colon 
\fH\cap X^H_{\ell}\in \{x_i\}_{i\in I\sm J},
\fH\sm X_{\cup[\ell]}^H\In \xX', |\fH\sm X_{\cup[\ell]}^H|= 1
\big\}
\big|.
\end{align}
\end{defin}

Let us describe Definition~\ref{def:pattern} in words.
Consider some fixed entry for $\ell\in[r]$.
Depending on $Z\in\{A,B\}$, we either count edges in $H_A$ or in $H_B$; note that the edges in $H_B=H_J-H$ always contain exactly one copied vertex in $X_{\cup J}^{H,B}$.
Further, $\xX\in X_{\sqcup I}^H$ determines the set $I$ and hence the definitions in~\eqref{eq:def pattern} and~\eqref{eq:def prepattern} display no additional dependence on $I$. 

We first describe the \fst-pattern entry $\bpZ(\xX,J)_\ell$ as defined in~\eqref{eq:def pattern}. 
The condition `$(\fH\cap X^H_{\ell})\sm\{x_i\}_{i\in I\sm J}\neq\emptyset$' means that $\fH$ has a non-empty intersection with $X_{\ell}^H$ which does not lie in~$\{x_i\}_{i\in I\sm J}$; in particular $\xX,\xX'\neq\fH$.
The last two conditions in~\eqref{eq:def pattern} mean that all vertices of $\fH\sm X_{\cup[\ell]}^H$ lie in~$\xX'$ and these are at least two vertices. In that sense, $\fH\cap X^H_{\ell}$ is the `last' vertex of $\fH$ not contained in~$\xX'$.\looseness=-1 

We now describe the \scd-pattern entry $\bppZ(\xX,J)_\ell$ as defined in~\eqref{eq:def prepattern}. 
The condition `$\fH\cap X^H_{\ell}\in \{x_i\}_{i\in I\sm J}$' means that $\fH$ has a non-empty intersection with $X_\ell^H$ that lies in $\{x_i\}_{i\in I\sm J}$; note that we may also count edges $\fH\in E(H_Z)$ with $\fH=\xX'$ if $\xX'\in E(H_Z)$.
The last two conditions in~\eqref{eq:def prepattern} mean that $k-1$ vertices of $\fH$ lie in $X_{\cup[\ell]}^H$ and the other vertex of $\fH$ lies in $\xX'$; note that the {last two conditions in~\eqref{eq:def prepattern} imply} that this must be the copied vertex $\fH\cap X_{\cup J}^{H,B}$ if $Z=B$.

\tikzstyle{hvertex}=[thick,circle,inner sep=0.cm, minimum size=2mm, fill=black, draw=black]
\tikzstyle{edge} = [fill,opacity=.3,fill opacity=0,line cap=round, line join=round, line width=14pt]
\begin{figure}[ht!]%
\vspace{-.4cm}
	\begin{center}
		\begin{tikzpicture}[scale = 0.4,every text node part/.style={align=center}]
		\def\ver{0.1} 
		
		\draw (9,4.7) node {$X_{i_1}^H$};		
		\draw (12,4.7) node {$X_{i_2}^H$};
		\draw (15,4.7) node {$X_{i_3}^H$};
		\draw (18,4.7) node {$X_{i_4}^H$};
		\draw (21,4.7) node {$X_{i_5}^H$};		

		\draw[very thick, fill=white] 
		(3,0) ellipse (1 and 4)
		(6,0) ellipse (1 and 4)
		(9,0) ellipse (1 and 4)
		(12,0) ellipse (1 and 4)
		(15,0) ellipse (1 and 4)
		(18,0) ellipse (1 and 4)
		(21,0) ellipse (1 and 4);		

		\coordinate (e0) at (9,0); 
		\coordinate (e1) at (12,0); 
		\coordinate (e2) at (15,0); 
		\coordinate (e3) at (18,0); 
		\coordinate (e4) at (21,0); 	
		\node[hvertex,fill=red,draw=red] (E0) at (e0){};					
		\node[hvertex,fill=red,draw=red] (E1) at (e1){};
		\node[hvertex,fill=red,draw=red] (E2) at (e2){};
		\node[hvertex,fill=red,draw=red] (E3) at (e3){};
		\node[hvertex,fill=red,draw=red] (E4) at (e4){};

		\coordinate (u) at (9,1.7); 
		\coordinate (v1) at (3,3); 
		\coordinate (v2) at (9,3); 
		\coordinate (v0) at (6,3); 
		\coordinate (v3) at (15,3); 
		\node[hvertex] (U) at (u){};				
		\node[hvertex] (V0) at (v0){};				
		\node[hvertex] (V1) at (v1){};
		\node[hvertex] (V2) at (v2){};
		\node[hvertex] (V3) at (v3){};

		\coordinate (a0) at (6,-2); 
		\coordinate (a1) at (9,-2); 
		\coordinate (a2) at (18,-2); 
		\node[hvertex] (A0) at (a0){};
		\node[hvertex] (A1) at (a1){};
		\node[hvertex] (A2) at (a2){};

		\coordinate (b1) at (12,2); 
		\node[hvertex] (B1) at (b1){};


		\coordinate (d1) at (3,0); 
		\coordinate (d2) at (6,0); 
		\node[hvertex] (D1) at (d1){};
		\node[hvertex] (D2) at (d2){};

		\coordinate (x) at (15,-2); 
		\node[hvertex] (X) at (x){};

\draw[edge,color=blue] (v1) -- (v3) -- (e4);
\draw[edge] (a0)--(a1) -- (e2) -- (a2) -- (e4);

\draw[edge] (v0)--(v2) -- (b1) -- (e3) -- (e4);

\draw[edge] (d1) -- (e2);	
\draw (4.5,0) node {$\fH$};			

\draw[edge] (u)--(e1) -- (x) -- (e3) -- (e4);		

		\node[hvertex,fill=red,draw=red] (E0) at (e0){};	
		\node[hvertex,fill=red,draw=red] (E1) at (12,0){};
		\node[hvertex,fill=red,draw=red] (E2) at (15,0){};
		\node[hvertex,fill=red,draw=red] (E3) at (18,0){};
		\node[hvertex,fill=red,draw=red] (E4) at (21,0){};

\node at (27,-5){};
\node at (-.2,-5){
$\bpA(\xX,\emptyset) =\big($};
\node at (22,-5){
$\big)$};

\node at (-.9,-6.5){
$\bppA(\xX,\emptyset) =\big($};
\node at (22,-6.5){
$\big)$};		

\def\y{-5}
\node at (3,\y){$0$};
\node at (6,\y){$1$};
\node at (9,\y){$0$};
\node at (12,\y){$1$};
\node at (15,\y){$1$};
\node at (18,\y){$0$};
\node at (21,\y){$0$};

\def\y{-5.2}
\node at (4.5,\y){$,$};
\node at (7.5,\y){$,$};
\node at (10.5,\y){$,$};
\node at (13.5,\y){$,$};
\node at (16.5,\y){$,$};
\node at (19.5,\y){$,$};

\def\y{-6.5}
\node at (3,\y){$0$};
\node at (6,\y){$0$};
\node at (9,\y){$0$};
\node at (12,\y){$1$};
\node at (15,\y){$0$};
\node at (18,\y){$2$};
\node at (21,\y){$0$};

\def\y{-6.7}
\node at (4.5,\y){$,$};
\node at (7.5,\y){$,$};
\node at (10.5,\y){$,$};
\node at (13.5,\y){$,$};
\node at (16.5,\y){$,$};
\node at (19.5,\y){$,$};

%
%
%
%
%

\node at (-1.4,-8){
$\bpB(\xX,\{i_1,i_2\}) =\big($};
\node at (22,-8){
$\big)$};

\node at (-2.1,-9.5){
$\bppB(\xX,\{i_1,i_2\}) =\big($};
\node at (22,-9.5){
$\big)$};		

\def\y{-8}
\node at (3,\y){$0$};
\node at (6,\y){$0$};
\node at (9,\y){\boxed{1}};
\node at (12,\y){\textcircled{\raisebox{-0.9pt}{1}}};
\node at (15,\y){$1$};
\node at (18,\y){$0$};
\node at (21,\y){$0$};

\def\y{-8.2}
\node at (4.5,\y){$,$};
\node at (7.5,\y){$,$};
\node at (10.5,\y){$,$};
\node at (13.5,\y){$,$};
\node at (16.5,\y){$,$};
\node at (19.5,\y){$,$};

\def\y{-9.5}
\node at (3,\y){$0$};
\node at (6,\y){$0$};
\node at (9,\y){$0$};
\node at (12,\y){$0$};
\node at (15,\y){\boxed{2}};
\node at (18,\y){$0$};
\node at (21,\y){$1$};

\def\y{-9.7}
\node at (4.5,\y){$,$};
\node at (7.5,\y){$,$};
\node at (10.5,\y){$,$};
\node at (13.5,\y){$,$};
\node at (16.5,\y){$,$};
\node at (19.5,\y){$,$};

\node at (-2.2,-11){
$\bpB(\xX,\{i_1,\ldots,i_4\}) =\big($};
\node at (22,-11){
$\big)$};

\node at (-2.9,-12.5){
$\bppB(\xX,\{i_1,\ldots,i_4\}) =\big($};
\node at (22,-12.5){
$\big)$};		

\def\y{-11}
\node at (3,\y){$0$};
\node at (6,\y){$0$};
\node at (9,\y){$0$};
\node at (12,\y){$1$};
\node at (15,\y){$1$};
\node at (18,\y){$2$};
\node at (21,\y){$0$};

\def\y{-11.2}
\node at (4.5,\y){$,$};
\node at (7.5,\y){$,$};
\node at (10.5,\y){$,$};
\node at (13.5,\y){$,$};
\node at (16.5,\y){$,$};
\node at (19.5,\y){$,$};

\def\y{-12.5}
\node at (3,\y){$0$};
\node at (6,\y){$0$};
\node at (9,\y){$0$};
\node at (12,\y){$0$};
\node at (15,\y){$0$};
\node at (18,\y){$0$};
\node at (21,\y){$4$};

\def\y{-12.7}
\node at (4.5,\y){$,$};
\node at (7.5,\y){$,$};
\node at (10.5,\y){$,$};
\node at (13.5,\y){$,$};
\node at (16.5,\y){$,$};
\node at (19.5,\y){$,$};

		\end{tikzpicture}
	\end{center}
	  \captionsetup{width=.97\linewidth}
\vspace{-.4cm}
	\caption{We have $k=5$, $I=\{i_1,\ldots,i_5\}$ in this figure and the set $\xX=\{x_{i_1},\ldots,x_{i_5} \}\in X_{\sqcup I}^H$, $\xX\notin E(H)$, consists of the five red vertices. Note that the blue hyperedge does not play a role for any of the displayed \fst-patterns or~\scd-patterns.
	{Let us explain the marked entries; all others can be checked similarly. 
	For $\bpB(\xX,\{i_1,i_2\})$ we consider~\eqref{eq:def pattern} for $Z=B$ and the displayed edge $\fH\in E(H)$.
	For $J=\{i_1,i_2\}$ we will obtain two copied edges $\fH_{i_1}:=(\fH\sm\{x_{i_1}\})\cup\{x_{i_1}'\}$ and $\fH_{i_2}:=(\fH\sm\{x_{i_2}\})\cup\{x_{i_2}'\}$ of $\fH$ in $E(H_B)=E(H_J)\sm E(H)$.
	By checking the conditions in~\eqref{eq:def pattern}, we note that $\fH_{i_2}$ accounts for the marked entry \boxed{1} and $\fH_{i_1}$ accounts for \textcircled{\raisebox{-0.9pt}{1}} of $\bpB(\xX,\{i_1,i_2\})$.
	Similarly, by checking the conditions in~\eqref{eq:def prepattern}, both edges $\fH_{i_1}$ and $\fH_{i_2}$ account for the marked entry \boxed{2} of $\bppB(\xX,\{i_1,i_2\})$.}
	}
	\label{fig:pattern}
\end{figure}
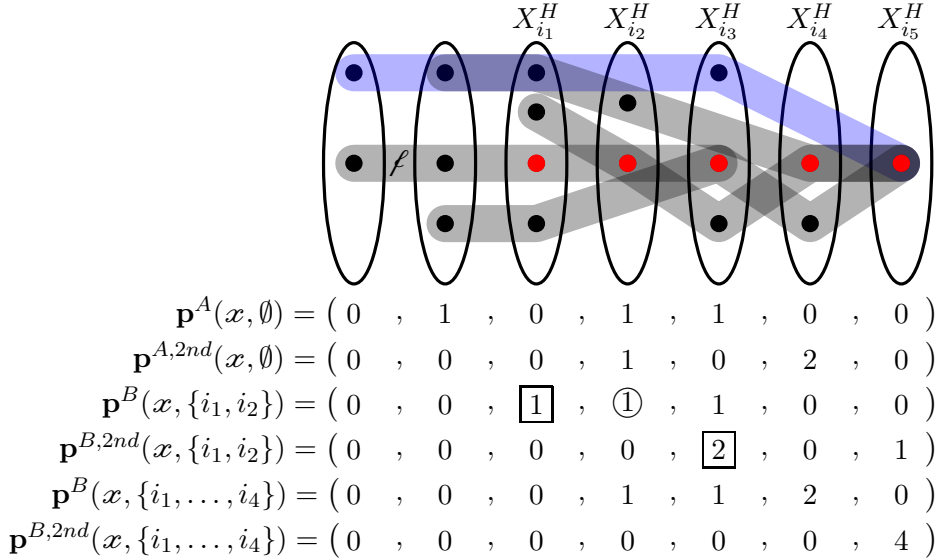

We make the following important observation concerning Definition~\ref{def:pattern}.
We claim that
\begin{align}\label{eq:norms equal}
\norm{\bpZ(\xX,J)}=\norm{\bppZ(\xX,J)}-\IND{\{\xX'\in E(H_Z) \}}.
\end{align}
To see that this is true, let us first assume that $\xX'\notin E(H_Z)$. Note that every $\fH$ that contributes to $\norm{\bpZ(\xX,J)}$,  has a `penultimate' vertex that lies in $\{x_i\}_{i\in I\sm J}$, {that is, there exists an index~$\ell$ such that the conditions in~\eqref{eq:def prepattern} are satisfied}, and thus~$\fH$ also contributes to $\norm{\bppZ(\xX,J)}$. Conversely, every~$\fH$ that contributes to $\norm{\bppZ(\xX,J)}$ has a `last' vertex not contained in $\{x_i\}_{i\in I\sm J}$ if $\xX'\notin E(H_Z)$, and thus $\fH$ also contributes to $\norm{\bpZ(\xX,J)}$. 
Hence, $\norm{\bpZ(\xX,J)}=\norm{\bppZ(\xX,J)}$ if $\xX'\notin E(H_Z)$. 
If $\xX'\in E(H_Z)$ (and thus $|I|=k$), then $\fH=\xX'$ additionally contributes to $\norm{\bppZ(\xX,J)}$ but not to $\norm{\bpZ(\xX,J)}$ because of the condition `$(\fH\cap X^H_{\ell})\sm\{x_i\}_{i\in I\sm J}\neq\emptyset$' in~\eqref{eq:def pattern}.
This implies~\eqref{eq:norms equal}.
{Further, note that $\xX'\in E(H_A)$ only if $J=\emptyset$, and $\xX'\in E(H_B)$ only if $|J|=1$.}

\medskip
We will also consider a set of vertex tuples that lie in an edge in~$\cH$ and have certain patterns.

\begin{defin}[$E_{\cH}(\pP,I,J)$]\label{def:e_H,p}
For all index sets $I\subseteq[r]$, and $\bp,\bp^{2nd}\in\bN_0^r$, let
\begin{align*}
E_{\cH}(\bp,\bp^{2nd},I):=\big\{
\xX\in \sX_{\sqcup I}\colon& \bpA(\xX,\emptyset)=\bp,\bppA(\xX,\emptyset)=\bp^{2nd}, \xX=\eH\cap\sX_{\sqcup I} \text{ for some $\eH\in E(\cH)$}
\big\}.
\end{align*}
\COMMENT{Note that $E_{\cH}(\bp,\bp^{2nd},I)\neq\emptyset$ only if $\norm{\bpp}-\norm{\bp}\in\{0,1\}$ by~\eqref{eq:norms equal}.}\COMMENT{Note that there might be multiple edges $\eH$ (in different reduced edges) that contain $\xX$, e.g. when $\xX$ is only a single vertex.}More generally, we allow to specify whether some vertices of $\xX\in \sX_{\sqcup I}$ lie in clusters with indices in $J\subseteq I$.
To this end, for all index sets $I\subseteq[r]$, all $J\subseteq I$, and $\pP=(\bpA,\bppA,\bpB,\bppB)\in\bN_0^r\times\bN_0^r\times\bN_0^r\times\bN_0^r=(\bN_0^r)^4$, let
\begin{align*}
E_{\cH}(\pP,I,J):=\big\{
\xX\in \sX_{\sqcup I}\colon& 
\bpZ(\xX,J)=\bpZ, \bppZ(\xX,J)=\bppZ \text{ for both $Z\in\{A,B\}$,}
\\&
 \xX=\eH\cap\sX_{\sqcup I} \text{ for some $\eH\in E(\cH)$}
\big\}.
\end{align*}
\end{defin}

\subsection{Edge testers}\label{sec:edge testers}
Recall that we consider some fixed $r_\circ\leq r$ and $\phi\colon \bigcup_{H\in \cH}\hX^H_{\cup \rc}\to V_{\cup \rc}$.
Let~$\cA$ be a collection of candidacy graphs $\cA_I=\bigcup_{H\in \cH} A_I^H\subseteq \cA_I(\phi)$ for all index sets $I\subseteq [r]\sm \rc$.\COMMENT{$\cA_I$ will be the candidacy graphs with respect to $\cH_+$.}
We will use weight functions on the edges of the candidacy graphs, which we also call edge testers, in order to track important quantities during our packing procedure.

We start with the definition of \emph{simple edge testers} (Definition~\ref{def:global edge tester}).
In Definition~\ref{def:general edge tester} we introduce more complex edge testers that include simple edge testers.
However, because we will frequently use weight functions on the candidacy graphs in form of simple edge testers, we include both definitions for the readers' convenience.

\medskip
To that end, given an index set $I\subseteq [r]$, a \fst-pattern vector $\bp$, a \scd-pattern vector~$\bpp$, vertices $\bc\in V_{\sqcup I}$ that we call \emph{centres}, and an \emph{(initial) weight function}
$\omega_\iota\colon\sX_{\sqcup {I}}\to\bR_{\ge 0}$ with $\supp(\omega_\iota)\subseteq E_\cH(\bp,\bpp,I)$, we define a \emph{simple edge tester $(\omega,\omega_\iota,\bc,\bp,\bpp)$ with respect to $\phi$ and the candidacy graphs in $\cA$} in the following Definition~\ref{def:global edge tester}. 
Ultimately, our aim is to track the $\omega_\iota$-weight of tuples in~$\sX_{\sqcup {I}}$ that are mapped onto the centres $\bc$.
We can think of $\omega\colon E(\cA_{I\sm \rc})\to \bR_{\ge 0}$ as an updated weight function that restricts the $\omega_\iota$-weight that still can be mapped onto the centres~$\bc$ with respect to~$\phi\colon \bigcup_{H\in \cH}\hX^H_{\cup \rc}\to V_{\cup \rc}$ and the candidacy graphs in~$\cA$.
Since $\supp(\omega_\iota)\subseteq E_\cH(\bp,\bpp,I)$,
we specify the \fst-pattern and \scd-pattern of the tuples in $\supp(\omega_\iota)$ and thus we know exactly how those tuples intersect with edges in $\cH$.
This will allow us to precisely control the weight of an edge tester during our packing procedure. 

\begin{defin}[Simple edge tester $(\omega,\omega_\iota,\bc,\bp,\bpp)$]\label{def:global edge tester}
For an index set $I\subseteq[r]$, $\omega_\iota\colon\sX_{\sqcup {I}}\to[0,s]$ for some $s\in\bR_{>0}$ with $\supp(\omega_\iota)\subseteq E_\cH(\bp,\bpp,I)$ for given $\bp, \bpp\in\bN_0^r$, and $\bc=\{c_i\}_{i\in I}\in V_{\sqcup I}$, let $\omega\colon E(\cA_{I\sm \rc})\to[0,s]$ be defined by
\begin{align}\label{eq:global edge tester}
\omega(\aA):=
\IND{\{\aA\cap V_{I\sm \rc}=\{c_i\}_{i\in I\sm \rc} \}}
\cdot\omega_\iota(\xX)
\end{align} 
for all $\aA\in E(A^H_{I\sm \rc})$ and $H\in\cH$ where $\xX=(\phi|_{V(H)})^{-1}(\{c_i\}_{i\in I\cap \rc})\cup (\aA\cap X^H_{\cup (I\sm \rc)})\in X_{\sqcup I}^H$.
If no such $\xX$ exists, we set $\omega(\aA):=0$.
We say $(\omega,\omega_\iota,\bc,\bp,\bpp)$ is a \emph{simple $s$-edge tester with respect to $(\omega_\iota,\bc,\bp,\bpp)$, $\phi$~and~$\cA$}.
\end{defin}

For the readers' convenience, let us discuss Definition~\ref{def:global edge tester} in detail.
Suppose we are given an (initial) weight function $\omega_\iota\colon\sX_{\sqcup {I}}\to[0,s]$ with centres $\bc=\{c_i\}_{i\in I}\in V_{\sqcup I}$ for an index set $I\subseteq[r]$ and $\supp(\omega_\iota)\subseteq E_\cH(\bp,\bpp,I)$ with \fst-pattern~$\bp$ and \scd-pattern $\bpp$. 
Recall that our aim is to control the $\omega_\iota$-weight of tuples in $\sX_{\sqcup {I}}$ that are mapped onto the centres $\bc$.
Therefore, for an edge $\aA\in E(\cA_{I\sm \rc})$, we put the weight $\omega_\iota(\xX)$ onto $\aA$ only if the following are satisfied: 
\begin{itemize}
\item $\aA$ contains the centres $\{c_i\}_{i\in I\sm \rc}$ of the not yet embedded clusters (which is incorporated by the indicator function in~\eqref{eq:global edge tester}), and
\item $\xX$ is such that the vertices of $\xX$ in $X^H_{\cup (I\sm \rc)}$ are contained in $\aA$ and $\phi$ maps the vertices in $\xX\cap X^H_{\cup (I\cap \rc)}$ onto $\{c_i\}_{i\in I\cap \rc}$.
\end{itemize}

\medskip
Let us now comment on the purpose of more complex edge testers as defined in Definition~\ref{def:general edge tester}.
Our partial packing procedure will only provide a packing $\phi\colon \bigcup_{H\in \cH}\hX^H_{\cup \rc}\to V_{\cup \rc}$ for $\rc\subseteq [r]$ that maps almost all vertices in $\bigcup_{H\in \cH}X^H_{\cup \rc}$ onto vertices in $V_{\cup \rc}$ and leaves the vertices $\bigcup_{H\in \cH}(X^H_{\cup \rc}\sm \hX^H_{\cup \rc})$ unembedded.
We will often call the vertices $\bigcup_{H\in \cH}(X^H_{\cup \rc}\sm \hX^H_{\cup \rc})$ \defn{unembedded by $\phi$} or simply the \defn{leftover} of the partial packing $\phi$. 
In the end, we will have to turn such a partial packing into a complete one. 
Therefore, we will utilize a second collection of candidacy graphs $\cB$ during the partial packing that tracks candidates that correspond to edges in $G$ that we reserved in the beginning for the completion step.\footnote{In fact, we partition the edge set of the host graph $G$ into two $k$-graphs $G_A$ and $G_B$, and $\cA$ will be a collection of candidacy graphs with respect to $\phi$ and $G_A$, and $\cB$ will be a collection of candidacy graphs with respect to $\phi$ and~$G_B$.}
In order for this to work, we have to take care that the leftover is well-behaved with respect to the candidacy graphs in $\cB$.
We achieve this by using weight functions on $2$-tuples consisting of one hyperedge of an $\cA$-candidacy graph and of a collection of edges within the $\cB$-candidacy graphs. 
That is, assume we are initially given a weight function  $\omega_\iota\colon\sX_{\sqcup I}\to[0,s]$ (that we therefore often call initial weight function) with centres $\bc=\{c_i\}_{i\in I}\in V_{\sqcup I}$ for an index set $I\subseteq[r]$ and recall that our overall aim is to track the $\omega_\iota$-weight of tuples in $\sX_{\sqcup {I}}$ that can be mapped onto the centres in $\bc$.
Now, in Definition~\ref{def:general edge tester} of our \emph{general edge testers} we allow to specify a set $J\subseteq I$ of indices where we track the $\omega_\iota$-weight of tuples $\xX=\{x_i\}_{i\in I}\in X_{\sqcup I}^H$ for all $H\in\cH$ such that {for each $j\in J$, the vertex $x_j$ can be mapped onto $c_j$ within the candidacy graph $B_j^H\in\cB$.}
That is, if the vertices $\{x_j\}_{j\in J}$ are left unembedded, then they can potentially still be mapped onto $\{c_j\}_{j\in J}$ during the completion process using the candidacy graphs $\cB$.
Further, we even allow to specify disjoint subsets $\JX$ and $\JV$ of $J$, where $\JX$ encodes that exactly the vertices $\{x_j\}_{j\in \JX}$ of $\xX$ are left unembedded, and $\JV$ encodes whether the tuple $\xX\in X_{\sqcup I}^H$ lies in a graph $H$ such that $\phi|_{V(H)}$ leaves the centres $\{c_j\}_{j\in \JV}$ uncovered.

Assume $\cB$ is a fixed collection of candidacy graphs $\cB_j=\bigcup_{H\in \cH} B_j^H\subseteq \cB_j(\phi)$ for all $j\in [r]$, defined as in Definition~\ref{def:candidacy graph}.
{To make our partial packing procedure more uniform, we will sometimes also treat vertices that are left unembedded by $\phi$ as embedded by some extension $\phi^+\colon \bigcup_{H\in \cH}X^H_{\cup \rc}\to V_{\cup \rc}$ of~$\phi$ (which only serves as a dummy extension and is not necessarily a packing).}

\begin{defin}[(General) edge tester $(\omega,\omega_\iota,J,\JX, \JV,\bc,\pP)$]\label{def:general edge tester}
For an index set $I\subseteq[r]$, $J\subseteq I$, disjoint sets $\JX, \JV\subseteq J$, $\omega_\iota\colon\sX_{\sqcup {I}}\to[0,s]$ for some $s\in\bR_{>0}$ with $\supp(\omega_\iota)\subseteq E_\cH(\pP,I,J)$ for given $\pP\in(\bN_0^r)^4$, and $\bc=\{c_i\}_{i\in I}\in V_{\sqcup I}$, let {$\omega\colon 
\bigcup_{H\in \cH}\big( 
E(A^H_{I\sm (\rc\cup J)})\sqcup \big(\bigsqcup_{j\in J}E(B^H_j)\big)
\big)
\to[0,s]$} be defined by
\begin{align}\label{eq:general edge tester}
\omega(\{\aA,\bA \}):=
\IND{\big\{(\aA\cup b_{\cup J})\cap V_{\cup((I\sm \rc)\cup J)}=\{c_i\}_{i\in (I\sm \rc)\cup J} \big\}}
\sum_{\xX \text{ is $\{\aA,\bA \}$-suitable}}
\omega_\iota(\xX)
\end{align} 
for all $\aA\in E(A^H_{I\sm (\rc\cup J)})$, {$\bA=\{b_j\}_{j\in J}\in \bigsqcup_{j\in J}E(B_j^H)$} and $H\in\cH$, where we say that $\xX\in X_{\sqcup I}^H$ is \defn{$\{\aA,\bA \}$-suitable} if
\begin{enumerate}[leftmargin=2cm,label={\rm (\roman*)\textsubscript{D\ref{def:general edge tester}}}]
\item\label{item:contained in ab} $(\aA\cup b_{\cup J} )\cap X_{\cup((I\sm \rc)\cup J)}^H=\xX\cap X_{\cup((I\sm \rc)\cup J)}^H$;

\item $\{c_i\}_{i\in (I\cap \rc)\sm J}\subseteq \phi(\xX\cap \hX_{\cup \rc}^H)
$;

\item\label{item:c-vertices left unembedded} $c_j\notin \phi(\hX_j^H)$ for all $j\in \JV \cap \rc$; \COMMENT{Note that this depends on $H$ for $\xX\in X_{\sqcup I}^H$.}

\item\label{item:x-vertices left unembedded} $\xX\cap(X_{\cup \rc}^H\sm \hX_{\cup \rc}^H)=\xX\cap X_{\cup(\JX\cap \rc)}^H$;

\item\label{item:B not on centers}
$\phi^+(
\xX\cap X_{\cup (J\cap[r_\circ])}^H
)
\cap\bc=\emptyset$.
\end{enumerate}
Note that for each $\{\aA,\bA\}$, there is at most one $\{\aA,\bA \}$-suitable tuple $\xX$.
If no $\{\aA,\bA \}$-suitable tuple~$\xX$ exists, we set $\omega(\{\aA,\bA \}):=0$.
We say $(\omega,\omega_\iota,J,\JX, \JV,\bc,\pP)$ is an \emph{$s$-edge tester with respect to $(\omega_\iota,J,\JX, \JV,\bc,\pP)$, {$(\phi,\phi^+)$}, $\cA$ and~$\cB$}.
\end{defin}

{We often write $\JXV$ for $\JX\cup\JV$ if $\JX$ and $\JV$ are fixed.}
Let us comment on Definition~\ref{def:general edge tester}.
Suppose we are given an initial weight function $\omega_\iota\colon\sX_{\sqcup {I}}\to[0,s]$ with centres $\bc=\{c_i\}_{i\in I}\in V_{\sqcup I}$ for an index set $I\subseteq[r]$ and $J\subseteq I$, and $\supp(\omega_\iota)\subseteq E_\cH(\pP,I,J)$.  
As in the case for simple edge testers, our aim is to control the $\omega_\iota$-weight of tuples $\xX$ in $\sX_{\sqcup {I}}$ that are mapped onto the centres~$\bc$.
The set $J\subseteq I$ allows us to specify some vertices $\xX\cap X_{\cup J}^H$ of such tuples $\xX$ that {are not yet embedded by $\phi^+$ onto their centers $\{c_i\}_{i\in J}$ and which} will potentially be embedded onto those during the completion process.
For the completion, we will use the candidacy graphs {$B_j^H$ in $\cB$ and therefore,  $(\xX\cap X_{j}^H)\cup \{c_j\}=b_j\in E(B_j^H)$ for each $j\in J$.}
Furthermore, the sets $\JX, \JV\subseteq J$ encode the situations that\looseness=-1
\begin{itemize}
\item only the vertices $\xX\cap X_{\cup(\JX\cap \rc)}^H$ of $\xX$ are left unembedded by $\phi$ (see \ref{item:x-vertices left unembedded}), or
\item the centres $\{c_j\}_{j\in \JV \cap \rc}$ are uncovered by $\phi|_{V(H)}$ (see \ref{item:c-vertices left unembedded}).
\end{itemize}
\COMMENT{The vertices $\{x_i\}_{i\in J\sm(\JX\cup \JX)}$ of a tuple $\xX=\{x_i\}_{i\in I}$ are vertices that are not left unembedded, and thus, in order to use them during the completion, they have to become activated.}\COMMENT{Figure~\ref{fig:leftover} might be a helpful illustration.}

Note that a (general) edge tester $(\omega,\omega_\iota,J=\emptyset,\JX=\emptyset,\JV=\emptyset,\bc,(\bp,\bpp,\mathbf{0},\mathbf{0}))$ is equivalent to a simple edge tester $(\omega,\omega_\iota,\bc,\bp,\bpp)$ with respect to $(\omega_\iota,\bc,\bp,\bpp)$, $\phi$~and~$\cA$.

\subsection{Sets of suitable $\cH$-edges}\label{sec:suitable H-edges}

Next, we define (sub)sets of $\cH$-edges that we track during our packing procedure.
Recall that we consider some fixed $r_\circ\leq r$ and $\phi\colon \bigcup_{H\in \cH}\hX^H_{\cup \rc}\to V_{\cup \rc}$.
In Definition~\ref{def:X_g,p,p,phi}, we define a set $\cX_{\gG,\bp,\bp^{2nd},\phi}(\cA)$ for every edge $\gG\in E(G)$ that contains the sets of vertices that are contained in an $\cH$-edge with \fst-pattern $\bp$ and \scd-pattern $\bp^{2nd}$, and that still could be mapped together onto $\gG$ with respect to~$\phi$ and the candidacy graphs in~$\cA$.
We can track the size of this set $\cX_{\gG,\bp,\bp^{2nd},\phi}(\cA)$ by using simple edge testers.

\begin{defin}[$\cX_{\gG,\bp,\bp^{2nd},\phi}(\cA)$]\label{def:X_g,p,p,phi}
For all $\gG=\gG_\circ\cup \gG_m\in E(G[V_{\cup\rR}])$ for some $\rR\in E(R)$ with $\gG_\circ\in \binom{V_{\cup \rc}}{k-m}$, $m\in[k]$, and with 
$I:=\rR\sm \rc$, $|I|=m$, $\gG_m\in V_{\sqcup I}$, and for all $\bp,\bp^{2nd}\in\bN_0^r$, 
let 
\begin{align}
\nonumber
\cX_{\gG,\bp,\bp^{2nd},\phi}(\cA)
:=\bigcup_{H\in \cH} \Big\{&
\xX_m\in X_{\sqcup I}^H\colon 
\\&\label{def:X_g,p,p,phi I}
\xX_{m}\cup\gG_m\in E(A_I^H),
\\&\label{def:X_g,p,p,phi II}
(\phi|_{V(H)})^{-1}(\gG_\circ)\cup\xX_m\in E_\cH(\bp,\bp^{2nd},\rR)
\Big\}.
\end{align}

Further, for $\omega_\iota\colon \sX_{\sqcup\rR}\to\{0,1\}$ with $\omega_\iota(\xX):=\IND{\{\xX\in E_\cH(\bp,\bpp,\rR) \} }$, we call the simple $1$-edge tester $(\omega,\omega_\iota,\gG,\bp,\bpp)$ with respect to $(\omega_\iota,\gG,\bp,\bpp)$, $\phi$ and $\cA$ (as defined in Definition~\ref{def:global edge tester}), the \defn{edge tester for $\cX_{\gG,\bp,\bp^{2nd},\phi}(\cA)$.}
\end{defin}

Let us describe Definition~\ref{def:X_g,p,p,phi} in words.
Suppose we are given an edge $\gG=\gG_\circ\cup \gG_m\in E(G[V_{\cup\rR}])$, where $\gG_\circ$ contains the vertices of $\gG$ that lie in clusters that are already embedded by $\phi$ and $\gG_m$ contains the remaining $m$ vertices of $\gG$ in the not yet embedded clusters with indices $I=\rR\sm \rc$.
For each $H\in\cH$, we track the set of vertices $\xX_{m}\in X_{\sqcup I}^H$ where 
\begin{itemize}
\item $\xX_{m}$ still could be mapped onto $\gG_m$ (that is, $\xX_{m}\cup\gG_m\in E(A_I^H)$ in~\eqref{def:X_g,p,p,phi I}), and
\item $\xX_{m}$ lies in an $H$-edge $\eH$ with \fst-pattern $\bp$ and \scd-pattern $\bpp$ such that if we map~$\xX_{m}$ onto~$\gG_m$, then $\eH$ is mapped onto~$\gG$ (that is, $(\phi|_{V(H)})^{-1}(\gG_\circ)\cup\xX_m\in E_\cH(\bp,\bp^{2nd},\rR)$ in~\eqref{def:X_g,p,p,phi II}).
\end{itemize}

Further, note that for a simple edge tester $(\omega,\omega_\iota,\gG,\bp,\bpp)$  for $\cX_{\gG,\bp,\bp^{2nd},\phi}(\cA)$ as in Definition~\ref{def:X_g,p,p,phi}, we have that $\omega(E(\cA_{I}))=|\cX_{\gG,\bp,\bp^{2nd},\phi}(\cA)|$ for $I=\rR\sm \rc$ because the indicator function in~\eqref{eq:global edge tester} corresponds to~\eqref{def:X_g,p,p,phi I}, and the choice of $\xX$ in~\eqref{eq:global edge tester} corresponds to~\eqref{def:X_g,p,p,phi II} by the definition of $\omega_\iota$.

\medskip
Next, we define in Definition~\ref{def:E_g,h,phi} a set $E_{\gG,\hG,\phi}(\cA)$ for all distinct $G$-edges $\gG,\hG$ with identical $G$-vertex in the last cluster such that $E_{\gG,\hG,\phi}(\cA)$ contains the tuples of $\cH$-edges $(\eH,\fH)$ with identical $\cH$-vertex in the last cluster, and $\eH$ and $\fH$ can still be mapped onto $\gG$ and $\hG$ with respect to $\phi$ and the candidacy graphs in~$\cA$. 
In this case, we ignore the patterns as we only aim for an upper bound on the number of these edges and have some room to spare.

\begin{defin}[$E_{\gG,\hG,\phi}(\cA)$]\label{def:E_g,h,phi}
For edges $\gG=\{v_{i_1},\ldots,v_{i_k} \}, \hG=\{w_{j_1},\ldots,w_{j_k} \}\in E(G)$ with $v_{i_k}=w_{j_k}$ and both $\zZ\in\{\gG,\hG\}$, let
\begin{align*}
E_{\zZ,\phi}(\cA) &:= E\left(\cH\left[
\phi^{-1}(\zZ\cap V_{\cup \rc})\cup 
\textstyle\bigcup_{i\in[r]\sm \rc}N_{\cA_i}(\zZ\cap V_i)
\right]\right);
\\E_{\gG,\hG,\phi}(\cA) &:=\left\{
(\eH,\fH) \in E_{\gG,\phi}(\cA)\times E_{\hG,\phi}(\cA) \colon \eH\cap \fH\cap \sX_{i_k}\neq \emptyset
\right\}.
\end{align*}
\end{defin}
\COMMENT{The last condition implies in particular that $i_k=j_k$ and $\eH$ and $\fH$ live in the same graph $H\in\cH$. Note that $(\eH,\fH) \in E_{\gG,\hG,\phi}(\cA)$ implies that $\eH\in E_{\gG,\phi}(\cA), \fH\in E_{\hG,\phi}(\cA)$.
}

\section{Approximate Packing Lemma}\label{sec:approx packing}

In this section we provide our `Approximate Packing Lemma' (Lemma~\ref{lem:packing lemma}).
Given a blow-up instance $(\cH,G,R,\cX,\cV)$, it allows us for one cluster to embed almost all vertices of $\bigcup_{H\in\cH}X_i^H$ into~$V_i$, while maintaining crucial properties for future embedding rounds of other clusters.
To describe this setup we define a \emph{packing instance} and collect some more notation.

\subsection{Packing instances}\label{sec:packing instance}
Our general understanding of a packing instance is as follows. 
Recall that we will consider the clusters of a blow-up instance one after another.
A packing instance arises from a blow-up instance where we have already embedded vertices of some clusters (which is given by a function $\phi_\circ$)
and focuses only on one particular cluster (denoted by $\bigcup_{H\in \cH}X_0^H$ and $V_0$) 
and all clusters that are close to the considered cluster (measured in the reduced graph $R$).
We track \emph{candidacy graphs} as defined in Definition~\ref{def:candidacy graph} and consider a collection of candidacy graphs~$\cA$ between the clusters in~$\cH$ and~$G$ that will be used for future embedding rounds.
In order to be able to turn a partial packing into a complete one in the end, we do not only track the collection of candidacy graphs in~$\cA$ but also a second collection of candidacy graphs~$\cB$, where the candidacy graphs in $\cB$ will be used for the completion step.
To that end, we also assume that the edges of~$G$ are partitioned into two parts $G_A$ and $G_B$ such that the edges in $G_A$ are used for the approximate packing and the edges in $G_B$ are reserved for the completion step. 
That is, we will think of the graphs in $\cA$ as candidacy graphs with respect to $\phi_\circ$ and $G_A$, and of the graphs in $\cB$ as candidacy graphs with respect to $\phi_\circ$ and $G_B$.

\medskip
We make this more precise.
Let $n,k,r,r_\circ\in \bN_0$.
We say $\sP=(\cH,G_A,G_B,R,\cA,\cB,\phi_\circ^-,\phi_\circ)$ is a \defn{packing instance of size $(n,k,r,r_\circ)$} if
\begin{itemize}
	\item $\cH$ is a collection of $k$-graphs, $G_A$ and $G_B$ are edge-disjoint $k$-graphs on the same vertex set, and $R$ is a $k$-graph where $V(R)=-[r_\circ]\cup[r]_0$;
	
	\item $\{X^H_i\}_{i\in \VR}$ is a vertex partition of $H\in\cH$ such that $|\eH\cap X_i^H|\leq 1$ for all $\eH\in E(H)$, $H\in\cH$;
	\item $\{V_i\}_{i\in \VR}$ is a vertex partition of $G_A$ as well as $G_B$;
	\item  $|X_i^H|=|V_i|=(1\pm1/2)n$ for each $i\in\VR$;
	\item for all $H\in\cH$, the hypergraph $H[X_{\cup\rR}^H]$ is a matching if $\rR\in E(R)$ and empty if $\rR\in\binom{V(R)}{k}\sm E(R)$;
		
	\item $\cA=\bigcup_{H\in\cH, I\subseteq[r]_0\colon \text{$I$ is an index set}}A^H_I$ is a union of candidacy graphs with respect to~$\phi_\circ$ and~$G_A$; in particular, $A^H_I$ is $2|I|$-uniform, and $A^H_i$ is a balanced bipartite $2$-graph with vertex partition $(X^H_i,V_i)$ for each $i\in[r]_0$;
\item $\cB=\bigcup_{H\in\cH, j\in V(R)}B^H_j$ is a union of candidacy graphs with respect to~$\phi_\circ$ and~$G_B$; in particular, $B^H_j$ is a balanced bipartite $2$-graph with vertex partition $(X^H_j,V_j)$ for each $j\in\VR$;

	\item $\phi_\circ^-\colon \bigcup_{H\in\cH} X^{H,-}_{\cup-[r_\circ]}\to V_{\cup-[r_\circ]}$ with $X^{H,-}_i\subseteq X_i^H$, $\phi_\circ^-(X_i^{H,-})\subseteq V_i$ and $\phi_\circ^-|_{X_i^{H,-}}$ is injective for all $H\in\cH, i\in-[r_\circ]$, and
$\phi_\circ\colon \bigcup_{H\in\cH} X^H_{\cup-[r_\circ]}\to V_{\cup-[r_\circ]}$ is an extension of $\phi_\circ^-$ with $\phi_\circ(X_i^H)= V_i$ and $\phi_\circ|_{X_i^H}$ is bijective for all $H\in\cH, i\in-[r_\circ]$.
\end{itemize}
For simplicity, we often write $G:=G_A\cup G_B$, $\sX_{i}:=\bigcup_{H\in\cH}X^H_{i}$, $\sX_{i}^-:=\bigcup_{H\in\cH}X^{H,-}_{i}$, $\sA_I:=\bigcup_{H\in\cH}A_I^H$, $\cB_i:=\bigcup_{H\in\cH}B_i^H$, $\cX_\circ:=\sX_{\cup -[r_\circ]}$, $\cX_\circ^-:=\sX_{\cup -[r_\circ]}^-$ and $\cV_\circ:=V_{\cup-[r_\circ]}$ for all $i\in V(R)$ and index sets $I\subseteq [r]_0$.
Note that the packing instance $\sP$ naturally corresponds to a blow-up instance $$(\cH,G,R,
\{X^H_i\}_{i\in \VR, H\in\cH}, \{V_i\}_{i\in \VR})$$ of size $(n,k,r+r_\circ+1)$. In particular, we also use the notation of Section~\ref{sec:notation blow-up}.
For the sake of a better readability, we stick to some conventions:

We will often use the letters $(Z,\cZ)\in\{(A,\cA),(B,\cB) \}$ as many arguments for candidacy graphs~$\cA_I$ with respect to $G_A$, and candidacy graphs $\cB_i$ with respect to $G_B$ are the same.
Whenever we write $xv\in E(\cZ_i)$ for some $i\in[r]$, we tacitly assume that $x\in\sX_i, v\in V_i$. We usually denote edges in~$\cZ_i$ by (non-calligraphic) letters $e, f$, and hyperedges in $\cA_I$ by $\aA$ and {a collection of edges from $\bigsqcup_{j\in J}E(\cB_j)$ by $\bA=\{b_j\}_{j\in J}$, where we allow to slightly abuse the notation and often treat $\bA$ as~$b_{\cup J}$.}
Whenever we write $\{v_{i_1},\ldots,v_{i_k} \}\in E(G_Z)$, we tacitly assume that $v_{i_\ell}\in V_{i_\ell}$ for all $\ell\in[k]$; analogously for $\{x_{i_1},\ldots,x_{i_k} \}\in E(H)$. 
We usually refer to hyperedges in~$G_Z$ with letters $\gG, \hG$, hyperedges in $\cH$ with letters $\eH, \fH$, and hyperedges in $R$ with $\rR$.

The aim of this section is to map almost all vertices of $\sX_0$ into $V_0$ by defining a function \mbox{$\sigma\colon \sX_0^\sigma\to V_0$} in~$\sA_0$ 
(that is, $x\sigma(x)\in E(\sA_0)$)
where $\sX_0^\sigma\subseteq \sX_0$, while maintaining several properties for the other candidacy graphs.
We identify such a function~$\sigma$ with its \emph{corresponding edge set}~$M(\sigma)$ defined as 
\begin{align}\label{eq:def M(sigma)}
M(\sigma):=\{xv\colon x\in \sX_0^\sigma, v\in V_0, \sigma(x)=v \}.
\end{align}

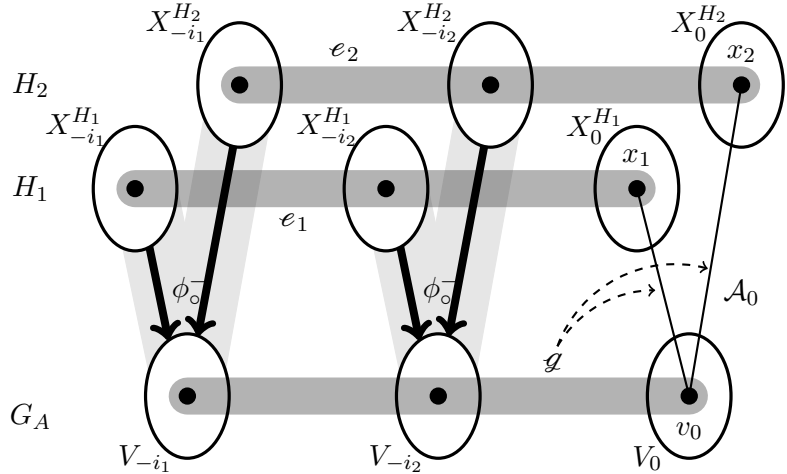
\begin{figure}[htb]
	\begin{center}
		\begin{tikzpicture}[scale = 0.55,every text node part/.style={align=center}]
		\def\ver{0.1} 

		\draw (-2.5,5) node {$H_1$};
		
		\draw (-1.4,6.5) node {$X_{-i_1}^{H_1}$};
		\draw (4.6,6.5) node {$X_{-i_2}^{H_1}$};
		\draw (11,6.5) node {$X_0^{H_1}$};

		\draw (-2.5,-0.5) node {$G_A$};
		
		\draw (.25,-1.5) node {$V_{-i_1}$};
		\draw (6.25,-1.5) node {$V_{-i_2}$};
		\draw (12.25,-1.5) node {$V_{0}$};
		 
		\draw (-2.5,7.5) node {$H_2$};
		
		\draw (1,9) node {$X_{-i_1}^{H_2}$};
		\draw (7,9) node {$X_{-i_2}^{H_2}$};
		\draw (13.5,9) node {$X_0^{H_2}$};

		\filldraw[fill=black!10, draw=black!10]
		(1.7,7)--(3.3,7) -- (2.05,0) -- (.45,0) -- cycle;

		\filldraw[fill=black!10, draw=black!10]
		(-0.8,6)--(0.8,6) -- (2.05,0) -- (.45,0) -- cycle;

		\draw[line width=3pt,->]
		(2.4,6.3) -- (1.5,1.45);

		\draw[line width=3pt,->]
		(.3,3.8) -- node[anchor=west] {$\phi_\circ^{-}$} (.8,1.35);

		\filldraw[fill=black!10, draw=black!10]
		(7.7,7)--(9.3,7) -- (8.05,0) -- (6.45,0) -- cycle;

		\filldraw[fill=black!10, draw=black!10]
		(5.2,6)--(6.8,6) -- (8.05,0) -- (6.45,0) -- cycle;

		\draw[line width=3pt,->]
		(8.4,6.3) -- (7.5,1.45);

		\draw[line width=3pt,->]
		(6.3,3.8) -- node[anchor=west] {$\phi_\circ^{-}$} (6.8,1.35);

		\draw[very thick, fill=white] 
		(1.25,0) ellipse (1 and 1.5)
		(7.25,0) ellipse (1 and 1.5)
		(13.25,0) ellipse (1 and 1.5)
		(0,5) ellipse (1 and 1.5)
		(6,5) ellipse (1 and 1.5)
		(12,5) ellipse (1 and 1.5)
		(2.5,7.5) ellipse (1 and 1.5)
		(8.5,7.5) ellipse (1 and 1.5)
		(14.5,7.5) ellipse (1 and 1.5);

		\coordinate (g1) at (1.25,0); 
		\node[hvertex] (G1) at (g1){};
		\node[hvertex] (g2) at (7.25,0){};
		\coordinate (g3) at (13.25,0); 		
		\node[hvertex,label={[label distance=.1cm]270:$v_0$}] (G3) at (g3){};
		\coordinate (x1) at (0,5); 		
		\node[hvertex] (X1) at (x1){};
		\node[hvertex] (x2) at (6,5){};
		\coordinate (x3) at (12,5); 		
		\node[hvertex,label={[label distance=.05cm]90:$x_1$}] (X3) at (x3){};
		\coordinate (y1) at (2.5,7.5); 
		\node[hvertex] (Y1) at (y1){};
		\node[hvertex] (y2) at (8.5,7.5){};
		\coordinate (y3) at (14.5,7.5); 		
		\node[hvertex,label={[label distance=.05cm]90:$x_2$}] (Y3) at (y3){};

		\draw[edge] (x1)--(x3);
		\draw[edge] (g1)--(g3);
		\draw[edge] (y1)--(y3);

		\draw[thick] (x3)--(g3);
		\draw[thick] (y3)--(g3);

		\draw (3.8,4.2) node {$\eH_1$};
		\draw (5,8.3) node {$\eH_2$};

		\draw (14.5,2.5) node {$\cA_0$};
		\draw (10,.8) node {$\gG$};
		\draw[thick,->, dashed] (10.1,1.2)  arc  (150:90:2.7cm);
		\draw[thick,->, dashed] (10.1,1.2)  arc  (160:75:3cm);
		\end{tikzpicture}
	\end{center}
		  \captionsetup{width=.9\linewidth}
\caption{We have $\cH=\{H_1,H_2\}$ and $k=3$.
If $\phi_\circ^{-}(\eH_1\sm\{x_1\})=\phi_\circ^{-}(\eH_2\sm\{x_2\})=\gG\sm\{v_0\}$, then we add the label $\gG$ to the edge set label $\psi(x_1v_0)$ and $\psi(x_2v_0)$ of the edges $x_1v_0$ and $x_2v_0$ of the candidacy graph~$\cA_0$. 
This captures the information that if $x_1$ or $x_2$ are mapped onto $v_0$ by $\sigma$ in $\cA_0$, then this embedding uses the edge $\gG\in E(G_A)$.}
	\label{fig:conflict-free packing}
\end{figure}

To incorporate that $\sigma$ has to be chosen such that each edge in $G_A$ is used at most once, we define an edge set labelling \defn{$\psi$ with respect to $\sP$} on $\sA_0$ as follows.
For every edge $xv\in E(\sA_0)$, we set
\begin{align}\label{eq:def edge set labelling}
\psi(xv):=\left\{
\phi_\circ^-(\eH\sm\{x\})\cup\{v\}\colon x\in\eH\in E(\cH) \text{ with } \eH\sm\{x\}\subseteq\cX_\circ^-
\right\}.
\end{align}
We defined the candidacy graphs $\sA_0$ in Definition~\ref{def:candidacy graph} such that $xv\in E(\sA_0)$ only if $\phi_\circ^-(\eH\sm\{x\})\cup\{v\}\in E(G_A)$ for each such edge $\eH$ as in~\eqref{eq:def edge set labelling}. 
That is, $\psi(xv)$ encodes the set of edges in~$G_A$ that are used for the packing when mapping $x$ onto $v$.
We say 
\begin{align}\label{eq:conflict free}
\begin{minipage}[c]{0.85\textwidth}\em
$\sigma\colon \sX_0^\sigma\to V_0$ is a \emph{conflict-free packing} if $\sigma|_{\sX_0^\sigma\cap X_0^H}$ is injective for all $H\in\cH$ and $\psi(e)\cap\psi(f)=\emptyset$ for all distinct $e,f\in M(\sigma)$.
\end{minipage}\ignorespacesafterend
\end{align}
Crucially note that the property that $\psi(e)\cap\psi(f)=\emptyset$ for all distinct $e,f\in M(\sigma)$ will guarantee that every edge in $G_A$ is used at most once. 
For an illustration, see Figure~\ref{fig:conflict-free packing}. 

Given a conflict-free packing $\sigma\colon \sX_0^\sigma\to V_0$ in $\sA_0$, we update the remaining candidacy graphs with respect to $\sigma$.
To account for the vertices in $\sX_0\sm\sX_0^\sigma$ that are left unembedded by $\sigma$, we will consider an extension $\sigma^+$ of~$\sigma$ such that $\sigma^+$ also maps every vertex $x_0\in \sX_0\sm \sX_0^\sigma$ to $V_0$ and $\sigma^+|_{X_0^H}$ is injective for all $H\in\cH$ (and hence bijective). 
We call such a~$\sigma^+$ a \defn{cluster-injective extension} of $\sigma$. 
The purpose of~$\sigma^+$ is that $\sigma^+|_{\sX_0\sm \sX_0^\sigma}$ will serve as a `dummy' extension resulting in an easier analysis of the packing process as $\sigma^+$ will impose further restriction that culminate in more consistent candidacy graphs.
Using Definition~\ref{def:candidacy graph}, we will consider the \defn{(updated) candidacy graphs} $A_I^H(\phi_\circ\cup \sigma^+)$ with respect to $\phi_\circ\cup \sigma^+$ and $G_A$ for index sets $I\subseteq[r]$, as well as the \defn{(updated) candidacy graphs} $B_j^H(\phi_\circ\cup \sigma^+)$ with respect to $\phi_\circ\cup \sigma^+$ and $G_B$ for $j\in V(R)$.\COMMENT{In Lemma~\ref{lem:packing lemma}, we consider subgraphs of $A_I^H(\phi_\circ\cup \sigma^+)$ by updating with respect to $H_+$, that is, we consider $A_I^{H_+}(\phi_\circ\cup \sigma^+)\subseteq A_I^H(\phi_\circ\cup \sigma^+)$.}

\medskip
To track our packing process, we carefully maintain quasirandom properties of the candidacy graphs throughout the entire procedure.
Our Approximate Packing Lemma will guarantee that we can find a conflict-free packing that behaves like an idealized random packing with respect to given sets of \emph{edge testers} (as defined in Definition~\ref{def:general edge tester}), and with respect to weight functions $\omega\colon\binom{E(\cA_0)}{\ell}^=\to[0,s]$ for $\ell\leq s$ that we will call \emph{local testers}.

To that end, assume we are given a packing instance $\sP=(\cH,G_A,G_B, R,\cA,\cB,\phi_\circ^-,\phi_\circ)$ of size $(n,k,r,r_\circ)$, and $\bd=(d_A,d_B, (d_{i,A})_{i\in [r]_0}, (d_{i,B})_{i\in \VR})$, 
and $q,t\in\bN$,
as well as a set $\cW_{edge}$ of edge testers.
We say $\sP$ is an \emph{$(\eps,q,t,\bd)$-packing instance with suitable edge testers $\cW_{edge}$}
if $|X_i^H|=|V_i|=(1\pm\eps)n$ for all $i\in\VR$, and the following properties are satisfied 
(recall~\eqref{eq:eps,t well-intersecting} and Definitions~\ref{def:general edge tester}--\ref{def:E_g,h,phi} for~\ref{item:P2}--\ref{item:P4}, respectively, {and that we write $\JXV$ for $\JX\cup \JV$}):
\begin{enumerate}[label={\rm (P\arabic*)}]
\item\label{item:P1} 
for all $i\in\VR$ and all pairs of disjoint sets $\cS_A,\cS_B\subseteq 
\bigcup_{\rR\in E(R)\colon i\in\rR} V_{\sqcup\rR\sm\{i\}}$ with $|\cS_A\cup\cS_B|\leq t$, we have 
$\big|V_i\cap N_{G_A}(\cS_A)\cap N_{G_B}(\cS_B)\big|
=(1\pm\eps)d_A^{|\cS_A|}d_B^{|\cS_B|}n;$

\item\label{item:P2} for all $H\in\cH,i\in[r]_0,j\in\VR$, we have that $A^H_i$ is $(\eps,d_{i,A})$-super-regular and $(\eps,q)$-well-intersecting with respect to $G_A$, and $B^H_j$ is $(\eps,d_{j,B})$-super-regular and $(\eps,q)$-well-intersecting with respect to $G_B$;

\item \label{item:P edge tester}
for every edge tester $(\omega,\omega_\iota,J,\JX, \JV,\bc,\pP)\in\cW_{edge}$ with respect to $(\omega_\iota,J,\JX, \JV,\bc,\pP)$, {$(\phi_\circ^-,\phi_\circ)$}, $\cA$ and $\cB$, with centres $\bc\in V_{\sqcup I}$ for $I\subseteq \VR$, $I_{r_0}:=(I\cap[r]_0)\sm J$, 
and patterns $\pP=(\bpA,\bppA,\bpB,\bppB)\in(\bN_0^{r_\circ+r+1})^4$, 
we have that
\begin{align*}
\omega\left(E(\sA_{I_{r_0}})\sqcup \textstyle\bigsqcup_{j\in J}E(\cB_j)\right)
=&
\Big(\IND\{\JXV\cap-[r_\circ]=\emptyset \}\pm\eps\Big)
\prod_{Z\in\{A,B\}}
d_Z^{\norm{\bpZ_{-[r_\circ]_0}}-\norm{\bppZ_{-[r_\circ]_0}}}
\\&\prod_{i\in I_{r_0}}d_{i,A}
\prod_{j\in  J}d_{j,B}
\frac{\omega_\iota(\sX_{\sqcup {I}})}{n^{|(I\cap-[r_\circ])\sm J|}}
\pm  n^{\eps};
\end{align*}

\item\label{item:P3.5} for all $\gG\in E(G_A[V_{\cup\rR}])$ for some $\rR\in E(R)$ with $\rR\cap[r]_0\neq\emptyset$, and all $\bp, \bpp\in\bN_0^{r_\circ+r+1}$, the set $\cW_{edge}$ contains the edge tester for $\cX_{\gG,\bp,\bp^{2nd},\phi_\circ^-}(\cA)$;

\item\label{item:P4} for all $\gG=\{v_{i_1},\ldots,v_{i_k}\}, \hG=\{w_{j_1},\ldots,w_{j_k}\}\in E(G_A)$ with $v_{i_k}=w_{j_k}$\COMMENT{Note that we implicitly assume that $i_k$ is the last cluster since we assume for an index set that $i_1\leq\ldots\leq i_k$.} and $I:=\{i_1,\ldots,i_k\}\neq \{j_1,\ldots,j_k\}=:J$, we have that
\begin{align*}
\left|E_{\gG,\hG,\phi_\circ^-}(\cA)\right|
\leq 
\max\left\{
n^{k-\left|(I\cup J)\cap -[r_\circ]\right|+\eps},
n^\eps
\right\}.
\end{align*}
\end{enumerate}
Note that~\ref{item:P1} also implies that $G_Z$ is $(3\eps,t,d_Z)$-typical with respect to $R$ for each $Z\in\{A,B\}$.

\medskip
Furthermore, we call a function $\omega\colon\binom{E(\cA_0)}{\ell}^=\to[0,s]$ for $\ell\leq s$ a \defn{local $s$-tester (for the $(\eps,q,t,\bd)$-packing instance $\sP$)} if $\normV{\omega}{\ell'}\leq n^{\ell-\ell'+\eps^2}$ for every $\ell'\in[\ell]$.
We introduce some more notation:
\begin{align}\label{eq:def b_i}
\begin{minipage}[c]{0.8\textwidth}\em
Let $E_i^R:=\{\rR\in E(R)\colon \{0,i\}=\rR\cap([r]_0 \cup\{i\}) \}$ and {$b_i:=|E_i^R|=\dg_{R[-[r_\circ]\cup\{0,i\}]}(\{0,i\})$} for each $i\in \VRsm$. 
For $I\subseteq \VR$, let $b_I:=\sum_{i\in I\sm\{0\}}b_i$.
For $i\in[r], j\in\VRsm$, let $d_{i,A}^{new}:=d_{i,A}d_A^{b_{i}}$ and $d_{j,B}^{new}:=d_{j,B}d_B^{b_{j}}$.
\end{minipage}
\ignorespacesafterend 
\end{align}
Note that for $i\in \VR\sm N_{R_\ast}[0]$, we have $b_i=0$ and thus $d_{i,Z}^{new}=d_{i,Z}$ for $Z\in\{A,B\}$.

\subsection{Approximate Packing Lemma}\label{sec:approx packing lem}
We now state the Approximate Packing Lemma, which is the key tool for the proof of our main result.

\begin{lemma}[Approximate Packing Lemma]\label{lem:packing lemma}
Let $1/n\ll\eps\ll\eps'\ll1/t\ll 1/k,1/q,1/r,1/(r_\circ+1),1/s.$
Suppose $(\cH,G_A,G_B,R,\cA,\cB,\phi_\circ^-,\phi_\circ)$ is an $(\eps,q,t,\bd)$-packing instance of size $(n,k,r,r_\circ)$ with $\bd\geq n^{-\eps}$ and suitable $s$-edge-testers $\cW_{edge}$. Suppose further that $\cW_{local}$ is a set of local $s$-testers, {$\cW_{0}$~is a set of tuples $(\omega,c)$ with $\omega\colon \sX_0\to[0,s]$, $c\in V_0$,} and $|\cW_{edge}|,$ $|\cW_{local}|,|\cW_{0}|\leq n^{4\log n}$, $|\cH|\leq n^{2k}$,\COMMENT{This is not restrictive and $|\cH|\leq n^{2k}$ is chosen arbitrarily.} as well as 
$e_{\cH}(\sX_{\sqcup\rR})\leq d_A n^k$ for all $\rR\in E(R)$.

Then there is a conflict-free packing $\sigma\colon \sX_0^\sigma\to V_0$ in $\sA_0$ and a cluster-injective extension $\sigma^+$ of $\sigma$ such that for all $H\in\cH$, we have $|\sX_0^\sigma\cap X_0^H|\geq(1-\eps')n$, and
for all index sets $I_A\subseteq[r]$, {$I_B\in \VR$} and $(Z,\cZ)\in\{(A,\cA),(B,\cB)\}$, there exist spanning subgraphs $Z_{I_Z}^{H,new}$ of 
the candidacy graphs $Z^{H}_{I_Z}(\phi_\circ\cup\sigma^+)$ with respect to $\phi_\circ\cup\sigma^+$ and $G_Z$
(where $\cZ_{I_Z}^{new}:=\bigcup_{H\in \cH}Z_{I_Z}^{H,new}$ and $\cZ^{new}$ is the collection of all $\cZ_{I_Z}^{new}$) such that

\begin{enumerate}[label={\rm (\Roman*)\textsubscript{L\ref{lem:packing lemma}}}]
\item\label{emb lem:sr} $Z^{H,new}_i$ is $(\eps',d_{i,Z}^{new})$-super-regular and $(\eps',q+\Delta(R))$-well-intersecting with respect to $G_Z$ for all $H\in\cH$, and all $i\in[r]$ if $Z=A$, and all $i\in \VRsm$ if $Z=B$;

\item \label{emb lem:general edge testers}
for every (general) $s$-edge tester $(\omega,\omega_\iota,J,\JX, \JV,\bc,\pP)\in\cW_{edge}$  with respect to \linebreak$(\omega_\iota,J,\JX, \JV,\bc,\pP)$, {$(\phi_\circ^-,\phi_\circ)$}, $\cA$ and $\cB$, with centres $\bc\in V_{\sqcup I}$ for $I\subseteq \VR$, $I_r:=(I\cap[r])\sm J$, $I_r\cup J\neq \emptyset$,\COMMENT{We only make a statement if there is either a relevant cluster in $J$ or in $I_r$ left.} and patterns $\pP=(\bpA,\bppA,\bpB,\bppB)\in(\bN_0^{r_\circ+r+1})^4$, 
the $s$-edge tester $(\omega^{new},\omega_\iota,J,\JX, \JV,\bc,\pP)$ defined as in Definition~\ref{def:general edge tester} with respect to $(\omega_\iota,J,\JX, \JV,\bc,\pP)$, {$(\phi_\circ^-\cup\sigma,\phi_\circ\cup\sigma^+)$}, $\cA^{new}$ and $\cB^{new}$ satisfies
\begin{align*}
\omega^{new}\left(E(\sA_{I_r}^{new})\sqcup \textstyle\bigsqcup_{j\in J}E(\cB_j^{new})\right)
=&
\left(\IND\{\JXV\cap-[r_\circ]_0=\emptyset \}\pm\eps'^2\right)
\prod_{Z\in\{A,B\}}
d_Z^{\norm{\bpZ_{-[r_\circ]_0}}-\norm{\bppZ_{-[r_\circ]_0}}}
\\&\prod_{i\in I_r}d_{i,A}^{new}
\prod_{j\in  J}d_{j,B}^{new}
\frac{\omega_\iota(\sX_{\sqcup {I}})}{n^{|(I\cap-[r_\circ]_0)\sm J|}}
\pm  n^{\eps'};
\end{align*}

\item\label{emb lem:local weight}
$\omega(M(\sigma))=(1\pm\eps'^2)\omega(E(\sA_0))/(d_{0,A}n)^\ell\pm n^\eps$ for every local $s$-tester $\omega\in\cW_{local}$ with \linebreak$\omega\colon\binom{E(\cA_0)}{\ell}^=\to[0,s]$;
\COMMENT{The error bound `$\pm n^\eps$' is correct, compare~\eqref{eq:w packing size}. We also make use of this, e.g. when establishing~\ind{q+1}\ref{vertex tester}, last display.}

\item\label{emb lem:hit tester}
{$\omega\left(
\{
x\in\sX_0\sm\sX_0^\sigma\colon \sigma^+(x)=c
\}
\right)
\leq \omega(\sX_0)/n^{1- \eps}+n^\eps$
for every $(\omega,c)\in\cW_{0}$.}
\end{enumerate}
\end{lemma}

Properties~\ref{emb lem:sr} and \ref{emb lem:general edge testers} ensure that~\ref{item:P2} and \ref{item:P edge tester} are also satisfied for the updated candidacy graphs $\cA^{new}$ and $\cB^{new}$, respectively.
Property~\ref{emb lem:local weight} states that $\sigma$ behaves like a random packing with respect to the local testers, which for instance can be used to establish~\ref{item:P4} for future packing rounds.
{Property~\ref{emb lem:hit tester} allows to control the weight on vertices that are not embedded by $\sigma$ but are nevertheless mapped onto a specific vertex $c$ by the extension $\sigma^+$.}

\begin{proof}
We split the proof into three parts. 
In Part~\ref{part:H_+} we construct an auxiliary supergraph $H_+$ of every $H\in\cH$ by adding some hyperedges to $H[X^H_{\cup {\rR}}]$ for every $\rR\in E(R)$ in order to make the packing procedure more uniform.
In Part~\ref{part:hypergraph} we construct an auxiliary hypergraph $\cH_{aux}$ for $\cA_0$ such that we can use Theorem~\ref{thm:hypermatching} to find a conflict-free packing in $\cA_0$. 
In order to be able to apply Theorem~\ref{thm:hypermatching}, we exploit~\ref{item:P4} as well as~\ref{item:P2} together with~\ref{item:P edge tester} to control $\Delta_2(\cH_{aux})$ and $\Delta(\cH_{aux})$, respectively.
In Part~\ref{part:weight functions} we define weight functions and employ the conclusion of Theorem~\ref{thm:hypermatching} to establish~\ref{emb lem:sr}--\ref{emb lem:local weight}.

\smallskip
Let $\Delta_R:=2\binom{r+r_\circ}{k-1}$. Note that $\Delta(R),\Delta(\cH),b_i\leq \Delta_R$ for all $i\in[r]$. For simplicity we write $d_0:=d_{0,A}$ and $\heps:=\eps^{1/2}$.
Further, we choose a new constant $\Delta$ such that $\eps'\ll1/\Delta\ll1/t$.

\begin{proofpart}\label{part:H_+}
Construction of $H_+$
\end{proofpart}
We construct an auxiliary supergraph $H_+$ of~$H$ by artificially adding some edges to $H$ for every $\rR\in E_i^R$ and $i\in\VRsm$. 
(Recall~\eqref{eq:def b_i} for the definition of~$E_i^R$ {and note that it can be that $\rR\in E_i^R\cap E_j^R$ for $i,j\in-[r_\circ]$}.)
For every $H\in\cH$, we proceed as follows. 
We obtain $H_+$ from~$H$ by adding {a minimal number of} hyperedges of size~$k$ subject to the conditions that for all $i\in\VRsm$ and $\rR \in E_i^R$, an $H_+$-edge
 meets every cluster in $H[X^H_{\cup\rR}]$ exactly once, and\looseness=-1
\begin{enumerate}[label={\rm (\alph*)}]
\item\label{item:everyone has backneighbour in H_+} {$\dg_{H_+[X_{\cup\rR}^H]}(x_i)\in\{1,2\}$}
for all $x_i\in X_i^H$, and $\dg_{H_+[X_{\cup\rR}^H]}(x_0)\leq 2$
for all $x_0\in X_0^H$; 
\item\label{item:no 2-overlap}
{for all $\eH\in E(H_+[X_{\cup\rR}^H])$, we have $|\{x\in \eH\colon \dg_{H_+[X_{\cup\rR}^H]}(x)=2 \}|\leq 1$;}
\item\label{item:pattern in H_+}
 for all $\{x_0,x_i\}\in X_0^H\sqcup X_i^H$, if
$\{x_0,x_i\}\subseteq \eH$ for some $\eH\in E(H_+)\sm E(H)$, then $\{x_0,x_i\}\nsubseteq \fH$ for all $\fH\in E(H)$.\COMMENT{Note that we only require that not both $x_0$ and $x_i$ are contained in $\fH$, that is, $\{x_0,x_i\}\nsubseteq \fH$ but it might well be that $\{x_0,x_i\}\cap\fH\neq\emptyset$.}\nopagebreak 
\end{enumerate}
\nopagebreak Note that~\ref{item:everyone has backneighbour in H_+}--\ref{item:pattern in H_+} can be met because $|X_i^H|=(1\pm\eps)n$ for all $i\in \VR$ {and $\Delta(R)\leq \Delta_R$}.
{For all $i\in\VRsm$, let~$H_+^i$ be an arbitrary but fixed $k$-graph $H\subseteq H_+^i\subseteq H_+$ such that for all $\rR\in E_i^R$, we have $\dg_{H_+^i[X_{\cup\rR}^H]}(x_i)=1$.
Observe that by the construction of $H_+$, we have $\dg_{H_+[X_{\cup\rR}^H]}(x_i)=1$ for all $x_i\in X_i^H$, $i\in[r]$, $\rR\in E_i^R$, and thus $H_+^i=H_+$ for all $i\in[r]$.}
We make some observations.\looseness=-1 \nopagebreak
\begin{align}
\allowdisplaybreaks
\label{eq:b_i backneighbours}
&\begin{minipage}[c]{0.9\textwidth}\em
For all $x\in X_i^H, i\in\VRsm$, we have $\sum_{\rR\in E_i^R}\dg_{H_+^i[X_{\cup\rR}^H]}(x)=b_i$.
\end{minipage}
\\\allowdisplaybreaks
\label{eq:pairs y,y'}
&\begin{minipage}[c]{0.9\textwidth}\em
By~\ref{item:everyone has backneighbour in H_+}, for every $x\in X_i^H, i\in\VRsm$, there are at most $\Delta_R$ vertices $x'\in X_i^H\sm\{x\}$ such that $x$ and $x'$ have a common neighbour in $X_0^H$ in $H_+[\cX_\circ,X_0^H,X_i^H]$, that is, $e_x\cap e_{x'}\cap X_0^H\neq\emptyset$ with $e_y\in E(H_+[\cX_\circ,X_0^H,X_i^H])$ and $y\in e_y$ for both $y\in\{x,x'\}$.
\end{minipage}
\\
\allowdisplaybreaks
\label{eq:pattern in H_+}
&\begin{minipage}[c]{0.9\textwidth}\em
If $\{x_0,x_i\}$ lies in an edge of $H_+-H$, then $\{x_0,x_i\}$ does not lie in an edge of $H$.
\end{minipage}
\\\label{eq:H_+ intersection cardinality 1}
&\begin{minipage}[c]{0.9\textwidth}\em
If $|\eH\cap X^H_{\cup {[r]}}|\geq 2$, then $\eH\in E(H)$ for all $\eH\in E(H_+)$.
\end{minipage}
\end{align}
\COMMENT{$H_+-H$ is defined in Section~\ref{sec:notation}.}\COMMENT{\eqref{eq:pairs y,y'} is used for the codegree condition in~\eqref{eq:codegree y,y'}.}\COMMENT{Observation for~\eqref{eq:pattern in H_+}:
We want to add edges to $\cH$ without changing the \scd-patterns of proper edges in $\cH$.}

We introduce some simpler notation how to denote edges in $H_+$ (respectively $H_+^i$) that contain a vertex~$x\in X_0^H$.
Let $\cH_+:=\bigcup_{H\in\cH}H_+$, {and for $i\in\VRsm$, let $\cH_+^i:=\bigcup_{H\in\cH}H_+^i$.} 
For all $x\in\sX_0$ and  $\by\in \sX_{\sqcup I}$ for some $I\subseteq \VR$,\COMMENT{For the definition it does not matter whether $I$ is an index set or not.} let
\begin{align}\label{eq:def E_x,y}
E_{x,\by}:=\big\{
\eH\in E(\cH_+[\cX_\circ,\sX_0,\sX_{\cup I}])\colon& 
x\in\eH, \eH\cap\sX_{\cup {I}}\subseteq\{x\}\cup\by\cup\cX_\circ,\eH\cap(\by\sm\sX_0)\neq \emptyset,
\\\nonumber&
\text{if $\eH\cap(\by\sm\sX_0)\in\sX_i$ for some $i\in\VRsm$, then $\eH\in E(\cH_+^i)$}
\big\}.
\end{align}
That is, $E_{x,\by}$ contains essentially all $\cH_+$-edges that contain $x\in\sX_0$ and a non-empty subset of~$\by$ and whose remaining vertices are already embedded and lie in $\cX_\circ$.
In particular, if $\by=y$ is a single vertex $y\in\sX_i$, then $E_{x,y}$ contains all {$\cH_+^i$-edges} that contain $x$ and $y$ and whose remaining vertices lie in~$\cX_\circ$.
Hence, note that by definition of $\cH_+$ and as observed in~\eqref{eq:b_i backneighbours}, we have that
\begin{align}\label{eq:size E_x,y}
\big|\textstyle\bigcup_{x\in\sX_0}E_{x,y}\big|=b_i \text{ for all $y\in\sX_i, i\in \VRsm$.}
\end{align}

\begin{proofpart}\label{part:hypergraph}
Applying Theorem~\ref{thm:hypermatching}
\end{proofpart}

Our strategy is to utilize Theorem~\ref{thm:hypermatching} to find the required conflict-free packing $\sigma$ in~$\sA_0$. 
To that end, we will define an auxiliary hypergraph $\cH_{aux}$ for $\cA_0$.
Let~$\psi\colon E(\sA_0)\to2^\cE$ be the edge set labelling with respect to the packing instance as defined in~\eqref{eq:def edge set labelling}.
For all $H\in\cH$, the hypergraph $H[X_{\cup\rR}^H]$ is a matching if $\rR\in E(R)$ and empty otherwise, and thus we have that $\norm{\psi}\leq \binom{r_\circ}{k-1}\leq\Delta_R$.
In the following, we may assume that $|\psi(e)|=\Delta_R$ for all $e\in E(\sA_0)$ as we may simply add distinct artificial dummy labels that we ignore afterwards again.

Further, let~$(V_0^H)_{H\in\cH}$ be disjoint copies of $V_0$, and for all $H\in\cH$ and $e=x_0v_0\in E(A_0^{H})$, let $e^H:=x_0v_0^H$ where $v_0^H$ is the copy of $v_0$ in $V_0^H$.
Let $\hG_e:=e^H\cup \psi(e)$ for each $e\in E(A_0^H)$, $H\in\cH$ and let $\cH_{aux}$ be the $(\Delta_R+2)$-graph with vertex set $\bigcup_{H\in\cH}(X_0^H\cup V_0^H)\cup \cE$
and edge set $\{\hG_e\colon e\in E(\sA_0) \}$.
A~key property of the construction of $\cH_{aux}$ is a bijection between conflict-free packings $\sigma$ in $\sA_0$ and matchings $\cM$ in $\cH_{aux}$ by assigning $\sigma$ to $\cM=\{\hG_e\colon e\in M(\sigma) \}$. (Recall that $M(\sigma)$ is the edge set corresponding to $\sigma$ as defined in~\eqref{eq:def M(sigma)}.\COMMENT{$M :=\{xv\colon x\in \sX_0^\sigma, v\in V_0, \sigma(x)=v \}$})

\begin{step}\label{step:hypergraph boundedness}
Estimating $\Delta(\cH_{aux})$ and $\Delta_2(\cH_{aux})$\nopagebreak
\end{step}

In order to apply Theorem~\ref{thm:hypermatching} to $\cH_{aux}$, we estimate $\Delta(\cH_{aux})$ and $\Delta_2(\cH_{aux})$.
We first claim that
\begin{align}\label{eq:Delta(H)}
\Delta(\cH_{aux})\leq(1+\eps^{2/3})d_0n=:\Delta_{aux}.
\end{align}
Since $A_0^{H}$ is $(\eps,d_0)$-super-regular and $|X_0^H|=|V_0|=(1\pm\eps)n$ for each $H\in\cH$, we have that an appropriate upper bound on $\Delta_\psi(\sA_0)$ immediately establishes~\eqref{eq:Delta(H)}.
In the following 
we derive such an upper bound on $\Delta_\psi(\sA_0)$ by employing property~\ref{item:P edge tester}.
For all $\rR\in E(R)$ with $0\in\rR$, $|\rR\cap-[r_\circ]|=k-1$, and all $\gG\in E(G_A[V_{\cup\rR}])$ with $\gG\sm\cV_\circ=\{v_0\}$, note that $\bigcup_{\bp,\bp^{2nd}}\cX_{\gG,\bp,\bp^{2nd},\phi_\circ^-}(\cA)$ contains by Definition~\ref{def:X_g,p,p,phi} all vertices $x_0\in N_{\sA_0}(v_0)$ (compare with~\eqref{def:X_g,p,p,phi I}) that are contained in an $\cH$-edge that could be mapped onto $\gG$ with respect to $\phi_\circ^-$ and $\cA$ (compare with~\eqref{def:X_g,p,p,phi II}).
Hence, by the definition of the edge set labelling~$\psi$ in~\eqref{eq:def edge set labelling}, 
$\gG$ appears as a label of $\psi$ on at most 
$\sum_{\bp,\bp^{2nd}}|\cX_{\gG,\bp,\bp^{2nd},\phi_\circ^-}(\cA)|$
edges of $\sA_0$.
Note that by Definition~\ref{def:e_H,p} of $E_\cH(\bp,\bp^{2nd},\rR)$, we obtain 
\begin{align}\label{eq:sum E(p,p,r)}
\sum_{\bp,\bp^{2nd}}|E_\cH(\bp,\bp^{2nd},\rR)|=e_\cH(\sX_{\sqcup \rR})\leq d_An^k,
\end{align}
where the last inequality holds by assumption of Lemma~\ref{lem:packing lemma}.
For all $x_0\in N_{A_0^H}(v_0)$ for some $H\in\cH$, 
let $\gG^{-1}_{x_0}:=(\phi_\circ^-|_{V(H)})^{-1}(\gG\sm\{v_0\})\cup\{x_0\}$.
If $\gG^{-1}_{x_0} \in \cX_{\gG,\bp,\bp^{2nd},\phi_\circ^-}(\cA)$ for  $\bp,\bp^{2nd}\in\bN_0^r$, then we have by~\eqref{def:X_g,p,p,phi II} that $\gG^{-1}_{x_0} \in E_\cH(\bp,\bpp,\rR)$, and thus
\begin{align*}
\bpA
\left(\gG^{-1}_{x_0},\emptyset\right)
=
\bp,
\quad
\bppA
\left(\gG^{-1}_{x_0},\emptyset\right)
=
\bp^{2nd},~~~\text{and}~~~
\norm{\bp}-\norm{\bp^{2nd}}
\stackrel{\eqref{eq:norms equal}}{=}
-1.
\end{align*}
By~\ref{item:P3.5}, the set $\cW_{edge}$ contains the (simple) edge tester $(\omega,\omega_\iota,\gG,\bp,\bp^{2nd})$ for $\cX_{\gG,\bp,\bp^{2nd},\phi_\circ^-}(\cA)$ (as defined in Definition~\ref{def:X_g,p,p,phi}) with $\omega(E(\cA_{0}))=|\cX_{\gG,\bp,\bp^{2nd},\phi_\circ^-}(\cA)|$
and $\omega_\iota(\sX_{\sqcup\rR})=|E_\cH(\bp,\bp^{2nd},\rR)|$.
Hence, by~\ref{item:P edge tester} (with $I=\rR$, $J=\JX=\JV=\emptyset$, $\bpA=\bp$, $\bppA=\bpp$, $\bpB=\bppB=\mathbf{0}$), we obtain\looseness=-1
\begin{align*}
\Delta_\psi(\sA_0) 
&\stackrel{\hphantom{{\eqref{eq:sum E(p,p,r)}}}}{\leq} \sum_{\bp,\bp^{2nd}}|\cX_{\gG,\bp,\bp^{2nd},\phi_\circ^-}(\cA)|
\\&
\stackrel{\ref{item:P edge tester}~}{\leq}
\sum_{\bp,\bp^{2nd}}\Big(
(1+\eps)d_A^{-1}d_0\frac{|E_\cH(\bp,\bp^{2nd},\rR)|}{n^{k-1}}+ n^{\eps}\Big)
\stackrel{\eqref{eq:sum E(p,p,r)}}{\leq}
(1+\eps^{2/3})d_0n.
\end{align*}
This establishes~\eqref{eq:Delta(H)}.

\medskip
Next, we claim that
\begin{align}\label{eq:Delta^c(H)}
\Delta_2(\cH_{aux})\leq n^\eps\leq \Delta_{aux}^{1-\eps^2}.
\end{align}
Note that the codegree in $\cH_{aux}$ of two vertices in $\bigcup_{H\in\cH}(X_0^H\cup V_0^H)$ is at most~1, and similarly, the codegree in $\cH_{aux}$ of a vertex in $\bigcup_{H\in\cH}(X_0^H\cup V_0^H)$ and a label in $\cE$ is at most~1 because $\Delta_\psi(A_0^H)\leq 1$ for all $H\in\cH$.
Hence, an appropriate upper bound on $\Delta_\psi^c(\sA_0)$ establishes~\eqref{eq:Delta^c(H)}. 
In the following we derive such an upper bound on $\Delta_\psi^c(\sA_0)$ by employing~\ref{item:P4}.
For all $\gG=\{v_{i_1},\ldots,v_{i_k}\}, \hG=\{w_{j_1},\ldots,w_{j_k}\}\in E(G_A)$ with $v_{i_k}=w_{j_k}\in V_0$, $(\gG\cup\hG)\sm\{v_{i_k}\}\subseteq \cV_\circ$, and $\{i_1,\ldots,i_k\}\neq \{j_1,\ldots,j_k\}$, note that $\gG$ and $\hG$ appear together as labels of $\psi$ on at most $|E_{\gG,\hG,\phi_\circ^-}(\cA)|$ edges of~$\sA_0$.
This follows immediately from Definition~\ref{def:E_g,h,phi} of $E_{\gG,\hG,\phi_\circ^-}(\cA)$.
Note further that $\left|\{i_1,\ldots,i_k,j_1,\ldots,j_k\}\cap -[r_\circ]\right|\geq k$ because $(\gG\cup\hG)\sm\{v_{i_k}\}\subseteq \cV_\circ$, and $\{i_1,\ldots,i_k\}\neq \{j_1,\ldots,j_k\}$.
Hence, by~\ref{item:P4}, we have 
\begin{align*}
\left|E_{\gG,\hG,\phi_\circ^-}(\cA)\right|
\leq 
\max\left\{
n^{k-\left|\{i_1,\ldots,i_k,j_1,\ldots,j_k\}\cap -[r_\circ]\right|+\eps},
n^\eps
\right\}= n^\eps,
\end{align*}
and thus, $\Delta_\psi^c(\sA_0)\leq n^\eps$, which establishes~\eqref{eq:Delta^c(H)}.

\begin{step}\label{step:hypergraph matching}
Applying Theorem~\ref{thm:hypermatching} to $\cH_{aux}$
\end{step}

Suppose $\cW=\bigcup_{\ell\in[\Delta]}\cW_\ell$ is a set of size at most $n^{4\log n}$ of given weight functions $\omega\in\cW_\ell$ for $\ell\in[\Delta]$ with $\omega\colon \binom{E(\sA_0)}{\ell}\to[0,\Delta]$ and 
\begin{align}\label{eq:normcondition}
\normV{\omega}{\ell'}\leq n^{\ell-\ell'+\eps^2} \text{ for every $\ell'\in[\ell]$.}
\end{align}
\COMMENT{All our considered weight functions fulfil the condition that $\normV{\omega}{\ell'}\leq n^{\ell-\ell'+\eps^2}$ for all $\ell'\in[\ell]$.
Therefore, the condition in Theorem~\ref{thm:hypermatching} that $\omega(E(\sA_0))\geq \normV{\omega}{\ell'}\Delta_{aux}^{\ell'+\eps^2}$ for all $\ell'\in[\ell]$ is already satisfied if $\omega(E(\sA_0))\geq n^{\ell+\eps/2}$.}
Note that every weight function $\omega\colon\binom{E(\sA_0)}{\ell}\to[0,\Delta]$ naturally corresponds to a weight function $\omega_{\cH_{aux}}\colon\binom{E(\cH_{aux})}{\ell}\to[0,\Delta]$ by defining $\omega_{\cH_{aux}}(\{\hG_{e_1},\ldots,\hG_{e_\ell}\}):=\omega(\{e_1,\ldots,e_\ell\})$.
We will explicitly specify~$\cW$ in Part~\ref{part:weight functions}, where every weight function $\omega\colon \binom{E(\sA_0)}{\ell}\to[0,\Delta]$ in $\cW_\ell$ for $\ell\in[\Delta]$ will be defined such that $\supp(\omega)\subseteq \bigcup_{H\in \cH}\binom{E(A_0^H)}{\ell}^=$ and in particular, such that the corresponding weight function $\omega_{\cH_{aux}}$ will also be clean, 
that is $\supp(\omega_{\cH_{aux}})\subseteq \binom{E(\cH_{aux})}{\ell}^=$.
Our main idea is to find a hypergraph matching in $\cH_{aux}$ that behaves like a typical random matching with respect to $\{\omega_{\cH_{aux}}\colon \omega\in\cW \}$ in order to establish~\ref{emb lem:sr}--\ref{emb lem:local weight}.

Suppose $\ell\in[\Delta]$ and $\omega\in\cW_{\ell}$. 
If $\omega(E(\sA_0))\geq n^{\ell+\eps/2}$, define $\tilde{\omega}:=\omega$.
Otherwise, arbitrarily choose $\tilde{\omega}\colon\binom{E(\sA_0)}{\ell}\to[0,\Delta]$ such that $\omega\leq\tilde{\omega}$, $\tilde{\omega}$ satisfies that $\supp(\tilde{\omega})\subseteq\bigcup_{H\in \cH} \binom{E(A_0^H)}{\ell}^\cl$ and in particular $\supp(\tilde{\omega}_{\cH_{aux}})\subseteq \binom{E(\cH_{aux})}{\ell}^\cl$, $\tilde{\omega}(E(\sA_0))=n^{\ell+\eps/2}$, and $\normV{\tilde{\omega}}{\ell'}\leq n^{\ell-\ell'+2\eps^2}$ for all $\ell'\in[\ell]$.
By~\eqref{eq:Delta(H)} and~\eqref{eq:Delta^c(H)}, we can apply Theorem~\ref{thm:hypermatching} 
(with $\Delta_{aux},\eps^2,\Delta_R+2,\Delta,\{\tilde{\omega}_{\cH_{aux}}\colon{\omega}\in\cW_\ell \}$ playing the roles of the parameters $\Delta,\delta,r,L,\cW_\ell$ of Theorem~\ref{thm:hypermatching}, respectively) to obtain a matching $\cM$ in $\cH_{aux}$ that corresponds to a conflict-free packing $\sigma\colon \sX_0^\sigma\to V_0$ in $\sA_0$ with its corresponding edge set $M(\sigma)$ that satisfies the following properties (where $\heps=\eps^{1/2}$):
\begin{align}
\label{eq:w packing size lower bound}
\omega(M(\sigma))&=(1\pm \eps)(1-\ell\eps^{2/3}) \frac{\omega(E(\sA_0))}{(d_0n)^\ell}\pm n^\eps
\\\label{eq:w packing size}
&=
(1\pm\heps)\frac{\omega(E(\sA_0))}{(d_0n)^\ell}\pm n^\eps \text{ for all } \omega\in\cW_\ell, \ell\in[\Delta].
\end{align}
\COMMENT{
\begin{align*}
(1\pm\Delta_{aux}^{-{\eps^2}/{50\Delta^2(\Delta_R+2)^{2} }})
\frac{\omega(E(\sA_0))}{((1+\eps^{2/3})d_0n)^\ell}
&=
(1\pm2\Delta_{aux}^{-{\eps^2}/{50\Delta^2(\Delta_R+2)^{2} }})
(1-\eps^{2/3})^\ell
\frac{\omega(E(\sA_0))}{(d_0n)^\ell}
\\&=
(1\pm3\Delta_{aux}^{-{\eps^2}/{50\Delta^2(\Delta_R+2)^{2} }})
(1-\ell\eps^{2/3})
\frac{\omega(E(\sA_0))}{(d_0n)^\ell}
\\&=
(1\pm\eps)(1-\ell\eps^{2/3})\frac{\omega(E(\sA_0))}{(d_0n)^\ell}
\end{align*}
}
\COMMENT{For $\omega\in\cW_\ell$ with $\omega(E(\sA_0))<n^{\ell+\eps/2}$, we have $\omega(M)\leq\tilde{\omega}(M)\leq (1+\heps)\frac{n^{\ell+\eps/2}}{(d_0n)^\ell}\leq n^\eps$.}

\begin{proofpart}\label{part:weight functions}
Employing weight functions to conclude~\ref{emb lem:sr}--\ref{emb lem:local weight}
\end{proofpart}

Let $\sigma\colon \sX_0^\sigma\to V_0$ be the conflict-free packing in $\cA_0$ as obtained in Part~\ref{part:hypergraph} 
and let $\sigma^+$ be a cluster-injective extension of $\sigma$ chosen uniformly and independently at random.
We will show that the random $\sigma^+$ satisfies with high probability the conclusions of the lemma and thus, there exists a suitable cluster-injective extension~$\sigma^+$ by picking one such extension deterministically.
We may assume that~\eqref{eq:w packing size} holds for a set of weight functions $\cW$. 
Each of these weight functions will only depend on our input parameters.
Hence, we could define them right away but for the sake of a cleaner presentation we postpone their definitions to the specific situations when we employ those weight functions to establish~\ref{emb lem:sr}--\ref{emb lem:local weight}.
We now define the candidacy graphs $Z_{I_Z}^{H,new}\subseteq Z_{I_Z}^{H}(\phi_\circ\cup\sigma^+)$ for $Z\in\{A,B\}$ and all index sets $I=I_A\subseteq [r]$, and {$i=I_B\in\VR$}. 
If $I_Z\cap\rR=\emptyset$ for all $\rR\in E(R)$ with $0\in\rR$, then we set $Z_{I_Z}^{H,new}:=Z_{I_Z}^H$.
Otherwise, let
\begin{align}\label{eq:def Z^H,new}
A_I^{H,new}:=A_I^{H_+}(\phi_\circ\cup\sigma^+),\quad\text{ and }\quad
B_i^{H,new}:=B_i^{H_+^i}(\phi_\circ\cup\sigma^+), 
\end{align}
with $A_I^{H_+}(\phi_\circ\cup\sigma^+)$  defined as in Definition~\ref{def:candidacy graph} with respect to $H_+$, $\phi_\circ\cup\sigma^+$ and~$G_A$, as well as $B_i^{H_+^i}(\phi_\circ\cup\sigma^+)$ defined with respect to $H_+^i$, $\phi_\circ\cup\sigma^+$ and~$G_B$. 
\COMMENT{
Note that by this definition, we do not require that $A_I^{H,new}$ is a subgraph of the previous candidacy graphs.
\\We recall Definition~\ref{def:candidacy graph} of $A_I^{H_+}(\phi_\circ\cup\sigma^+)$:
For all $H_+$ and an index set $I\subseteq [r]$, let $A_I^{H_+}(\phi_\circ\cup\sigma^+)$ be the $2|I|$-graph with vertex set $X_{\cup I}^H\cup V_{\cup I}$ and $\bigcup_{i\in I}\{x_i,v_i\}\in E(A_I^{H_+}(\phi_\circ\cup\sigma^+))$ for $\{x_i,v_i\}\in X_i^H\sqcup V_i$ (we could put here the condition that $\bigcup_{i\in I}\{x_i,v_i\}\in E(A_I^H)$)
 if all $\eH=\eH_\circ\cup \eH_{m'}\in E(H_+[X^H_{\cup -[r_\circ]_0},X_{i_1}^H,\ldots,X_{i_{m'}}^H])$ with $\eH_\circ\subseteq\binom{ X^H_{\cup -[r_\circ]_0}}{k-{m'}}$, $\eH_{m'}=\{x_{i_1},\ldots,x_{i_{m'}}\}\subseteq \{x_i\}_{i\in I}$ and ${m'}\in[|I|]$,
satisfy 
\begin{align*}
(\phi_\circ\cup\sigma^+)(\eH_\circ)\cup \{v_{i_1},\ldots,v_{i_{m'}}\} \in E(G[\cV_\circ,V_{i_1},\ldots,V_{i_{m'}}]).
\end{align*}}

\smallskip
Before we establish~\ref{emb lem:sr}--\ref{emb lem:local weight}, we first estimate $|\sX_0^\sigma\cap X_0^H|$ for each $H\in\cH$. 
We define a weight function $\omega_{H}\colon E(A_0^H)\to\{0,1\}$ for each $H\in\cH$ by $\omega_{H}(e):=\IND{\{e\in E(A_0^H) \}}$, and add $\omega_{H}$ to~$\cW$.
Note that $\omega_{H}(E(A_0^H))=(1\pm 3\eps)d_0n^2$ because $A_0^H$ is $(\eps,d_0)$-super-regular by~\ref{item:P2}. 
By~\eqref{eq:w packing size lower bound}, we obtain
\begin{align*}
|\sX_0^\sigma\cap X_0^H|
=\omega_{H}(M(\sigma))
&=(1\pm5\eps)(1-\eps^{2/3})\frac{d_0n^2}{d_0n}\pm n^\eps,
\end{align*}
and thus,
\begin{align}\label{eq:packing size}
(1-\eps^{2/3}/2)n
\geq
|\sX_0^\sigma\cap X_0^H|
\geq
(1-\heps)n.
\end{align}

We first prove~\ref{emb lem:general edge testers}, as we can use this for establishing~\ref{emb lem:sr}.

\begin{step}\label{step:edge testers}
Preparation for checking~\ref{emb lem:general edge testers}\nopagebreak
\end{step}

We will even show that~\ref{emb lem:general edge testers} holds for edge testers in $\cW_{edge}\cup\cW_{edge}'$, where $\cW_{edge}'$ is a set of suitable edge testers satisfying~\ref{item:P edge tester} that we will explicitly specify in Step~\ref{step:sr} when establishing~\ref{emb lem:sr}.
Throughout Steps~\ref{step:edge testers}--\ref{step:conclusion II 3} let $(\omega,\omega_\iota,J,\JX, \JV,\bc,\pP)\in \cW_{edge}\cup\cW_{edge}'$ be fixed.
That is, we fix an index set $I\subseteq \VR$, $J\subseteq I$, disjoint sets $\JX, \JV\subseteq J$, and let $I_{r_0}:=(I\cap[r]_0)\sm J$, $I_r:=(I\cap[r])\sm J$, $\JXV:=\JX\cup \JV$, and we fix 
{$\omega\colon E(\cA_{I_{r_0}})\sqcup \bigsqcup_{j\in J}E(\cB_j)\to[0,s]$}, $\omega_\iota\colon \sX_{\sqcup I}\to[0,s]$, $\bc=\{c_i\}_{i\in I}\in V_{\sqcup I}$, and $\pP=(\bpA,\bppA,\bpB,\bppB)\in(\bN_0^{r_\circ+r+1})^4$ with $\bp^Z=(p_i^Z)_{i\in\VR},\bp^{Z,2nd}=(p_i^{Z,2nd})_{i\in\VR}\in\bN_0^{r_\circ+r+1}$ for $Z\in\{A,B\}$, and we may assume that $I\cap\rR\neq\emptyset$ for some $\rR\in E(R)$ with $0\in\rR$ because otherwise we do not update the weight of the edge tester. 
Recall that the statement~\ref{emb lem:general edge testers} concerns the weight of the edge tester $(\omega^{new},\omega_\iota,J,\JX, \JV,\bc,\pP)$ defined as in Definition~\ref{def:general edge tester} with respect to $(\omega_\iota,J,\JX, \JV,\bc,\pP)$, {$(\phi_\circ^-\cup\sigma,\phi_\circ\cup\sigma^+)$}, $\cA^{new}$ and $\cB^{new}$.

We will consider three different cases depending on whether $0\in (V(R)\sm I)\cup (J\sm \JXV)$, $0\in I\sm J$, and $0\in \JXV$.
Even though we proceed similarly in each of these cases, the effects on~\ref{emb lem:general edge testers} are quite different in each scenario as we try to illude in the following.
Recall that~\ref{emb lem:general edge testers} ensures that~\ref{item:P edge tester} is also satisfied for the updated candidacy graphs. 
If $0\in (V(R)\sm I)\cup (J\sm \JXV)$,
we have to update the density factors whereas the magnitude of ${\omega_\iota(\sX_{\sqcup {I}})}/{n^{|(I\cap-[r_\circ])\sm J|}}$ in~\ref{item:P edge tester} equals the magnitude of ${\omega_\iota(\sX_{\sqcup {I}})}/{n^{|(I\cap-[r_\circ]_0)\sm J|}}$ in~\ref{emb lem:general edge testers}. 
In contrast, if $0\in I\sm J$, we additionally have to ensure that the magnitude of ${\omega_\iota(\sX_{\sqcup {I}})}/{n^{|(I\cap-[r_\circ])\sm J|}}$ in~\ref{item:P edge tester} will be updated by a factor of $n^{-1}$ to obtain the the magnitude of ${\omega_\iota(\sX_{\sqcup {I}})}/{n^{|(I\cap-[r_\circ]_0)\sm J|}}$ in~\ref{emb lem:general edge testers}.
If $0\in \JXV$ the magnitudes are again equal, but besides updating the densities we additionally have to consider that $0\in \JXV$ and thus~\ref{item:P edge tester} will potentially be updated by the factor
$(\IND\{\JXV\cap-[r_\circ]_0=\emptyset \}\pm\eps'^2)=0\pm\eps'^2$.\COMMENT{In this case, we use the conclusion for the weight of the edge tester as an upper bound.}

\medskip
We collect some common notation that will be used to establish~\ref{emb lem:general edge testers}.
Recall that the centres~$\bc$ are fixed and for all $H\in\cH$ and $\aA\bA=\{\aA,\bA\}\in E(\cA_{I_{r_0}})\sqcup \textstyle\bigsqcup_{j\in J}E(\cB_j)$ with $\bA=\{b_j\}_{j\in J}$, we have $\omega(\aA\bA)>0$ only if $\{c_i\}_{i\in I_{r_0}\cup J}=(\aA\cup \bA)\cap V_{\cup {(I_{r_0}\cup J)}}$ by Definition~\ref{def:general edge tester} of an edge tester in~\eqref{eq:general edge tester}.
{(Recall that we allow to treat $\bA=\{b_j\}_{j\in J}\in\bigsqcup_{j\in J}E(\cB_j) $ as $b_{\cup J}$.)}
To that end, for all $H\in\cH$ and $\aA\bA\in E(A_{I_{r_0}}^H)\sqcup \textstyle\bigsqcup_{j\in J}E(B_j^H)$ with $\{c_i\}_{i\in I_{r_0}\cup J}=(\aA\cup\bA)\cap V_{\cup {(I_{r_0}\cup J)}}$, let $\by_{\aA\bA}=\{y_i\}_{i\in I_{r_0}\cup J}:=(\aA\cup\bA)\cap X^H_{\cup {(I_{r_0}\cup J)}}$. 
Our overall strategy in all three cases is to define for vertices $x$ in some set $\sX_0^{{\aA\bA}}\subseteq X_0^H$ a target set $T_{x,{\aA\bA}}$ of suitable images for $x$ such that if all 
$x\in \sX_0^{{\aA\bA}}$ are embedded into $T_{x,{\aA\bA}}$, then $\ab$ (or $\ab\sm\{c_0,y_0\}\}$ if $0\in I\sm J$) is an element in $E(A_{I_{r}}^{H,new})\sqcup \textstyle\bigsqcup_{j\in J}E(B_j^{H,new})$.
Hence for all $H\in\cH$, $x\in X_0^H$ and $\aA\bA\in E(A_{I_{r_0}}^H)\sqcup \textstyle\bigsqcup_{j\in J}E(B_j^H)$ with $\omega(\aA\bA)>0$, we define the following sets.
We give more motivation for these definitions in the subsequent paragraph.
{For an edge $\eH\in E(H)$, let $\rR_\eH\in E(R)$ be such that $\eH\in X_{\cup\rR_\eH}^H$.}
\begin{align*}
\cS_{x,{\aA\bA},\bar{J}}&:=
\big\{\phi_\circ(\eH)\cup \{c_i\}_{i\in \rR_\eH\cap I_r} \colon \eH\in E_{x,\{y_i\}_{i\in I_r}}
\big\};
\\
\cS_{x,{\aA\bA},J}&:=
\big\{S=\phi_\circ(\eH\sm \{y_j\})\cup \{c_i\}_{i\in \rR_\eH\cap(I_r\cup \{j\})} \colon 
|S|=k-1, \eH\in E_{x,\by_{\aA\bA}}, j\in J\sm\{0\}, y_j\in\eH  
\big\};
\\
V_{x,\ab}&:=V_0\cap N_{G_A}(\cS_{x,{\aA\bA},\bar{J}})
\cap N_{G_B}(\cS_{x,{\aA\bA},J});
\\T_{x,{\aA\bA}}&:=V_{x,\ab}\cap N_{A_0^H}(x).
\end{align*}
That is, $\cS_{x,{\aA\bA},\bar{J}}$ and $\cS_{x,{\aA\bA},J}$ are sets of $(k-1)$-sets.
In general, these $(k-1)$-sets consist of the image $\phi_\circ(\eH)=\phi_\circ(\eH\cap\cX_\circ)$ of an $H_+$-edge $\eH\in E_{x,\by_{\aA\bA}}$ together with the centres corresponding to the clusters that $\eH$ intersects. 
The set $\cS_{x,{\aA\bA},\bar{J}}$ contains all $(k-1)$-sets that only intersect with clusters of $I_{r}$, whereas $\cS_{x,{\aA\bA},J}$ contains $(k-1)$-sets that intersect with a cluster of $J\sm\{0\}$.
Consequently, $T_{x,{\aA\bA}}$ is the intersection of the $A_0^H$-neighbourhood of $x$ in~$V_0$ with the common neighbourhood $V_{x,\ab}$ in $G_A$ and $G_B$  of all these $(k-1)$-sets in $\cS_{x,{\aA\bA},\bar{J}}$ and $\cS_{x,{\aA\bA},J}$ (see also Figures~\ref{fig:target sets 0notinI} and~\ref{fig:target sets 0inI}).
Note that $\cS_{x,{\aA\bA},\bar{J}}\cup \cS_{x,{\aA\bA},J}=\emptyset$ if $E_{x,\by_{\aA\bA}}=\emptyset$. 
{Further, since $\omega(\ab)>0$ and it is required for the edge tester~$\omega$ that~$\phi_\circ$ does not map any vertices $\{y_j\}_{j\in J\cap-[r_\circ]}$ onto its centres by~\ref{item:B not on centers}, we have that~$\cS_{x,{\aA\bA},\bar{J}}$ and~$\cS_{x,{\aA\bA},J}$ are disjoint sets.} 
We estimate the sizes of $T_{x,\ab}$ and $V_{x,\ab}$ in Step~\ref{step:target sets}.
Further, for $\ab=\{\aA,\bA\}\in E(\cA_{I_{r_0}})\sqcup \textstyle\bigsqcup_{j\in J}E(\cB_j)$, let $x_\aA:=\aA\cap X_0^H$, $x_\bA:=\bA\cap X_0^H$ (note that $x_\aA$ or $x_\bA$ might be empty), and
\begin{align}
\label{eq:def X_0^y} 
\sX_0^{{\aA\bA}}&:= 
\left\{x\in X_0^H\sm\{x_\aA\}\colon |\cS_{x,{\aA\bA},\bar{J}}\cup \cS_{x,\aA\bA,J}|\geq 1 \right\}.
\end{align}

\begin{substep}
Weight functions to establish~\ref{emb lem:general edge testers}
\end{substep}

We emphasize again that the general strategy for establishing~\ref{emb lem:general edge testers} is to define tuple weight functions for the edges between the vertices $x\in \sX_0^\ab$ and their corresponding target sets $T_{x,\ab}$, which we will do in this step depending on the three cases whether $0\in (V(R)\sm I)\cup (J\sm \JXV)$, $0\in I\sm J$, and $0\in \JXV$.

For all $H\in\cH$, $\ab=\{\aA,\bA\}\in E(A_{I_{r_0}}^H)\sqcup \textstyle\bigsqcup_{j\in J}E(B_j^H)$ with $\omega(\ab)>0$, and $\sX_0^\ab$ as defined in~\eqref{eq:def X_0^y}, we make the following definition. (For notational convenience, we treat $\{\emptyset\}$ as $\emptyset$ in the following definition.)
\begin{align}\label{eq:def E_ab}
E_{{\aA\bA}}:=
\bigg\{&
\{\aA\cap(X_0^H\cup V_0)\}\cup
\{e_x\}_{x\in \sX_0^{{\aA\bA}} } \in \binom{E(A_0^H)}{\IND\{0\in I\sm J\}+|\sX_0^{{\aA\bA}}|}^\cl \colon
\\\nonumber&
 e_x\in E\left(A_0^H[\{x\},T_{x,{\aA\bA}}]\right) \text{ for all $x\in \sX_0^{{\aA\bA}}$}  
 \bigg\}.
\end{align}
Let us explain the definition of $E_\ab$.
If $0\in (V(R)\sm I)\cup (J\sm \JXV)$, then $\aA\cap (X_0^H\cup V_0)=\emptyset$. 
Thus, $E_\ab$ is the set of clean $|\sX_0^\ab|$-tuples of edges in $E(A_0^H)$ between a vertex $x\in\sX_0^\ab$ and its target set $T_{x,\ab}$.
If $0\in I\sm J$, then $\aA\cap (X_0^H\cup V_0)=\{x_\aA,c_0\}$.\COMMENT{Because $\ab$ is such that $\omega(\ab)>0$ and thus $\aA$ must contain the centre $c_0$.} 
(Recall that $x_\aA=\aA\cap X_0^H$.)
Thus, $E_\ab$ is the set of clean $(1+|\sX_0^\ab|)$-tuples of edges in $E(A_0^H)$ where we additionally require that the tuple contains the edge $x_\aA c_0$. 
If $0\in \JXV$, we will not make use of the definition of $E_\ab$. 

\medskip
With $E_\ab$ we can define the following weight function $\omega_{\aA\bA}\colon
\binom{E(\sA_0)}{\IND\{0\in I\sm J\}+|\sX_0^{\aA\bA}|}\to[0,s]$ 
by 
\begin{align*}
\omega_{\aA\bA}(\be):=\omega({\aA\bA})\cdot\IND{\{\be\in E_{\aA\bA}\}}.
\end{align*}
The motivation behind this is the following observation for the two cases $0\in (V(R)\sm I)\cup (J\sm \JXV)$, and $0\in I\sm J$.
{We claim that} for $\aA\bA=\{\aA,\bA\}\in E(A_{I_{r_0}}^H)\sqcup\textstyle\bigsqcup_{j\in J} E(B_j^H)$ with $\omega(\aA\bA)>0$, if $\omega_{\aA\bA}(M(\sigma))>0$ {and $\sigma^+(x_\bA)\neq c_0$ if $0\in J$},\COMMENT{Note that the weight can take values in $[0,s]$.} then {$\ab^{new}:=\{\aA\sm(\aA\cap(X_0^H\cup V_0)),\bA\}\in E(A_{I_{r}}^{H,new})\sqcup \textstyle\bigsqcup_{j\in J}E(B_j^{H,new})$.
To see that this is true, note that if $\omega_{\aA\bA}(M(\sigma))>0$, then the definition of the target sets $T_{x,\ab}$ implies that~\eqref{eq:condition updating} of Definition~\ref{def:candidacy graph} for the updated candidacy graphs $A^{H,new}_{I_r}=A^{H_+}_{I_{r}}(\phi_\circ\cup\sigma^+)$ and {$B_j^{H,new}=B_j^{H_+^j}(\phi_\circ\cup\sigma^+)$} is satisfied.
Hence in this case, for the edge tester~$\omega^{new}$ as defined in the statement of~\ref{emb lem:general edge testers}, we obtain by~\eqref{eq:general edge tester} of Definition~\ref{def:general edge tester} of~$\omega^{new}$ that $\omega^{new}\big(\ab^{new}\big)=\omega(\aA\bA)$ {requiring that $\sigma^+(x_\bA)\neq c_0$ if $0\in J$ so that \ref{item:B not on centers} of Definition~\ref{def:general edge tester} is satisfied}.
Note that $\ab^{new}=\ab$ if $0\in (V(R)\sm I)\cup (J\sm \JXV)$.}

In order to ensure that $\sigma^+(x_\bA)\neq c_0$ if $0\in J$, we apply an inclusion-exclusion principle and introduce another weight function $\omega_\ab^-$ that accounts for the weight in the case that $0\in J$ and \mbox{$\sigma^+(x_\bA)= c_0$}.
To that end, similarly as in~\eqref{eq:def E_ab}, 
for all $H\in\cH$, $\ab=\{\aA,\bA\}\in E(A_{I_{r_0}}^H)\sqcup \textstyle\bigsqcup_{j\in J}E(B_j^H)$ with $\omega(\ab)>0$, and $\sX_0^\ab$ as defined in~\eqref{eq:def X_0^y}, we define the following set of~edge~tuples\looseness=-1
\begin{align*}
E_{{\aA\bA}}^-:=
\bigg\{&
\{\bA\cap(X_0^H\cup V_0)\}\cup
\{e_x\}_{x\in \sX_0^{{\aA\bA}}\sm\{x_\bA\} } \in \binom{E(A_0^H)}{\IND\{0\in  J\}+|\sX_0^{{\aA\bA}}\sm\{x_\bA\}|}^\cl \colon
\\\nonumber&
 e_x\in E\left(A_0^H[\{x\},T_{x,{\aA\bA}}]\right) \text{ for all $x\in \sX_0^{{\aA\bA}}\sm\{x_\bA\}$}  
 \bigg\}.
\end{align*}
Analogously to $\omega_\ab$, we define the weight function $\omega_{\aA\bA}^-\colon
\binom{E(\sA_0)}{\IND\{0\in J\}+|\sX_0^{\aA\bA}\sm\{x_\bA\}|}\to[0,s]$ 
by 
$\omega_{\aA\bA}^-(\be):=\omega({\aA\bA})\cdot\IND{\{\be\in E_{\aA\bA}^-\}}$.

The size of the tuple weight functions depends on the cardinality of~$\sX_0^\ab$.
To that end, let 
\begin{align}\label{eq:def b_max}
b_{\max}:=\max\{|\sX_0^{\aA\bA}|\colon \aA\bA\in E(\cA_{I_{r_0}})\sqcup \textstyle\bigsqcup_{j\in J}E(\cB_j), \omega(\aA\bA)>0 \},
\end{align}
and we will group the tuple functions for all  $\ab$ with $\omega(\aA\bA)>0$ into all possible $(\IND\{0\in I\sm J\}+b)$-tuple weight functions for $b\in[b_{\max}]_0$.
To this end, for each $b\in[b_{\max}]_0$, we set 
\begin{align}
\label{eq:def N_b}
\Omega_b&:=\left\{  \aA\bA\in E(\cA_{I_{r_0}})\sqcup \textstyle\bigsqcup_{j\in J}E(\cB_j) \colon \omega(\aA\bA)>0, |\sX_0^{\aA\bA}|=b \right\}, 
\\\nonumber
\Omega_b^-&:=\left\{  \aA\bA\in E(\cA_{I_{r_0}})\sqcup \textstyle\bigsqcup_{j\in J}E(\cB_j) \colon \omega(\aA\bA)>0, |\sX_0^{\aA\bA}\sm\{x_\bA\}|=b \right\}, 
\end{align}
as well as 
\begin{align*}
\omega_{b}:=\sum_{\aA\bA\in \Omega_b}\omega_{\aA\bA}, \quad\text{ and }\quad 
\omega_{b}^-:=\sum_{\aA\bA\in \Omega_b^-}\omega_{\aA\bA}^-.
\end{align*}
In the two cases when $0\in (V(R)\sm I)\cup (J\sm \JXV)$, and $0\in I\sm J$, we will see that $\sum_{b\in[b_{\max}]_0}\omega_b(M(\sigma))$ is the major contribution to $\omega^{new}\big(
E(\sA_{I_r}^{new})\sqcup \bigsqcup_{j\in J}E(\cB_j^{new})
\big)$.
However, since $\sX_0^\sigma$ is a proper subset of $\sX_0$, we additionally need to consider those $\ab\in E(\cA_{I_{r_0}})\sqcup \bigsqcup_{j\in J}E(\cB_j)$ for which a relevant vertex in~$\sX_0^\ab$ has not been embedded by~$\sigma$ as this might also contribute to the weight of $\omega^{new}$.
This is also the case when $0\in \JXV$, because then we require that either $x_\bA=\bA\cap X_0^H$ is not embedded or no $H$-vertex is mapped onto the centre $c_0$.

We collect some notation.
For all $\aA\bA\in E(\cA_{I_{r_0}})\sqcup \textstyle\bigsqcup_{j\in J}E(\cB_j)$ with $\omega(\aA\bA)>0$, let $H_\ab\in\cH$ be such that $\aA\bA\in E(A_{I_{r_0}}^{H_\ab})\sqcup \bigsqcup_{j\in J}E(B_j^{H_\ab})$, and let $\sX_0^{\ab,\bars}:=\sX_0^\ab\sm\sX_0^\sigma$.
Recall that $\JX, \JV\subseteq J$ are disjoint sets and $\JXV=\JX\cup \JV$.
Let $z_\bA:= x_\bA$ if $0\in \JX$, and let $z_\bA:= c_0$ if $0\in \JV$.
We can now define the following set of edge tuples that we describe in detail below.
For $b\in[b_{\max}]_0$, $\ell\in[b]_0$, $m_A,m_B\in[\Delta]_0$, let
\begin{align}\label{eq:def Gamma substep 0 notin I}
\Gamma_b(\ell,m_A,m_B)
:=\Big\{&
\ab\in \Omega_b\colon
|\sX_0^{\ab,\bars}|=\ell, 
\sum_{x\in\sX_0^{\ab,\bars}}|\cS_{x,\ab,\bar{J}}|=m_A,
\sum_{x\in\sX_0^{\ab,\bars}}|\cS_{x,\ab,J}|=m_B,
\\&\nonumber
\sigma(x)\in T_{x,\ab} \text{ for all $x\in \sX_0^{\ab}\cap\sX_0^\sigma$, }
\text{if $0\in I\sm J$ then $\sigma(x_\aA)=\{c_0\}$,}
\\&\nonumber
\text{if $0\in \JXV$ then $z_\bA\in(\sX_0\sm\sX_0^\sigma)\cup(V_0\sm \sigma(X_0^{H_\ab}\cap \sX_0^\sigma))$}
\Big\}.
\end{align}
That is, $\Gamma_b(\ell,m_A,m_B)$ is the set of edges $\ab\in \Omega_b$ such that there exists an $\ell$-set $\sX_0^{\ab,\bars}$ of vertices in $\sX_0^\ab$ which are not embedded by $\sigma$, and these $\ell$ vertices in $\sX_0^{\ab,\bars}$ contribute $m_A$ and $m_B$ many $(k-1)$-sets, and all remaining $b-\ell$ vertices $x\in\sX_0^{\ab}\cap\sX_0^\sigma$ are embedded onto their target set~$T_{x,\ab}$.
Additionally, if $0\in I\sm J$, then we require that $x_\aA$ is mapped onto $c_0$ by $\sigma$, and if $0\in \JX$, then we require that $x_\bA$ is not embedded by $\sigma$, and if $0\in \JV$, then we require that no vertex of~$X_0^{H_\ab}$ is mapped onto $c_0$ by $\sigma$.
Further, let
\begin{align}\label{eq:def Gamma hit substep 0 notin I}
\Gamma_b^{hit}(\ell,m_A,m_B):=
\Big\{
\ab\in\Gamma_b(\ell,m_A,m_B)
\colon& 
\sigma^+(x)\in V_{x,\ab} \text{ for all $x\in\sX_0^{\ab,\bars}$,}
\\\nonumber&
\text{if $x_\bA\in \sX_0^{\ab,\bars}$ then $\sigma^+(x_\bA)\neq c_0$}
\Big\}.
\end{align}
That is, $\Gamma_b^{hit}(\ell,m_A,m_B)\subseteq \Gamma_b(\ell,m_A,m_B)$ contains those edges $\ab\in \Omega_b$, where the not embedded vertices in $\sX_0^\ab$ are nevertheless mapped onto their target set $V_{x,\ab}$ by the random cluster-injective extension $\sigma^+$ of $\sigma$, and $\sigma^+$ does not map $x_\bA$ onto $c_0$.
Thus, in such a case the weight of $\ab$ will `accidentally' be taken into account in addition to the `real' contribution given by~$\sigma$ (compare also with~\eqref{eq:new total weight} below).

Crucially note that in the two cases when $0\in (V(R)\sm I)\cup (J\sm \JXV)$ and $0\in I\sm J$, we know for $\ab=\{\aA,\bA\}\in \Omega_b$ that $\ab^{new}=\{\aA\sm(\aA\cap(X_0^H\cup V_0)),\bA\}\in E(\cA_{I_{r}}^{new})\sqcup \bigsqcup_{j\in J}E(\cB_j^{new})$ 
only if 
$\sigma^+(x_\bA)\neq c_0$ if $0\in J$
and
either $\omega_{\aA\bA}(M(\sigma))=\omega(\aA\bA)>0$ or
$\ab\in\Gamma_b^{hit}(\ell,m_A,m_B)$ for some $\ell\in[b]$, $m_A,m_B\in[\Delta]_0$.\COMMENT{In this case it must hold that $b>0$.} (Recall that $\ab^{new}=\ab$ if $0\in (V(R)\sm I)\cup (J\sm \JXV)$.)
This holds by~\eqref{eq:condition updating} of Definition~\ref{def:candidacy graph} and since we defined $\cA_{I_{r}}^{new}$ and $\cB_j^{new}$ in~\eqref{eq:def Z^H,new} as updated candidacy graphs with respect to $\phi_\circ\cup\sigma^+$.
(Note that since we choose $\sigma^+$ as a `dummy' enlargement, we do not require that $\sigma^+(x)\in N_{\cA_0}(x)$, which is the reason why $\sigma^+(x)\in V_{x,{\aA\bA}}$ in~\eqref{eq:def Gamma hit substep 0 notin I} instead of $\sigma^+(x)\in T_{x,{\aA\bA}}$.)
Hence, for the two cases when $0\in (V(R)\sm I)\cup (J\sm \JXV)$ and $0\in I\sm J$, we make the following key observation:\looseness=-1
\begin{align}
\nonumber
&\omega^{new}\left(
E(\sA_{I_r}^{new})\sqcup \textstyle\bigsqcup_{j\in J}E(\cB_j^{new})
\right)
\\\label{eq:new total weight}&= 
\sum_{b\in[b_{\max}]_0}\left(\omega_{b}(M(\sigma))
-
{\IND\{0\in J\sm\JXV \}\omega_{b}^-(M(\sigma))}
\right)
+ 
\sum_{\substack{b\in[b_{\max}],\ell\in[b],\\m_A,m_B\in[\Delta]_0}}
\omega\left(\Gamma_b^{hit}(\ell,m_A,m_B)\right).
\end{align}

In the case that $0\in \JXV$, it suffices in view of the statement to establish an upper bound for 
$\omega^{new}\left(
E(\sA_{I_r}^{new})\sqcup \textstyle\bigsqcup_{j\in J}E(\cB_j^{new})
\right)$.
Similarly as in~\eqref{eq:new total weight}, we have in this case for $\ab\in\Omega_b$ and $b\in[b_{\max}]_0$ that 
$\aA\bA\in E(\cA_{I_{r}}^{new})\sqcup \textstyle\bigsqcup_{j\in J}E(\cB_j^{new})$ 
only if 
$\aA\bA\in\Gamma_b^{hit}(\ell,m_A,m_B)$ for some $\ell\in[b]_0$, $m_A,m_B\in[\Delta]_0$.
This holds by~\eqref{eq:condition updating} of Definition~\ref{def:candidacy graph} and since we defined $\cA_{I_{r}}^{new}$ and $\cB_j^{new}$ in~\eqref{eq:def Z^H,new} as updated candidacy graphs with respect to $\phi_\circ\cup\sigma^+$.
Hence, for the case that $0\in \JXV$, we make the following key observation:
\begin{align}
\label{eq:new total weight 0 in J}
\omega^{new}\left(
E(\sA_{I_r}^{new})\sqcup \textstyle\bigsqcup_{j\in J}E(\cB_j^{new})
\right)
\leq
\sum_{\substack{b\in[b_{\max}]_0,\ell\in[b]_0,\\m_A,m_B\in[\Delta]_0}}
\omega\left(\Gamma_b^{hit}(\ell,m_A,m_B)\right).
\end{align}

In Steps~\ref{step:main contribution}--\ref{step:bound w(Gamma^hit)}, we will estimate the weight of the contributing terms in~\eqref{eq:new total weight} and~\eqref{eq:new total weight 0 in J}.
To do so, we first determine the sizes of the target sets $T_{x,\ab}$ and $V_{x,\ab}$ in the next step.

\begin{substep}\label{step:target sets}
Size of the target sets $T_{x,\ab}$ and $V_{x,\ab}$\nopagebreak
\end{substep}

Let ${b\in[b_{\max}]}$ and $\ab\in \Omega_b$ be fixed.
We observe that
\begin{align}
\label{eq:size V_x,ab}
|V_{x,{\aA\bA}}|=(1\pm 3\eps)d_A^{|\cS_{x,{\aA\bA},\bar{J}}|}d_B^{|\cS_{x,{\aA\bA},J}|}n,
\quad\text{ and }\quad
|T_{x,{\aA\bA}}|=(1\pm 3\eps)d_A^{|\cS_{x,{\aA\bA},\bar{J}}|}d_B^{|\cS_{x,{\aA\bA},J}|}d_0n,
\end{align}
for $x\in\sX_0^\ab$,
where we used~\ref{item:P1} and that $A_0^H$ is $(\eps,d_0)$-super-regular and $(\eps,q)$-well-intersecting for each $H\in\cH$. 
For an illustration of the sets $\cS_{x,{\aA\bA},\bar{J}}$ and $T_{x,{\aA\bA}}$ in the case that $0\notin I$ and $J=\emptyset$, see Figure~\ref{fig:target sets 0notinI}.

\tikzstyle{edge} = [fill,opacity=.3,fill opacity=0,line cap=round, line join=round, line width=14pt]
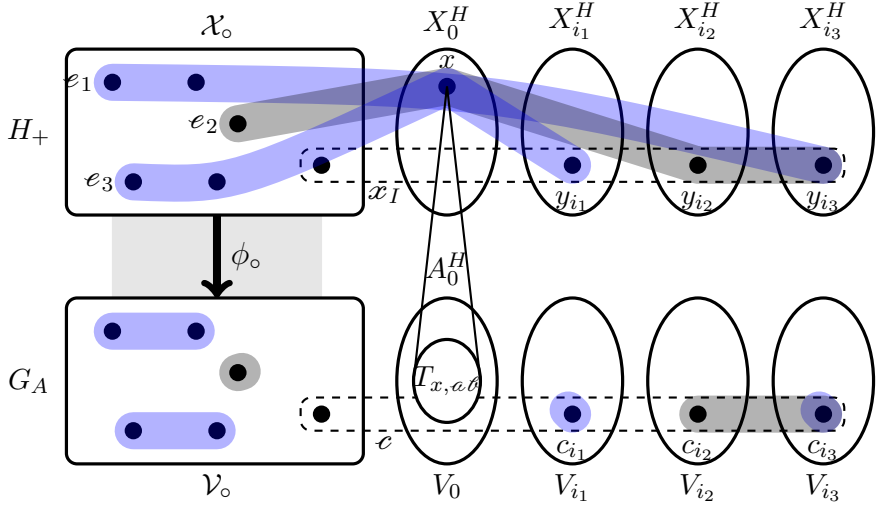
\begin{figure}[htb!]
\vspace{-.2cm}
	\begin{center}
		\begin{tikzpicture}[scale = 0.55,every text node part/.style={align=center}]
		\def\ver{0.1} 

		\draw (-1,0) node {$H_+$};			 
		\draw (-1,-6) node {$G_A$};

		\draw[fill=black!10, draw=black!10]
		(1,-2) rectangle (6,-4);

		\draw[line width=3pt,->]
		(3.5,-2) -- node[anchor=west] {\large$\phi_\circ$} (3.5,-4);

		\draw[very thick, fill=white,rounded corners] (-.1,2) rectangle (7,-2);
		\node[above] at (3.5,2) {$\cX_\circ$};
		
		\draw[very thick, fill=white] 
		(9,0) ellipse (1.2 and 2)
		(12,0) ellipse (1.2 and 2)
		(15,0) ellipse (1.2 and 2)
		(18,0) ellipse (1.2 and 2);
		\node[above] at (9,2) {$X_0^H$};
		\node[above] at (12,2) {$X_{i_1}^H$};
		\node[above] at (15,2) {$X_{i_2}^H$};
		\node[above] at (18,2) {$X_{i_3}^H$};

		\draw[very thick, fill=white,rounded corners] (-.1,-4) rectangle (7,-8);
		\node[below] at (3.5,-8) {$\cV_\circ$};

		\draw[very thick, fill=white] 
		(9,-6) ellipse (1.2 and 2)
		(12,-6) ellipse (1.2 and 2)
		(15,-6) ellipse (1.2 and 2)
		(18,-6) ellipse (1.2 and 2);
		\node[below] at (9,-8) {$V_0$};
		\node[below] at (12,-8) {$V_{i_1}$};
		\node[below] at (15,-8) {$V_{i_2}$};
		\node[below] at (18,-8) {$V_{i_3}$};

		\node[hvertex,label={[label distance=.1cm]270:$y_{i_1}$}] (y1) at (12,-.8){};
		\coordinate (y2) at (15,-.8); 		
		\node[hvertex,label={[label distance=.1cm]270:$y_{i_2}$}] (Y2) at (y2){};
		\node[hvertex,label={[label distance=.1cm]270:$y_{i_3}$}] (y3) at (18,-.8){};
		\coordinate (x) at (9,1.1); 
		\node[hvertex,label=above:$x$] (X) at (x){};

		\node[hvertex] (xI) at (6,-.8){};
		\draw[thick, fill=white,rounded corners,dashed,fill opacity=0] (5.5,-.4) rectangle (18.5,-1.2);
		\draw (7.5,-1.5) node {$\xX_I$};

		\draw[thick, fill=white,rounded corners,dashed,fill opacity=0] (5.5,-6.4) rectangle (18.5,-7.2);
		\draw (7.5,-7.5) node {$\bc$};

		\node[hvertex] (cI) at (6,-6.8){};
		\coordinate (c1) at (12,-6.8); 
		\node[hvertex,label={[label distance=.1cm]270:$c_{i_1}$}] (C1) at (c1){};
		\coordinate (c2) at (15,-6.8); 
		\coordinate (c3) at (18,-6.8); 				
		\node[hvertex,label={[label distance=.1cm]270:$c_{i_2}$}] (C2) at (c2){};
		\node[hvertex,label={[label distance=.1cm]270:$c_{i_3}$}] (C3) at (c3){};
		\coordinate (e0) at (1,1.2); 
		\coordinate (e1) at (3,1.2); 
		\node[hvertex,label=left:$\eH_1$] (E0) at (e0){};
		\node[hvertex] (E1) at (e1){};

		\coordinate (f0) at (4,.2); 
		\node[hvertex,label=left:$\eH_2$] (F0) at (f0){};

		\coordinate (g0) at (1.5,-1.2); 
		\node[hvertex,label=left:$\eH_3$] (G0) at (g0){};
		\node[hvertex] (g1) at (3.5,-1.2){};

		\coordinate (ve0) at (1,-4.8); 
		\coordinate (ve1) at (3,-4.8); 		
		\node[hvertex] (Ve0) at (ve0){};
		\node[hvertex] (Ve1) at (ve1){};

		\coordinate (vf0) at (4,-5.8); 
		\node[hvertex] (Vf0) at (vf0){};

		\coordinate (vg0) at (1.5,-7.2); 
		\coordinate (vg1) at (3.5,-7.2); 		
		\node[hvertex] (Vg0) at (vg0){};
		\node[hvertex] (Vg1) at (vg1){};

		\draw[edge,blue] (18,-.8)  .. controls (10,1.1) .. (e0);
		\draw[edge] (18,-.8) -- (y2) --(x) --(f0);
		\draw[edge,blue] (12,-.8)--(x) .. controls (3.9,-1.3) ..(g0);

		\draw[edge,blue] (ve0)--(ve1);
		\draw[edge] (vf0)--(4.1,-5.77);
		\draw[edge,blue] (vg0)--(vg1);

		\draw[edge,blue] (c1)--(11.9,-6.7);
		\draw[edge,blue] (c3)--(17.9,-6.7);
		\draw[edge] (c2)--(c3);

		\draw[very thick, fill=white] 
		(9,-6) ellipse (.8 and 1);

		\draw[thick] (x)--(8.2,-6);
		\draw[thick] (x)--(9.8,-6);

		\draw (9,-6) node {$T_{x,\ab}$};
		\draw (9,-3.3) node {$A_0^H$};

		\end{tikzpicture}
	\end{center}
	  \captionsetup{width=.95\linewidth}
\caption{This illustrates the case that $0\notin I$ and for simplicity $J=\emptyset$.
That is, we consider the edge $\ab=\{c_{i_1},c_{i_2},c_{i_3},y_{i_1},y_{i_2},y_{i_3} \}\in E(A_{I_{r_0}}^H)$ with $\by_{\ab}=\{y_{i_1},y_{i_2},y_{i_3} \}\subseteq\xX_I$.
Note that the edges $\eH_1,\eH_2,\eH_3$ in~$H_+$ belong to $E_{x,\by_{\ab}}$. 
For the set $\cS_{x,\ab,\bar{J}}$, we have $\cS_{x,\ab,\bar{J}}=\{ \phi_\circ(\eH_1)\cup\{c_{i_1}\}, \phi_\circ(\eH_2)\cup\{c_{i_2},c_{i_3}\},
\phi_\circ(\eH_3)\cup\{c_{i_3}\} \}.$
Accordingly, $T_{x,\ab}$ is the intersection in $V_0$ of the $G_A$-neighbourhoods of these $(k-1)$-sets in $\cS_{x,\ab,\bar{J}}$ and the neighbourhood of $x$ in~$A_0^H$.
Note that the blue edges $\eH_1$ and $\eH_3$ in~$H_+$ satisfy that $|\eH_1\cap\by_\ab|=|\eH_3\cap\by_\ab|=1$, and thus they do not account for the \fst-pattern $\bpA(\xX_I,\emptyset)$ of $\xX_I$ by~\eqref{eq:def pattern} of Definition~\ref{def:pattern}. 
By considering all possible $x\in X_0^H$, there are in total $b_{I_r}=\sum_{i\in I_r}b_i$ many such blue edges in~$H_+$.
Further, note that the grey edge~$\eH_2$ satisfies that $|\eH_2\cap\by_\ab|=2$ and thus, $\eH_2$ belongs to $H$ by~\eqref{eq:H_+ intersection cardinality 1} and accounts for the \fst-pattern $\bpA(\xX_I,\emptyset)$ of $\xX_I$.
Again, by considering all possible $x\in X_0^H$, there are in total $\bpA(\xX_I,\emptyset)_0$ many such grey edges in~$H_+$.\looseness=-2
}
	\label{fig:target sets 0notinI}
\end{figure}
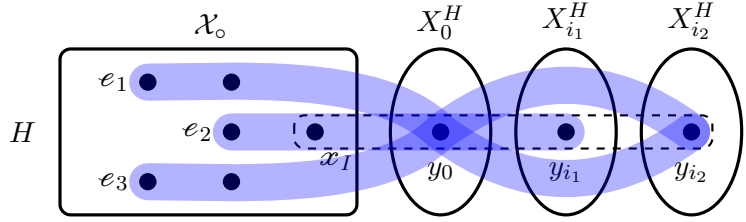
\begin{figure}[htb!]
\vspace{-.5cm}
	\begin{center}
		\begin{tikzpicture}[scale = 0.55,every text node part/.style={align=center}]
		\def\ver{0.1} 

		\draw (-1,0) node {$H$};			 

%

		\draw[very thick, fill=white,rounded corners] (-.1,2) rectangle (7,-2);
		\node[above] at (3.5,2) {$\cX_\circ$};
		
		\draw[very thick, fill=white] 
		(9,0) ellipse (1.2 and 2)
		(12,0) ellipse (1.2 and 2)
		(15,0) ellipse (1.2 and 2);
		\node[above] at (9,2) {$X_0^H$};
		\node[above] at (12,2) {$X_{i_1}^H$};
		\node[above] at (15,2) {$X_{i_2}^H$};



		\coordinate (y1) at (12,0); 
		\node[hvertex,label={[label distance=.15cm]270:$y_{0}$}] (y0) at (9,0){};
		\node[hvertex,label={[label distance=.15cm]270:$y_{i_1}$}] (Y1) at (y1){};
		\node[hvertex,label={[label distance=.15cm]270:$y_{i_2}$}] (y2) at (15,0){};

		\node[hvertex] (xI) at (6,0){};
		\draw[thick, fill=white,rounded corners,dashed,fill opacity=0] (5.5,.4) rectangle (15.5,-.4);
		\draw (6.5,-.7) node {$\xX_I$};


		\coordinate (e0) at (2,1.2); 
		\node[hvertex,label=left:$\eH_1$] (E0) at (e0){};
		\node[hvertex] (e1) at (4,1.2){};

		\coordinate (f0) at (4,0); 
		\node[hvertex,label=left:$\eH_2$] (F0) at (f0){};

		\coordinate (g0) at (2,-1.2); 
		\node[hvertex,label=left:$\eH_3$] (G0) at (g0){};
		\node[hvertex] (g1) at (4,-1.2){};
		
		\draw[edge,blue]
		(15,0)  .. controls (13,-1.5) and (11,-1.5) .. (9,0) .. controls (7.5,1.5) and (e1).. (e0);
		\draw[edge,blue] (y1)--(f0);
		\draw[edge,blue] (15,0) .. controls (13,1.5) and (11,1.5) .. (9,0) .. controls (7.5,-1.5) and (g1).. (g0);

%
%

%
%
%
%
%

		\end{tikzpicture}
	\end{center}
  \captionsetup{width=.95\linewidth}
\caption{This illustrates the case that $0\in I$
and $\by_{\ab}=\{y_{0},y_{i_1},y_{i_2} \}\subseteq\xX_I$, and we assume that $I_{r_0}=\{0,i_1,i_2\}$ and $J=\emptyset$. Note that $y_0=x_\aA$.
By Definition~\ref{def:global edge tester} of an edge tester we require that $\xX_I$ lies in an $H$-edge if we assigned positive weight to the tuple $\xX_I$. 
Hence, by the construction of $H_+$ (see~\eqref{eq:pattern in H_+}), we have that $\eH_1,\eH_2,\eH_3$ are edges in $H$.
Note that the edges $\eH_1,\eH_2,\eH_3$ in~$H$ belong to $E_{y_0,y_{i_1}}\cup E_{y_0,y_{i_2}}$. 
By definition, $y_0=x_\aA\notin \sX_0^\ab$ and we thus do not consider the possible target set $T_{y_0,\ab}$ because $y_0$ has to be mapped onto the centre $c_0$ in~$G$.
Hence, the edges $\eH_1, \eH_2, \eH_3$ do not account for $\sum_{x\in\sX_0^\ab}|\cS_{x,\ab,\bar{J}}|$. 
Since $\by_{\ab}\subseteq\xX_I$, there are $\bppA(\xX_I,\emptyset)_0$ many such blue edges by the definition of the \scd-pattern $\bppA(\xX_I,\emptyset)$ in Definition~\ref{def:pattern}.
}
	\label{fig:target sets 0inI}
\end{figure}

For $b_{I_r}=\sum_{i\in I_r}b_i$ and $b_{J}=\sum_{j\in J}b_j$ as defined in~\eqref{eq:def b_i}, we claim that
\begin{align}\label{eq:size S_x,ab,barJ}
\sum_{x\in\sX_0^{{\aA\bA}}}|\cS_{x,{\aA\bA},\bar{J}}|=b_{I_r}+p_0^A-p_0^{A,2nd},
\quad\text{ and }\quad
\sum_{x\in\sX_0^{{\aA\bA}}}|\cS_{x,\aA\bA,J}|=b_{J}+p_0^B-p_0^{B,2nd}.
\end{align}

{We establish the first equation in~\eqref{eq:size S_x,ab,barJ}; the second one then follows similarly.
At this part of the proof it is crucial to refresh Definition~\ref{def:pattern} because we make use of all the details of the pattern definitions. 
Further, recall that $\by_{\aA\bA}=\{y_i\}_{i\in I_{r_0}\cup J}:=(\aA\cup\bA)\cap X^H_{\cup {(I_{r_0}\cup J)}}$ and thus $x_\aA=y_0$ if $0\in I_{r_0}$; otherwise, $x_\aA=\emptyset$.
{Note that in order to compute $|\cS_{x,{\aA\bA},\bar{J}}|$ it is equivalent to count $|E_{x,\{y_i\}_{i\in I_r}}|$.}
By definition of~$\cH_+^i$, we have that $|\bigcup_{x\in\sX_0}E_{x,y_i}|=b_i$ for all $i\in (I_{r}\cup J)\sm\{0\}$ (see~\eqref{eq:size E_x,y}). 
That is, there are $b_{I_r}=\sum_{i\in I_r}b_i$ edges $\eH\in \bigcup_{x\in\sX_0}E_{x,\{y_i \}_{i\in I_r} }$ with $|\eH\cap \{y_i \}_{i\in I_r}|=1$.
{Out of these $b_{I_r}$ edges, we claim that $p_0^{A,2nd}$ edges $\eH$ satisfy that $x_\aA\in \eH$, that is, $\sum_{i\in I_r}|E_{x_\aA,y_i}|= p_0^{A,2nd}$.} 
For an illustration of these edges in $\bigcup_{i\in I_r}E_{x_\aA,y_i}$, see Figure~\ref{fig:target sets 0inI}.
Indeed, since $\omega(\ab)>0$, we have by Definition~\ref{def:global edge tester} of an edge tester that $\by_{\ab}\subseteq \xX_I$ for some $\xX_I\in\sX_{\sqcup I}$ with $\xX_I\in E_\cH(\pP,I,J)$.
{Thus, by Definition~\ref{def:e_H,p} of $E_{\cH}(\pP,I,J)$, we have that $\bpA(\xX_I,J)=\bpA=(p_i^A)_{i\in \VR}$ and $\bppA(\xX_I,J)=\bppA=(p_i^{A,2nd})_{i\in \VR}$.}
Hence, by the definition of a \scd-pattern in Definition~\ref{def:pattern} and because $\cH_+$ is constructed such that each subset $\{x_\aA,y_i\}$ only lies in proper edges of $\cH$ due to~\ref{item:pattern in H_+} and~\eqref{eq:pattern in H_+},\COMMENT{Here we crucially use that $\xX_I$ is contained in an $\cH$-edge.} 
we have $p_0^{A,2nd}=\sum_{i\in I_r}|E_{x_\aA,y_i}|$.
By the definition of $\sX_0^\ab$ in~\eqref{eq:def X_0^y}  which excludes $x_\aA$, this accounts for the term `$b_{I_r}-p_0^{A,2nd}$' in~\eqref{eq:size S_x,ab,barJ}.
{It is worth pointing out that $p_0^{A,2nd}=p_0^{B,2nd}=0$ if $0\notin I_{r_0}$ because then $x_\aA=\emptyset$.}}

{Further, we claim that there are $p_0^A$ edges $\eH\in \bigcup_{x\in\sX_0}E_{x,\{y_i \}_{i\in I_r}}$ with $|\eH\cap \{y_i \}_{i\in I_r}|\geq2$ but $x_\aA\notin \eH$.
Since $\by_{\ab}\subseteq \xX_I$ for some $\xX_I\in\sX_{\sqcup I}$ with $\xX_I\in E_{\cH}(\pP,I,J)$, we have that $\bpA(\xX_I, J)=\bpA=(p_i^A)_{i\in \VR}$.
Hence, by the definition of a \fst-pattern in Definition~\ref{def:pattern} and by~\eqref{eq:H_+ intersection cardinality 1}, there are~$p_0^A$ edges $\eH\in \bigcup_{x\in\sX_0}E_{x,\{y_i \}_{i\in I_r}}$ with $|\eH\cap \{y_i \}_{i\in I_r}|\geq2$ because of the last two conditions in~\eqref{eq:def pattern}, and all of these edges $\eH$ satisfy that $x_\aA\notin \eH$ due to the condition `$(\fH\cap X^H_{\ell})\sm\{x_i\}_{i\in I\sm J}\neq\emptyset$' in~\eqref{eq:def pattern}.
Altogether, this implies~\eqref{eq:size S_x,ab,barJ}.}
\TASK{
In order to compute $|\cS_{x,{\aA\bA},J}|$, we count tuples $(\eH,j)$ for edges $\eH$ and indices $j\in J\sm\{0\}$ that satisfy the conditions in the definition of $\cS_{x,{\aA\bA},J}$. 
\textbf{This counting equivalence is possible since we require~\ref{item:B not on centers} for the edge tester, and because an $\cH_+$-edge $\eH$ satisfies~\ref{item:no 2-overlap}.}
By definition of~$\cH_+^j$, we claim for all $j\in J\sm\{0\}$ that 
\begin{align*}
\sum_{x\in\sX_0}\left|
\left\{
\eH\in E_{x,\by_{\ab}}\colon
y_j\in\eH,
|\{y_j\}\cup(\eH\cap X_{\cup (I\cap[r])}^H)|= 1
\right\}
\right|
=
b_j
\end{align*}
Observe that each such edge $\eH$ yields a $(k-1)$-set $S\in\cS_{x,\ab,J}$ for $j\in J\sm\{0\}$, and these are the edges where only one `free' vertex remains. (We could also write $\eH\in E_{x,\{y_j\}_{j\in J\sm\{0\}}}$ in the display above.)
In particular, these edges satisfy the conditions in the definition of $\cS_{x,{\aA\bA},J}$. 
That is, there are $b_{J}=\sum_{j\in J\sm\{0\}}b_j$ such edges $\eH\in \bigcup_{x\in\sX_0}E_{x,\by_{\ab}}$ and entries $j\in J\sm\{0\}$ (some edges might account multiple times to different $j,j'$) with one `free' vertex.
Let us call these combinations/tuples of an edge $\eH\in \bigcup_{x\in\sX_0}E_{x,\by_{\ab}}$ and an entry $j\in J\sm\{0\}$ with one `free' vertex a \defn{co-edge $(\eH,j)$}.
Out of these $b_{J}$ co-edges, we claim that $p_0^{B,2nd}$ co-edges $(\eH,j)$ satisfy that $x_\aA\in \eH$.
Indeed, since $\omega(\ab)>0$, we have by Definition~\ref{def:global edge tester} of an edge tester that $\by_{\ab}\subseteq \xX_I$ for some $\xX_I\in\sX_{\sqcup I}$ with $\xX_I\in E_\cH(\pP,I,J)$.
{Thus, by Definition~\ref{def:e_H,p} of $E_{\cH}(\pP,I,J)$, we have that $\bpB(\xX_I,J)=\bpB=(p_i^B)_{i\in \VR}$ and $\bppB(\xX_I,J)=\bppB=(p_i^{B,2nd})_{i\in \VR}$.}
Hence, by the definition of a \scd-pattern in Definition~\ref{def:pattern} and because $\cH_+$ is constructed such that each subset $\{x_\aA,y_j\}$ only lies in proper edges of $\cH$ due to~\ref{item:pattern in H_+} and~\eqref{eq:pattern in H_+},\COMMENT{Here we crucially use that $\xX_I$ is contained in an $\cH$-edge.} 
we have $p_0^{B,2nd}$ co-edges that contain $x_\aA$.
To that end, note that every co-edge $(\eH,j)$ that contains $x_\aA$ accounts to $p_0^{B,2nd}$, where we imagine that the $j$-th vertex is copied to the right.
Vice versa every edge that accounts to $p_0^{B,2nd}$ corresponds to a co-edge that contains $x_\aA$.
By the definition of $\sX_0^\ab$ in~\eqref{eq:def X_0^y}  which excludes $x_\aA$, this accounts for the term `$b_{J}-p_0^{B,2nd}$' in~\eqref{eq:size S_x,ab,barJ}.
{It is worth pointing out that $p_0^{A,2nd}=p_0^{B,2nd}=0$ if $0\notin I_{r_0}$ because then $x_\aA=\emptyset$.}
\\
Further, we claim that
\begin{align*}
\sum_{j\in J\sm\{0\},x\in\sX_0^\ab}
\left|
\left\{
\eH\in E_{x,\by_{\ab}}\colon
y_j\in\eH,
|\{y_j\}\cup(\eH\cap X_{\cup (J\cap[r])}^H)|=1\leq |\eH\cap X_{\cup I_r}^H)|
\right\}
\right|
=
p_0^B.
\end{align*}
This follows by~\eqref{eq:def pattern} of Definition~\ref{def:pattern}. Altogether, this implies~\eqref{eq:size S_x,ab,barJ}.
}

\medskip
Hence, for $\ab\in \Omega_b$, we have by~\eqref{eq:size V_x,ab} and~\eqref{eq:size S_x,ab,barJ} that
\begin{align}\label{eq:prod T_x,ab}
\prod_{x\in\sX_0^{\aA\bA}}|T_{x,\aA\bA}|=(1\pm\heps)d_A^{b_{I_r}+p_0^A-p_0^{A,2nd}}d_B^{b_{J}+p_0^B-p_0^{B,2nd}}d_0^bn^b.
\end{align}

\begin{step}\label{step:main contribution}
Estimating $\sum_{b\in[b_{\max}]_0}\omega_{b}(M(\sigma))$ in~\eqref{eq:new total weight}\nopagebreak
\end{step}

In this step we estimate the contribution of the first term $\sum_{b\in[b_{\max}]_0}\omega_{b}(M(\sigma))$ in~\eqref{eq:new total weight}.
Throughout this step, let us consider the case that $0\in (V(R)\sm I)\cup (J\sm \JXV)$.
At the end of the step we explain how the estimate of $\sum_{b\in[b_{\max}]_0}\omega_{b}(M(\sigma))$ changes if $0\in I\sm J$.
Note, if $0\in (V(R)\sm I)\cup (J\sm \JXV)$, then $\omega_0$ is the empty function and thus,
$\sum_{b\in[b_{\max}]_0}\omega_{b}(M(\sigma))=\sum_{b\in[b_{\max}]}\omega_{b}(M(\sigma))$. (That is, $b=0$ is only relevant if $0\in I\sm J$.) 

We first consider $\omega_{b}(M(\sigma))$ for some ${b\in[b_{\max}]}$.
By~\eqref{eq:prod T_x,ab} and the definition of $\omega_b$, we have
\begin{align}\label{eq:weight w_b}
\omega_b(E(\sA_0))=(1\pm2\heps)
d_A^{b_{I_r}+p_0^A-p_0^{A,2nd}}
d_B^{b_{J}+p_0^B-p_0^{B,2nd}}
 d_0^b  n^b\sum_{\aA\bA\in \Omega_b}\omega(\aA\bA).
\end{align}

We verify that $\omega_b$ satisfies~\eqref{eq:normcondition}.
For all $\{e_1,\ldots,e_{b'}\}\in\binom{E(\sA_0)}{b'}, b'\in[b]$, the number of edges $\{e_{b'+1},\ldots,e_{b} \}$ such that $\be=\{e_1,\ldots,e_{b}\}\in\binom{E(\sA_0)}{b}$ with $\omega_b(\be)>0$ is at most $\Delta n^{b-b'}$ (recall that we have chosen $\Delta$ such that $\eps'\ll1/\Delta\ll1/t\ll 1/k,1/r,1/r_\circ,1/s$),\COMMENT{For fixed $\{e_1,\ldots,e_{b'}\}\in\binom{E(\sA_0)}{b'}$ there are at most $\Delta_R^{2b'}$ possible vertices in $\sX_0$ that can be centres of a star weight function for $\{e_{b'+1},\ldots,e_{b} \}$. Out of these possible centres we choose $b-b'$, that is, we have at most
$$\binom{\Delta_R^{2b'}}{b-b'}\leq \left(\Delta_R^{2b'}\right)^{b-b'}\leq \left(\Delta_R^{2\Delta_R}\right)^{\Delta_R}\leq \Delta_R^{2\Delta_R^2}.$$ 
We take $\Delta_R^{3\Delta_R^2}$ to account for different sizes of the clusters, that is, $|V_0|=(1\pm\eps)n$.}
implying that 
$\normV{\omega_b}{b'}\leq \Delta^2 n^{b-b'}\leq n^{b-b'+\eps^2}.$
Hence, by adding $\omega_b$ to $\cW$,~\eqref{eq:w packing size} implies that
\begin{align}
\nonumber\omega_b(M(\sigma))
&\stackrel{\hphantom{\eqref{eq:weight w_b}}}{=}
(1\pm\heps)\frac{\omega_b(E(\sA_0))}{(d_0n)^b}\pm n^\eps
\\&\label{eq:bound w_p(M)}
\stackrel{\eqref{eq:weight w_b}}{=}
(1\pm\heps^{1/2})
d_A^{b_{I_r}+p_0^A-p_0^{A,2nd}}
d_B^{b_{J}+p_0^B-p_0^{B,2nd}}
\sum_{\aA\bA\in \Omega_b}\omega(\aA\bA)\pm  n^\eps.
\end{align}
Finally, observe that
\begin{align*}
\sum_{b\in[b_{\max}]}\sum_{\aA\bA\in \Omega_b} \omega(\aA\bA)
=
\sum_{\aA\bA\in E(\cA_{I_{r_0}})\sqcup \bigsqcup_{j\in J}E(\cB_j)}
\omega(\aA\bA)
=
\omega\left(
E(\cA_{I_{r_0}})
\sqcup \textstyle\bigsqcup_{j\in J}E(\cB_j)
\right),
\end{align*}
and thus~\eqref{eq:bound w_p(M)} implies that
\begin{align}\label{eq:bound sum w_p(M)}
\sum_{b\in[b_{\max}]}\omega_b(M(\sigma))
=(1\pm\heps^{1/3})
d_A^{b_{I_r}+p_0^A-p_0^{A,2nd}}
d_B^{b_{J}+p_0^B-p_0^{B,2nd}}
\omega\left(
E(\cA_{I_{r_0}})\sqcup \textstyle\bigsqcup_{j\in J}E(\cB_j)
\right)
\pm n^{2\eps},
\end{align}
which is the desired estimate of $\sum_{b\in[b_{\max}]}\omega_b(M(\sigma))=\sum_{b\in[b_{\max}]_0}\omega_b(M(\sigma))$ in the case that $0\in (V(R)\sm I)\cup (J\sm \JXV)$.

\medskip
Let us now assume that $0\in I\sm J$ and we explain how the estimate of $\sum_{b\in[b_{\max}]_0}\omega_b(M(\sigma))$ changes.
(Note that we allow for $b=0$.)
The intuition is that we additionally require that $x_\aA$ is mapped onto $c_0$ which we would expect to happen in an idealized random setting with probability roughly $(d_0n)^{-1}$.
In fact,~\eqref{eq:weight w_b} is still true in the case that $0\in I\sm J$, but note that $\omega_b$ is now a $(1+b)$-tuple weight function which yields an additional factor of $(d_0n)^{-1}$ in~\eqref{eq:bound w_p(M)} and thus also in~\eqref{eq:bound sum w_p(M)}.
Hence, we obtain~\eqref{eq:bound sum w_p(M)} with an additional factor of $(d_0n)^{-1}$ as the desired estimate in the case that $0\in I\sm J$.

\begin{step}\label{step:main -contribution}
Estimating $\sum_{b\in[b_{\max}]_0}\omega_{b}^-(M(\sigma))$ in~\eqref{eq:new total weight} if $0 \in J\sm\JXV$
\end{step}

In this step we establish an upper bound on the contribution of the minuend $\sum_{b\in[b_{\max}]_0}\omega_{b}^-(M(\sigma))$ in~\eqref{eq:new total weight} and therefore, suppose that $0 \in J\sm\JXV$.

We first consider $\omega_{b}^-(M(\sigma))$ for some ${b\in[b_{\max}]_0}$.
It suffices to establish only a rough upper bound which follows directly by the definition of $\omega_b^-$:
\begin{align}\label{eq:weight w_b^-}
\omega_b^-(E(\sA_0))\leq
2n^b
\sum_{\ab\in \Omega_b^-}\omega(\ab).
\end{align}
It is easy to verify that $\omega_b^-$ satisfies~\eqref{eq:normcondition}.
Further, note that $\omega_b^-$ is a $(b+1)$-tuple weight function-
Hence, by adding $\omega_b$ to $\cW$,~\eqref{eq:w packing size} implies that
\begin{align}\label{eq:bound w_p(M)^-}
\omega_b^-(M(\sigma))
\leq
(1+\heps)\frac{\omega_b^-(E(\sA_0))}{(d_0n)^{b+1}}+ n^\eps
\stackrel{\eqref{eq:weight w_b^-}}{\leq}
n^{\heps-1}
\sum_{\ab\in \Omega_b^-}\omega(\ab)+  n^\eps.
\end{align}
Finally, observe that
\begin{align*}
\sum_{b\in[b_{\max}]_0}\sum_{\aA\bA\in \Omega_b^-} \omega(\aA\bA)
\leq
\omega\left(
E(\cA_{I_{r_0}})
\sqcup \textstyle\bigsqcup_{j\in J}E(\cB_j)
\right),
\end{align*}
and thus~\eqref{eq:bound w_p(M)^-} implies that
\begin{align}\label{eq:bound sum w_p(M)^-}
\sum_{b\in[b_{\max}]_0}\omega_b^-(M(\sigma))
\leq 
n^{\heps-1}
\omega\left(
E(\cA_{I_{r_0}})\sqcup \textstyle\bigsqcup_{j\in J}E(\cB_j)
\right)
+ n^{2\eps}.
\end{align}
Hence, combining~\eqref{eq:bound sum w_p(M)} and~\eqref{eq:bound sum w_p(M)^-} in the case that $0\in J\sm\JXV$, yields that
\begin{align}\label{eq:bound totalsum w_p(M)}
\sum_{b\in[b_{\max}]_0}\left(\omega_b(M(\sigma))
-
\omega_b^-(M(\sigma))
\right)
=&(1\pm2\heps^{1/3})
d_A^{b_{I_r}+p_0^A-p_0^{A,2nd}}
d_B^{b_{J}+p_0^B-p_0^{B,2nd}}
\\\nonumber&
\omega\left(
E(\cA_{I_{r_0}})\sqcup \textstyle\bigsqcup_{j\in J}E(\cB_j)
\right)
\pm 2n^{2\eps}.
\end{align}

\begin{step}\label{step:bound w(Gamma)}
Estimating $\omega(\Gamma_b(\ell,m_A,m_B))$\nopagebreak
\end{step}

In this step we derive an upper bound for
$\omega(\Gamma_b(\ell,m_A,m_B))$ for fixed $b\in[b_{\max}]_0$, $\ell\in[b]_0$, $m_A,m_B\in[\Delta]_0$, and $\Gamma_b(\ell,m_A,m_B)\supseteq \Gamma_b^{hit}(\ell,m_A,m_B)$ as defined in~\eqref{eq:def Gamma substep 0 notin I}. 
We will use this bound in the subsequent Step~\ref{step:bound w(Gamma^hit)} to derive an upper bound for $\omega(\Gamma_b^{hit}(\ell,m_A,m_B))$ as in~\eqref{eq:new total weight} and~\eqref{eq:new total weight 0 in J}.
Throughout this step, let us again consider the case that $0\in (V(R)\sm I)\cup (J\sm \JXV)$, and thus $b,\ell>0$. At the end of the step, we explain how the estimate changes if $0\in I\sm J$ or $0\in \JXV$.

Our general strategy is based on the following inclusion-exclusion principle.
For every $\ab\in \Omega_b$, we estimate the
\begin{align}\label{eq:estimate b-l hit}
\begin{minipage}[c]{0.8\textwidth}\em
$\omega$-weight $\omega(\ab)$ with tuples $\xX\in \binom{\sX_0^\ab}{b}$ such that $b-\ell$ vertices $x$ of $\xX$ are mapped onto their target set $T_{x,\ab}$ and the remaining $\ell$ vertices $x_1,\ldots,x_\ell$ of $\xX$ satisfy 
\begin{enumerate}
\item[{\rm ($\ast$)}]
\begin{center} $\quad\displaystyle\sum_{i\in[\ell]}|\cS_{x_i,\ab,\bar{J}}|=m_A,\quad\quad
\sum_{i\in[\ell]}|\cS_{x_i,\ab,J}|=m_B,$
\end{center}
\end{enumerate}
\end{minipage}
\\\label{eq:estimate b-l hit, l embedded}
\begin{minipage}[c]{0.8\textwidth}\em
$\omega$-weight $\omega(\ab)$ as in~\eqref{eq:estimate b-l hit} with tuples $\xX$ where we additionally require that the remaining $\ell$ vertices vertices $x_1,\ldots,x_\ell$ of $\xX$ are embedded by $\sigma$ and  satisfy~{\rm ($\ast$)}.
\end{minipage}
\ignorespacesafterend
\end{align}
Now, subtracting~\eqref{eq:estimate b-l hit, l embedded} from~\eqref{eq:estimate b-l hit} yields the 
\begin{align}\label{eq:estimate substract}
\begin{minipage}[c]{0.8\textwidth}\em
$\omega$-weight $\omega(\ab)$ as in~\eqref{eq:estimate b-l hit} with tuples $\xX$ where at least one of the remaining $\ell$ vertices $x_1,\ldots,x_\ell$ of $\xX$ is not embedded by $\sigma$, and $x_1,\ldots,x_\ell$ satisfy {\rm ($\ast$)}.
\end{minipage}
\ignorespacesafterend
\end{align}
Hence, summing over the $\omega$-weight as in~\eqref{eq:estimate substract} for all $\ab\in\Omega_b$ yields an upper bound for \linebreak$\omega(\Gamma_b(\ell,m_A,m_B))$ as defined in~\eqref{eq:def Gamma substep 0 notin I} when $0\in (V(R)\sm I)\cup (J\sm \JXV)$.

\medskip
First, in order to estimate~\eqref{eq:estimate b-l hit} and~\eqref{eq:estimate b-l hit, l embedded}, we define the following sets of tuples of edges in~$\cA_0$.
For all $H\in\cH$ and $\aA\bA\in E(A_{I_{r_0}}^H)\sqcup \textstyle\bigsqcup_{j\in J}E(B_j^H)$ with $\ab\in \Omega_b$, let
\begin{align}
\allowdisplaybreaks
\nonumber
\sX_0^{\ab,\bars,m_A,m_B}&:=
\bigg\{
\sX \in\binom{\sX_0^\ab}{\ell}
\colon \sum_{x\in\sX}|\cS_{x,\ab,\bar{J}}|=m_A,
\sum_{x\in\sX}|\cS_{x,\ab,J}|=m_B
\bigg\};
\allowdisplaybreaks
\\\label{eq:def E_a,b b-l hit}
E_{\ab}^\eqref{eq:estimate b-l hit}
&:=
\bigcup_{\sX \in\sX_0^{\ab,\bars,m_A,m_B}}
\bigg\{ 
\{e_x\}_{x\in  \sX_0^\ab\sm\sX}
\in\binom{E(A_0^H)}{b-\ell}^\cl \colon 
\\\nonumber&\hphantom{:=
\bigcup_{\sX\in\sX_0^{\ab,\bars,m_A,m_B}}
\bigg\{ }
e_x\in E(A_0^H[\{x\},T_{x,\ab}]) \text{ for all $x\in\sX_0^\ab\sm\sX$}
\bigg\};
\allowdisplaybreaks
\\
\label{eq:def E_a,b b-l hit, l embedded}
E_{\ab}^\eqref{eq:estimate b-l hit, l embedded}
&:=
\bigcup_{\sX\in\sX_0^{\ab,\bars,m_A,m_B}}
\bigg\{ 
\{e_x\}_{x\in\sX_0^{\ab}}
\in\binom{E(A_0^H)}{b}^\cl \colon 
\\\nonumber&\hphantom{:=
\bigg\{ }
x\in e_x \text{ for all $x\in\sX_0^\ab$, }
e_x\in E(A_0^H[\{x\},T_{x,\ab}]) \text{ for all $x\in\sX_0^\ab\sm\sX$}
\bigg\}.
\end{align}
Note that the edges in $E_{\ab}^\eqref{eq:estimate b-l hit}$ and $E_{\ab}^\eqref{eq:estimate b-l hit, l embedded}$ correspond to the described situations in~\eqref{eq:estimate b-l hit} and~\eqref{eq:estimate b-l hit, l embedded}, respectively.
We define weight functions $\omega_{\ab}^\eqref{eq:estimate b-l hit}\colon \binom{E(A_0^H)}{b-\ell}\to[0,s]$ and $\omega_{\ab}^\eqref{eq:estimate b-l hit, l embedded}\colon \binom{E(A_0^H)}{b}\to[0,s]$ by 
\begin{align*}
\omega_{\ab}^\eqref{eq:estimate b-l hit}(\be):=\IND\{\be\in E_{\ab}^\eqref{eq:estimate b-l hit}\}
\cdot\omega(\ab),
\quad\text{ and }\quad
\omega_{\ab}^\eqref{eq:estimate b-l hit, l embedded}(\be):=\IND\{\be\in E_{\ab}^\eqref{eq:estimate b-l hit, l embedded}\}
\cdot\omega(\ab).
\end{align*}
Let 
\begin{align*}
\omega_{\Gamma}^\eqref{eq:estimate b-l hit}(\be):=\sum_{{\aA\bA}\in \Omega_b}\omega_{\ab}^\eqref{eq:estimate b-l hit}(\be),
\quad\text{ and }\quad
\omega_{\Gamma}^\eqref{eq:estimate b-l hit, l embedded}(\be):=\sum_{{\aA\bA}\in \Omega_b}\omega_{\ab}^\eqref{eq:estimate b-l hit, l embedded}(\be).
\end{align*}

We estimate $\omega_{\Gamma}^\eqref{eq:estimate b-l hit}(E(\cA_0))$ and  $\omega_{\Gamma}^\eqref{eq:estimate b-l hit, l embedded}(E(\cA_0))$.
By~\eqref{eq:prod T_x,ab} and the definition of  $E_{\ab}^\eqref{eq:estimate b-l hit}$ and $ E_{\ab}^\eqref{eq:estimate b-l hit, l embedded}$ in~\eqref{eq:def E_a,b b-l hit} and~\eqref{eq:def E_a,b b-l hit, l embedded}, respectively, we obtain
\begin{align}
\label{eq:weight w_b b-l hit}
\omega_{\Gamma}^\eqref{eq:estimate b-l hit}(E(\cA_0))
&=(1\pm2\heps)
d_A^{b_{I_r}+p_0^A-p_0^{A,2nd}-m_A}
d_B^{b_{J}+p_0^B-p_0^{B,2nd}-m_B}
(d_0n)^{b-\ell}\sum_{\aA\bA\in \Omega_b}\omega(\aA\bA);
\\\label{eq:weight w_b b-l hit, l embedded}
\omega_{\Gamma}^\eqref{eq:estimate b-l hit, l embedded}(E(\cA_0))
&=(1\pm2\heps)
d_A^{b_{I_r}+p_0^A-p_0^{A,2nd}-m_A}
d_B^{b_{J}+p_0^B-p_0^{B,2nd}-m_B}
 (d_0n)^b\sum_{\aA\bA\in \Omega_b}\omega(\aA\bA).
\end{align}

Again, we can add $\omega_{\Gamma}^\eqref{eq:estimate b-l hit}$ and $\omega_{\Gamma}^\eqref{eq:estimate b-l hit, l embedded}$ to $\cW$ and employ property~\eqref{eq:w packing size}.
This yields that
\begin{align}
\nonumber
\omega_{\Gamma}^\eqref{eq:estimate b-l hit}(M(\sigma))
&
\stackrel{\hphantom{\eqref{eq:weight w_b b-l hit}}}{=}
(1\pm\heps)\frac{\omega_{\Gamma}^\eqref{eq:estimate b-l hit}(E(\cA_0))}{(d_0n)^{b-\ell}}\pm n^\eps
\\&\label{eq:M-weight b-l hit}
\stackrel{\eqref{eq:weight w_b b-l hit}}{=}
(1\pm4\heps)
d_A^{b_{I_r}+p_0^A-p_0^{A,2nd}-m_A}
d_B^{b_{J}+p_0^B-p_0^{B,2nd}-m_B}
\sum_{\aA\bA\in \Omega_b}\omega(\aA\bA)\pm n^\eps
\end{align}
and
\begin{align}
\nonumber
\omega_{\Gamma}^\eqref{eq:estimate b-l hit, l embedded}(M(\sigma))
&
\stackrel{\hphantom{\eqref{eq:weight w_b b-l hit, l embedded}}}{=}
(1\pm\heps)\frac{\omega_{\Gamma}^\eqref{eq:estimate b-l hit, l embedded}(E(\cA_0))}{(d_0n)^b}\pm n^\eps
\\&\label{eq:M-weight b-l hit, l embedded}
\stackrel{\eqref{eq:weight w_b b-l hit, l embedded}}{=}
(1\pm4\heps)
d_A^{b_{I_r}+p_0^A-p_0^{A,2nd}-m_A}
d_B^{b_{J}+p_0^B-p_0^{B,2nd}-m_B}
\sum_{\aA\bA\in \Omega_b}\omega(\aA\bA)\pm n^\eps.
\end{align}
Finally, as observed in~\eqref{eq:estimate substract}, subtracting~\eqref{eq:M-weight b-l hit, l embedded} from~\eqref{eq:M-weight b-l hit} gives us an upper bound on $\omega(\Gamma_b(\ell,m_A,m_B))$.
We obtain
\begin{align}\nonumber
\omega(\Gamma_b(\ell,m_A,m_B))
&\leq \omega_{\Gamma}^\eqref{eq:estimate b-l hit}(M(\sigma))-\omega_{\Gamma}^\eqref{eq:estimate b-l hit, l embedded}(M(\sigma))
\\&\label{eq:bound w-Gamma_b(l,m_A,m_B)}
\leq
8\heps
d_A^{b_{I_r}+p_0^A-p_0^{A,2nd}-m_A}
d_B^{b_{J}+p_0^B-p_0^{B,2nd}-m_B}
\sum_{\aA\bA\in \Omega_b}\omega(\aA\bA)+ 2n^\eps.
\end{align}

\medskip
Let us now first assume that $0\in I\sm J$ and we explain how the estimate on $\omega(\Gamma_b(\ell,m_A,m_B))$ changes.
If $0\in I\sm J$, then by the definition of $\Gamma_b(\ell,m_A,m_B)$ in~\eqref{eq:def Gamma substep 0 notin I}, we additionally require that $x_\aA$ is mapped onto $c_0$ by $\sigma$ which we would again expect to happen with probability roughly $(d_0n)^{-1}$.
That is, we have to modify the definitions in~\eqref{eq:def E_a,b b-l hit} and~\eqref{eq:def E_a,b b-l hit, l embedded} by additionally adding the edge $x_\aA c_0$ to the tuples. 
Again, the estimates for the total weights in~\eqref{eq:weight w_b b-l hit} and~\eqref{eq:weight w_b b-l hit, l embedded} are still true but we obtain an additional factor of $(d_0n)^{-1}$ in~\eqref{eq:M-weight b-l hit} and~\eqref{eq:M-weight b-l hit, l embedded} as the sizes of the tuple functions increased by $1$.
Thus, we also obtain~\eqref{eq:bound w-Gamma_b(l,m_A,m_B)} with an additional factor of $(d_0n)^{-1}$ which will be our desired estimate in the case that $0\in I\sm J$.

Finally, let us assume that $0\in \JXV$ and we explain how the estimate on $\omega(\Gamma_b(\ell,m_A,m_B))$ changes.
If $0\in \JXV$, then by the definition of $\Gamma_b(\ell,m_A,m_B)$ in~\eqref{eq:def Gamma substep 0 notin I}, we additionally require that either $x_\bA$ is left unembedded by $\sigma$, or no $H^\ab$-vertex is mapped onto $c_0$.
This can be achieved by modifying the definition in~\eqref{eq:def E_a,b b-l hit, l embedded} such that the edge tuples are increased by adding the $A_0^H$-edges $e_{z_\bA}$ such that $z_\bA\in e_{z_\bA}$.\COMMENT{
\begin{align*}
E_{\ab,b+1}
&:=
\bigcup_{\sX\in\sX_0^{\ab,\bars,m_A,m_B}}
\bigg\{ 
\{e_{z_\bA}\}\cup\{e_x\}_{x\in\sX}
\cup
\{e_x\}_{x\in  \sX_0^\ab\sm\sX}
\in\binom{E(A_0^H)}{b+1}^\cl \colon 
\\\nonumber&\hphantom{:=\bigg\{ }
x\in e_x \text{ for all $x\in\sX_0^\ab\cup\{z_\bA\}$, }
e_x\in E(A_0^H[\{x\},T_{x,\ab}]) \text{ for all $x\in\sX_0^\ab\sm\sX$}
\bigg\};
\end{align*}}
This ensures that for $z_\bA\in\{x_\bA,c_0\}$, we either have $x_\bA$ is left unembedded by~$\sigma$, or no $H^\ab$-vertex is mapped onto $c_0$.
The modification adds another factor of $d_0n$ to the total weight in~\eqref{eq:weight w_b b-l hit, l embedded} but also another factor of $(d_0n)^{-1}$ to~\eqref{eq:M-weight b-l hit, l embedded}.
Thus,~\eqref{eq:bound w-Gamma_b(l,m_A,m_B)} will also be our desired estimate in the case that $0\in \JXV$. 

\begin{step}\label{step:bound w(Gamma^hit)}
Estimating $\omega(\Gamma_b^{hit}(\ell,m_A,m_B))$ in~\eqref{eq:new total weight} and~\eqref{eq:new total weight 0 in J}\nopagebreak
\end{step}

We will use the bounds on $\omega(\Gamma_b(\ell,m_A,m_B))$ of Step~\ref{step:bound w(Gamma)} to derive an upper bound for the $\omega$-weight $\omega(\Gamma^{hit}_b(\ell,m_A,m_B))$.
Let us first assume the two cases that $0\in (V(R)\sm I)\cup (J\sm \JXV)$, and $0\in \JXV$, since both cases yield the same bound on $\omega(\Gamma_b(\ell,m_A,m_B))$ in~\eqref{eq:bound w-Gamma_b(l,m_A,m_B)}.
The general idea is that we will obtain additional factors `$d_A^{m_A}$' and `$d_B^{m_B}$' in~\eqref{eq:bound w-Gamma_b(l,m_A,m_B)} when we extend $\sigma$ to $\sigma^+$, that is, when we embed the~$\ell$ unembedded $\cH_+$-neighbours of each $\ab$ contributing to $\omega(\Gamma_b(\ell,m_A,m_B))$. 
Only if these $\ell$ vertices are mapped onto their target set $V_{x,\ab}$, then $\ab$ also contributes to $\omega(\Gamma^{hit}_b(\ell,m_A,m_B))$; that is, if $\sigma^+(x)\in V_{x,\ab}$ {for all $x\in\sX_0^\ab\sm\sX_0^\sigma$}. (Recall the definition of $\Gamma^{hit}_b(\ell,m_A,m_B)$ in~\eqref{eq:def Gamma hit substep 0 notin I}.)
This happens roughly with probability $d_A^{m_A}d_B^{m_B}$.
We proceed with the details.

\medskip
Note that $\Gamma_b^{hit}(\ell,m_A,m_B)=\Gamma_b(\ell,m_A,m_B)$ for $\ell=0$, and thus we may consider fixed
$b\in[b_{\max}]$, $\ell\in[b]$, $m_A,m_B\in[\Delta]_0$ with $m_A+m_B>0$.
Note that we extend $\sigma|_{V(H)}$ to $\sigma^+|_{V(H)}$ for every $H\in\cH$ by choosing a bijective mapping of $X_0^H\sm\sX_0^\sigma$ into $V_0\sm\sigma(X_0^H\cap\sX_0^\sigma)$ uniformly and independently at random. 
To that end, for $\aA\bA\in E(A_{I_{r_0}}^H)\sqcup \bigsqcup_{j\in J} E(B_j^H)$, $H\in\cH$, and $x\in\sX_0^{\ab}$,
let $V_{x,\ab,\bar{\sigma}}:=V_{x,\ab}\sm\sigma(X_0^H\cap\sX_0^\sigma)$.
We first estimate $|V_{x,\ab,\bars}|$.
To that end, let $\aA\bA\in E(A_{I_{r_0}}^H)\sqcup \textstyle\bigsqcup_{j\in J}E(B_j^H)$, $H\in\cH$, and $x\in\sX_0^{\ab}$ be fixed.
For every $v\in V_{x,\ab}$, we define a weight function $\omega_v\colon E(A_0^H)\to\{0,1\}$ by $\omega_v(e):=\IND\{v\in e \}$ and let $\omega_{V_{x,\ab}}:=\sum_{v\in V_{x,\ab}}\omega_v$.
Observe that $\omega_{V_{x,\ab}}(M(\sigma))$ counts the vertices $v\in V_{x,\ab}\sm V_{x,\ab,\bar{\sigma}}$.
Hence,
\begin{align}\label{eq:size V_x,ab,bars}
|V_{x,\ab,\bar{\sigma}}|=|V_{x,\ab}|-\omega_{V_{x,\ab}}(M(\sigma)).
\end{align}
Since $A_0^H$ is $(\eps,d_0)$ super-regular, we have
\begin{align*}
\omega_{V_{x,\ab}}(E(A_0^H))=
(1\pm 3\eps)d_0n|V_{x,\ab}|.
\end{align*}
Adding $\omega_{V_{x,\ab}}$ to $\cW$ and employing~\eqref{eq:w packing size lower bound} yields that
\begin{align*}
\omega_{V_{x,\ab}}(M(\sigma))=
(1\pm\eps)(1-\eps^{2/3})\frac{\omega_{V_{x,\ab}}(E(A_0^H))}{d_0n}\pm n^\eps
=(1\pm 5\eps)(1-\eps^{2/3})|V_{x,\ab}|.
\end{align*}
We conclude that
\begin{align}\label{eq:M-size V_x,ab,bars}
|V_{x,\ab,\bar{\sigma}}|
\stackrel{\eqref{eq:size V_x,ab,bars}}{\leq}
2\eps^{2/3}|V_{x,\ab}|
\stackrel{\eqref{eq:size V_x,ab}}{\leq}
3\eps^{2/3} d_A^{|\cS_{x,{\aA\bA},\bar{J}}|}d_B^{|\cS_{x,{\aA\bA},J}|}n.
\end{align}

For $\ab\in \Gamma_b(\ell,m_A,m_B)$, let $\sX_0^{\ab,\bars}:=\sX_0^\ab\sm\sX_0^\sigma$. By the definition of $\Gamma_b(\ell,m_A,m_B)$ in~\eqref{eq:def Gamma substep 0 notin I}, we have\looseness=-1 
\begin{align*}
\sum_{x\in\sX_0^{\ab,\bars}}|\cS_{x,\ab,\bar{J}}|=m_A,\quad
\sum_{x\in\sX_0^{\ab,\bars}}|\cS_{x,\ab,J}|=m_B.
\end{align*}
Hence, by~\eqref{eq:M-size V_x,ab,bars} and because $|\sX_0^{\ab,\bars}|=\ell$, we obtain
\begin{align}\label{eq:size of target sets sigma^+}
\prod_{x\in\sX_0^{\ab,\bars}}|V_{x,\ab,\bars}|
\leq
d_A^{m_A}d_B^{m_B} 
(3\eps^{2/3}n)^\ell.
\end{align}
By~\eqref{eq:packing size}, we have that $V_0\sm \sigma(X_0^H\cap\sX_0^\sigma)\geq \eps^{2/3}n/3$. 
Now, with~\eqref{eq:size of target sets sigma^+} we obtain that the probability that all $\ell$ vertices $x\in\sX_0^{\ab,\bars}$ are mapped onto their target set $V_{x,\ab,\bars}$ and if $x_\bA\in\sX_0^{\ab,\bars}$ then it is not mapped onto $c_0$ --- that is, the probability that $\ab\in \Gamma_b(\ell,m_A,m_B)$ is also contained in $\Gamma_b^{hit}(\ell,m_A,m_B)$ --- is at most
\begin{align*}
\frac{d_A^{m_A}d_B^{m_B} 
(3\eps^{2/3}n)^\ell}
{\eps^{2/3}n\cdot(\eps^{2/3}n-1)\cdots(\eps^{2/3}n-\ell+1)/3^\ell}
\leq 
10^\ell d_A^{m_A}d_B^{m_B}.
\end{align*}

Finally, we can derive an upper bound for the expected value of $\omega(\Gamma_b^{hit}(\ell,m_A,m_B))$.
\begin{align*}
\expn{\omega\left(\Gamma_b^{hit}(\ell,m_A,m_B)\right)}
&
\stackrel{\hphantom{\eqref{eq:bound w-Gamma_b(l,m_A,m_B)}}}{\leq}
\omega\left(\Gamma_b(\ell,m_A,m_B)\right)10^\ell d_A^{m_A}d_B^{m_B}
\\&
\stackrel{\eqref{eq:bound w-Gamma_b(l,m_A,m_B)}}{\leq}
\eps^{1/4}
d_A^{b_{I_r}+p_0^A-p_0^{A,2nd}}
d_B^{b_{J}+p_0^B-p_0^{B,2nd}}
\sum_{\aA\bA\in \Omega_b}\omega(\aA\bA)+ n^{2\eps}.
\end{align*} 

By using Theorem~\ref{thm:McDiarmid} and a union bound, we can establish concentration with probability, say, at least {$1-\eul^{-n^{\eps}}$.}
Thus, we conclude
\begin{align}\label{eq:bound w(Gamma^hit)}
\omega\left(\Gamma_b^{hit}(\ell,m_A,m_B)\right)
\leq
2\eps^{1/4}
d_A^{b_{I_r}+p_0^A-p_0^{A,2nd}}
d_B^{b_{J}+p_0^B-p_0^{B,2nd}}
\sum_{\aA\bA\in \Omega_b}\omega(\aA\bA)+ 2n^{2\eps}
\end{align} 
for all $b\in[b_{\max}]$, $\ell\in[b]$, $m_A,m_B\in[\Delta]_0$ with $m_A+m_B>0$.
By summing over all values of $b,\ell,m_A,m_B$, we obtain
\begin{align}
&\nonumber
\sum_{\substack{b\in[b_{\max}],\ell\in[b],\\m_A,m_B\in[\Delta]_0}}
\omega\left(\Gamma_b^{hit}(\ell,m_A,m_B)\right)
\\\label{eq:bound sum w(Gamma^hit)}
\stackrel{\eqref{eq:bound w(Gamma^hit)}}{\leq}
&\eps^{1/5}
d_A^{b_{I_r}+p_0^A-p_0^{A,2nd}}
d_B^{b_{J}+p_0^B-p_0^{B,2nd}}
\omega\left(
E(\cA_{I_{r_0}})\sqcup \textstyle\bigsqcup_{j\in J}E(\cB_j)
\right)
+n^{3\eps},
\end{align}
which is the desired estimate for the second term in~\eqref{eq:new total weight} for the case that $0\in (V(R)\sm I)\cup (J\sm \JXV)$, and it is the desired estimate for the right hand side in~\eqref{eq:new total weight 0 in J} for the case that $0\in \JXV$ (where we also allow $b,\ell=0$ in the summation).

In the case that $0\in I\sm J$, we obtain~\eqref{eq:bound sum w(Gamma^hit)} with an additional factor of $(d_0n)^{-1}$ since the estimate on $\omega\left(\Gamma_b(\ell,m_A,m_B)\right)$ from the previous step yields an additional factor of $(d_0n)^{-1}$.

\begin{step}\label{step:conclusion II 1}
Concluding~\ref{emb lem:general edge testers} if $0\in (V(R)\sm I)\cup (J\sm \JXV)$\nopagebreak
\end{step}

Finally, we can establish~\ref{emb lem:general edge testers} if $0\in (V(R)\sm I)\cup (J\sm \JXV)$ by the estimates derived in the Steps~\ref{step:main contribution} and~\ref{step:bound w(Gamma^hit)}.
So, let us assume that $0\in (V(R)\sm I)\cup (J\sm \JXV)$.
Using~\eqref{eq:bound sum w_p(M)} (respectively~\eqref{eq:bound totalsum w_p(M)} if $0\in J\sm\JXV$) and~\eqref{eq:bound sum w(Gamma^hit)} in our key observation~\eqref{eq:new total weight} yields that
\begin{align}
\nonumber
\omega^{new}\left(
E(\sA_{I_r}^{new})\sqcup \textstyle\bigsqcup_{j\in J}E(\cB_j^{new})
\right)
=&
(1\pm3\heps^{1/3})
d_A^{b_{I_r}+p_0^A-p_0^{A,2nd}}
d_B^{b_{J}+p_0^B-p_0^{B,2nd}}
\\&\label{eq:w^new P3 0 notin I}
\omega\left(
E(\cA_{I_{r_0}})\sqcup \textstyle\bigsqcup_{j\in J}E(\cB_j)
\right)
\pm2n^{3\eps}.
\end{align}
We will use~\ref{item:P edge tester} for $\omega\left(
E(\cA_{I_{r_0}})\sqcup \textstyle\bigsqcup_{j\in J}E(\cB_j)
\right)$ in order to obtain~\ref{emb lem:general edge testers} from~\eqref{eq:w^new P3 0 notin I}.
For $0\notin I$ or $0\in J\sm \JXV$, we have that 
\begin{align}\label{eq:d_A^new}
d_A^{b_{I_r}}\prod_{i\in I_{r_0}}d_{i,A}=\prod_{i\in I_r}d_A^{b_i}d_{i,A}=\prod_{i\in I_r}d_{i,A}^{new}
\end{align} for $d_{i,A}^{new}$ defined as in~\eqref{eq:def b_i} because $I_{r_0}=I_r$, and similarly
\begin{align}\label{eq:d_B^new}
d_B^{b_{J}}\prod_{j\in J}d_{j,B}
=\prod_{j\in J}d_B^{b_j}d_{j,B}
=\prod_{j\in J}d_{j,B}^{new},
\end{align}
because by the definition in~\eqref{eq:def b_i}, we have $b_0=0$ and $d_B^{b_0}=d_B^0=1$ and $d_{0,B}=d_{0,B}^{new}$.
Hence, altogether using~\ref{item:P edge tester} together with~\eqref{eq:w^new P3 0 notin I}, we obtain\COMMENT{ $\norm{\bp^Z_{-[r_\circ]_0}}=\norm{\bp^Z_{-[r_\circ]}}+p_0^Z$, and $\norm{\bp^{Z,2nd}_{-[r_\circ]_0}}=\norm{\bp^{Z,2nd}_{-[r_\circ]}}+p_0^{Z,2nd}$.}
\begin{align*}
&\omega^{new}\left(
E(\sA_{I_r}^{new})\sqcup\textstyle\bigsqcup_{j\in J} E(\cB_j^{new})
\right)
\\
&\stackrel{\hphantom{\eqref{eq:d_A^new},\eqref{eq:d_B^new}}}{=}
(1\pm3\heps^{1/3})
d_A^{b_{I_r}+p_0^A-p_0^{A,2nd}}
d_B^{b_{J}+p_0^B-p_0^{B,2nd}}
\\&\hphantom{\stackrel{\hphantom{\eqref{eq:d_A^new},\eqref{eq:d_B^new}}}{=}}
\Bigg(
\Big(\IND\{\JXV\cap-[r_\circ]=\emptyset \}\pm\eps\Big)
d_A^{\norm{\bpA_{-[r_\circ]}}-\norm{\bppA_{-[r_\circ]}}}
d_B^{\norm{\bpB_{-[r_\circ]}}-\norm{\bppB_{-[r_\circ]}}}
\\&
\hphantom{\stackrel{\hphantom{\eqref{eq:d_A^new},\eqref{eq:d_B^new}}}{=}}
\prod_{i\in I_{r_0}}d_{i,A}
\prod_{j\in  J}d_{j,B}
\frac{\omega_\iota(\sX_{\sqcup {I}})}{n^{|(I\cap-[r_\circ])\sm J|}}
\pm  n^{\eps}
\Bigg)\pm2n^{3\eps}
\\
&\stackrel{\eqref{eq:d_A^new},\eqref{eq:d_B^new}}{=}
\Big(\IND\{\JXV\cap-[r_\circ]_0=\emptyset \}\pm\eps'^2\Big)
d_A^{\norm{\bpA_{-[r_\circ]_0}}-\norm{\bppA_{-[r_\circ]_0}}}
d_B^{\norm{\bpB_{-[r_\circ]_0}}-\norm{\bppB_{-[r_\circ]_0}}}
\\&
\hphantom{\stackrel{\hphantom{\eqref{eq:d_A^new},\eqref{eq:d_B^new}}}{=}}
\prod_{i\in I_r}d_{i,A}^{new}
\prod_{j\in  J}d_{j,B}^{new}
\frac{\omega_\iota(\sX_{\sqcup {I}})}{n^{|(I\cap-[r_\circ]_0)\sm J|}}
\pm  n^{\eps'},
\end{align*}
which establishes~\ref{emb lem:general edge testers} in the case of $0\notin I$ or $0\in J\sm \JXV$.

\begin{step}\label{step:conclusion II 2}
Concluding~\ref{emb lem:general edge testers} if $0\in I\sm J$
\end{step}

Similar as in Step~\ref{step:conclusion II 1}, we can establish~\ref{emb lem:general edge testers} if $0\in I\sm J$ by the estimates derived in the Steps~\ref{step:main contribution} and~\ref{step:bound w(Gamma^hit)}.
So, let us assume that $0\in I\sm J$, and recall that we obtained an additional factor of $(d_0n)^{-1}$ in~\eqref{eq:bound sum w_p(M)} and~\eqref{eq:bound sum w(Gamma^hit)}.
Together with our key observation~\eqref{eq:new total weight} this yields that
\begin{align}
\nonumber
&\omega^{new}\left(
E(\sA_{I_r}^{new})\sqcup \textstyle\bigsqcup_{j\in J}E(\cB_j^{new})
\right)
\\\label{eq:w^new P3 0 in I}
&=\quad
(1\pm2\heps^{1/3})
d_A^{b_{I_r}+p_0^A-p_0^{A,2nd}}d_B^{b_{J}+p_0^B-p_0^{B,2nd}}
(d_0n)^{-1}
\omega\left(
E(\cA_{I_{r_0}})\sqcup \textstyle\bigsqcup_{j\in J}E(\cB_j)
\right)
\pm2n^{2\eps}.
\end{align}
We will use~\ref{item:P edge tester} for $\omega\left(
E(\cA_{I_{r_0}})\sqcup \textstyle\bigsqcup_{j\in J}E(\cB_j)
\right)$ in order to obtain~\ref{emb lem:general edge testers} from~\eqref{eq:w^new P3 0 in I}.
For $0\in I\sm J$, we have that 
\begin{align}\label{eq:d_A^new 0 in I}
d_0^{-1}d_A^{b_{I_r}}\prod_{i\in I_{r_0}}d_{i,A}=\prod_{i\in I_r}d_A^{b_i}d_{i,A}=\prod_{i\in I_r}d_{i,A}^{new}
\end{align} for $d_{i,A}^{new}$ defined as in~\eqref{eq:def b_i} because $d_{0}=d_{0,A}$, and similarly
\begin{align}\label{eq:d_B^new 0 in I}
d_B^{b_{J}}\prod_{j\in J}d_{j,B}=\prod_{j\in J}d_{j,B}^{new}.
\end{align}
Hence, altogether using~\ref{item:P edge tester} together with~\eqref{eq:w^new P3 0 in I}, we obtain\COMMENT{ $\norm{\bp^Z_{-[r_\circ]_0}}=\norm{\bp^Z_{-[r_\circ]}}+p_0^Z$, and $\norm{\bp^{Z,2nd}_{-[r_\circ]_0}}=\norm{\bp^{Z,2nd}_{-[r_\circ]}}+p_0^{Z,2nd}$.}
\begin{align*}
&\omega^{new}\left(
E(\sA_{I_r}^{new})\sqcup \textstyle\bigsqcup_{j\in J}E(\cB_j^{new})
\right)
\\
&\stackrel{\hphantom{\eqref{eq:d_A^new},\eqref{eq:d_B^new}}}{=}
(1\pm\heps^{1/3})
d_A^{b_{I_r}+p_0^A-p_0^{A,2nd}}d_B^{b_{J}+p_0^B-p_0^{B,2nd}}
(d_0n)^{-1}
\\&\hphantom{\stackrel{\hphantom{\eqref{eq:d_A^new},\eqref{eq:d_B^new}}}{=}}
\Bigg(
\Big(\IND\{\JXV\cap-[r_\circ]=\emptyset \}\pm\eps\Big)
d_A^{\norm{\bpA_{-[r_\circ]}}-\norm{\bppA_{-[r_\circ]}}}
d_B^{\norm{\bpB_{-[r_\circ]}}-\norm{\bppB_{-[r_\circ]}}}
\\&
\hphantom{\stackrel{\hphantom{\eqref{eq:d_A^new},\eqref{eq:d_B^new}}}{=}\Bigg(}
\prod_{i\in I_{r_0}}d_{i,A}
\prod_{j\in  J}d_{j,B}
\frac{\omega_\iota(\sX_{\sqcup {I}})}{n^{|(I\cap-[r_\circ])\sm J|}}
\pm  n^{\eps}
\Bigg)\pm2n^{2\eps}
\\
&\stackrel{\eqref{eq:d_A^new 0 in I},\eqref{eq:d_B^new 0 in I}}{=}
\Big(\IND\{\JXV\cap-[r_\circ]_0=\emptyset \}\pm\eps'^2\Big)
d_A^{\norm{\bpA_{-[r_\circ]_0}}-\norm{\bppA_{-[r_\circ]_0}}}
d_B^{\norm{\bpB_{-[r_\circ]_0}}-\norm{\bppB_{-[r_\circ]_0}}}
\\
&\hphantom{\stackrel{\hphantom{\eqref{eq:d_A^new},\eqref{eq:d_B^new}}}{=}\Big(}
\prod_{i\in I_r}d_{i,A}^{new}
\prod_{j\in  J}d_{j,B}^{new}
\frac{\omega_\iota(\sX_{\sqcup {I}})}{n^{|(I\cap-[r_\circ]_0)\sm J|}}
\pm  n^{\eps'},
\end{align*}
which establishes~\ref{emb lem:general edge testers} in the case of $0\in I\sm J$.

\begin{step}\label{step:conclusion II 3}
Concluding~\ref{emb lem:general edge testers} if $0\in \JXV$\nopagebreak
\end{step}

For the last case that $0\in \JXV$, we employ the estimate derived in Step~\ref{step:bound w(Gamma^hit)} in our key observation~\ref{eq:new total weight 0 in J}.
For $0\in \JXV$, we have that 
\begin{align}\label{eq:d_A^new 0 in J}
d_A^{b_{I_r}}\prod_{i\in I_{r_0}}d_{i,A}
=\prod_{i\in I_r}d_A^{b_i}d_{i,A}
=\prod_{i\in I_r}d_{i,A}^{new}
\end{align} for $d_{i,A}^{new}$ defined as in~\eqref{eq:def b_i} because $I_{r_0}=I_r$, and
\begin{align}\label{eq:d_B^new 0 in J}
d_B^{b_{J}}\prod_{j\in J}d_{j,B}
=\prod_{j\in J}d_B^{b_j}d_{j,B}
=\prod_{j\in J}d_{j,B}^{new},
\end{align}
because by the definition in~\eqref{eq:def b_i}, we have $b_0=0$ and $d_B^{b_0}=d_B^0=1$ and $d_{0,B}=d_{0,B}^{new}$.
Hence, altogether using~\ref{item:P edge tester} together with~\eqref{eq:bound sum w(Gamma^hit)} in our key observation~\eqref{eq:new total weight 0 in J}, we obtain\COMMENT{ $\norm{\bp^Z_{-[r_\circ]_0}}=\norm{\bp^Z_{-[r_\circ]}}+p_0^Z$, and $\norm{\bp^{Z,2nd}_{-[r_\circ]_0}}=\norm{\bp^{Z,2nd}_{-[r_\circ]}}+p_0^{Z,2nd}$.}
\begin{align*}
&\omega^{new}\left(
E(\sA_{I_r}^{new})\sqcup \textstyle\bigsqcup_{j\in J}E(\cB_j^{new})
\right)
\\
&\stackrel{\hphantom{\eqref{eq:d_A^new},\eqref{eq:d_B^new}}}{\leq}
\eps^{1/5}
d_A^{b_{I_r}+p_0^A-p_0^{A,2nd}}
d_B^{b_{J}+p_0^B-p_0^{B,2nd}}
\\&\hphantom{\stackrel{\hphantom{\eqref{eq:d_A^new},\eqref{eq:d_B^new}}}{=}}
\Bigg(
\Big(\IND\{\JXV\cap-[r_\circ]=\emptyset \}+\eps\Big)
d_A^{\norm{\bpA_{-[r_\circ]}}-\norm{\bppA_{-[r_\circ]}}}
d_B^{\norm{\bpB_{-[r_\circ]}}-\norm{\bppB_{-[r_\circ]}}}
\\&
\hphantom{\stackrel{\hphantom{\eqref{eq:d_A^new},\eqref{eq:d_B^new}}}{=}\Bigg(}
\prod_{i\in I_{r_0}}d_{i,A}
\prod_{j\in  J}d_{j,B}
\frac{\omega_\iota(\sX_{\sqcup {I}})}{n^{|(I\cap-[r_\circ])\sm J|}}
+  n^{\eps}
\Bigg)+2n^{2\eps}
\\
&\stackrel{\eqref{eq:d_A^new 0 in J},\eqref{eq:d_B^new 0 in J}}{\leq}
\eps'^2
d_A^{\norm{\bpA_{-[r_\circ]_0}}-\norm{\bppA_{-[r_\circ]_0}}}
d_B^{\norm{\bpB_{-[r_\circ]_0}}-\norm{\bppB_{-[r_\circ]_0}}}
\prod_{i\in I_r}d_{i,A}^{new}
\prod_{j\in  J}d_{j,B}^{new}
\frac{\omega_\iota(\sX_{\sqcup {I}})}{n^{|(I\cap-[r_\circ]_0)\sm J|}}
+  n^{\eps'},
\end{align*}
which establishes~\ref{emb lem:general edge testers} in the case of $0\in \JXV$ and concludes the proof of~\ref{emb lem:general edge testers}.

\begin{step}\label{step:sr}
Checking~\ref{emb lem:sr}
\end{step}

In order to establish~\ref{emb lem:sr}, we fix $H\in\cH$, $Z\in\{A,B\}$ and $i\in[r]$ if $Z=A$, $i\in\VRsm$ if $Z=B$, and we may assume that $i\in N_{R_\ast}(0)$ otherwise $Z_i^{H,new}=Z_i^{H}$.
We will first show that $Z_i^{H,new}$ as defined in~\eqref{eq:def Z^H,new} is $(\eps',d_i^{new})$-regular by employing Theorem~\ref{thm:new almost quasirandom}. 
In order to do so, we verify that every vertex in $X_i^H$ has the appropriate degree in $Z_i^{H,new}$ and that most pairs of vertices in $X_i^H$ have the appropriate common neighbourhood in $Z_i^{H,new}$. These properties follow easily due to the typicality of $G_Z$ and because $Z_i^H$ is $(\eps,q)$-well-intersecting. 
Finally, we show that also each vertex in $V_i$ has the correct degree in $Z_i^{H,new}$ by employing~\ref{emb lem:general edge testers}. Altogether this will imply that $Z_i^{H,new}$ is $(\eps',d_i^{new})$-super-regular.
Since we basically obtain $Z_i^{H,new}$ from~$Z_i^{H}$ by restricting the neighbourhood of every vertex by $b_i\leq\Delta(R)$ additional $(k-1)$-sets in $G_Z$ (see~\eqref{eq:1 S_y} and~\eqref{eq:degree y} below), we will obtain directly from~\ref{item:P1} that $Z_i^{H,new}$ is $(\eps',q+\Delta(R))$-well-intersecting.
We proceed with the details.\looseness=-1

For every vertex $y\in X_i^H$, let
\begin{align}
\label{eq:1 S_y}
\cS_{y}:=&\left\{\phi_\circ(\eH)\cup \sigma^+(x)\colon \eH\in E_{x,y}, x\in\sX_0 \right\}
\end{align}
with $E_{x,y}$ defined as in~\eqref{eq:def E_x,y}.
Note that $\cS_y\subseteq 
\bigcup_{\rR\in E(R)\colon i\in\rR} V_{\sqcup\rR\sm\{i\}}$ and $|\cS_y|=b_i$.
Since $G_Z$ is $(\eps,t,d_Z)$ typical with respect to $R$ by~\ref{item:P1}, and since $Z_i^H$ is $(\eps,d_{i,Z})$-super-regular and $(\eps,q)$-well-intersecting by~\ref{item:P2}, we conclude by the definition of $Z_i^{H,new}$ as in~\eqref{eq:def Z^H,new} that for all $y\in X_i^H$, we have
\begin{align}\label{eq:degree y}
\deg_{Z^{H,new}_i}(y)=|N_{Z^H_i}(y)\cap
N_{G_Z}(\cS_y)|=(1\pm\heps)d_{i,Z}d_{Z}^{b_i}|V_i|.
\end{align}
Note that~\eqref{eq:degree y} implies in particular that the density of $X_i^H$ and $V_i$ in $Z_i^{H,new}$ is $d_{Z_i^{H,new}}(X_i^H,V_i)=(1\pm\heps)d_{i,Z}d_{Z}^{b_i}$.\COMMENT{We defined the density in Section\ref{sec:graph regularity}.}

We can proceed similarly as for the conclusion~\eqref{eq:degree y} and obtain that all but at most $n^{3/2}$ pairs $\{y,y'\}\in \binom{X_i^H}{2}$ satisfy
\begin{align}\label{eq:codegree y,y'}
|N_{Z^{H,new}_i}(y\land y')|=(1\pm\heps)(d_{i,Z}d_Z^{b_i})^2|V_i|.
\end{align}
To see~\eqref{eq:codegree y,y'}, note that 
\begin{align*}
N_{Z^{H,new}_i}(y\land y')
=
N_{Z^H_i}(y\land y')\cap 
N_{G_Z}(\cS_y)\cap 
N_{G_Z}(\cS_{y'}),
\end{align*}
and for all but at most $2\Delta_Rn$ pairs $\{y,y'\}\in \binom{X_i^H}{2}$, we have $|\cS_{y}\cup \cS_{y'}|=2b_i$ by~\eqref{eq:pairs y,y'}. 
By employing again~\ref{item:P1} and~\ref{item:P2}, we obtain~\eqref{eq:codegree y,y'}, where we used~\eqref{eq:eps,t well-intersecting} that for all but at most $2n\cdot n^{1/4+\eps}$ pairs $\{y,y'\}\in \binom{X_i^H}{2}$ the sets of $(k-1)$-sets corresponding to $N_{Z_i^H}(y)$ and $N_{Z_i^H}(y')$ are disjoint.
Hence, all but at most $2n^{5/4+\eps}+2\Delta_Rn\leq n^{3/2}$ pairs $\{y,y'\}\in \binom{X_i^H}{2}$ satisfy \eqref{eq:codegree y,y'}.

We can now easily derive an upper bound for the number of $4$-cycles in $Z_i^{H,new}$ by~\eqref{eq:codegree y,y'}. 
To that end, note that every pair $\{y,y'\}\in \binom{X_i^H}{2}$ together with a pair of common neighbours in $N_{Z^{H,new}_i}(y\land y')$ forms a $4$-cycle in $Z_i^{H,new}$.
Hence, by~\eqref{eq:codegree y,y'}, the number of $4$-cycles in $Z_i^{H,new}$ is at most
\begin{align*}
C_4(Z_i^{H,new})
&\leq
\frac{|X_i^H|^2}{2}\cdot (1+2\heps)\frac{(d_{i,Z}d_Z^{b_i})^4|V_i|^2}{2}
+
n^{3/2}\cdot n^2
\\
&\leq 
(1+3\heps)\frac{(d_{i,Z}d_Z^{b_i})^4|X_i^H|^2|V_i|^2}{4}.
\end{align*}
Thus, we can apply Theorem~\ref{thm:new almost quasirandom} and obtain that 
\begin{align}\label{eq:A_i^H,* regular}
Z_i^{H,new}\text{ is $(\eps',d_{i,Z}d_Z^{b_i})$-regular.}
\end{align}

For every $v\in V_i$, in order to control the degree of $v$ in $Z^{H,new}_i$, we define weight functions $\omega_v\colon E(Z_i^H)\to\{0,1\}$ by $\omega_v(x_iv_i):=\IND{\{v_i=v\}}$, 
and $\omega_\iota\colon\sX_{i}\to\{0,1\}$ by $\omega_\iota:=\IND{\{x\in X_i^H \}}$,
and add the $1$-edge tester $$(\omega_v,\omega_\iota,J=J_Z,\JX=\emptyset,\JV=\emptyset,\bc=\{v\},\pP=(\mathbf{0},\mathbf{0},\mathbf{0},\mathbf{0}))$$ 
to $\cW_{edge}'$ for $J_A:=\emptyset$ and $J_B:=\{i\}$.
This is indeed a (general) $1$-edge tester satisfying Definition~\ref{def:global edge tester}.
In particular, $\bpA(x)=\bppA(x)=\bpB(x)=\bppB(x)=\mathbf{0}$ for every $x\in X_i^H$ by Definition~\ref{def:pattern} because the \fst-pattern and \scd-pattern of a single vertex is always $\mathbf{0}$.
Since $Z_i^H$ is $(\eps,d_{i,Z})$-super-regular, we have that
\begin{align*}
\omega_v(E(Z_i^H))=(1\pm\eps)d_{i,Z}|X_i^H|.
\end{align*}
Hence in particular, this general edge tester satisfies~\ref{item:P edge tester}.
By~\ref{emb lem:general edge testers}, we obtain\COMMENT{We apply \ref{emb lem:general edge testers} for the edge tester $$(\omega_v,\omega_\iota,J=J_Z,\JX=\emptyset,\JV=\emptyset,\bc=\{v\},\pP=(\mathbf{0},\mathbf{0},\mathbf{0},\mathbf{0})).$$
Note that
\begin{center}
\begin{tabular}{l|l}
parameter & plays the role of 
\\\hline
$\omega_v^{new}(E(Z_i^H))$ &$\omega^{new}\left(E(\sA_{I_r}^{new})\sqcup E(\cB_J^{new})\right)$  
\\
$(1\pm \eps'^2)$ & $\Big(\IND\{\JXV\cap-[r_\circ]_0=\emptyset \}\pm\eps'^2\Big)$
\\
$1=d_A^0d_B^0$ & $d_A^{\norm{\bpA_{-[r_\circ]_0}}-\norm{\bppA_{-[r_\circ]_0}}}
d_B^{\norm{\bpB_{-[r_\circ]_0}}-\norm{\bppB_{-[r_\circ]_0}}}$
\\
$d_{i,Z}^{new}$ & $\prod_{i\in I_r}d_{i,A}^{new}
\prod_{j\in  J}d_{j,B}^{new}$
\\
$|X_i^H|$ & $\omega_\iota(\sX_{\sqcup {I}})$
\\
$1=n^0$ & $n^{|(I\cap-[r_\circ]_0)\sm J|}$.
\end{tabular}
\end{center}
  }
\begin{align*}
\dg_{Z^{H,new}_i}(v)
=
\omega_v^{new}(E(Z^{H,new}_i))
\stackrel{\ref{emb lem:general edge testers}}{=}
(1\pm \eps'^2)d_{i,Z}^{new}|X_i^H|\pm n^{\eps'}.
\end{align*}
Together with~\eqref{eq:degree y} and~\eqref{eq:A_i^H,* regular}, this implies that $Z_i^{H,new}$ is $(\eps',d_{i,Z}^{new})$-super-regular because $d_{i,Z}^{new}=d_{i,Z}d_Z^{b_i}$ (see~\eqref{eq:def b_i}).
Further, since the neighbourhood of every vertex $y\in X_i^H$ in $Z_i^{H,new}$ is the intersection of a set $\cS_y^{new}$ of $(k-1)$-sets in $G_Z$ (see~\eqref{eq:degree y}), and every $y\in X_i^H$ is contained in at most $n^{1/4+\eps}+\Delta_R$\COMMENT{Because $Z_i^H$ is $(\eps,q)$-well-intersecting, see~\eqref{eq:eps,t well-intersecting}, and because of~\eqref{eq:pairs y,y'}.} pairs $\{y,y'\}\in\binom{X_i^H}{2}$ such that $\cS_y^{new}\cap \cS_{y'}^{new}\neq\emptyset$, 
we also obtain that $Z_i^{H,new}$ is $(\eps',q+\Delta(R))$-well-intersecting as defined in~\eqref{eq:eps,t well-intersecting}.
This establishes~\ref{emb lem:sr}.

\begin{step}\label{step:local weight}
Checking~\ref{emb lem:local weight}\nopagebreak
\end{step}

For every $\omega\in\cW_{local}$ with $\omega\colon\binom{E(\cA_0)}{\ell}^=\to[0,s]$, we add $\omega$ to $\cW$. 
Hence,~\eqref{eq:w packing size} yields~\ref{emb lem:local weight}.

\begin{step}\label{step:hit tester}
Checking~\ref{emb lem:hit tester}
\end{step}

In order to establish~\ref{emb lem:hit tester}, we fix $(\omega,c)\in\cW_{0}$ with $\omega\colon\sX_0\to[0,s]$ and $c\in V_0$.
By~\eqref{eq:packing size}, we have that $V_0\sm \sigma(X_0^H\cap\sX_0^\sigma)\geq \eps^{2/3}n/3$ for every $H\in\cH$, and thus, the probability for a vertex $x\in \sX_0\sm\sX_0^\sigma$ to be mapped onto $c$ is at most $3/\eps^{2/3}n$.
We therefore expect that $\omega\left(
\{
x\in\sX_0\sm\sX_0^\sigma\colon \sigma^+(x)=c
\}
\right)
\leq \omega(\sX_0)/n^{1- 2\eps}$.
By an application of Theorem~\ref{thm:McDiarmid} and a union bound, we can establish concentration with probability, say at least $1-\eul^{-n^{\eps}}$.
This establishes~\ref{emb lem:hit tester} and completes the proof of Lemma~\ref{lem:packing lemma}.
\end{proof}

\section{Iterative packing}\label{sec:mainpart}

{
In this section we essentially prove our main result, Theorem~\ref{thm:main_new}.
We prove the following lemma whose statement is very similar
because we only require additionally that for every graph $H\in\cH$ and every reduced edge $\rR\in E(R)$, the graph $H[X_{\cup\rR}^H]$ is a matching.
This reduction can be achieved by an application of Lemma~\ref{lem:refining partitions}  and simplifies several arguments;
it is presented in the proof of Theorem~\ref{thm:main_new} in Section~\ref{sec:mainproofs}.
}

\begin{lemma}\label{lem:main matching}
Let $1/n\ll\eps\ll1/t\ll \alpha,1/k$ and $r\leq n^{2\log n}$ as well as $d\geq n^{-\eps}$. 
Suppose $(\cH,G,R,\cX,\cV)$ is an  $(\eps,t,d)$-typical and $\alpha^{-1}$-bounded blow-up instance of size $(n,k,r)$ 
and \mbox{$|\cH|\leq  n^{2k}$}.
Suppose that $e_\cH(\sX_{\sqcup \rR}) \leq (1-\alpha)dn^k$ for all $\rR \in E(R)$,
and $H[X_{\cup\rR}^H]$ is a matching
if $\rR\in E(R)$ and empty if $\rR\in\binom{[r]}{k}\sm E(R)$ for each $H\in\cH$.
Suppose $\cW_{set},\cW_{ver}$ are sets of  $\alpha^{-1}$-set testers and $\alpha^{-1}$-vertex testers of size at most $n^{3\log n}$, respectively.
Then there is a packing $\phi$ of $\cH$ into $G$ such that\looseness=-1 
\begin{enumerate}[label={\rm (\roman*)}]
	\item\label{item:partition1} $\phi(X_i^H)=V_i$ for all $i\in [r]$ and $H\in \cH$;
	\item\label{item:set testers1} $|W\cap \bigcap_{j\in [\ell]}\phi(Y_j)|= |W||Y_1|\cdots |Y_\ell|/n^\ell \pm \alpha n$ for all $(W,Y_1,\ldots,Y_\ell)\in \cW_{set}$;
	\item\label{item: vertex tester1} $\omega(\phi^{-1}(\bc))=(1\pm\alpha)\omega(\sX_{\sqcup I})/n^{|I|}\pm n^\alpha$ for all $(\omega,\bc)\in \cW_{ver}$ with centres $\bc$ in~$I$.
\end{enumerate}
\end{lemma}

\lateproof{Lemma~\ref{lem:main matching}}
We split the proof into five steps.
In Step~\ref{step:colouring}, we define a vertex colouring of the reduced graph which will incorporate in which order we consider the clusters in turn.
In Step~\ref{step:partitioning}, we partition $G$ into two edge-disjoint subgraphs $G_A$ and $G_B$.
In Step~\ref{step:candidacy}, we introduce candidacy graphs and edge testers that we track for the partial packing in Step~\ref{step:induction}, where we iteratively apply Lemma~\ref{lem:packing lemma} and consider the clusters in turn with respect to the ordering of the clusters given by the colouring obtained in Step~\ref{step:colouring}. We only use the edges of $G_A$ for the partial packing in Step~\ref{step:induction} such that we can complete the packing in Step~\ref{step:completion} using the edges of $G_B$.

\begin{step}\label{step:colouring}
Notation and colouring of the reduced graph
\end{step}

We will proceed cluster by cluster in Step~\ref{step:induction} to find a function that packs almost all vertices of~$\cH$ into $G$. 
Since we allow $r$ to grow with $n$ and only require that $r\leq n^{\log n}$, we need to carefully control the growth of the error term. 
Recall that $R_\ast$ is the $2$-graph with vertex set $V(R)$ and edge set $\bigcup_{\rR\in E(R)}\binom{\rR}{2}$.\COMMENT{Also defined in Section~\ref{sec:notation}.}
Let $c\colon V(R)\to [T]$ be a proper vertex colouring of $R_\ast^3$ where $T:=k^3\alpha^{-3}$.
The colouring naturally yields an order in which we consider the clusters in turn. 
To this end, we simply relabel the cluster indices such that the colour values are non-decreasing; that is, 
$c(1)\leq\cdots\leq c(r)$.
Note that the sets $(c^{-1}(j))_{j\in[T]}$ are independent in $R_\ast^3$. 
We choose new constants $\eps_0,\eps_1,\ldots,\eps_T,\mu,\gamma$ such that\looseness=-1  $$\eps\ll\eps_0\ll\eps_1\ll\cdots\ll\eps_T\ll\mu\ll\gamma\ll1/t\ll\alpha,1/k.$$

\medskip
For $i,q\in[r]$ and $I\subseteq[r]$, we define counters $c_i(q), c_I(q), m_i(q)$ (see \eqref{eq:def c_i(q)}--\eqref{eq:def m_i(q)} below).
Our intuition is the following:
If we think of $[q]$ as the indices of clusters that have already been embedded, then $c_i(q)$ is the largest colour of an already embedded cluster in the closed neighbourhood of $i$ in~$R_\ast$.
That is, $c_i(q)$ is the largest colour that is relevant to $i$ after embedding the first $q$ clusters, and~$c_i(q)$ will incorporate how to update the error term.

To be more precise, for $i,q\in[r]$ and an index set $I\subseteq[r]$ (that is, $I\subseteq\rR\in E(R)$), let
\begin{align}\label{eq:def c_i(q)}
c_i(q)&:=\max\left\{
\{0\}\cup\{c(j)\colon j\in N_{R_\ast}[i]\cap[q]\}
\right\};
\\\label{eq:def max c_i(q)}
c_I(q)&:=
\max_{i\in I}\left\{0, c_i(q)
\right\}.
\end{align}

\COMMENT{$c_{\emptyset}(q)=0$}Similarly as $c_i(q)$, we define $m_i(q)$ as the number of edges in $R$ that contain $i$ and where $k-1$ clusters excluding $i$ have already been embedded. 
That is, for $i,q\in[r]$, let
\begin{align}
\label{eq:def m_i(q)}
m_i(q)&:=
\left|
\left\{
\rR\in E(R)\colon i\in\rR,|\rR\cap[q]\sm\{i\}|=k-1
\right\}
\right|.
\end{align}
{Further, for all $i\in[r]$ and index sets $I\subseteq [r]$, we set $c_i(0)=c_I(0)=m_i(0):=0$.}

For $q\in[r]$, recall that $\sX_{q}=\textstyle\bigcup_{H\in\cH}X^H_{q}$, and we set 
\begin{align*} 
 \cX_q& :=\textstyle\bigcup_{\ell\in[q]}\sX_{\ell}, 
& \cV_q&:=\textstyle\bigcup_{\ell\in [q]}V_{\ell}. 
\end{align*}

\begin{step}\label{step:partitioning}
Partitioning the edges of $G$\nopagebreak
\end{step}

In order to reserve an exclusive set of edges for the completion in Step~\ref{step:completion}, we partition the edges of $G$ into two subgraphs $G_A$ and $G_B$.
For each edge $\gG$ of~$G$ independently, we add $\gG$ to $G_B$ with probability $\gamma$ and otherwise to $G_A$. 
Let $d_A:=(1-\gamma)d$ and $d_B:=\gamma d$.
Using Chernoff's inequality {and a union bound}, we can easily conclude that with probability at least $1-1/n$ it holds that
\begin{align}\label{eq:G_Z sr}
\begin{minipage}[c]{0.9\textwidth}\em
for all $i\in[r]$ and all pairs of disjoint sets $\cS_A,\cS_B\subseteq 
\bigcup_{\rR\in E(R)\colon i\in\rR} V_{\sqcup\rR\sm\{i\}}$ with $|\cS_A\cup\cS_B|\leq t$, we have 
$\big|V_i\cap N_{G_A}(\cS_A)\cap N_{G_B}(\cS_B)\big|
=(1\pm\eps_0)d_A^{|\cS_A|}d_B^{|\cS_B|}n.$
\end{minipage}
\ignorespacesafterend
\end{align}
Hence, we may assume that $G$ is partitioned into $G_A$ and $G_B$ such that~\eqref{eq:G_Z sr} holds.
In particular, \eqref{eq:G_Z sr} implies that $G_Z$ is $(2\eps_0,t,d_Z)$-typical with respect to $R$ for both $Z\in\{A,B\}$.

\begin{step}\label{step:candidacy}
Partial packings, candidacy graphs and edge testers
\end{step}

For $q\in[r]$, we call $\phi\colon \bigcup_{H\in\cH,i\in[q]}\hX^H_i\to\cV_q$ a \emph{$q$-partial packing} if $\hX^H_i\subseteq X^H_i$ and $\phi(\hX^H_i)\subseteq V_i$ for all $H\in\cH, i\in[q]$ such that $\phi$ is a packing of $(H[\hX^H_1\cup\ldots\cup \hX^H_q])_{H\in\cH}$ into $G_A[\cV_q]$. Note that $\phi|_{\hX^H_i}$ is injective for all $H\in\cH$ and $i\in[q]$.
For convenience, we often write
\begin{align*}
\cX_q^{\phi} :=\textstyle\bigcup_{H\in\cH,i\in[q]}\hX^H_i. 
\end{align*}
Further, we call $\phi^+\colon\cX_q\to\cV_q$ a \defn{cluster-injective extension} of $\phi$ if $\phi^+$ is an extension of $\phi$ such that $\phi^+|_{X_i^H}$ is injective (and thus bijective) and $\phi^+(X_i^H)=V_i$ for all $H\in\cH, i\in[q]$.
Note that we do not even require that $\phi^+$ is an embedding of $H[X^H_1\cup\ldots\cup X^H_q]$ for $H\in\cH$.

Suppose $q\in[r]$ and $\phi_q\colon \cX_q^{\phi_q}\to \cV_q$ is a $q$-partial packing with a cluster-injective extension $\phi_q^+$.
We consider two kinds of candidacy graphs as in Definition~\ref{def:candidacy graph}:
Candidacy graphs $A^H_{I}(\phi_q^+)$ with respect to $\phi_q^+$ and $G_A$ for all index sets $I\subseteq[r]\sm[q]$, 
and candidacy graphs $B^H_{j}(\phi_q^+)$ with respect to $\phi_q^+$ and $G_B$ for all $j\in[r]$.
The candidacy graphs $A^H_{I}(\phi_q^+)$ will be used to extend the $q$-partial packing $\phi_q$ to a $(q+1)$-partial packing $\phi_{q+1}$ via Lemma~\ref{lem:packing lemma} in Step~\ref{step:induction}, whereas the candidacy graphs $B^H_{j}(\phi_r^+)$ will be used for the completion in Step~\ref{step:completion}.

\medskip
Given $\phi_q$, $\phi_q^+$ and a collection $\cA^q$ and $\cB^q$ of candidacy graphs $A_{I}^{H,q}\subseteq A_{I}^H(\phi_q^+)$ and $B_{j}^{H,q}\subseteq B_{j}^H(\phi_q^+)$ for all index sets $I\subseteq [r]\sm[q]$ and $j\in[r]$, we introduce (general) edge testers as in Definition~\ref{def:general edge tester} to track several quantities during our packing procedure.  

To that end, we first define a set $\cW_{initial}$ of tuples $(\omega_\iota,J,\JX,\JV,\bc,\pP)$. {We also define a superset~$\cW_{hit}$ of $\cW_{ver}$ containing tuples $(\omega,\bc)$.}
For every vertex tester $(\omega_{ver},\bc)\in\cW_{ver}$ as in the assumptions of Lemma~\ref{lem:main matching} with centres $\bc\in V_{\sqcup I}$ for an index set $I\subseteq [r]$, and for all $J\subseteq I$\COMMENT{$J$ can also be empty.} and pairs of disjoint sets $\JX, \JV\subseteq J$,
all $\bpA, \bppA,\bpB, \bppB\in[k\alpha^{-1}]^r_0$, we define a tuple 
$(\omega_\iota,J,\JX, \JV,\bc,\pP)$ with $\omega_\iota\colon\sX_{\sqcup I}\to[0,\alpha^{-1}]$, $\pP:=(\bpA,\bppA,\bpB,\bppB)$, by 
\begin{align}\label{eq:def initial weight vertex tester}
\omega_\iota(\xX):=\IND{\{\xX\in E_\cH(\pP,I,J) \}} \cdot\omega_{ver}(\xX),
\end{align}
and we add this tuple to $\cW_{initial}$. (Recall Definition~\ref{def:e_H,p} for $E_\cH(\pP,I,J)$.) We also add $(\omega_{ver},\bc)$~to~$\cW_{hit}$.\looseness=-1

Similarly, for every $\rR\in E(R)$, all $J\subseteq \rR$ and pairs of disjoint sets $\JX, \JV\subseteq J$, $\gG\in  V_{\sqcup \rR}$, and $\bpA, \bppA, \bpB, \bppB\in[k\alpha^{-1}]^r_0$, we define a tuple $(\omega_\iota,J,\JX, \JV,\gG,\pP)$ with $\omega_\iota\colon\sX_{\sqcup \rR}\to\{0,1\}$, $\pP:=(\bpA,\bppA,\bpB,\bppB)$, by 
\begin{align}\label{eq:def initial weight boundedness}
\omega_\iota(\xX):=\IND{\{\xX\in E_\cH(\pP,\rR,J) \}},
\end{align}
and we add this tuple to $\cW_{initial}$, and
we define a tuple $(\omega,\gG)$ with $\omega\colon\sX_{\sqcup {\rR}}\to\{0,1\}$ by $\omega(\xX):=\IND\{\xX\in E(\cH) \}$ and add $(\omega,\gG)$ to $\cW_{hit}$.

To control the number of unembedded $H$-vertices in one graph $H$ that could potentially be mapped onto a fixed vertex $v$ during the completion, we define for all $j\in[r]$, $H\in\cH$, and $v\in V_j$, a tuple $(\omega_{\iota,H},J=\{j\},\JX=\{j\},\JV=\emptyset,\bc=\{v\},(\mathbf{0},\mathbf{0},\mathbf{0},\mathbf{0}))$ with $\omega_{\iota,H}\colon\sX_{j}\to\{0,1\}$ by 
\begin{align}\label{eq:def edge tester leftover}
\omega_{\iota,H}(\xX):=\IND{\{\xX\in X_j^H \}},
\end{align}
and we add this tuple to $\cW_{initial}$.\COMMENT{
This is an edge tester to control the leftover of a single vertex $v$ in one graph $H$. Note that a single vertex has always \fst-pattern $\mathbf{0}$.
}

For one single graph $H\in \cH$,
we also consider tuples with only two centres.
That is, for all $H\in\cH$, $\rR\in E(R)$, distinct $j,j_X\in\rR$, $v\in V_j$, $w\in V_{j_X}$, and $\bpA= \bppA=\mathbf{0}, \bpB, \bppB\in[k\alpha^{-1}]^r_0$, we define $I=J:=\{j,j_X\}$, $\JX:=\{j_X\}$, $\JV:=\emptyset$ and a tuple $(\omega_\iota,J,\JX, \JV,\{v,w\},\pP)$ with $\omega_\iota\colon\sX_{\sqcup I}\to\{0,1\}$, $\pP:=(\bpA,\bppA,\bpB,\bppB)$, by
\begin{align}\label{eq:def initial weight H-edge leftover}
\omega_\iota(\xX):=\IND{\{\xX\in E_\cH(\pP,I,J), \xX\subseteq V(H) \}},
\end{align}
and we add this tuple to $\cW_{initial}$.

\medskip
We now define a set $\cW_{edge}^q=\cW_{edge}^q(\phi_q,\phi_q^+,\cA^q,\cB^q)$ of edge testers with respect to the elements in $\cW_{initial}$. 
That is, for every $(\omega_\iota,J,\JX, \JV,\bc,\pP)\in\cW_{initial}$, let $(\omega_q,\omega_\iota,J,\JX, \JV,\bc,\pP)$ be the edge tester with respect to $(\omega_\iota,J,\JX, \JV,\bc,\pP)$, $(\phi_q,\phi_q^+)$, $\cA^q$ and $\cB^q$ as in Definition~\ref{def:general edge tester}, and we add $(\omega_q,\omega_\iota,J,\JX, \JV,\bc,\pP)$ to $\cW_{edge}^q$.

\begin{step}\label{step:induction}
Induction
\end{step}

We inductively prove that the following statement~\ind{q} holds for all $q\in[r]_0$, which will provide a partial packing of $\cH$ into $G_A$. 
\begin{itemize}
\item[\ind{q}.] 
For all $H\in\cH$, there exists a $q$-partial packing $\phi_q\colon \cX_q^{\phi_q}\to \cV_q$ with $|\cX_q^{\phi_q}\cap X^H_i|\geq (1-\eps_{c_i(q)})n$ for all $i\in [q]$, and with a cluster-injective extension $\phi_q^+$ of~$\phi_q$, 
and for all index sets $I_A\subseteq [r]\sm[q]$, $I_B\in[r]$, and $(Z,\cZ)\in\{(A,\cA),(B,\cB)\}$, there exist subgraphs $Z^{H,q}_{I_Z}$ of the candidacy graphs $Z^H_{I_Z}(\phi_q^+)$ with respect to $\phi_q^+$ and $G_Z$ (where $\cZ_{I_Z}^{q}:=\bigcup_{H\in\cH}Z_{I_Z}^{H,q}$ and $\cZ^q$ is the collection of all $\cZ_{I_Z}^{q}$) such that
\begin{enumerate}[label=\text{(\alph*)}]
\item\label{Z candidacy superregular} $Z^{H,q}_{i_Z}$ is $(\eps_{c_{i_Z}(q)},d_Z^{m_{i_Z}(q)})$-super-regular and $(\eps_{c_{i_Z}(q)},c_{i_Z}(q)t^{1/2})$-well-intersecting with respect to $G_Z$ for all $i_A\in[r]\sm[q]$, $i_B\in[r]$ and $Z\in\{A,B\}$;\COMMENT{Recall the definitions of $c_i(q)$ and $m_i(q)$ in~\eqref{eq:def c_i(q)} and \eqref{eq:def m_i(q)}.}

\item\label{edge weight}
for every edge tester $(\omega_q,\omega_\iota,J,\JX, \JV,\bc,\pP)\in\cW^q_{edge}(\phi_{q},\cA^q,\cB^q)$ with centres $\bc\in V_{\sqcup I}$ for $I\subseteq [r]$, non-empty $I_{q}:=(I\sm[q])\cup J$, patterns~$\pP=(\bpA,\bppA,\bpB,\bppB)$, and with $\JXV:=\JX\cup\JV$, we have that
\begin{align*}
\omega_q\left(E(\cA_{I_{q}\sm J}^{q})\sqcup \textstyle\bigsqcup_{j\in J} E(\cB_j^q)\right)
=
&\Big(\IND\{\JXV\cap[q]=\emptyset \}\pm \eps_{c_{I_q}(q)}\Big)
\prod_{Z\in\{A,B\}}
d_Z^{\norm{\bpZ_{[q]}}-\norm{\bppZ_{[q]}}}
\\&\prod_{i\in I_q\sm J}d_{A}^{m_i(q)}
\prod_{j\in  J}d_{B}^{m_j(q)}
\frac{\omega_\iota(\sX_{\sqcup {I}})}{n^{|(I\cap[q])\sm J|}}
\pm  n^{\eps_{c_{I_q}(q)}};
\end{align*}

\item\label{e,f suitable} 
for all $\gG=\{v_{i_1},\ldots,v_{i_k}\}, \hG=\{w_{j_1},\ldots,w_{j_k}\}\in E(G_A)$ with $v_{i_k}=w_{j_k}$, $I:=\{i_1,\ldots,i_k\}\neq \{j_1,\ldots,j_k\}=:J$, 
and $\eps^\ast:=\max\{\eps_{c_{I\sm[q]}(q)},\eps_{c_{J\sm[q]}(q)}\}$, we have that
$$
\big|E_{\gG,\hG,\phi_q}(\cA^q)\big|
\leq 
\max\left\{
n^{k-\left|(I\cup J)\cap [q]\right|+\eps^\ast},
n^{\eps^\ast}
\right\};
$$

\item\label{treffer auf g} 
for all $(\omega,\bc)\in\cW_{hit}$ with $\bc\in V_{\sqcup I}$ and all non-empty $J\subseteq I$, we have that
\begin{align*}
&\omega\left(
\{
\xX\in \sX_{\sqcup {I}}\colon \xX\subseteq \eH\in E(\cH), 
\phi_q^+(\xX\cap\sX_{\cup (J\cap[q])})\subseteq \bc
\}
\right)
\\&
\leq 
\omega(\sX_{\sqcup {I}})/n^{|J\cap [q]|-\eps_{c_{J\sm[q]}(q)}}+n^{\eps_{c_{J\sm[q]}(q)}};
\end{align*}

\item \label{G_B leftover for typicality}
$\big|
(V_i\sm \phi_q(X_i^H)) \cap  N_{G_B}(\cS)
\big|
\leq \eps_T \big|
V_i \cap  N_{G_B}(\cS)
\big|$ for all $H\in\cH$, $i\in[q]$ and $\cS\subseteq \bigcup_{\rR\in E(R)\colon i\in\rR} V_{\sqcup\rR\sm\{i\}}$ with $|\cS|\leq t$;

\item\label{set tester} 
$|W\cap\bigcap_{j\in[\ell]} \phi_q(Y_j)|=|W||Y_1|\cdots|Y_\ell|/n^\ell\pm\alpha^2 n$ for all $(W,Y_1,\ldots,Y_\ell)\in\cW_{set}$ with $W\subseteq V_i$, $i\in[q]$;

\item\label{vertex tester}
$\omega(\phi_{q}^{-1}(\bc))=(1\pm\eps_T)\omega(\sX_{\sqcup {I}})/n^{|I|}\pm n^{\eps_T}$ for all $(\omega,\bc)\in\cW_{ver}$ with centres $\bc\in V_{\sqcup I}$ for $I\subseteq[q]$.
\end{enumerate}
\end{itemize}

Let us first explain~\ind{q}\ref{Z candidacy superregular}--\ref{vertex tester}.
Properties~\ind{q}\ref{Z candidacy superregular}--\ref{e,f suitable} are used to establish~\ind{q+1} by applying Lemma~\ref{lem:packing lemma}. 
In particular, these properties will imply that the assumptions~\ref{item:P2}, \ref{item:P edge tester} and \ref{item:P4} are satisfied, respectively, in order to apply Lemma~\ref{lem:packing lemma}.
Properties~\ind{q}\ref{edge weight} and~\ref{G_B leftover for typicality} can be used to control the leftover for the completion in Step~\ref{step:completion}.
Properties~\ind{r}\ref{set tester} and~\ref{vertex tester} will imply the conclusions~\ref{item:set testers1} and~\ref{item: vertex tester1} of Lemma~\ref{lem:main matching} as we merely modify the $r$-partial packing~$\phi_r$ during the completion in Step~\ref{step:completion}.

\bigskip
We now inductively prove that~\ind{q} holds for all $q\in[r]_0$.
The statement~\ind{0} holds for $\phi_0$ and $\phi_0^+$ being the empty function and $Z_{I_Z}^{H,0}$ being complete multipartite $2|I_Z|$-graphs:
Clearly, for all $Z\in\{A,B\}$, and all index sets $I_A\subseteq[r]$ and $I_B\in[r]$, the candidacy graph $Z_{I_Z}^H(\phi_0^+)$ is complete $2|I_Z|$-partite.
For~\ind{0}\ref{edge weight}, consider an edge tester $(\omega,\omega_\iota,J,\JX,\JX,\bc,\pP)\in\cW_{edge}^0$ with centres $\bc\in V_{\sqcup I}$ and note that, by the definition of an edge tester (see Definition~\ref{def:general edge tester}), we have $\omega(E(\cA_{I\sm J})\sqcup \textstyle\bigsqcup_{j\in J}E(\cB_{j}))=\omega_\iota(\sX_{\sqcup {I}})$.
For~\ind{0}\ref{e,f suitable}, we observe that for all $\gG=\{v_{i_1},\ldots,v_{i_k}\}, \hG=\{w_{j_1},\ldots,w_{j_k}\}\in E(G_A)$ with $v_{i_k}=w_{j_k}$, we have that $|E_{\gG,\hG,\phi_0}(\cA^0)|\leq \alpha^{-1} e_\cH(\sX_{i_1},\ldots,\sX_{i_k})\leq \alpha^{-1}n^k\leq n^{k+\eps_0}$ because $e_\cH(\sX_{\sqcup {\rR}})\leq (1-\alpha)dn^k$ and $\Delta(H)\leq \alpha^{-1}$ for each $H\in\cH$ by assumption.
(Recall Definition~\ref{def:E_g,h,phi} of~$E_{\gG,\hG,\phi_0}(\cA^0)$.)
\ind{0}\ref{treffer auf g}--\ref{vertex tester} are vacuously true.

\medskip
Hence, we assume the truth of~\ind{q} for some $q\in[r-1]_0$ and let $\phi_q\colon \cX_q^{\phi_q}\to \cV_q$, $\phi_q^+$, and $\cA^q$ and $\cB^q$ be as in~\ind{q}.
Any function $\sigma\colon \sX_{q+1}^\sigma\to V_{q+1}$ with $\sX_{q+1}^\sigma\subseteq\sX_{q+1}$ naturally extends $\phi_q$ to a function $\phi_{q+1}:=\phi_q\cup\sigma$ with $\phi_{q+1}\colon \cX_q^{\phi_q}\cup\sX_{q+1}^\sigma\to \cV_{q+1}$.\COMMENT{That is,
\begin{align*}
\phi_{q+1}(x):=\begin{cases}
\phi_q(x) & \text{if $x\in \cX_q^{\phi_q}$,}\\
\sigma(x) & \text{if $x\in \sX_{q+1}^\sigma$.}
\end{cases}
\end{align*}
}

We now make a key observation based on the Definition~\ref{def:candidacy graph} of candidacy graphs:
By definition of the candidacy graphs $\sA^q_{q+1}=\bigcup_{H\in\cH}A_{q+1}^{H,q}$ where $A_{q+1}^{H,q}\subseteq A_{q+1}^{H}(\phi_q^+)$, if~$\sigma$ is a conflict-free packing in $\sA^q_{q+1}$ as defined in~\eqref{eq:conflict free}, then $\phi_{q+1}$ is a $(q+1)$-partial packing. See also Figure~\ref{fig:conflict-free packing} in Section~\ref{sec:packing instance}.

We aim to apply Lemma~\ref{lem:packing lemma} in order to obtain a conflict-free packing $\sigma$ in~$\sA_{q+1}^q$.
To this end, we consider subgraphs $\cH_{q+1}, G_{A,q+1}, G_{B,q+1}, R_{q+1}$ of $\cH, G_A, G_B, R$, respectively, that consist only of the  `relevant' clusters when finding a conflict-free packing in $\sA_{q+1}^q$. 
Note that all relevant clusters lie in $N_{R_\ast^3}[q+1]$.
That is, by considering all clusters in $N_{R_\ast^3}[q+1]$, we also account for hyperedges $\rR\in E(R)$ and all $R$-edges that intersect $\rR$ where $q+1\notin\rR$ but $\rR\cap \rR_{q+1}\neq\emptyset$ for some $\rR_{q+1}\in E(R)$ with $q+1\in\rR_{q+1}$.\COMMENT{To consider such edges it would suffice to consider $R_\ast^2$. 
However, for the pattern definitions~\ref{def:pattern} and~\ref{def:penpattern}, we also want to consider all hyperedges that intersect such edges and these edges may live in $R_\ast^3$.}
Let $\cQ:=[q]\cap N_{R_\ast^3}(q+1)$ and for each $Z\in\{A,B\}$, let
\begin{align*}
\cH_{q+1} &:=\bigcup_{H\in\cH}H\left[
\textstyle\bigcup_{i\in N_{R_\ast^3}[q+1]}X^H_i \right];
\\G_{Z,q+1} &:=G_Z\left[
\textstyle\bigcup_{i\in N_{R_\ast^3}[q+1]}V_i
\right];
\\R_{q+1}&:=R\left[N_{R_\ast^3}[q+1]\right].
\end{align*}
\COMMENT{Or equivalent: $\cH_{q+1}=\cH\left[\sX_{\cup{N_{R_\ast^3}[q+1]}} \right]$ and $G_{q+1}=G_A\left[ V_{\cup N_{R_\ast^3}[q+1]}\right]$.}
Correspondingly, for $(Z,\cZ)\in\{(A,\cA),(B,\cB)\}$, we also define a subset $\cZ^q[R_{q+1}]$ of $\cZ^q$. 
Let $\cZ^q[R_{q+1}]$ be the collection of all candidacy graphs  $\cZ_{I_Z}^q$ for index sets $I_A\subseteq N_{R_\ast^3}[q+1]\sm[q]$, $I_B\in N_{R_\ast^3}[q+1]$.\COMMENT{In fact, we only have to update candidacy graphs in $R_\ast^2$. $R_\ast^3$ is only needed to account for all patterns.}

Following the definition of a packing instance in Section~\ref{sec:packing instance}, we observe that 
$$\sP:=(\cH_{q+1}, G_{A,q+1}, G_{B,q+1}, R_{q+1}, \cA^q[R_{q+1}], \cB^q[R_{q+1}],\phi_q|_{\sX_{\cup \cQ}},\phi_q^+|_{\sX_{\cup \cQ}})$$ is a packing instance of size $$(n,k,|N_{R_\ast^3}(q+1)\sm[q]|, |\cQ|).$$ 
Further, we claim that $\sP$ is an $(\eps_{c(q+1)-1},(c(q+1)-1)t^{1/2},t,\bd)$-packing instance with suitable edge testers $\cW^q_{edge}$, 
where $$\bd=(d_A,d_B,(d_A^{m_i(q)})_{i\in N_{R_\ast^3}[q+1]\sm[q]},(d_B^{m_i(q)})_{i\in N_{R_\ast^3}[q+1]}).$$
To establish this claim, we first make some important observations.
By the definition of $c_i(q)$ in~\eqref{eq:def c_i(q)} and $c_I(q)$ in~\eqref{eq:def max c_i(q)}, we have:
\begin{align}\label{eq:change c_i}
\begin{minipage}[c]{0.9\textwidth}\em
If $i\in N_{R_\ast}(q+1)$, then $c(q+1)=c_i(q+1)>c_I(q)$ for every index set $I\subseteq[r]$ with $i\in I$.
\end{minipage}\ignorespacesafterend
\end{align}
Note that for the inequality in~\eqref{eq:change c_i} we used that an index set $I$ is contained in some hyperedge $\rR\in E(R)$, and no vertex of a hyperedge $\rR$ in $R$ has two neighbours in $R_\ast$ that are coloured alike as we have chosen the vertex colouring as a colouring in~$R_\ast^3$.\COMMENT{Note that we really need a colouring of $R_\ast^3$ for the inequality. 
}
In particular, we infer from~\eqref{eq:change c_i} that
\begin{align}\label{eq:eps changes}
\begin{minipage}[c]{0.9\textwidth}\em
for all $i\in N_{R_\ast}(q+1)$, we have $\eps_{c(q+1)-1}=\eps_{c_i(q+1)-1}\geq \eps_{c_i(q)}$ and $\eps_{c(q+1)}=\eps_{c_i(q+1)}$.
For every index set $I\subseteq[r]\sm[q+1]$ with $I\cap\rR\neq\emptyset$ for some $\rR\in E(R)$ and $q+1\in\rR$, we have $\eps_{c(q+1)-1}=\eps_{c_I(q+1)-1}\geq \eps_{c_I(q)}$ and $\eps_{c(q+1)}=\eps_{c_I(q+1)}$.
\end{minipage}\ignorespacesafterend
\end{align}
Similar, by the definition of $m_i(q)$ in~\eqref{eq:def m_i(q)}, we have:
\begin{align}\label{eq:change m_i}
\begin{minipage}[c]{0.9\textwidth}\em
If $i\in N_{R_\ast}(q+1)$, then $$m_i(q+1)=m_i(q)+
\Big|\Big\{\rR\in E(R)\colon \{q+1,i\}=\rR\cap\big(([r]\sm[q])\cup\{i\}\big) \Big\}\Big|.$$
If $i\in [r]\sm N_{R_\ast}(q+1)$, then $m_i(q+1)=m_i(q)$.
\end{minipage}\ignorespacesafterend
\end{align}

Hence to see that $\sP$ is an $(\eps_{c(q+1)-1},(c(q+1)-1)t^{1/2},t,\bd)$-packing instance, note that \ref{item:P1} follows from~\eqref{eq:G_Z sr}, 
property~\ref{item:P2} follows from~\ind{q}\ref{Z candidacy superregular}, 
property~\ref{item:P edge tester} follows from~\ind{q}\ref{edge weight},
property~\ref{item:P3.5} holds by the definition of the edge testers in~\eqref{eq:def initial weight boundedness},
and~\ref{item:P4} follows from~\ind{q}\ref{e,f suitable}.

\medskip
Observe further that by assumption, we have $|\cH|\leq n^{2k}$, and $e_\cH(\sX_{\sqcup \rR})\leq (1-\alpha)dn^k\leq d_An^k$ for all $\rR\in E(R_{q+1})$.
Hence, we can apply Lemma~\ref{lem:packing lemma} to $\sP$ with\COMMENT{I checked that again.}
\begin{center}
\begin{tabular}{r||c|c|c|c|c|c|c|c}
parameter& $n$  & $\eps_{c(q+1)-1}$ & $\eps_{c(q+1)}$  & $t$ &  $(c(q+1)-1)t^{1/2}$ & $\alpha^{-1}$ & $|N_{R_\ast^3}(q+1)\sm[q]|$ & $|\cQ|$ 
\\ \hline
plays the role of & $n$  & $\eps$ & $\eps'$ & $t$ &  $q$ & $s$ & $r$ & $r_\circ$ 
\vphantom{\Big(}
\end{tabular}
\end{center}

\medskip\noindent
and with 
\begin{itemize}
\item local $\alpha^{-1}$-testers in $\cW_{local}$ that we will define explicitly in Steps \ref{step:ind e,f suitable}--\ref{step:ind vertex tester} when establishing~\ind{q+1}\ref{e,f suitable}--\ref{vertex tester};
\item tuples $(\omega,c)$ in $\cW_0$ that we will define explicitly in Step~\ref{step:ind treffer auf G} when establishing~\ind{q+1}\ref{treffer auf g};
\item edge testers in $\cW_{edge}^{q}$.
\end{itemize}

Let $\sigma\colon\sX_{q+1}^\sigma\to V_{q+1}$ be the conflict-free packing in $\cA_{q+1}^q$ obtained from Lemma~\ref{lem:packing lemma} with $|\sX_{q+1}^\sigma\cap X_{q+1}^H|\geq (1-\eps_{c(q+1)})n$ for all $H\in\cH$, which extends $\phi_q$ to $\phi_{q+1}=\phi_q\cup\sigma$ with $\phi_{q+1}\colon \cX_q^{\phi_q}\cup\sX_{q+1}^\sigma\to \cV_{q+1}$, and let $M=M(\sigma)$ be the corresponding edge set to $\sigma$ defined as in~\eqref{eq:def M(sigma)}.
Further, let $\sigma^+$ be the cluster-injective extension of $\sigma$ obtained from Lemma~\ref{lem:packing lemma}.
Analogously, $\sigma^+$ extends $\phi_q^+$ to a cluster-injective extension $\phi_{q+1}^+:=\phi_q^+\cup\sigma^+$ of $\phi_{q+1}$ with $\phi_{q+1}^+\colon \cX_{q+1}\to\cV_{q+1}$.\COMMENT{That is,
\begin{align*}
\phi_{q+1}^+(x):=\begin{cases}
\phi_q^+(x) & \text{if $x\in \cX_q$,}\\
\sigma^+(x) & \text{if $x\in \sX_{q+1}$.}
\end{cases}
\end{align*}
}

For all $H\in\cH$, $(Z,\cZ)\in\{(A,\cA),(B,\cB) \}$, and all index sets $I_A\subseteq [r]\sm[q+1]$, {$I_B\in[r]$} with $I_Z\cap\rR\neq\emptyset$ for some $\rR\in E(R)$ with $q+1\in\rR$, let $Z_{I_Z}^{H,q+1}\subseteq Z_{I_Z}^H(\phi_q^+|_{\sX_{\cup \cQ}}\cup \sigma^+)=Z_{I_Z}^H(\phi_q^+\cup \sigma^+)$ be the candidacy graph $Z_{I_Z}^{H,new}=:Z_{I_Z}^{H,q+1}$ obtained from Lemma~\ref{lem:packing lemma} satisfying~\ref{emb lem:sr}--\ref{emb lem:local weight}.

For all  $H\in\cH$, $(Z,\cZ)\in\{(A,\cA),(B,\cB) \}$, and all index sets $I_A\subseteq [r]\sm[q+1]$, {$I_B\in[r]$} with $I_Z\cap\rR=\emptyset$ for all $\rR\in E(R)$ with $q+1\in\rR$, note that $m_i(q)=m_i(q+1)$ for all $i\in I_Z$ and $Z_{I_Z}^H(\phi_q^+)=Z_{I_Z}^H(\phi_{q+1}^+)$. 
Thus, in such a case we set $Z_{I_Z}^{H,q+1}:=Z_{I_Z}^{H,q}$. 
Let $\cZ_{I_Z}^{q+1}:=\bigcup_{H\in \cH}Z_{I_Z}^{H,q+1}$ and let $\cZ^{q+1}$ be the collection of all $\cZ_{I_Z}^{q+1}$ for all index sets $I_A\subseteq [r]\sm[q+1]$, {$I_B\in[r]$}.
We will employ Lemma~\ref{lem:packing lemma}\ref{emb lem:sr}--\ref{emb lem:local weight} to establish~\ind{q+1}\ref{Z candidacy superregular}--\ref{vertex tester}.

\begin{substep}\label{step:ind A candidacy sr}
Checking~\ind{q+1}{\rm \ref{Z candidacy superregular}} 
\end{substep}

We fix some $H\in\cH$, $Z\in\{A,B\}$ and establish~\ind{q}\ref{Z candidacy superregular} for the candidacy graph $Z_{i_Z}^{H,q+1}$.
For all $i_A\in N_{R_\ast}(q+1)\sm[q]$ and $i_B\in N_{R_\ast}(q+1)$, we have by our observation~\eqref{eq:change m_i} for $m_{i_Z}(q+1)$ that $d_A^{m_{i_Z}(q+1)}=d_{Z,{i_Z}}^{new}$ for $d_{Z,{i_Z}}^{new}$ in~\ref{emb lem:sr} as defined in~\eqref{eq:def b_i}.
Hence with~\ref{emb lem:sr},~\eqref{eq:change c_i} and~\eqref{eq:change m_i}, we obtain that $Z_{i_Z}^{H,q+1}$ is $(\eps_{c_{i_Z}(q+1)},d_Z^{m_{i_Z}(q+1)})$-super-regular and $(\eps_{c_{i_Z}(q+1)},c_{i_Z}(q+1)t^{1/2})$-well-intersecting for all $i_A\in N_{R_\ast}(q+1)\sm[q]$ and $i_B\in N_{R_\ast}(q+1)$.
For all $i_A\in[r]\sm (N_{R_\ast}(q+1)\cup[q+1])$ and $i_B\in[r]\sm N_{R_\ast}[q+1]$, we have $m_{i_Z}(q)=m_{i_Z}(q+1)$ and $Z^{H,q+1}_{i_Z}=Z^{H,q}_{i_Z}$. 
Hence with~\ind{q}\ref{Z candidacy superregular}, we also obtain in this case that $Z^{H,q+1}_{i_Z}$ is $(\eps_{c_{i_Z}(q+1)},d_Z^{m_{i_Z}(q+1)})$-super-regular and $(\eps_{c_{i_Z}(q+1)},c_{i_Z}(q+1)t^{1/2})$-well-intersecting.
This establishes~\ind{q+1}{\rm \ref{Z candidacy superregular}}.

\begin{substep}\label{step:ind e suitable}
Checking~\ind{q+1}{\rm \ref{edge weight}}\nopagebreak
\end{substep}

In order to establish~\ind{q+1}\ref{edge weight}, we fix $(\omega_{q+1},\omega_\iota,J,\JX,\JV,\bc,\pP)\in\cW^{q+1}_{edge}(\phi_{q+1},\phi_{q+1}^+,\cA^{q+1},\cB^{q+1})$ with centres $\bc\in V_{\sqcup I}$ for $I\subseteq [r]$, non-empty $I_{q+1}:=(I\sm[q+1])\cup J$, and patterns~$\pP=(\bpA,\bppA,\bpB,\bppB)$.

If $I_{q+1}\cap\rR=\emptyset$ for all $\rR\in E(R)$ with $q+1\in\rR$, {then the conflict-free packing $\sigma$ does not have an effect at all on the considered edge tester;} that is, $(\omega_{q},\omega_\iota,J,\JX,\JV,\bc,\pP)=(\omega_{q+1},\omega_\iota,J,\JX,\JV,\bc,\pP)$ by Definition~\ref{def:general edge tester}, and thus, \ind{q+1}\ref{edge weight} holds by \ind{q}\ref{edge weight}.
In this case, note in particular that if~$\pP$ is such that $\omega_\iota(\sX_{\sqcup I})>0$, then $\norm{\bpZ_{[q]}}=\norm{\bpZ_{[q+1]}}$ and $\norm{\bppZ_{[q]}}=\norm{\bppZ_{[q+1]}}$ for all $Z\in\{A,B\}$ because by Definition~\ref{def:pattern} of a \fst-pattern and \scd-pattern we have for the $(q+1)$th entries that 
$\bpZ_{q+1}=\bppZ_{q+1}=0$ as $I_{q+1}\cap\rR=\emptyset$ for all $\rR\in E(R)$ with $q+1\in\rR$.

Hence, we may assume that $I_{q+1}\cap\rR\neq\emptyset$ for some $\rR\in E(R)$ with $q+1\in\rR$.
{It is important to note that the edge tester $(\omega_{q+1},\omega_\iota,J,\JX,\JV,\bc,\pP)\in\cW^{q+1}_{edge}(\phi_{q+1},\phi_{q+1}^+,\cA^{q+1},\cB^{q+1})$ is defined in the same way as the edge tester $(\omega^{new},\omega_\iota,J,\JX, \JV,\bc,\pP)$ that we obtain from~\ref{emb lem:general edge testers}, and thus, they are identical.}
As in Step~\ref{step:ind A candidacy sr}, for $i\in I_{q+1}$, we have by our observation~\eqref{eq:change m_i} for $m_i(q+1)$ that $d_Z^{m_i(q+1)}=d_{i,Z}^{new}$ for $Z\in\{A,B\}$ and $d_{i,Z}^{new}$ in~\ref{emb lem:sr} as defined in~\eqref{eq:def b_i}.
Hence with~\eqref{eq:eps changes}, \eqref{eq:change m_i} and~\ref{emb lem:general edge testers}, we obtain that the edge tester $(\omega_{q+1},\omega_\iota,J,\JX,\JV,\bc,\pP)\in\cW^{q+1}_{edge}(\phi_{q+1},\phi_{q+1}^+,\cA^{q+1},\cB^{q+1})$ with respect to $(\omega_\iota,J,\JX,\JV,\bc,\pP)$, $\phi_{q+1}$, $\cA^{q+1}$ and $\cB^{q+1}$ satisfies \ind{q+1}\ref{edge weight}.

\begin{substep}\label{step:ind e,f suitable}
Checking~\ind{q+1}{\rm \ref{e,f suitable}}\nopagebreak
\end{substep}

In order to establish~\ind{q+1}\ref{e,f suitable}, we fix $\gG=\{v_{i_1},\ldots,v_{i_k}\}, \hG=\{w_{j_1},\ldots,w_{j_k}\}\in E(G)$ with $v_{i_k}=w_{j_k}$, $I:=\{i_1,\ldots,i_k\}\neq \{j_1,\ldots,j_k\}=:J$,  
and $\eps^\ast_q:=\max\{\eps_{c_{I\sm[q]}(q)},\eps_{c_{J\sm[q]}(q)}\}$.

If $q+1\notin I\cup J$, we have $|(I\cup J)\cap [q]|=|(I\cup J)\cap [q+1]|$, and thus~\ind{q+1}\ref{e,f suitable} holds by~\ind{q}\ref{e,f suitable}.

Hence, by symmetry, we may assume that $q+1=i_\ell$ for some $\ell\in[k]$.\COMMENT{Then~\eqref{eq:change c_i} implies that $\max\{ \eps_{c_{I\sm[q+1]}(q+1)},\eps_{c_{J\sm[q+1]}(q+1)} \}=\eps_{c(q+1)}$.}
Our strategy is to define a weight function that bounds from above the number of elements in $E_{\gG,\hG,\phi_{q}}(\cA^{q})$ that still can be present in $E_{\gG,\hG,\phi_{q+1}}(\cA^{q+1})$ by employing~\ref{emb lem:local weight}.\COMMENT{\label{codegree comment}We use that $\cA^{q+1}$ consists of subgraphs of $\cA^q$.}
For $(\eH,\fH)\in E_{\gG,\hG,\phi_q}(\cA^q)$, let $\{x_{i_\ell}\}:=\eH\cap \sX_{i_\ell}$, and we define a weight function $\omega_{\eH,\fH}\colon E(\sA_{q+1}^q)\to\{0,1\}$ by $\omega_{\eH,\fH}(xv):=\IND_{\{xv=x_{i_\ell}v_{i_\ell}\}}$. 
Note that $(\eH,\fH)\in E_{\gG,\hG,\phi_{q+1}}(\cA^{q+1})$ only if $\omega_{\eH,\fH}(M)=1$ as it is necessary that $\sigma$ embeds $x_{i_\ell}$ onto $v_{i_\ell}$.
Let $\omega_{\gG,\hG}:=\sum_{(\eH,\fH) \in E_{\gG,\hG,\phi_q}(\cA^q)} \omega_{\eH,\fH}.$
A~key observation is that\COMMENT{Note that we do not have to include a ``$\Gamma$''-term as in previous steps since in $E_{\gG,\hG,\phi_{q+1}}(\cA^{new})$ we only consider edges where all vertices are embedded by $\phi_{q+1}$ and thus, none of its vertices are left unembedded, see Definition~\ref{def:E_g,h,phi}.}
\begin{align}\label{eq:updated suitable e,f}
\left|E_{\gG,\hG,\phi_{q+1}}(\cA^{q+1})\right|\leq
\omega_{\gG,\hG}(M).
\end{align}
By the definition of $\omega_{\gG,\hG}$, we have that
\begin{align}\label{eq:w_b weight suitable e,f}
\omega_{\gG,\hG}(E(\sA_{q+1}^q))=|E_{\gG,\hG,\phi_q}(\cA^q)|.
\end{align}
By adding $\omega_{\gG,\hG}$ to $\cW_{local}$ and by employing~\ref{emb lem:local weight}, we obtain with~\eqref{eq:w_b weight suitable e,f} that
\begin{align}\label{eq:M-weight suitable e,f}
\omega_{\gG,\hG}(M)
= (1\pm\eps_{c(q+1)}^2)(d_A^{m_{q+1}(q)}n)^{-1}|E_{\gG,\hG,\phi_q}(\cA^q)| \pm n^{\eps_{c(q+1)-1}}.
\end{align}
We further observe that
\begin{align*}
|E_{\gG,\hG,\phi_q}(\cA^q)|
\stackrel{\text{\ind{q}\ref{e,f suitable}}}
\leq 
\max\left\{
n^{k-\left| (I\cup J)\cap [q]\right|+\eps_q^\ast},
n^{\eps_q^\ast}
\right\},
\end{align*}
and thus by~\eqref{eq:eps changes}, \eqref{eq:M-weight suitable e,f} and because $d\geq n^{-\eps}$, we finally obtain\COMMENT{$(d_A^{m_{q+1}(q)}n)^{-1}\leq 1/n^{1-\eps^{1/2}}$.}
\begin{align*}
\omega_{\gG,\hG}(M)
\leq
\max\left\{
n^{k-\left| (I\cup J)\cap [q]\right|-1+\eps_{c(q+1)}},
n^{\eps_{c(q+1)}}
\right\}.
\end{align*}
By our key observation~\eqref{eq:updated suitable e,f}, this establishes~\ind{q+1}\ref{e,f suitable}.\COMMENT{Recall that
\eqref{eq:change c_i} implies that $\max\{ \eps_{c_{I\sm[q+1]}(q+1)},\eps_{c_{J\sm[q+1]}(q+1)} \}=\eps_{c(q+1)}$ because $q+1=i_\ell$ for some $\ell\in[k]$.
}

\begin{substep}\label{step:ind treffer auf G}
Checking~\ind{q+1}{\rm \ref{treffer auf g}}\nopagebreak
\end{substep}

In order to establish~\ind{q+1}{\rm \ref{treffer auf g}}, we fix $(\omega,\bc)\in\cW_{hit}$ with $\bc\in V_{\sqcup I}$ and $J\subseteq I$.
In view of the statement, we may assume that $q+1\in J$, otherwise \ind{q+1}{\rm \ref{treffer auf g}} holds by \ind{q}{\rm \ref{treffer auf g}}.
Our general strategy is to define two weight functions and employ~\ref{emb lem:local weight} and~\ref{emb lem:hit tester} to derive the desired upper bound.

For $p\in\{q,q+1\}$, let $W_p:=\{
\xX\in \sX_{\sqcup {I}}\colon \xX\subseteq \eH \text{ for some $\eH\in E(\cH)$, }
\phi_p^+(\xX\cap\sX_{\cup (J\cap[p])})\subseteq\bc
\}$, and let \mbox{$c:=\bc\cap V_{q+1}$}.
We define a local tester $\omega_\sigma\colon E(\cA^{q+1})\to[0,\alpha^{-1}]$ by $$\omega_\sigma(uv):=
\sum_{\xX\in W_q\colon u\in \xX}\IND\{v=c \}\omega(\xX),$$ and we add $\omega_\sigma$ to $\cW_{local}$.
We also define a tuple $(\omega_{\sigma^+},c)$ with $\omega_{\sigma^+}\colon\sX_{q+1}\to[0,\alpha^{-1}]$ by $$\omega_{\sigma^+}(u):= \sum_{\xX\in W_q\colon u\in \xX} \omega(\xX),$$ and we add $(\omega_{\sigma^+},c)$ to~$\cW_0$.
We make the following observation
\begin{align}\label{eq:upperbound hit}
\omega(W_{q+1})
= 
\omega_\sigma(M)
+
\omega_{\sigma^+}(\{
x\in\sX_{q+1}\sm\sX_{q+1}^\sigma\colon \sigma^+(x)=c
\}).
\end{align}
We first employ~\ref{emb lem:local weight} to derive an upper bound on $\omega_\sigma(M)$ 
\begin{align}
\nonumber
\omega_\sigma(M)&
\stackrel{\rm\ref{emb lem:local weight}}{\leq}
2 (d_A^{m_{q+1}(q)}n)^{-1}\omega_\sigma(E(\cA^{q+1}))+n^{\eps_{c(q+1)-1}}
\\\nonumber
&\stackrel{\text{~\ind{q}\ref{treffer auf g}~}}{\leq}
2 (d_A^{m_{q+1}(q)}n)^{-1}
\frac{\omega(\sX_{\sqcup {I}})}{n^{|J\cap[q]|-\eps_{c(q+1)-1}}}+2n^{\eps_{c(q+1)-1}}
\\\label{eq:hit 1st bound}
&
\stackrel{\hphantom{\rm\ref{emb lem:local weight}}}{\leq}
\frac{\omega(\sX_{\sqcup {I}})}{n^{|J\cap[q+1]|-\eps^{1/2}_{c(q+1)-1}}}+2n^{\eps_{c(q+1)-1}}.
\end{align}
Next, we employ~\ref{emb lem:hit tester} to derive an upper bound for the last term of~\eqref{eq:upperbound hit}
\begin{align}\nonumber
\omega_{\sigma^+}(\{
x\in\sX_{q+1}\sm\sX_{q+1}^\sigma\colon \sigma^+(x)=c
\})
&\stackrel{\ref{emb lem:hit tester}}{\leq}
\omega_{\sigma^+}(\sX_{q+1})/n^{1-\eps_{c(q+1)-1}}+n^{\eps_{c(q+1)-1}}
\\\label{eq:hit 2nd bound}
&\stackrel{\text{~\ind{q}\ref{treffer auf g}~}}{\leq}
\frac{\omega(\sX_{\sqcup {I}})}{n^{|J\cap[q+1]|-2\eps_{c(q+1)-1}}}+2n^{\eps_{c(q+1)-1}}.
\end{align}

Finally, plugging~\eqref{eq:hit 1st bound} and~\eqref{eq:hit 2nd bound} into~\eqref{eq:upperbound hit}, yields that
\begin{align*}
\omega(W_{q+1})\leq 
\frac{\omega(\sX_{\sqcup {I}})}{n^{|J\cap[q+1]|-\eps_{c(q+1)}}}+n^{\eps_{c(q+1)}}.
\end{align*}
Together with~\eqref{eq:eps changes}, this establishes~\ind{q+1}{\rm \ref{treffer auf g}}.

\begin{substep}\label{step:ind G_B leftover for typicality}
Checking~\ind{q+1}{\rm \ref{G_B leftover for typicality}}\nopagebreak
\end{substep}

In order to establish~\ind{q+1}{\ref{G_B leftover for typicality}}, we fix $H\in\cH$ and $\cS\subseteq 
\bigcup_{\rR\in E(R)\colon q+1\in\rR} V_{\sqcup\rR\sm\{q+1\}}$ with $|\cS|\leq t$.
Let $W:=V_{q+1} \cap N_{G_B}(\cS)$.
Our general strategy is to define a weight function that estimates the number of vertices in $W\cap \sigma(X_{q+1}^H)$ from which we can derive an upper bound for $|W\sm \sigma(X_{q+1}^H)|$.

Hence, let $\omega_W\colon E(\cA_{q+1}^q)\to\{0,1\}$ be defined by $\omega_W(e):=\IND\{w\in e, w\in W\}$.
A key observation is that
\begin{align}\label{eq:G_B-leftover key observation}
|W\sm \sigma(X_{q+1}^H)|\leq |W|-\omega_W(M).
\end{align}
By the definition of $\omega_W$ and because $A_{q+1}^{H,q}$ is $(\eps_{c_{q+1}(q)},d_A^{m_{q+1}(q)})$-super-regular for every $H\in\cH$ by \ind{q}\ref{Z candidacy superregular}, we have that
\begin{align}\label{eq:weight omega_W}
\omega_W(E(\cA_{q+1}^q))
=
(1\pm 3\eps_{c_{q+1}(q)})d_A^{m_{q+1}(q)}n|W|.
\end{align}
By adding $\omega_{W}$ to $\cW_{local}$ and by employing~\ref{emb lem:local weight}, we obtain with~\eqref{eq:weight omega_W} that
\begin{align*}
\omega_W(M)
&\stackrel{\ref{emb lem:local weight}}{=}
(1\pm\eps_{c(q+1)}^2)
\frac{\omega_{W}(E(\cA_{q+1}^q))}{d_A^{m_{q+1}(q)} n}\pm n^{\eps_{c_{q+1}(q)}}
\stackrel{\eqref{eq:weight omega_W}}{=}
(1\pm\eps_T)|W|.
\end{align*}
Now, this together with~\eqref{eq:G_B-leftover key observation} implies that $|W\sm \sigma(X_{q+1}^H)|
\leq 
\eps_T|W|$.
Together with~\ind{q}\ref{G_B leftover for typicality} this establishes~\ind{q+1}\ref{G_B leftover for typicality}.

\begin{substep}\label{step:ind set tester}
Checking~\ind{q+1}{\rm \ref{set tester}}\nopagebreak
\end{substep}

Let $(W,Y_1,\ldots,Y_\ell)\in\cW_{set}$ be a set tester with $W\subseteq V_{q+1}$ and $Y_j\subseteq X_{q+1}^{H_j}$ for all $j\in[\ell]$.
We define 
\begin{align*}
E_{(W,Y_1,\ldots,Y_\ell)}:=
\bigg\{
\{e_1,\ldots,e_\ell\}\in\bigsqcup_{j\in[\ell]} E\left(A_{q+1}^{H_j,q}[W,Y_j]\right)
\colon \bigcap_{j\in [\ell]}e_j\neq\emptyset
\bigg\}
\end{align*}
and a weight function $\omega_{{(W,Y_1,\ldots,Y_\ell)}}\colon\binom{E(\cA_{q+1}^q)}{\ell}\to\{0,1\}$ by $\omega_{{(W,Y_1,\ldots,Y_\ell)}}(\be):=\IND{\{\be\in E_{(W,Y_1,\ldots,Y_\ell)}  \}}$.
Note that 
\begin{align}\label{eq:weight set tester}
\textstyle\omega_{{(W,Y_1,\ldots,Y_\ell)}}(M)=\big|W\cap\bigcap_{j\in[\ell]}\sigma(Y_j)\big|.
\end{align}
In view of the statement, we may assume that $|W|,|Y_j|\geq\eps_{c(q+1)}n$ for all $j\in[\ell]$.\COMMENT{Otherwise $|W\cap \bigcap_{j\in [\ell]}\sigma(Y_j\cap \sX_0^\sigma)|\leq \eps_{c(q+1)} n$ and $\eps_{c(q+1)} n= |W||Y_1|\cdots |Y_\ell|/n^\ell \pm \alpha n/2$.}
Since $A_{q+1}^{H,q}$ is $(\eps_{c_{q+1}(q)},d_A^{m_{q+1}(q)})$-super-regular for every $H\in\cH$ by \ind{q}\ref{Z candidacy superregular}, we obtain by Fact~\ref{fact:regularity} and because $\ell\leq \alpha^{-1}$ that there are at most $\eps_{c_{q+1}(q)}^{1/2}n$ vertices in~$W$ that do not have $(1\pm\eps_{c_{q+1}(q)})d_A^{m_{q+1}(q)}|Y_j|$ neighbours in $Y_j$ for every $j\in[\ell]$.
Hence, we obtain that
\begin{align}\label{eq:weight W_set}
\omega_{(W,Y_1,\ldots,Y_\ell)}(E(\sA_{q+1}^q))=|E_{(W,Y_1,\ldots,Y_\ell)}|=(1\pm\eps_{c(q+1)})(d_A^{m_{q+1}(q)})^\ell|W||Y_1|\cdots|Y_\ell|.
\end{align}
We check that $\omega_{(W,Y_1,\ldots,Y_\ell)}$ is a local tester: 
For all  $\{e_1,\ldots,e_{\ell'}\}\in \binom{E(\cA_{q+1}^q)}{\ell'}$, $\ell'\in[\ell]$, 
the number of edges $\{e_{\ell'+1},\ldots,e_\ell \}$ such that $\be:=\{e_j\}_{j\in[\ell]}\in\binom{E(\cA_{q+1}^q)}{\ell}$ with $\omega_{(W,Y_1,\ldots,Y_\ell)}(\be)$ is at most $(2n)^{\ell-\ell'}$,\COMMENT{Let $\{e_1,\ldots,e_{\ell'}\}$ be $\ell'$ edges all containing one specific $x\in W$.
This also fixes the endpoint $x$ for the edges $\{e_{\ell'+1},\ldots,e_\ell\}$. For each of these $\ell-\ell'$ edges, there are at most $2n$ choices for its other endpoint.} implying that $\norm{\omega_{(W,Y_1,\ldots,Y_\ell)}}_{\ell'}\leq n^{\ell-\ell'+\eps_{c_{q+1}(q)}^2}$.

Hence, $\omega_{(W,Y_1,\ldots,Y_\ell)}$ is a local tester and we can add $\omega_{E_{(W,Y_1,\ldots,Y_\ell)}}$ to $\cW_{local}$.
By~\ref{emb lem:local weight}, we conclude that  
\begin{align*}
\omega_{(W,Y_1,\ldots,Y_\ell)}(M)
=(1\pm\eps_{c(q+1)}^2)
\frac{\omega_{(W,Y_1,\ldots,Y_\ell)}(E(\sA_{q+1}^q))}
{(d_A^{m_{q+1}(q)}n)^\ell}
\pm n^{\eps_{c(q+1)}}
\stackrel{\eqref{eq:weight W_set}}{=}\frac{|W||Y_1|\cdots|Y_\ell|}{n^\ell}\pm \alpha^2 n,
\end{align*}
which establishes~\ind{q+1}\ref{set tester} by~\eqref{eq:weight set tester}.

\begin{substep}\label{step:ind vertex tester}
Checking~\ind{q+1}{\rm \ref{vertex tester}}\nopagebreak
\end{substep}

In order to establish~\ind{q+1}{\rm \ref{vertex tester}}, let $(\omega_{ver},\bc)\in\cW_{ver}$ be an $\alpha^{-1}$-vertex tester with centres $\bc=\{c_i\}_{i\in I}\in V_{\sqcup I}$ where $I\subseteq [q+1]$ and $q+1\in I$.
By~\eqref{eq:def initial weight vertex tester}, we defined in particular for all $\bpA,\bppA\in[k\alpha^{-1}]_0^r$ and $\pP:=(\bpA,\bppA,\mathbf{0},\mathbf{0})$ a tuple $(\omega_{\iota},J=\emptyset,\JX=\emptyset,\JV=\emptyset,\bc,\pP)$ with initial weight function $\omega_\iota=:\omega_{\iota,\pP}$ corresponding to $(\omega_{ver},\bc)$. That is, by~\eqref{eq:def initial weight vertex tester}, we have
\begin{align}\label{eq:relation weight w_0 w_ver}
\omega_{ver}(\sX_{\sqcup {I}})
=
\sum_{\pP\colon\bpA,\bppA\in[k\alpha^{-1}]_0^r} \omega_{\iota,\pP}(\sX_{\sqcup {I}}).
\end{align}
Note that $\omega_{\iota,\pP}(\sX_{\sqcup {I}})>0$ only if {$\norm{\bpA_{[q]}}=\norm{\bppA_{[q]}}$} by~\eqref{eq:norms equal} since $|I|\leq k-1$\COMMENT{$|I|\leq k-1$ implies that the indicator function in~\eqref{eq:norms equal} equals $0$ and thus, $\norm{\bpA_{[q]}}=\norm{\bppA_{[q]}}$.
We defined vertex testers such that $|I|\leq k-1$, see comment~\ref{comment size I vertex tester}.}\COMMENT{Note that $q+1\in I$. 
Nevertheless we have $\norm{\bpA_{[q]}}=\norm{\bppA_{[q]}}$ because $\bpA_{q+1}=\bppA_{q+1}=0$ since the last entry of a \fst-pattern and \scd-pattern is always $0$. Thus, $d_A^{\norm{\bpA_{[q]}}-\norm{\bppA_{[q]}}}=d_A^0=1$}.
For $\bpA,\bppA\in[k\alpha^{-1}]_0^r$ and $\pP:=(\bpA,\bppA,\mathbf{0},\mathbf{0})$ with {$\norm{\bpA_{[q]}}=\norm{\bppA_{[q]}}$}, let
$\tau_\pP:=(\omega_{q},\omega_{\iota},J=\emptyset,\JX=\emptyset,\JV=\emptyset,\bc,\pP)$ be the edge tester with respect to $(\omega_{\iota},J=\emptyset,\JX=\emptyset,\JV=\emptyset,\bc,\pP)$, $(\phi_{q},\phi_q^+)$, $\cA^q$ and $\cB^q$, which is contained in $\cW_{edge}^q$, and $\omega_{q}\colon E(\cA_{q+1}^q)\to[0,\alpha^{-1}]$,
and let $\omega_{q,\pP}:=\omega_q$.
By~\ind{q}\ref{edge weight}, we obtain\COMMENT{Since $q+1\in I$ and $I\subseteq[q+1]$, we have that $I_q=(I\sm[q])\cup J=\{q+1\}$ in \ind{q}\ref{edge weight}, and thus
$$\prod_{i\in I_q\sm J}d_{A}^{m_i(q)}=d_A^{m_{q+1}(q)}.$$}
\begin{align}\label{eq:weight w_q,p,p}
\omega_{q,\pP}(E(\cA_{q+1}^q))=
(1\pm\eps_{c_{q+1}(q)})
d_A^{m_{q+1}(q)}
\frac{\omega_{\iota,\pP}(\sX_{\sqcup {I}})}{n^{|I|-1}}\pm n^{\eps_{c_{q+1}(q)}} .
\end{align}
A key observation is that
\begin{align}\label{eq:key vertex tester}
\omega_{ver}(\phi_{q+1}^{-1}(\bc))=
\sum_{\pP\colon\bpA,\bppA\in[k\alpha^{-1}]_0^r}
 \omega_{q,\pP}(M),
\end{align}
which follows from the definition of the edge tester $\tau_\pP$ for $\pP=(\bpA,\bppA,\mathbf{0},\mathbf{0})$ as in Definition~\ref{def:general edge tester}.\COMMENT{In particular, this is a simple edge tester because $J=\emptyset,\JX=\emptyset,\JV=\emptyset$.}
Note that $\cA_{q+1}^q$ is $(\eps_{c_{q+1}(q)},d_A^{m_{q+1}(q)})$-super-regular by \ind{q}\ref{Z candidacy superregular}.
For all $\bpA,\bppA\in[k\alpha^{-1}]_0^r$ with {$\norm{\bpA_{[q]}}=\norm{\bppA_{[q]}}$} and $\pP:=(\bpA,\bppA,\mathbf{0},\mathbf{0})$, we add $\omega_{q,\pP}\colon E(\cA_{q+1}^q)\to[0,\alpha^{-1}]$ to $\cW_{local}$ and obtain by~\ref{emb lem:local weight} that
\begin{align*}
\omega_{q,\pP}(M)
&
\stackrel{\ref{emb lem:local weight}}{=}
(1\pm\eps_{c(q+1)}^2)
\frac{\omega_{q,\pP}(E(\cA_{q+1}^q))}{d_A^{m_{q+1}(q)} n}\pm n^{\eps_{c_{q+1}(q)}}
\\&
\stackrel{~\eqref{eq:weight w_q,p,p}~~}{=}
(1\pm\eps_{c(q+1)})\frac{\omega_{\iota,\pP}(\sX_{\sqcup {I}})}{n^{|I|}}\pm 2n^{\eps_{c_{q+1}(q)}}.
\end{align*}
Together with~\eqref{eq:relation weight w_0 w_ver} and~\eqref{eq:key vertex tester}, this implies that
\begin{align*}
\omega_{ver}( \phi_{q+1}^{-1}(\bc))=
(1\pm\eps_T)\omega_{ver}(\sX_{\sqcup {I}})/n^{|I|}\pm n^{\eps_{T}},
\end{align*}
which establishes~\ind{q+1}\ref{vertex tester}.

\begin{step}\label{step:completion}
Completion
\end{step}

Let $\phi_r\colon \bigcup_{H\in \cH, i\in[r]}\hX_i^H\to\cV_r$ be an $r$-partial packing satisfying~\ind{r} with $(\eps_T,d_i)$-super-regular and $(\eps_T,t^{2/3})$-well-intersecting\COMMENT{By~\ind{r}\ref{Z candidacy superregular} we have that $B_i^H$ is $(\eps_T,d_i)$-super-regular and $(\eps_T,Tt^{1/2})$-well-intersecting. } candidacy graphs $B_i^H:=B_i^{H,r}\subseteq B_i^H(\phi_r^+)$ where $d_i:=d_B^{\dg_R(i)}=d_B^{m_i(r)}$ for all $i\in[r]$ and $\cX_r^{\phi_r}=\bigcup_{H\in \cH, i\in[r]}\hX_i^H$. 
We will apply a random packing procedure in order to complete the partial packing~$\phi_r$ using the edges in~$G_B$.
Recall that $\eps_T\ll\mu\ll\gamma\ll1/t\ll\alpha,1/k$ and we often call the vertices $\bigcup_{H\in \cH, i\in[r]}(X_i^H\sm\hX_i^H)$ \emph{unembedded (by $\phi_r$)} or the \emph{leftover (of $\phi_r$)}. 
Our general strategy is as follows.
For every $H\in\cH$ in turn, we choose a set $Y_i^H\subseteq \hX_i^H$ for all $i\in[r]$ of size roughly~$\mu n$ by selecting every vertex uniformly at random with probability $\mu$ and adding $X_i^H\sm\hX_i^H$ deterministically. 
Afterwards we apply a random matching argument to pack $H[Y^H_{\cup[r]}]$ into~$G_B$, which together with~$\phi_r$ yields a complete packing of~$H$ into~$G_A\cup G_B$.
Before we proceed with the details of our random packing procedure in Claim~\ref{claim:random packing} (Steps~\ref{step:random packing procedure}--\ref{step:few activated weight on vertex testers}), we verify in Claim~\ref{claim:matching} (Steps~\ref{step:matching argument}--\ref{step:claim step2}) that we can indeed pack a subgraph of one single $H\in\cH$ into~$G_B$ using another random embedding argument as long as our random packing procedure does not deviate too much from its expected behaviour.
To that end, we collect some more notation in Step~\ref{step:notation for completion}, and establish several important leftover conditions in Step~\ref{step:leftover conditions}.

\begin{substep}\label{step:notation for completion}
Notation for the completion
\end{substep}

We introduce some more notation.
We arbitrarily enumerate the graphs in~$\cH$ and write $\cH=\{H_1,\ldots,H_{|\cH|} \}$.
For $G^\circ\subseteq G_B$ and $B_i^H$ with $H\in\cH$, $i\in[r]$,
\begin{align}\label{eq:def B^G}
\begin{minipage}[c]{0.91\textwidth}\em
let $(B_i^H)^{G^\circ}$ be the subgraph of $B_i^H$ where $N_{(B_i^H)^{G^\circ}}(x)=V_i\cap N_{G_B-G^\circ}(\cS_x)$ for every $x\in X_i^H$,
\end{minipage}\ignorespacesafterend
\end{align}
and $\cS_x\subseteq \binom{V(G)}{k-1}$ is the set such that $N_{B_i^H}(x)=N_{G_B}(\cS_x)$. 
Note that $\cS_x$ exists because $B_i^H$ is $(\eps_T,t^{2/3})$-well-intersecting.
This implies in particular that we removed every edge $xv$ from $B_i^H$ for which there exists an edge $\eH\in E(H)$ such that $\phi_r(\eH\sm\{x\})\cup\{v\}\in E(G^\circ)$.\COMMENT{Indeed it suffices to only considers $\phi_r$. Note that we also remove the candidates used by `bad' edges in a second step, that is, the edges with respect to $\phi_q^{bad}$, see~\eqref{eq:def B-candidacy}. In the end, we will consider a subgraph of $(B_{q+1}^H)^{G^\circ}(\phi_q^{bad})$, see display after~\indC{q}.}
We may think of $E(G^\circ)$ as the edge set in $G_B$ that we have already used in our completion step for packing some other graphs of $\cH$ into $G_A\cup G_B$. 
Consequently, $(B_i^H)^{G^\circ}$ is the subgraph of the candidacy graph $B_i^H$ that only contains an edge $xv$ if we do not use an edge in $G^\circ$ when we would map $x$ onto~$v$.
To count the number of removed edges incident to a vertex $v\in V_i$ in $B_i^H$, we define
\begin{align}\label{eq:def 1 degree}
\rho^H_{G^\circ}(v):=\big|\big\{
x\in N_{B_i^H}(v)\colon S \cup\{v\}\in E(G^\circ)\text{ for some } S\in\cS_x 
 \big\}\big|,
\end{align}
where $\cS_x\subseteq \binom{V(G)}{k-1}$ is the set such that $N_{B_i^H}(x)=N_{G_B}(\cS_x)$.
Note that
\begin{align}\label{eq:v removed vertices B^circ}
\dg_{B_i^H}(v)-\dg_{(B_i^H)^{G^\circ}}(v)
= \rho^H_{G^\circ}(v) \text{ for every $v\in V_i$, $i\in[r]$.}
\end{align}

During our random packing procedure we will guarantee that $\rho^H_{G^\circ}(v)$ is negligibly small (Step~\ref{step:codegree bounds2}) and that the probability for a $G_B$-edge to be `used' during the completion is appropriately small (at most $\mu^{3/4}$, see Claim~\ref{claim:prob to use G_B edge}).
To that end, we control certain conditions for the leftover of $\phi_r$ in the next step.\looseness=-2

\begin{substep}\label{step:leftover conditions}
Controlling the leftover for the completion \nopagebreak
\end{substep}

In this step we make some important observations for the completion. 
In general, we will employ~\ind{r}\ref{edge weight} multiple times in order to suitably control the leftover, that is, the structure of the vertices that are left unembedded by $\phi_r$.

\medskip
We start with an observation how to control the number of neighbours of a vertex $v$ in $B_j^H$ that are left unembedded by $\phi_r$.
\begin{claim}\label{claim:leftover obs}
$|N_{B_j^H}(v)\cap (X_j^H\sm \hX_j^H)|
\leq 2\eps_T d_jn$ for all $H\in\cH$, $j\in[r]$ and $v\in V_j$.
\end{claim}

\claimproof
Recall that we defined edge testers in~\eqref{eq:def edge tester leftover} for all $H\in\cH$, $j\in[r]$ and $v\in V_j$ to count the number of neighbours of $v$ in $B_j^H$ that are left unembedded by $\phi_r$.
Let $(\omega,\omega_{\iota,H},J=\{j\},\JX=\{j\},\JV=\emptyset,\bc=\{v\},(\mathbf{0},\mathbf{0},\mathbf{0},\mathbf{0}))$ be the edge tester in $\cW_{edge}^r(\phi_r,\phi_r^+,\cA^r,\cB^r)$ that we obtain from~\ind{r}\ref{edge weight}.
Definition~\ref{def:general edge tester} of this edge tester implies that
\begin{align*}
|N_{B_j^H}(v)\cap (X_j^H\sm \hX_j^H)|
=
\omega(E(\cB_j^r))
=
\omega(E(B_j^H))
\stackrel{{\text{\ind{r}\ref{edge weight}}} 
}{\leq}
2\eps_T d_jn.
\end{align*}
This establishes Claim~\ref{claim:leftover obs}.
\endclaimproof

Next, we control the number of neighbours of a vertex $v$ in $B_i^H$ that are embedded but lie in an $H$-edge that contains unembedded vertices.\COMMENT{Recall that $H_\ast$ is the underlying $2$-graph of $H$.}
\begin{claim}\label{claim:leftover obs H-neighbours of v}
$|N_{B_i^H}(v)\cap \cX_r^{\phi_r}\cap N_{H_\ast}(\cX_r\sm \cX_r^{\phi_r})|
\leq \eps_T^{1/2}d_in$ for all $H\in\cH$, $i\in[r]$ and $v\in V_i$.
\end{claim}

\claimproof
Our general strategy is as follows.
We do not only consider a single vertex $v\in V_i$ but also a second vertex $w\in V_{j}$ for some $j\in N_{R_\ast}(i)$.
We can use our defined edge testers and employ~\ind{r}\ref{edge weight} to count for a fixed vertex $w$ the number of $2$-sets $\{x_i,x_{j}\}$ where $x_i\in N_{B_i^H}(v)$ and $x_{j}\in X_j^H$ is left unembedded.
Hence, by summing over all possible choices of $w$, we count all such $2$-sets but multiple times. 
Hence, by a double counting argument, for one fixed $j$ and all choices of $w\in V_{j}$, we can establish an upper bound for $|N_{B_i^H}(v)\cap \cX_r^{\phi_r}\cap N_{H_\ast}(X_{j}^H \sm \cX_r^{\phi_r})|$. 
In the end, this implies Claim~\ref{claim:leftover obs H-neighbours of v} as there are at most $k\alpha^{-1}$ choices for $j$.
We proceed with the details.

We fix  $H\in\cH$, $i\in[r]$, $v\in V_i$, $j\in N_{R_\ast}(i)$ and consider $w\in V_{j}$.
For all $\pP=(\mathbf{0},\mathbf{0},\bpB,\bppB)$, we defined a tuple $(\omega_\iota,J=\{i,j\}, J_X=\{j\}, J_V=\emptyset, \{v,w\}, \pP)$ in~\eqref{eq:def initial weight H-edge leftover}. 
Let $(\omega,\omega_{\iota},J,\JX, \JV,\{v,w\},\pP)$ be the edge tester in $\cW_{edge}^r(\phi_r,\phi_r^+,\cA^r,\cB^r)$ with respect to $(\omega_{\iota},J,\JX, \JV,\{v,w\},\pP)$, {$(\phi_r,\phi_r^+)$}, $\cA^r$, and $\cB^r$ that we obtain from~\ind{r}\ref{edge weight}. 
In order to be able to distinguish these edge testers according to the pattern vector $\pP$, we write $\omega_{\iota,\pP}:=\omega_\iota$ and $\omega_{\pP}:=\omega$.
Note that $\sum_{\pP}\omega_{\iota,\pP}(X_i^H\sqcup X_j^H)\leq 2\alpha^{-1}n$ by~\eqref{eq:def initial weight H-edge leftover}.
Thus, by Definition~\ref{def:general edge tester} of an edge tester and by summing over all patterns $\pP=(\mathbf{0},\mathbf{0},\bpB,\bppB)$ with $\bpB,\bppB\in[k\alpha^{-1}]^r_0$, we can employ conclusion~\ind{r}\ref{edge weight} to count the tuples $\{x_i,x_{j}\}\in X_i^H\sqcup X_{j}^H$ where $x_{j}$ is left unembedded but $\{x_i,x_{j}\}$ could still be mapped onto $\{v,w\}$; that is, {$\{\{x_i,v\},\{x_{j},w\} \}\in E(B_i^H)\sqcup E(B_j^H)$}.
Note that we count each such tuple $\{x_i,x_{j}\}$ multiple times, namely, for every $w\in V_{j}$ such that {$\{\{x_i,v\},\{x_{j},w\} \}\in E(B_i^H)\sqcup E(B_j^H)$}. 
By the Definition~\ref{def:candidacy graph} of the candidacy graphs and because~$B_{j}^H$ is $(\eps_T,d_{j})$-super-regular, there are $|N_{B_{j}^H}(x_{j})|=(1\pm2\eps_T)d_{j}n$ choices for $w\in V_{j}$ such that {$\{\{x_i,v\},\{x_{j},w\} \}\in E(B_i^H)\sqcup E(B_j^H)$}.\COMMENT{In order that {$\{\{x_i,v\},\{x_{j},w\} \}\in E(B_i^H)\sqcup E(B_j^H)$} we need by~\eqref{eq:condition updating} that $\phi_r^+(\eH_\circ)\cup\{v,w\}\in E(G_B)$ for every relevant $(k-2)$-set $\phi_r^+(\eH_\circ)$ (which might even include dummy sets).
Since $B_j^H$ is $(\eps_T,t^{2/3})$-well-intersecting, there is a set $\cS_{x_j}$ such that $N_{B_j^H}(x_j)=V_j\cap N_{G_B}(\cS_{x_j})$. We can use $\cS_{x_j}$ to obtain the number of $w$'s in $V_j$ such that $\phi_r^+(\eH_\circ)\cup\{v,w\}\in E(G_B)$ as follows:
Modifying every $S\in\cS_{x_j}$ to $S':=(S\sm\phi_r(x_i))\cup \{v\}$ yields a set of $(k-1)$-sets $\cS_{x_j}'$ and $V_j\cap N_{G_B}(\cS_{x_j})$ contains all $w$'s in $V_j$ such that $\phi_r^+(\eH_\circ)\cup\{v,w\}\in E(G_B)$. 
Since $|\cS_{x_j}|=|\cS_{x_j}'|$, this yields that there are $|N_{B_{j}^H}(x_{j})|=|V_j\cap N_{G_B}(\cS_{x_j})|=(1\pm2\eps_T)d_{j}n$ choices for $w$.} 
Altogether, this implies that
\begin{align*}
\left|
N_{B_i}(v)\cap \cX_r^{\phi_r}\cap N_{H_\ast}(X_{j}^H\sm\cX_r^{\phi_r})
\right| 
&\stackrel{\hphantom{{\text{\ind{r}\ref{edge weight}}}}
}{\leq} 
\left((1-2\eps_T)d_{j}n\right)^{-1}
\sum_{w\in V_{j}}\sum_{\pP} \omega_\pP\left(E(\cB_{i}^r)\sqcup E(\cB_{j}^r)\right) 
\\&
\stackrel{{\text{\ind{r}\ref{edge weight}}} 
}{\leq}
\left((1-2\eps_T)d_{j}n\right)^{-1}
2n \sum_{\pP}\left( \eps_Td_id_{j}\omega_{\iota,\pP}(\sX_{\sqcup J})+n^{\eps_T}\right) 
\\&
\stackrel{\hphantom{{\text{\ind{r}\ref{edge weight}}}}
}{\leq} 
4d_{j}^{-1}\cdot \left(2\alpha^{-1}\eps_Td_id_{j}n+n^{\eps_T}\right)
\\&
\stackrel{\hphantom{{\text{\ind{r}\ref{edge weight}}}}
}{\leq} \eps_T^{2/3}d_in.
\end{align*} 
\COMMENT{Note that we used in the second line that $d_B^{\norm{\bpB_{[r]}}-\norm{\bppB_{[r]}}}=0$ because $\omega_{\iota,\pP}>0$ only if $\norm{\bpB_{[r]}}=\norm{\bppB_{[r]}}$ by~\eqref{eq:norms equal} since we only consider two clusters $i,j$.
\\However, there is one technical detail: if $k=2$, then we obtain that $d_B^{\norm{\bpB_{[r]}}-\norm{\bppB_{[r]}}}=d_B^{-1}$, which adds a multiplicative factor of $d_B^{-1}$ in line 2. 
But if $k=2$, then we must also take the summation $\sum_{w\in N_{G_B}(v)}$ instead of  $\sum_{w\in V_{j}}$, and $\sum_{w\in N_{G_B}(v)\cap V_j}1\leq 2d_Bn$. 
I do not want to mention and confuse the reader with this...
}Summing over all $j\in N_{R_\ast}(i)$ establishes Claim~\ref{claim:leftover obs H-neighbours of v}.
\endclaimproof

The last observation for controlling the leftover concerns the number of $\cH$-edges that contain unembedded vertices and could still be mapped onto an edge $\gG\in E(G_B)$.
We define the following set in a slightly more general way as we will use this definition again in Step~\ref{step:few activated weight on vertex testers}, {where we also consider subsets of edges}.
For all $\rR\in E(R)$, $I\subseteq \rR$, $\bc\in V_{\sqcup I}$, non-empty $J\subseteq I$, and all pairs of disjoint sets $\JX, \JV\subseteq J$, let
\begin{align}\label{eq:def leftover g-edge}
E_{\phi_r}(\bc,J,\JX, \JV)
:=
\bigcup_{H\in \cH}
\Big\{
&\xX\in X_{\sqcup I}^H\colon 
\xX\subseteq \eH \text{ for some $\eH\in E(\cH)$,}
\bc\cap V_{\cup(I\sm J)}\subseteq\phi_r(\xX),
\\&\nonumber
 \phi_r(\xX)\neq\bc,
\bc\cap V_{\cup \JV}\subseteq\bc\sm \phi_r(\hX_{\cup I}^H),
\xX\sm \cX_r^{\phi_r} = \xX\cap X_{\cup \JX}^H, 
\\&\nonumber
\{\xX\cap X_j^H,\bc\cap V_j\}\in E(B_j^H) 
\text{ for all $j\in J$}
\Big\}.
\end{align}
Let us first explain this definition in words for the following more special case. 
For $I=\rR\in E(R)$, $\bc=\gG\in E(G_B)$, $\JX, \JV, J$ as above, $E_{\phi_r}(\gG,J,\JX, \JV)$ is the set of all $H$-edges $\eH$ for all $H\in\cH$ such that
\begin{enumerate}[label={\rm (\roman*)\textsubscript{\eqref{eq:def leftover g-edge}}}]
\item\label{item:def1 leftover} $\{\eH\cap X_j^H,\gG\cap V_j\}$ is an edge in $B_j^H$ for every $j\in J$,
\item the $k-|J|$ vertices $\eH\cap X^H_{\cup {(\rR\sm J)}}$ are mapped onto $\gG\cap V_{\cup (\rR\sm J)}$, 
\item  $\gG\cap V_{\cup \JV}$ is a subset of the vertices of $\gG$ onto which no $H$-vertex is embedded by $\phi_r$, and
\item\label{item:def4 leftover} all vertices in $\eH$ but the $|\JX|$ vertices $\eH\cap X^H_{\cup \JX}$ of $\eH$ are embedded by $\phi_r$.
\end{enumerate} 
This means, that if we modified $\phi_r$ and allowed the vertices $\eH\cap X_{\cup (J\sm \JX)}^H$ to be embedded somewhere else, we could potentially embed $\eH$ onto $\gG$.
Note that for an edge $\gG\in E(G_B)$, we naturally have that $\phi_r(\eH)\neq\gG$ because $\phi_r$ maps $\cH$-edges onto $G_A$-edges.

\begin{claim}\label{claim:leftover for g}
$\left|E_{\phi_r}(\gG,J,\JX, \JV)\right|
\leq
\left(\IND\{\JXV=\emptyset \}+2\eps_T\right)\gamma^{-1}
 n^{|J|}
\prod_{j\in  J} d_j$,
 for all $\rR\in E(R)$, $\gG\in E(G_B[V_{\cup\rR}])$, non-empty $J\subseteq\rR$, all pairs of disjoint sets $\JX, \JV\subseteq J$, and $\JXV=\JX\cup\JV$.
Note that $\eps_T\gamma^{-1}\leq \eps_T^{1/2}$. Hence 
$\left|E_{\phi_r}(\gG,J,\JX, \JV)\right|
\leq
\eps_T^{1/2}
 n^{|J|}
\prod_{j\in  J} d_j$
if $\JXV\neq\emptyset$.
\end{claim}

\claimproof
We fix $\rR=I$, $\gG$, $J$, $\JX$ and $\JV$ as in the statement of Claim~\ref{claim:leftover for g}, and recall that in~\eqref{eq:def initial weight boundedness} we defined a tuple $(\omega_{\iota},J,\JX, \JV,\gG,\pP)$ for each $\pP\in([k\alpha^{-1}]_0^r)^4$.
Let $(\omega,\omega_{\iota},J,\JX, \JV,\gG,\pP)$ be the edge tester in $\cW_{edge}^r(\phi_r,\phi_r^+,\cA^r,\cB^r)$ with respect to $(\omega_{\iota},J,\JX, \JV,\gG,\pP)$, $(\phi_r,\phi_r^+)$, $\cA^r$, and $\cB^r$ that we obtain from~\ind{r}\ref{edge weight}. 
In order to be able to distinguish these edge testers according to the patterns $\pP\in([k\alpha^{-1}]_0^r)^4$, we write $\omega_{\iota,\pP}:=\omega_\iota$ and $\omega_{\pP}:=\omega$.
Note that for $\pP=(\bpA,\bppA,\bpB,\bppB)$ such that $\omega_{\iota,\pP}(\sX_{\sqcup {\rR}})>0$, we have that $\supp(\omega_{\iota,\pP})\subseteq E_\cH(\pP,\rR,J)$ by Definition~\ref{def:general edge tester} and thus, by~\eqref{eq:norms equal} and because $I=\rR\in E(R)$ and $J\subseteq\rR$ is non-empty, it holds that {$\norm{\bpA}-\norm{\bppA}=0$ and $\norm{\bpB}-\norm{\bppB}\in\{-1,0\}$.}
By Definition~\ref{def:general edge tester} of an edge tester and by summing over all patterns $\pP\in([k\alpha^{-1}]_0^r)^4$, we can utilize our general edge testers to count the edges in $E_{\phi_r}(\gG,J,\JX, \JV)$.
Note that the properties~\ref{item:def1 leftover}--\ref{item:def4 leftover} of the definition of $E_{\phi_r}(\gG,J,\JX, \JV)$ in~\eqref{eq:def leftover g-edge} correspond to the properties~\ref{item:contained in ab}--\ref{item:x-vertices left unembedded} of Definition~\ref{def:general edge tester} for a general edge tester, respectively.
Due to~\ref{item:B not on centers} of Definition~\ref{def:general edge tester}, the edge tester $(\omega_\pP,\omega_{\iota,\pP},J,\JX,\JV,\pP)$ additionally requires for an element $\eH\in E_\cH(\pP,I,J)$ that $\phi_r^+(\eH\cap \sX_{\cup J})\cap \gG=\emptyset$.
That is, elements $\eH\in E_{\phi_r}(\gG,J,\JX, \JV)$ where some vertices $\eH\cap \sX_{\cup J}$ are already embedded by $\phi_r^+$ onto $\gG$ during the partial packing procedure, are not counted by any edge tester. 
However, our intuition is that the number of such edges only yields a minor order contribution, and in fact, we can employ~\ind{r}\ref{treffer auf g} to also account for those edges.
We make the following observation
\begin{align}
\nonumber&\left|E_{\phi_r}(\gG,J,\JX, \JV)\right|
\\\label{eq:sum claim3}&~\leq~
\sum_{\pP}\omega_{\pP}\left(\textstyle\bigsqcup_{j\in J}E(\cB_j^r)\right)
+\sum_{j\in J}\left|
\left\{
\eH\in E(\cH[\sX_{\cup {\rR}}])\colon
\phi_r^+(\eH\cap \sX_{\cup (\{j\}\cup(\rR\sm J) )})\subseteq\gG
\right\}
\right|.
\end{align}
We obtain an upper bound on the first term in~\eqref{eq:sum claim3} by employing~\ind{r}\ref{edge weight} as follows
\begin{align}
\nonumber\sum_{\pP}\omega_{\pP}\left(\textstyle\bigsqcup_{j\in J}E(\cB_j^r)\right)
&\stackrel{{\text{\ind{r}\ref{edge weight}}} 
}{\leq}
\sum_{\pP}
\left(
\left(\IND\{\JXV=\emptyset \}+\eps_T\right)
 d_B^{-1}
 \prod_{j\in  J} d_j
 \cdot
\frac{\omega_{\iota,\pP}(\sX_{\sqcup {\rR}})}{n^{k-|J|}}
+n^{\eps_T}
\right)
\\&\nonumber
\stackrel{\hphantom{{\text{\ind{r}\ref{edge weight}}}}
}{\leq}
\left(\IND\{\JXV=\emptyset \}+\eps_T\right) d_B^{-1}
\prod_{j\in  J} d_j
\cdot
\frac{\sum_{\pP}|E_{\cH}(\pP,\rR,J)|}{n^{k-|J|}}
+n^{2\eps_T}
\\&\label{eq:bound1 sum claim3}
\stackrel{\hphantom{{\text{\ind{r}\ref{edge weight}}}}
}{\leq} 
\left(\IND\{\JXV=\emptyset \}+\eps_T\right)
\gamma^{-1}
n^{|J|}\prod_{j\in  J} d_j,
\end{align}
where we used that $\sum_{\pP}|E_{\cH}(\pP,\rR,J)|\leq d_An^k$ and $d_B^{-1}d_A=(\gamma d)^{-1}(1-\gamma)d\leq \gamma^{-1}$.

An upper bound on the second term in~\eqref{eq:sum claim3} can be obtained by employing~\ind{r}\ref{treffer auf g}
\begin{align}
\label{eq:bound2 sum claim3}
\sum_{j\in J}\left|
\left\{
\eH\in E(\cH[\sX_{\cup {\rR}}])\colon
\phi_r^+(\eH\cap \sX_{\cup (\{j\}\cup(\rR\sm J) )})\subseteq\gG
\right\}
\right|
\leq 
\sum_{j\in J}
n^{k-|\{j\}\cup (\rR\sm J)|+2\eps_T}
\leq n^{|J|-1+3\eps_T}.
\end{align}
Substituting~\eqref{eq:bound1 sum claim3} and~\eqref{eq:bound2 sum claim3} in~\eqref{eq:sum claim3} establishes Claim~\ref{claim:leftover for g}.\endclaimproof
\COMMENT{
It is worth pointing out (at least for ourselves) that Claim~\ref{claim:leftover for g} is one of the crucial goals that we need to allow for vanishing density.
It is rather easy to obtain an upper bound `$\leq \eps_T^{1/2}n^{|J|}$'. 
Hence, a major reason for the complex definition of edge testers is to obtain this more general form for controlling the leftover in Claim~\ref{claim:leftover for g}.
}

\begin{substep}\label{step:matching argument}
Embedding one $H\in\cH$ by a random argument -- Claim~\ref{claim:matching}\nopagebreak
\end{substep}

We proceed with our argument for embedding one subgraph of a fixed graph $H\in\cH$.
Suppose we are given $Y_i^H\subseteq X_i^H$ for all $i\in[r]$. 
For all $\rR\in E(R)$ and $i\in[r]$, let 
\begin{align*}
E_\rR^{bad}&:=\big\{
\eH\in E(H[X_{\cup\rR}^H])\colon |\eH\cap Y_{\cup\rR}^H|\geq 2
\big\};
&\hspace{-.3cm}\cE^{bad}:=\textstyle\bigcup_{\rR\in E(R)}E_\rR^{bad};~
\quad
Y_i^{bad}&:=\{\eH\cap Y_i^H\colon \eH\in \cE^{bad} \};
\\
\cE^{good}&:=\{\eH\in E(H)\colon |\eH\cap Y_{\cup[r]}^H|=1 \};
&\hspace{-.35cm}~~~Y_i^{good}:= Y_i^H\sm Y_i^{bad}. \quad\quad~~
\hphantom{Y_i^{bad}}&\hphantom{:=\{\eH\cap Y_i^H\colon \eH\in \cE^{bad} \};}
\end{align*}

\begin{claim}\label{claim:matching}
Suppose $G^\circ\subseteq G_B$, $H\in\cH$, and $Y_i^H\subseteq X_i^H$ for all $i\in[r]$ such that the following hold for all $i\in[r]$, where  $W_i^H:=V_i\sm\phi_r(\hX_i^H\sm Y_i^H)$:
\begin{enumerate}[label={\rm(\Alph*)\textsubscript{C\ref{claim:matching}}}]
	\item\label{item:size X} 
$X_i^H\sm\hX_i^H\subseteq Y_i^H\subseteq X_i^H$ and $|Y_i^H|=|W_i^H|= (1\pm \eps_T^{1/2}) \mu n$;

	\item\label{item:few bad activated}
 $|E_\rR^{bad}|\leq \mu^{3/2}n$  for all $\rR\in E(R)$;

	\item\label{item:few bad privates} 
$|N_{B_i^H}(v)\cap Y_i^{bad}|\leq \mu^{3/2}d_in$ for all $i\in[r]$, $v\in W_i^H$;

	\item \label{item:B sr}
$(B_i^H)^{G^\circ}$ is $(\mu^{1/5},d_i)$-super-regular 
for all $i\in[r]$;

	\item \label{item:B subgraph sr}
 $(B_i^H)^{G^\circ}[Y_i^H,W_i^H]$ is $(\mu^{1/6},d_i)$-super-regular
   for all $i\in[r]$;

	\item \label{item:G_B-Gcirc typical} $G_B-G^\circ$ is $(\mu^{1/2},t,d_B)$-typical with respect to $R$;

	\item \label{item:mu well-intersecting}
for all $\cS\subseteq
\bigcup_{\rR\in E(R)\colon i\in\rR} V_{\sqcup\rR\sm\{i\}}$ with $|\cS|\leq t$, we have
$$\big|W_i^H\cap N_{G_B-G^\circ}(\cS)\big| = (1\pm\eps_T^{1/2}) \mu \big|V_i\cap N_{G_B-G^\circ}(\cS)\big|.$$
\end{enumerate}
Then there exists a probability distribution of the embeddings $\phi^{\tH}$ of $\tH:=H[Y^H_{\cup[r]}]$ into $\tG:= G_B[W^H_{\cup[r]}]-G^\circ$ with $\phi^H:=\phi^{\tH} \cup \phi_r|_{V(H)\sm V(\tH)}$ such that
\begin{enumerate}[label={\rm(\Alph*)\textsubscript{C\ref{claim:matching}}}]
\setcounter{enumi}{7}
\item\label{item:phi' packing} $\phi^H$ is an embedding of $H$ into $G$
where $\phi^H(X_i^H)= V_i$ for all $i\in[r]$;
\setcounter{enumi}{9}
\item\label{item:H-edges in G_B} $\phi^H$ embeds all $H$-edges that contain a vertex in $Y_{\cup[r]}^H$ onto an edge in $G_B-G^\circ$;
\item\label{item:random packing}
$\displaystyle\prob{\phi^H(\eH\cap X_{\cup I}^H)=\{v_i\}_{i\in I}}\leq {\prod_{i\in I\colon \eH\cap Y_i^H\neq\emptyset}2(\mu d_in)^{-1}}$
\\for all $m\in[k]$, index sets $I\in\binom{[r]}{m}$, $\{v_i\}_{i\in I}\in V_{\sqcup I}$, and $H$-edges $\eH$ that contain a vertex in $Y_{\cup I}^H$.
\end{enumerate}
\end{claim}
\claimproof
We split the proof of Claim~\ref{claim:matching} into two parts, Step~\ref{step:claim step1} and~\ref{step:claim step2}.
In Step~\ref{step:claim step1}, we greedily embed all the $H$-edges of $\cE^{bad}$ into $G_B-G^\circ$ by considering the clusters in turn.
Afterwards, in Step~\ref{step:claim step2}, we are left with the $H$-edges in $\cE^{good}$; that is, where only a single vertex is not yet embedded. 
The assumptions~\ref{item:few bad activated} and~\ref{item:few bad privates} will guarantee that we only used few edges in $G_B-G^\circ$ in Step~\ref{step:claim step1}, such that we merely have to modify our candidacy graphs. 
Hence, we can easily find a perfect matching in each of these candidacy graphs to embed the $H$-edges of $\cE^{good}$, which will complete the embedding of $H$.
This matching procedure can be performed independently for each cluster as the $H$-edges in $\cE^{good}$ only contain a single vertex that is not yet embedded.
Clearly, this approach will establish~\ref{item:phi' packing} and~\ref{item:H-edges in G_B}.

In both steps, we will embed the vertices of $H$ by a random procedure in order to establish~\ref{item:random packing}. 
We will do so by making several random choices sequentially which naturally yields a probability distribution.
Since some of these choices lead to instances that do not yield a valid or good embedding, we will discard some of these instances and say the random procedure terminates with failure in these cases.
We will show that the proportion of choices with failure instances is exponentially small, that is, the probability that the random procedure terminates with failure is exponentially small in $n$.
This allows us to discard some of these choices / instances and to restrict our probability space to the remaining `nice' outcomes. 
Since the failure probability is exponentially small in $n$, this does not have a significant effect on the probability in~\ref{item:random packing}.

\begin{substep}\label{step:claim step1}
Proof of Claim~\ref{claim:matching} -- Embedding the $H$-edges in $\cE^{bad}$\nopagebreak
\end{substep}

Suppose $\phi_q^{bad}\colon Y_{\cup[q]}^{bad}\to W^H_{\cup[q]}$ is an injective function. 
Similarly as we defined candidacy graphs in Definition~\ref{def:candidacy graph}, we will also define updated candidacy graphs of $B_i^H$ with respect to $\phi_q^{bad}$ for $i\in[r]\sm[q]$.
To that end, for all $i\in[r]\sm[q]$ and $B\subseteq B_i^H$, let $B(\phi_q^{bad})$ be the spanning subgraph of~$B$, where we keep the edge $xv$ of $B$ in $B(\phi_q^{bad})$ if all $\eH\in\cE^{bad}$ with $\eH\cap Y_{\cup([r]\sm[q])}^{bad}=\{x\}$\COMMENT{This also implies that $\eH\cap Y_{\cup[q]}^{bad}\neq\emptyset$.} satisfy that 
\begin{align}\label{eq:def B-candidacy}
&\phi_r|_{V(H)\sm V(\tH)}(\eH\sm V(\tH)) \cup \phi_q^{bad}(\eH\cap Y_{\cup[q]}^{bad})\cup\{v\} \in E(G_B-G^\circ).
\end{align} 
Observe that~\eqref{eq:def B-candidacy} is very similar to~\eqref{eq:condition updating} in Definition~\ref{def:candidacy graph}.
For all $q\in[r]$, $i\in[r]\sm[q]$, $y\in Y_i^{bad}$,~let
\begin{align*}
b_q(y):=\left|
\left\{
\eH\in \cE^{bad}\colon \eH\cap Y_{\cup([r]\sm[q])}^{bad}=\{y\}
\right\}
\right|.
\end{align*} 
That is, $b_q(y)$ is the number of edges $\eH$ in $\cE^{bad}$ containing $y$ whose $(k-1)$-set $\eH\sm\{y\}$ has already been embedded by $\phi_r|_{V(H)\sm V(\tH)}\cup\phi_q^{bad}$.

We make one more observation and fix $y\in Y_{q+1}^{bad}$.
Since $B_{q+1}^H$ is $(\eps_T,d_B^{\deg_{R}(q+1)})$-super-regular and $(\eps_T,t^{2/3})$-well-intersecting by~\ind{r}\ref{Z candidacy superregular}, 
we have that for all $\rR\in E(R)$ and $\eH\in E_\rR^{bad}$ with $\eH\cap Y_{\cup([r]\sm[q])}^{bad}=\{y\}$, there exists
a $(k-1)$-set $S_\rR=\phi_r^+(\eH\sm\{y\})\in V_{\sqcup(\rR\sm\{q+1\})}$, as well as there exists a set $\cS_y$ of $(k-1)$-sets with $S_\rR\in \cS_y$ and $|\cS_y|=\dg_R(q+1)$ such that
$N_{B_{q+1}^H}(y)
= 
V_{q+1}\cap N_{G_B}(\cS_y)$
and $|N_{B_{q+1}^H}(y)|=
(1\pm\eps_T)d_B^{|\cS_y|}n$.
(See~\eqref{eq:eps,t well-intersecting} for the definition of well-intersecting.)
Note that~$\phi_r^+$ embeds $\eH\sm\{y\}$ onto $S_\rR$ but $\phi_r|_{V(H)\sm V(\tH)}$ does not. 
Hence, $S_\rR$ only serves as a `dummy' $(k-1)$-set for updating the candidacy graph $B_{q+1}^H$ conveniently and to artificially restrict the candidates for~$y$. 
That is, $V_{q+1}\cap N_{G_B}(\cS_y\sm\{S_\rR\})$ are also suitable candidates where we could embed~$y$, assuming that~{$\eH\sm\{y\}$} has not been embedded yet.
Let $\cS_y^{dummy}$ be the set of all these $(k-1)$-sets $S_\rR$ for~$y$, and note that $|\cS_y^{dummy}|=b_q(y)$.

During our process of embedding the $H$-edges in~$\cE^{bad}$, we will have to drop this artificial restriction of the candidate sets. 
To that end, suppose we are given a spanning subgraph $B\subseteq B_{q+1}^H$,  $y\in Y_{q+1}^{bad}$, and there exists a set $\cS'_y$ of $(k-1)$-sets such that we can write
\begin{align*}
N_B(y)=V_{q+1}\cap N_{G_B-G^\circ}(\cS_y').
\end{align*}
Then let $B^{\ominus dummy}$ be the spanning supergraph of $B$ where the neighbourhood of each vertex $y\in Y_{q+1}^{bad}$ is given by
\begin{align}\label{eq:add dummy}
N_{B^{\ominus dummy}}(y)=V_{q+1}\cap N_{G_B-G^\circ}(\cS_y'\sm\cS_y^{dummy}).
\end{align}

\medskip
We inductively prove that the following statement~\indC{q} holds for all $q\in[r]_0$, which will extend $\phi_r|_{V(H)\sm V(\tH)}$ by embedding the edges in $\cE^{bad}$ into $G_B-G^\circ$.
To that end, we define a set of good pairs of vertices.
\begin{align}\label{eq:Y-good pairs}
\begin{minipage}[c]{0.9\textwidth}\em
For every $i\in[r]$, let $Y_i^{good~pairs}\subseteq\binom{Y_i^{good}}{2}$ be the set containing all pairs $\{y,y'\}$ with $\cS_y\cap \cS_{y'}=\emptyset$ where $\cS_x\subseteq 
\bigcup_{\rR\in E(R)\colon i\in\rR} V_{\sqcup\rR\sm\{i\}}$ is such that $N_{B_i^H}(x)=V_i\cap  N_{G_B}(\cS_x)$ for each $x\in\{y,y'\}$.
\end{minipage}\ignorespacesafterend
\end{align}
\COMMENT{That is, $Y_i^{good~pairs}$ is the set of pairs in $Y_i^{good}$ with appropriate codegree.}
We note for future reference that 
\begin{align}\label{eq:number of bad Y-pairs}
\bigg|
\binom{Y_i^{good}}{2}\sm Y_i^{good~pairs}
\bigg|
\leq 2n\cdot n^{1/4+\eps_T}\leq n^{4/3},
\end{align}
since $B_i^H$ is $(\eps_T,t^{2/3})$-well-intersecting as defined in~\eqref{eq:eps,t well-intersecting}.
\medskip

\begin{itemize}
\item[\indC{q}.] 
There exists a probability distribution of the injective functions $\phi_q^{bad}\colon Y^{bad}_{\cup[q]}\to W^H_{\cup[q]}$ with $\phi'_q:=\phi_q^{bad}\cup \phi_r|_{V(H)\sm V(\tH)}$ such that
\begin{enumerate}[label={\rm(\Roman*)\textsubscript{C\ref{claim:matching}}}]
\item\label{C1I} $\phi'_q$ is an embedding of $H'_q:=H[Y_{\cup[q]}^{bad}\cup (V(H)\sm V(\tH))]$ into $G$;
\item\label{C1II} all edges in $H[Y_{\cup[q]}^{bad}\cup (V(H)\sm V(\tH))]$ that contain a vertex in $Y_{\cup[q]}^{bad}$ are embedded on an edge in $G_B-G^\circ$;
\item\label{C1II.5}
for every vertex $x\in Y_i^{good}$ and $i\in[q]$, there are at most $\mu^{4/3}d_in$ vertices $y\in Y_i^{bad}$ with $\phi_q^{bad}(y)\in N_{B_i^H}(x)$;

\item\label{C1III} 
$\displaystyle\prob{
\phi_q'(\eH\cap X_{\cup I}^H)=\{v_i\}_{i\in I}
}\leq 
{\prod_{i\in I\colon\eH\cap Y_i^H\neq\emptyset}2(\mu d_in)^{-1}}$
\\for all $\rR\in E(R)$,  $m\in[k]$, $I\in\binom{\rR\cap[q]}{m}$, $\{v_i\}_{i\in I}\in V_{\sqcup I}$ and $\eH\in E_\rR^{bad}$.
\end{enumerate}
\end{itemize}

\medskip
The statement~\indC{0} clearly holds for $\phi_0^{bad}$ being the empty function.
Hence, we assume the truth of~\indC{q} for some $q\in[r-1]_0$. 
Our strategy to establish~\indC{q+1} is as follows.
We extend the probability space given in~\indC{q} by making further random choices.
For $\phi_q^{bad}$ as in~\indC{q},
we aim to find a matching $\sigma_{q+1}^{bad}\colon Y_{q+1}^{bad}\to W^H_{q+1}$ in
a suitable candidacy graph between $Y_{q+1}^{bad}$ and $W^H_{q+1}$
that extends $\phi_q^{bad}$ to $\phi_{q+1}^{bad}:=\phi_q^{bad}\cup\sigma_{q+1}^{bad}$.
If we can find such a matching $\sigma_{q+1}^{bad}$, then~\indC{q+1}\ref{C1I} and~\ref{C1II} will hold by the definition of this suitable candidacy graph.
In particular, we will find $\sigma_{q+1}^{bad}$ by a random procedure to also ensure~\ref{C1III}. 
We will discard an exponentially small proportion of random choices during this procedure in order to satisfy~\ref{C1II.5} and to obtain a suitable embedding~$\sigma_{q+1}^{bad}$.

Let us describe this suitable candidacy graph. 
We will choose $\sigma_{q+1}^{bad}$ randomly in 
\begin{align*}
\tB:=(B_{q+1}^H)^{G^\circ}(\phi_q^{bad})^{\ominus dummy}[Y_{q+1}^{bad},W^H_{q+1}].
\end{align*}
That is, $\tB$ arises from $B_{q+1}^H$ as follows.
\begin{itemize}
\item First, we restrict the candidate sets in $B_{q+1}^H$ to those edges whose corresponding edges in~$G_B$ have not been used in~$G^\circ$ for packing graphs $H_1,\ldots, H_h$ in previous rounds.
\item Second, we restrict the candidate sets with respect to the packing $\phi_q^{bad}$ of the vertices in~$Y_{\cup[q]}^{bad}$ according to~\eqref{eq:def B-candidacy}.
\item Third, we drop the restriction of the candidate sets to the dummy $(k-1)$-sets as in~\eqref{eq:add dummy}.
\item In the end, we consider the induced subgraph of this candidacy graph on the bad vertices~$Y_{q+1}^{bad}$ and the vertices $W_{q+1}^H$ that can be used for the completion.
\end{itemize}
For the sake of a better readability, let $B:=(B_{q+1}^H)^{G^\circ}$.

To guarantee the existence of $\sigma_{q+1}^{bad}$, we will show that the degree of every vertex $y\in Y_{q+1}^{bad}$ is sufficiently large in $B$ and also in $\tB$.
Let $y\in Y_{q+1}^{bad}$ be fixed.
Since $B_{q+1}^H$ is $(\eps_T,t^{2/3})$-well-intersecting by~\ind{r}\ref{Z candidacy superregular}, there exists a set $\cS_y\subseteq
\bigcup_{\rR\in E(R)\colon q+1\in\rR} V_{\sqcup\rR\sm\{q+1\}}$ of $(k-1)$-sets with $|\cS_y|\leq t^{2/3}$ such that $N_{B_{q+1}^H}(y)=V_{q+1}\cap N_{G_B}(\cS_y)$. 
Hence, by the definition of $B=(B_{q+1}^H)^{G^\circ}$ in~\eqref{eq:def B^G}, we have
\begin{align}\label{eq:B^circ-neighbourhood y}
N_{B}(y)=V_{q+1}\cap N_{G_B-G^\circ}(\cS_y),
\quad \text{ and }\quad
\dg_{B}(y)
=(1\pm 2\mu^{1/5})d_{q+1}n,
\end{align}
because $B$ is $(\mu^{1/5},d_{q+1})$-super-regular by~\ref{item:B sr}.
Of course we have to restrict the potential images of $y$ according to the vertices we already embedded by $\phi_q^{bad}$.
To this end, let 
\begin{align*}
\cS_y^{bad}:=\left\{ \phi_q^{bad}(\eH\cap Y_{\cup[q]}^{bad})\cup \phi_r|_{V(H)\sm V(\tH)}(\eH\sm V(\tH)) \colon \eH\in\cE^{bad}, \eH\cap Y_{\cup([r]\sm[q])}^{bad}=\{y\} \right\} 
\end{align*}
and note that $|\cS_y^{bad}|=b_q(y)$.
By the definition of the candidacy graph $B(\phi_q^{bad})$ in~\eqref{eq:def B-candidacy} and by~\eqref{eq:B^circ-neighbourhood y}, we obtain 
\begin{align*}
N_{B(\phi_q^{bad})}(y)=V_{q+1}\cap 
N_{G_B-G^\circ}(\cS_y\cup\cS_y^{bad}),
\end{align*}
and thus,
together with~\ref{item:G_B-Gcirc typical} and~\eqref{eq:B^circ-neighbourhood y}, we have that
\begin{align*}
\dg_{B(\phi_q^{bad})}(y)
=(1\pm \mu^{1/6})d_{q+1}d_B^{b_q(y)}n.
\end{align*}
Since $|\cS_y^{dummy}|=b_q(y)$, this implies
\begin{align}\label{eq:B-degree add dummy y}
\dg_{B(\phi_q^{bad})^{\ominus dummy}}(y)= \Big|
V_{q+1}\cap 
N_{G_B-G^\circ}\big((\cS_y\sm\cS_y^{dummy})\cup\cS_y^{bad}\big)
\Big|
=(1\pm \mu^{1/6})d_{q+1}n.
\end{align}
\COMMENT{Note that $(\cS_y\sm\cS_y^{dummy})\cap\cS_y^{bad}=\emptyset$.}Now, we obtain
\begin{align}\label{eq:B'-neighbourhood y}
N_{\tB}(y)= 
W^H_{q+1}\cap N_{G_B-G^\circ}\big((\cS_y\sm\cS_y^{dummy})\cup\cS_y^{bad}\big)
\end{align}
and thus, by~\ref{item:mu well-intersecting} and~\eqref{eq:B-degree add dummy y}, we have that
\begin{align}\label{eq:B'-degree y}
\dg_{\tB}(y)
=(1\pm \mu^{1/7})d_{q+1}\mu n.
\end{align}
In order to guarantee~\ref{C1III}, we find $\sigma_{q+1}^{bad}$ via the following random procedure:
\begin{itemize}
\item for every vertex $y\in Y_{q+1}^{bad}$ in turn, we choose a neighbour in $N_{\tB}(y)$ uniformly at random among all neighbours that have not been chosen in previous turns;
\item we terminate the random procedure with failure at some step of the procedure, say at the turn of some vertex $y\in Y_{q+1}^{bad}$, 
\begin{itemize}
\item  if we have less than $2(\mu d_{q+1} n)/3$ choices to select an image for $y$ in $N_{\tB}(y)$, or
\item if there is a vertex $x\in Y_{q+1}^{good}$ such that~\ref{C1II.5} is violated, that is, we mapped in previous turns already more than $\mu^{4/3}d_{q+1}n$ vertices of $Y_{q+1}^{bad}$ into $N_{B_{q+1}^H}(x)$.
\end{itemize}
\end{itemize}
We show in the following claim that this random procedure terminates with  failure only with exponentially small probability.
If the procedure does not terminate with failure, we obtain a random $Y_{q+1}^{bad}$-saturating matching $\sigma_{q+1}^{bad}\colon Y_{q+1}^{bad}\to W^H_{q+1}$ in $\tB$, which by definition of the candidacy graph $B_{q+1}^H$ and $\tB$ implies~\indC{q+1}\ref{C1I}--\ref{C1II.5}.
Further, $\sigma_{q+1}^{bad}$ satisfies the following.
\begin{align}\label{eq:bad random perfect matching}
\begin{minipage}[c]{0.9\textwidth}\em
For all $y\in Y_{q+1}^{bad}$, $w\in W^H_{q+1}$, we have that $\sigma_{q+1}^{bad}(y)=w$ with probability at most $2(\mu d_{q+1}n)^{-1}$.
\end{minipage}\ignorespacesafterend
\end{align}
\COMMENT{The upper bound on the probability holds for any choice of $\sigma_i^{bad}$ for $i\in[q]$.}
Hence, the following claim together with~\indC{q}\ref{C1III} establishes~\indC{q+1}\ref{C1III}.
\begin{claim}\label{claim choices for bad edges}
The random procedure for computing $\sigma^{bad}_{q+1}$ terminates with failure with probability at most $\eul^{-n^{1/2}}$.\looseness=-1
\end{claim}

\claimproof
Let $Y_{q+1}^{bad}=\{y_1,\ldots, y_m\}$ and we consider every vertex in turn, that is, $y_{\ell+1}$ will be treated after $y_\ell$.
For all $x\in Y_{q+1}^H$ and $\ell\in[m]$, let $\xi_\ell(x)$ be the random variable that counts the number of covered neighbours so far, that is, the number of vertices $v\in N_{B_{q+1}^H}(x)$ such that $\sigma_{q+1}^{bad}(y)=v$ for some $y\in\{y_i\}_{i\in[\ell]}$.
We say the random procedure fails at step $\ell\in[m]$, if $\ell$ is the smallest integer such that 
\begin{itemize}
\item $\xi_{\ell}(z)>\mu^{2/5}d_{q+1}\mu n$ for some vertex $z\in Y_{q+1}^{good}\cup \{y_i\}_{i\in [m]\sm[\ell]}$.
\end{itemize}
We show that the random procedure fails at some step $\ell\in[m]$ with probability at most $\eul^{-n^{2/3}}$. 
A union bound then establishes Claim~\ref{claim choices for bad edges}.

We fix $\ell\in[m]$ and a vertex $z\in Y_{q+1}^{good}\cup \{y_i\}_{i\in [m]\sm[\ell]}$.
By employing~\eqref{eq:B'-neighbourhood y}, \eqref{eq:B'-degree y} and that $B_{q+1}^H$ is $(\eps_T,t^{2/3})$-well-intersecting, we obtain
\COMMENT{
Let $$\cS_{y,z}:=
\left((\cS_{y}\sm\cS_{y}^{dummy})\cup\cS_{y}^{bad}\right)
\cup 
\left((\cS_{z}\sm\cS_{z}^{dummy})\cup\cS_{z}^{bad}\right),$$
and note that for all but at most $n^{1/4+\eps_T}\leq n^{1/3}$ vertices $y\in\{y_i\}_{i\in[\ell]}$, we have $|\cS_{y,z}|=2\deg_{R}(q+1)$ because it holds that
$$\left((\cS_{y}\sm\cS_{y}^{dummy})\cup\cS_{y}^{bad}\right)
\cap 
\left((\cS_{z}\sm\cS_{z}^{dummy})\cup\cS_{z}^{bad}\right)=\emptyset
$$
for all but at most $n^{1/4+\eps_T}$ vertices $y\in\{y_i\}_{i\in[\ell]}$ since $B_{q+1}^H$ is $(\eps_T,t^{2/3})$-well-intersecting as defined in~\eqref{eq:eps,t well-intersecting}.
Hence, 
\begin{align*}
\left|N_{\tB}(y)\cap N_{\tB}(z)\right|
&
=
\left|
W^H_{q+1}\cap N_{G_B-G^\circ}(\cS_{y,z})
\right|
=
(1\pm \mu^{1/7})d_{q+1}^2\mu n.
\end{align*}
for all but at most $n^{1/3}$ vertices $y\in\{y_i\}_{i\in[\ell]}$.
}
\begin{align}\label{eq:common choices}
|N_{\tB}(y\land z)|= (1\pm \mu^{1/7})d_{q+1}^2\mu n \text{ for all but at most $n^{1/3}$ vertices $y\in\{y_i\}_{i\in[\ell]}$.} 
\end{align} 
Further, at the turn of each $y_i$, $i\in[\ell]$, we have at least 
\begin{align}\label{eq:all choices}
|N_{\tB}(y_i)|-\mu^{2/5}d_{q+1}\mu n 
\stackrel{\eqref{eq:B'-degree y}}{\geq}
(1-2\mu^{1/7})d_{q+1}\mu n
\end{align}
choices for the embedding of $y_i$.\COMMENT{For the embedding of $y_i$, there might also be some neighbours in $N_{\tB}(y_i)$ that are already utilized by the previous choices for $y_1,\ldots, y_{i-1}$. 
However, if we haven't terminated with failure yet, these are at most $\mu^{2/5}d_{q+1}\mu n$ neighbours.} 
Hence, by~\eqref{eq:common choices} and~\eqref{eq:all choices}, the probability that a vertex $y_i$, $i\in[\ell]$ which satisfies~\eqref{eq:common choices} is mapped into $N_{\tB}(z)$ is at most
\begin{align*}
\frac{(1+ \mu^{1/7})d_{q+1}^2\mu n}{(1- 2\mu^{1/7})d_{q+1}\mu n}
\leq 2 d_{q+1}.
\end{align*}
This implies that
\begin{align*}
\expn{\xi_{\ell}(z)}
\leq 2d_{q+1}\ell + n^{1/3}
\leq 2d_{q+1}|Y_{q+1}^{bad}| + n^{1/3}
\leq 3d_{q+1} \alpha^{-1}\mu^{3/2}n,
\end{align*}
where we used that $|Y_{q+1}^{bad}|\leq\alpha^{-1}\mu^{3/2}n$ by~\ref{item:few bad activated}. 
Hence, Theorem~\ref{thm:McDiarmid} implies that $\xi_{\ell}(z)>\mu^{2/5}d_{q+1}\mu n$ with probability, say, at most $\eul^{-n^{2/3}}$ for some $z\in Y_{q+1}^{good}\cup \{y_i\}_{i\in [m]\sm[\ell]}$.
Thus, the random procedure fails at step $\ell$ with probability at most~$\eul^{-n^{2/3}}$.
A simple union bound completes the proof of Claim~\ref{claim choices for bad edges}.
\endclaimproof

\begin{substep}\label{step:claim step2}
Proof of Claim~\ref{claim:matching} -- Embedding the $H$-edges in $\cE^{good}$\nopagebreak
\end{substep}

Let $\phi_r^{bad}\colon Y^{bad}_{\cup[r]}\to W^H_{\cup[r]}$ and $\phi_r'=\phi_r^{bad}\cup\phi_r|_{V(H)\sm V(\tH)}$ be as in~\indC{r} obtained in Step~\ref{step:claim step1}.
Recall that for all $i\in[r]$, we have $Y_i^{good}=Y_i^H\sm Y_i^{bad}$ and let $W_i^{good}:=W_i^H\sm\phi_r^{bad}(Y_i^{bad})$, and thus clearly, $|Y_i^{good}|=|W_i^{good}|$.
For every $i\in[r]$, we aim to embed the vertices $Y_i^{good}$ onto $W_i^{good}$ by finding a perfect matching in $(B_i^H)^{G^\circ}$.

Let $i\in[r]$ be fixed.
By~\ref{item:B subgraph sr}, we have that $\hB_i:=(B_i^H)^{G^\circ}[Y_i^H,W_i^H]$ is $(\mu^{1/6},d_i)$-super-regular. We show that for every vertex in $Y_i^{good}\cup W_i^{good}$, we only removed few incident edges when take the subgraph $B_i^{good}:=\hB_i[Y_i^{good},W_i^{good}]$ of~$\hB_i$.
For a vertex $v\in W_i^{good}$, we removed at most $\mu^{3/2}d_in$ incident edges by~\ref{item:few bad privates}.
For a vertex $x\in Y_i^{good}$, we removed at most $\mu^{4/3}d_in$ incident edges by~\ref{C1II.5}.
Hence, by employing Fact~\ref{fact:regularity robust}, we obtain that $B_i^{good}$ is $(\mu^{1/19},d_i)$-super-regular for every $i\in[r]$.\looseness=-1 

Our strategy is to apply Lemma~\ref{lem:m-factor} that allows us to find a regular spanning subgraph of $B_i^{good}$ from which we can easily take a random perfect matching. 
In order to satisfy the assumptions of Lemma~\ref{lem:m-factor}, we show that the common neighbourhood of most of the pairs in $Y_i^{good}$ is also not too large in $B_i^{good}$.
To that end, we fix a pair of vertices $\{y,y'\}\in Y_i^{good~pairs}$  as defined in~\eqref{eq:Y-good pairs} and let $\cS_y, \cS_{y'}$ be the sets of $(k-1)$-sets such that $N_{B_i^H}(x)=V_i\cap N_{G_B}(\cS_x)$ for each $x\in\{y,y'\}$.
We have that\looseness=-1
\begin{align}\label{eq:common neighbourhood tB_i}
N_{\hB_i}(y\land y') 
= 
W_i^H \cap N_{(B_i^H)^{G^\circ}}(y\land y')
= W_i^H \cap
N_{G_B-G^\circ}(\cS_y\cup \cS_{y'}).
\end{align}
Hence, we obtain
\begin{align}\nonumber
|N_{B_i^{good}}(y\land y')|
&\stackrel{\hphantom{\eqref{eq:common neighbourhood tB_i},\ref{item:mu well-intersecting}}}{\leq}
 |N_{\hB_i}(y\land y')|
\\\nonumber&\stackrel{\eqref{eq:common neighbourhood tB_i},\ref{item:mu well-intersecting}}{\leq}
(1+\eps_T^{1/2})\mu 
\big|
V_i\cap  N_{G_B-G^\circ}(\cS_y\cup \cS_{y'})
\big|
\\
\label{eq:codegree in B_i^good}
&\stackrel{~~~~\ref{item:G_B-Gcirc typical}~~~}{\leq}
(1+ 2\mu^{1/2})\mu d_i^2n,
\end{align}
where we used for the last equality that $\cS_y\cap \cS_{y'}=\emptyset$ since $\{y,y'\}\in Y_i^{good~pairs}$, and thus, $d_B^{|\cS_y\cup \cS_{y'}|}=d_i^2$. 
Hence, \eqref{eq:number of bad Y-pairs} implies that all but at most $n^{4/3}$ pairs of vertices in $Y_i^{good}$ satisfy~\eqref{eq:codegree in B_i^good}.

Finally, we can apply Lemma~\ref{lem:m-factor} and obtain a spanning $\mu d_in/2$-regular subgraph of $B_i^{good}$. 
In particular, $B_i^{good}$ contains $\mu d_in/2$ edge-disjoint perfect matchings, from which we choose one perfect matching $\sigma_i^{good}\colon Y_i^{good}\to W_i^{good}$ for each $i\in[r]$ uniformly and independently at random.
We crucially observe:
\begin{align}\label{eq:good random perfect matching}
\begin{minipage}[c]{0.9\textwidth}\em
For all $i\in[r]$, $y_i\in Y_i^{good}$, $w_i\in W_i^{good}$, we have that $\sigma_i^{good}(y_i)=w_i$ with probability at most $2(\mu d_i n)^{-1}$.
\end{minipage}\ignorespacesafterend
\end{align} 
Further, since all the vertices in $H_\ast[Y_{\cup[r]}^{good}]$ are isolated, we have that $\phi^{good}:=\bigcup_{i\in[r]}\sigma_i^{good}$ is an injective function $\phi^{good}\colon Y_{\cup[r]}^{good}\to W_{\cup[r]}^{good}$ such that $\phi_r|_{V(H)\sm V(\tH)}\cup\phi^{good}$ is a random packing of $H[Y_{\cup[q]}^{good}\cup (V(H)\sm V(\tH))]$ into $G$ that embeds all edges in $H[Y_{\cup[r]}^{good}\cup (V(H)\sm V(\tH))]$ that contain a vertex in $Y_{\cup[r]}^{good}$ on an edge in $G_B-G^\circ$. 
In particular, since we find the perfect matchings in the candidacy graphs $(B_i^H)^{G^\circ}$ induced on $Y_i^{good}\cup W_i^{good}$ for each $i\in[r]$, we have for $\phi^H:=\phi_r|_{V(H)\sm V(\tH)}\cup\phi_r^{bad}\cup\phi^{good}$ that 
\begin{itemize}
\item $\phi^H$ is an embedding of $H$ into $G$
where $\phi^H(X_i^H)= V_i$ for all $i\in[r]$, which establishes~\ref{item:phi' packing};
\item all $H$-edges that contain a vertex in $Y_{\cup[r]}^H$ are embedded on an edge in $G_B-G^\circ$, which establishes~\ref{item:H-edges in G_B};
\item for all $m\in[k]$, index sets $I\in\binom{[r]}{m}$, $\{v_i\}_{i\in I}\in V_{\sqcup I}$, and $H$-edges $\eH$ that contain a vertex in $Y_{\cup I}^H$, we have that $\phi^H(\eH\cap X_{\cup I}^H)=\{v_i\}_{i\in I}$ with probability at most ${\prod_{i\in I\colon \eH\cap Y_i^H\neq\emptyset}2(\mu d_in)^{-1}}$  by~\indC{r}\ref{C1III} and \eqref{eq:good random perfect matching}, which establishes~\ref{item:random packing}. 
\end{itemize}
This completes the proof of Claim~\ref{claim:matching}.
\endclaimproof

\begin{substep}\label{step:random packing procedure}
The random packing procedure -- Claim~\ref{claim:random packing}\nopagebreak
\end{substep}

We now proceed to our Random Packing Procedure (RPP). 
Let $G_0^\circ$ be the edgeless graph on $V(G)$ and let $\phi^0$ be the empty function. 
We perform the following random procedure.\COMMENT{Previously we activated every vertex in $\hX_i^H$ with probability $\mu n/|\hX_i^H|$.}

\medskip
\fbox{
\begin{minipage}[c]{.9\textwidth}\em
\noindent{\it\underline{Random Packing Procedure (RPP)}}

{\medskip\noindent\it For $h=1,\ldots,|\cH|$ do:}
\begin{itemize}
\it
\item Set $H:=H_{h}$.
For all $i\in[r]$, independently activate every vertex in $\hX_i^H$ with probability $\mu$ and let $Y_i^H$ be the union of $X_i^H\sm\hX_i^H$ and all activated vertices in~$\hX_i^H$.
\item If the assumptions of Claim~\ref{claim:matching} are satisfied, apply Claim~\ref{claim:matching} and obtain a random packing $\phi^H$ that satisfies~\ref{item:phi' packing}--\ref{item:random packing}; otherwise terminate with failure.
\item Set $\phi^{h}:=\phi^{h-1}\cup \phi^H$ and $G_{h}^\circ:=G_{h-1}^\circ \cup (\phi^H(H)\cap G_B)$.
\end{itemize}
\end{minipage}\ignorespacesafterend
}

\medskip
\begin{claim}\label{claim:random packing}
With probability at least $1-1/n$, the RPP terminates without failure and satisfies conclusion~\ref{item: vertex tester1} of Lemma~\ref{lem:main matching}.\COMMENT{Conclusions~\ref{item:partition1} and~\ref{item:set testers1} are satisfied by a deterministic conclusion.}
\end{claim}
\claimproof
We prove Claim~\ref{claim:random packing} in Steps~\ref{step:random packing step1}--\ref{step:few activated weight on vertex testers}.
Our general strategy is to guarantee that we can apply Claim~\ref{claim:matching} in each turn of the procedure.
To that end, we will introduce a collection of random variables that we call \emph{identifiers}.
Such an identifier indicates an unlikely event and if this event happens, we say the identifier \emph{detects alarm} and we simply terminate the RPP with failure. 
This means, if an identifier detects alarm at some turn $h\in[|\cH|]$, we terminate the RPP and deactivate all further identifiers; that is, the probability that they detect alarm is set to $0$.
We show that the probability that an individual identifier detects alarm is exponentially small in $n$. 
In the end, a union bound over all identifiers will imply that with probability at least $1-1/n$, none of the identifiers will detect alarm and thus, the RPP terminates without failure.

For most of the identifiers, it follows by a standard application of  Theorem~\ref{thm:McDiarmid} that the probability to detect alarm is exponentially small (in fact, often Chernoff's inequality suffices).
To that end, in the subsequent Steps~\ref{step:random packing step1}--\ref{step:few activated weight on vertex testers}, we often describe only the random variables to which we apply Theorem~\ref{thm:McDiarmid}.

\medskip
First, let us observe that~\ind{r} implies that
\begin{align}\label{eq:small leftover}
\big|
X_i^{H}\sm\hX_i^{H}
\big|
=
\big|
V_i\sm \phi_r(\hX_i^{H})
\big|
\leq 2\eps_Tn \text{ for all $H\in\cH$, $i\in[r]$.}
\end{align}
Further, for $Y_i^H$ as in the RPP, let $W_i^H:= V_i\sm \phi_r(\hX_i^H\sm Y_i^H)$.
For convenience, we also call a vertex $w\in W_i^H\cap \phi_r(\hX_i^H)$ \emph{activated} by $H\in\cH$. (Recall that we only activate vertices in $\hX_i^H$.) 

\begin{substep}\label{step:random packing step1}
Proof of Claim~\ref{claim:random packing} -- Establishing \ref{item:size X}--\ref{item:mu well-intersecting}\nopagebreak
\end{substep}

In this step, in order to establish \ref{item:size X}--\ref{item:mu well-intersecting} at each turn of the RPP, we consider several random variables for which we individually introduce an identifier that detects alarm if the considered random variable is not within a factor of $(1\pm\eps)$ of its expectation.\COMMENT{In view of the statements, we may assume that the expectation is not too small; in particular, when we establish~\ref{item:few bad activated} and~\ref{item:few bad privates}.}
As mentioned above, for each identifier a standard application of Theorem~\ref{thm:McDiarmid} implies that the probability to detect alarm is exponentially small, say, $\eul^{-n^{1/2}}$.
Let us only describe the random variables that we consider for establishing \ref{item:size X}--\ref{item:mu well-intersecting}.

\medskip
To establish~\ref{item:size X}, for each $H\in\cH$, $i\in[r]$, we consider the sum of indicator variables which indicate whether a vertex is activated in $X_i^H$. 
Together with~\eqref{eq:small leftover}, this implies~\ref{item:size X}. 

\medskip
To establish~\ref{item:few bad activated}, for each $H\in\cH$ and $\rR\in E(R)$, we consider the random variable that counts how many $H$-edges $\eH$ lie in $X_{\cup\rR}^H$ where at least two vertices of $\eH$ are activated; in view of the statement, we may assume that $e_H(X_{\cup\rR}^H)\geq \mu^{3/2}n$. 
Note that the probability for an edge $\eH$ that at least two vertices are activated is at most $k^2\mu^2$.
Together with~\eqref{eq:small leftover},\COMMENT{There are at most $2 k\eps_Tn$ edges $\eH$ in $X_{\cup\rR}^H$ where at most one vertex is activated and at least one vertex is left unembedded.} this implies~\ref{item:few bad activated}.

\medskip
To establish~\ref{item:few bad privates}, for all $H\in\cH$, $i\in[r]$, and $v\in V_i$, we consider the random variable $\xi$ that counts how many $B_i^H$-neighbours of $v$ in $X_i^H$ are activated and lie in an $H$-edge $\eH$ where a second vertex in $\eH$ is either activated or left unembedded.
With Claim~\ref{claim:leftover obs H-neighbours of v} we have that 
$\expn{\xi}\leq 2k^2\mu^2\alpha^{-1}d_in+ 2\mu \eps_T^{1/2}d_in\leq \mu^{5/3}d_in$. 
\COMMENT{\new{We have $|N_{B_i^H}(v)|\leq 2d_in$ and each such neighbour is contained in at most $\alpha^{-1}$ edges $\eH$ in $H$.
The probability that two vertices of $\eH$ are activated is at most $k^2\mu^2$.
This accounts for the expected number $2k^2\mu^2\alpha^{-1}d_in$ of $\cH$-edges where at least two vertices are activated.
By Claim~\ref{claim:leftover obs H-neighbours of v}, there are at most $\eps_T^{1/2}d_in$ neighbours $x$ of $v$ in $B_i^H$ that lie in an edge where some vertex is left unembedded. 
The probability that $x$ is activated is $\mu$.
This accounts for the expected number $2\mu \eps_T^{1/2}d_in$ of $B_i^H$-neighbours $x$ of $v$ such that $x$ becomes activated and lies in an edge with an unembedded vertex.
In the end, by Claim~\ref{claim:leftover obs}, we might also have to remove $2\eps_Td_in$ $B_i^H$-neighbours $x$ of $v$ where $x$ is left unembedded (and lies in an edge where a second vertex is either unembedded or activated, but we do not care).}
}
Together with Claim~\ref{claim:leftover obs}, this implies~\ref{item:few bad privates}.

\medskip
To establish~\ref{item:mu well-intersecting}, for all $i\in[r]$, $\cS\subseteq
\bigcup_{\rR\in E(R)\colon i\in\rR} V_{\sqcup\rR\sm\{i\}}$, $|\cS|\leq t$, we consider the sum of indicator variables which each indicates whether a vertex in $\phi_r(\hX_i^H)\cap N_{G_B-G_h^\circ}(\cS)$ is activated.
By~\ind{r}\ref{G_B leftover for typicality}, we further have that 
\begin{align}\label{eq:leftover intersection k-1 sets}
\big|
(V_i\sm \phi_r(\hX_i^H)) \cap N_{G_B}(\cS)
\big|
\leq \eps_T \big|
V_i \cap N_{G_B}(\cS)
\big|.
\end{align}
Altogether, this implies~\ref{item:mu well-intersecting}.

\medskip
In order to establish~\ref{item:B sr}--\ref{item:G_B-Gcirc typical}, we claim that if no identifier detected alarm until turn $h\in[|\cH|]$, then for all $i\in[r]$
\begin{align}\label{eq:bound used k-1 sets}
\left|
 \textstyle\bigcup_{S\in\cS}
\left(
V_i \cap N_{G^\circ_h}(S) \cap N_{G_B}(\cS)
\right)
\right| 
&\leq 
\mu^{2/3}d_B^{|\cS|}n 
\text{ for all $\cS\subseteq
\textstyle\bigcup_{\rR\in E(R)\colon i\in\rR} V_{\sqcup\rR\sm\{i\}}$, $|\cS|\leq t$;}
\\\label{eq:bound 1 degree}
\rho^H_{G^\circ_h}(v)&\leq \mu^{2/3}d_in \text{ for all $H\in\cH$, $v\in V_i$,}
\end{align}
where $\rho^H_{G^\circ_h}(v)$ is defined as in~\eqref{eq:def 1 degree}.
We verify~\eqref{eq:bound used k-1 sets} and \eqref{eq:bound 1 degree} in the subsequent Steps~\ref{step:bound used k-1 sets} and \ref{step:codegree bounds2}, respectively, and first establish~\ref{item:B sr}--\ref{item:G_B-Gcirc typical} assuming~\eqref{eq:bound used k-1 sets} and \eqref{eq:bound 1 degree}.

\medskip
To establish~\ref{item:G_B-Gcirc typical}, recall that $G_B$ is $(2\eps_0,t,d_B)$-typical with respect to $R$ by~\eqref{eq:G_Z sr}.
Hence, we obtain from~\eqref{eq:bound used k-1 sets} and the definition of typicality that $G_B-G^\circ_h$ is $(\mu^{1/2},t,d_B)$-typical with respect to $R$, which implies~\ref{item:G_B-Gcirc typical}.

\medskip
To establish~\ref{item:B sr} for $H=H_{h+1}$, recall that $B_i^H$ is $(\eps_T,d_i)$-super-regular and $(\eps_T,t^{2/3})$-well-intersecting with respect to $G_B$ for all $i\in[r]$. 
Observe that there exists a set $\cS_x$ for every $x\in X_i^H$ with $|\cS_x|\leq t^{2/3}$ and $N_{B_i^H}(x)=V_i\cap N_{G_B}(\cS_x)$.
By the definition of $(B_i^H)^{G_h^\circ}$ in~\eqref{eq:def B^G}, we have that 
$N_{(B_i^H)^{G^\circ_h}}(x)=V_i\cap N_{G_B-G^\circ_h}(\cS_x)$.
This together with~\eqref{eq:bound used k-1 sets}, and~\eqref{eq:bound 1 degree} together with~\eqref{eq:v removed vertices B^circ}, implies that we removed at most $\mu^{2/3}d_in$ edges incident to every vertex to obtain $(B_i^H)^{G_h^\circ}$ from $B_i^H$. 
Now Fact~\ref{fact:regularity robust} yields that $(B_i^H)^{G_h^\circ}$ is $(\mu^{1/5},d_i)$-super-regular.
This implies~\ref{item:B sr}.

\medskip
To establish~\ref{item:B subgraph sr} for $H=H_{h+1}$,
we exploit that $(B_i^H)^{G_h^\circ}$ is $(\mu^{1/5},d_i)$-super-regular and only have to show that every vertex in $(B_i^H)^{G_h^\circ}[Y_i^H,W_i^H]$ has the appropriate degree.
To that end, we consider the following random variables for all $x\in X_i^H$, $v\in V_i$, $i\in[r]$.
For $x$, we consider the random variable that counts the number of activated vertices by $H$ in $N_{(B_i^H)^{G_h^\circ}}(x)$, and by employing~\eqref{eq:leftover intersection k-1 sets} and that $B_i^H$ is $(\eps_T,t)$-well-intersecting with respect to $G_B$, we expect that $N_{(B_i^H)^{G_h^\circ}}(x)\cap W_i^H$ has size $\mu |N_{(B_i^H)^{G_h^\circ}}(x)| \pm 2\eps_Td_i n$.
This yields the appropriate degree for $x$. 
For~$v$, we consider the random variable that counts the number of activated vertices in $N_{(B_i^H)^{G_h^\circ}}(v)$, and by employing Claim~\ref{claim:leftover obs}, we expect that $N_{(B_i^H)^{G_h^\circ}}(v)\cap Y_i^H$ has size $\mu |N_{(B_i^H)^{G_h^\circ}}(v)| \pm 2\eps_Td_i n$.
This yields the appropriate degree for $v$. 
Altogether this implies~\ref{item:B subgraph sr}.

\begin{substep}\label{step:prob to use G_B edge}
Probability to use a $G_B$-edge during the completion\nopagebreak
\end{substep}

We say an edge $\gG\in E(G_B)$ is \defn{used} during the RPP if there exists a graph $H_h\in\cH$ and an edge $\eH\in H_h$ such that $\phi^h(\eH)=\gG$.
In this step we show the following claim that we will apply to establish~\eqref{eq:bound used k-1 sets} and \eqref{eq:bound 1 degree} in Steps~\ref{step:bound used k-1 sets} and \ref{step:codegree bounds2}.

\begin{claim}\label{claim:prob to use G_B edge}
For every edge $\gG\in E(G_B)$, the probability that $\gG$ is used during the RPP is at most~$\mu^{3/4}$.
\end{claim}
\claimproof
Let $\gG\in E(G_B)$ and $\rR\in E(R)$ with $\gG=\{v_i\}_{i\in \rR}\in E(G_B[V_{\cup\rR}])$ be fixed.
We consider different cases and sets of $\cH$-edges that could potentially be embedded onto $\gG$, say in each case with probability at must $\mu^{4/5}$.
In the end, a union bound will establish Claim~\ref{claim:prob to use G_B edge}.
Therefore, for all $m\in[k]$, $J\in\binom{\rR}{m}$, we consider different sets of edges in $\cH[\sX_{\cup {\rR}}]$ where the $m$ vertices corresponding to the clusters in $J$ are either unembedded or activated and can potentially be mapped onto $\{v_i\}_{i\in J}$, and where the remaining $k-m$ vertices corresponding to the clusters in $\rR\sm J$ have already been embedded onto $\{v_i\}_{i\in\rR\sm J}$.
Let $m\in[k]$ and $J\in\binom{\rR}{m}$ be fixed.

\medskip
We first consider the set of $H$-edges $\eH=\{x_i\}_{i\in\rR}$ for all $H\in\cH$ where no vertex of $\eH$ is left unembedded by $\phi_r$ and for every vertex of $\gG$, there is an $H$-vertex that is mapped onto $\gG$.
That is,
the $m$ vertices $\{x_i\}_{i\in J}$ as well as the $m$ vertices $\{v_i\}_{i\in J}$ are activated, the vertices $\{x_i\}_{i\in J}$ can potentially be mapped onto $\{v_i\}_{i\in J}$, and $\phi_r(\{x_i\}_{i\in\rR\sm J})=\{v_i\}_{i\in\rR\sm J}$.
To that end, we first consider the set $\widetilde{E}:=E_{\phi_r}(
\gG, J,\JX=\emptyset,\JV=\emptyset)\sm E_{\phi_r}(
\gG, J,\JX=\emptyset,\JV\neq\emptyset)$
as defined in~\eqref{eq:def leftover g-edge}.
(Note that we defined this set in~\eqref{eq:def leftover g-edge} only in the more convenient way that $\gG\cap V_{\cup \JV}\subseteq \gG\sm \phi_r(X_{\cup\rR}^H)$, which is the reason why we remove the set $E_{\phi_r}(
\gG, J,\JX=\emptyset,\JV\neq\emptyset)$ from the current consideration.)
For an edge $\eH=\{x_i\}_{i\in\rR}\in \widetilde{E}$, 
in order that $\{x_i\}_{i\in J}$ can be mapped onto $\{v_i\}_{i\in J}$ during the completion, it must hold that 
the vertices $\{x_i\}_{i\in J}$ and the vertices $\{v_i\}_{i\in J}$ must become activated.
Since $\phi_r(\xX)\neq\gG$ and $\phi_r(\{x_i\}_{i\in\rR\sm J})=\{v_i\}_{i\in\rR\sm J}$, it holds that $\phi_r(\{x_i\}_{i\in J})\neq\{v_i\}_{i\in J}$.
Hence, the probability that $\{x_i\}_{i\in J}$ and $\{v_i\}_{i\in J}$ become activated is at most $\mu^{m+1}$ because every vertex in $\hX_i^H$ for all $i\in[r]$, $H\in\cH$ is activated independently with probability $\mu$.
Further, by~\ref{item:random packing}, activated vertices $\{x_i\}_{i\in J}$ are mapped onto $\{v_i\}_{i\in J}$ with probability at most $2^{m}\prod_{i\in J}(\mu d_i n)^{-1}$.
Altogether, this implies that that the probability that some edge in $\widetilde{E}$ is mapped onto~$\gG$ is at most 
\begin{align}
\nonumber&\mu^{m+1}2^{m}\prod_{i\in J}(\mu d_i n)^{-1}
\big|E_{\phi_r}\left(
\gG, J,\JX=\emptyset,\JV=\emptyset
\right)\big|
\\\label{eq:prob g is used by only activated vertices}
\stackrel{\text{Claim~\ref{claim:leftover for g}}}{\leq}
&\mu^{m+1}2^{m+1}\prod_{i\in J}(\mu d_i n)^{-1}
\cdot
\gamma^{-1}
n^{m}\prod_{i\in J} d_i
\leq
 \mu^{4/5},
\end{align}
where we used for the last inequality that $\mu\ll\gamma$.

Next, we consider the sets $E_{\phi_r}(
\gG, J,\JX,\JV)$
where $\JXV=\JX\cup \JV$ is non-empty for disjoint $\JX,\JV\subseteq J$.
For an edge $\eH=\{x_i\}_{i\in\rR}\in E_{\phi_r}(
\gG, J,\JX,\JV)$,
the vertices $\{x_i\}_{i\in J}$ are mapped onto $\{v_i\}_{i\in J}$ with probability at most $2^m\prod_{i\in J}(\mu d_i n)^{-1}$ by~\ref{item:random packing}.\COMMENT{Some vertices of $\{x_i\}_{i\in J}$ or $\{v_i\}_{i\in J}$ might have to become activated but we simply do not care. (It happens with probability at most $1$.)}
This implies that the probability that some edge in $E_{\phi_r}(
\gG, J,\JX=\emptyset,\JV)$ is mapped onto~$\gG$ is at most 
\begin{align}
\nonumber&2^{m}\prod_{i\in J}(\mu d_i n)^{-1}
\big|E_{\phi_r}\left(
\gG, J,\JX,\JV
\right)\big|
\\\label{eq:prob g is used v unembedded}
\stackrel{\text{Claim~\ref{claim:leftover for g}}}{\leq}
&2^{m}\prod_{i\in J}(\mu d_i n)^{-1}
\cdot
\eps_T^{1/2}
n^{m}\prod_{i\in J} d_i
\leq \eps_T^{1/3}.
\end{align}

Now, Claim~\ref{claim:prob to use G_B edge} is established by a union bound over all $m\in[k]$, $J\in\binom{\rR}{m}$, and all possible sets $\JX, \JV$ that we considered to be fixed in~\eqref{eq:prob g is used by only activated vertices} and \eqref{eq:prob g is used v unembedded}.
\endclaimproof

\begin{substep}\label{step:bound used k-1 sets}
Proof of Claim~\ref{claim:random packing} -- Bound in~\eqref{eq:bound used k-1 sets}\nopagebreak
\end{substep}

We use Claim~\ref{claim:prob to use G_B edge} to verify the claimed bound in~\eqref{eq:bound used k-1 sets}.
We fix $i\in[r]$, $\cS\subseteq
\bigcup_{\rR\in E(R)\colon i\in\rR} V_{\sqcup\rR\sm\{i\}}$
 with $|\cS|\leq t$ and $S\in\cS$.
By Claim~\ref{claim:prob to use G_B edge}, each $G_B$-edge is used with probability at most $\mu^{3/4}$ during the RPP, and thus, by an application of Lemma~\ref{lem:convenient freedman}\COMMENT{We consider the set of edges $\cE$ in $G_B$ that correspond to $V_i\cap N_{G_B}(\cS)$ and contain the set $S$.
We introduce a random variable $X_e$ for every $e\in\cE$ that indicates whether $e$ is used during the RPP or not.
Hence, $X_e\in\{0,1\}$ and since every edge is used with probability at most $\mu^{3/4}$ during the RPP, we have
$\sum_{e\in\cE}\mathbb{E}'[X_e]\leq \mu^{3/4}|V_i\cap N_{G_B}(S)|$. 
Now we can apply Lemma~\ref{lem:convenient freedman}.
}, we have that 
$|V_i \cap N_{G^\circ_h}(S)\cap N_{G_B}(\cS)| 
\leq 
\mu^{7/10}d_B^{|\cS|}n$ with probability at least, say,  $1-\eul^{-n^{3/4}}$.
Otherwise we detect alarm. 
Together with a union bound this implies the claimed bound in~\eqref{eq:bound used k-1 sets} because $t\cdot\mu^{7/10}d_B^{|\cS|}n\leq \mu^{2/3}d_B^{|\cS|}n$.

\begin{substep}\label{step:codegree bounds2}
Proof of Claim~\ref{claim:random packing} -- Bound for $\rho^H_{G^\circ_h}(v)$\nopagebreak
\end{substep}

We also use Claim~\ref{claim:prob to use G_B edge} to verify the claimed bound for $\rho^H_{G^\circ_h}(v)$ in~\eqref{eq:bound 1 degree}.
(Recall the definition of $\rho_{G^\circ_h}^H(v)$ in~\eqref{eq:def 1 degree}.)
For all $H\in\cH$, $i\in[r]$ and $v\in V_i$, we have that~$v$ has at most $(1+2\eps_T)d_in$ neighbours in~$B_i^H$.
For each such neighbour $x\in N_{B_i^H}(v)$, there exists a set $\cS_x\subseteq \binom{V(G)}{k-1}$ with $|\cS_x|\leq t^{2/3}$ such that $N_{B_i^H}(x)=N_{G_B}(\cS_x)$.
Hence, there exist at most $2t^{2/3}d_in$ edges $\gG=S\cup\{v\}$ in $G_B$ for some $S\in\bigcup_{x\in N_{B_i^H}(v)}\cS_x$.
By Claim~\ref{claim:prob to use G_B edge}, each such $G_B$-edge~$\gG$ is used with probability at most~$\mu^{3/4}$ during the entire RPP.
We therefore expect that for each $h\in[|\cH|]$ at most $\mu^{3/4}2t^{2/3}d_in$ edges incident to $v$ in $B_i^H$ have to be removed when we obtain $(B_i^H)^{G^\circ_h}$.
Hence, by an application of  Lemma~\ref{lem:convenient freedman} and a union bound, we obtain that $\rho_{G^\circ_h}^H(v)\leq \mu^{2/3}d_in$  for all $h\in[|\cH|]$, $H\in\cH$, $i\in[r]$, $v\in V_i$  with probability at least, say, $1-\eul^{-n^{3/4}}$.
Otherwise we detect alarm. 
This implies the claimed bound in~\eqref{eq:bound 1 degree}.

\begin{substep}\label{step:few activated weight on vertex testers}
Proof of Claim~\ref{claim:random packing} -- Establishing~\ref{item: vertex tester1} of Lemma~\ref{lem:main matching}
\end{substep}

In order to establish conclusion~\ref{item: vertex tester1} of Lemma~\ref{lem:main matching}, let $(\omega,\bc)\in\cW_{ver}$ with centres $\bc=\{c_i\}_{i\in I}$ in~$I$ be fixed.
We split the proof into two parts, depending on whether we activate a vertex of a tuple $\xX\in\sX_{\sqcup {I}}$ that $\phi_r$ already mapped onto $\bc$, or whether a tuple $\xX\in\sX_{\sqcup {I}}$ contains vertices in $\cX_r\sm \cX_r^{\phi_r}$ and $\xX$ is mapped onto $\bc$ during the completion.

Let us first consider tuples $\xX\in\sX_{\sqcup {I}}$ that we have already embedded onto $\bc$, that is, $\phi_r(\xX)=\bc$, and where some vertices of $\xX$ become activated.
We claim that during the entire RPP not too many such tuples become activated.
That is, we claim that
\begin{align}\label{eq:few activated weight on vertex testers}
\omega\left(
\bigcup_{H\in \cH}\left\{
\xX\in X_{\sqcup I}^H\colon 
\phi_r(\xX)=\bc,
\xX\cap Y_{\cup I}^H\neq\emptyset
\right\}
\right)
\leq \mu^{1/2}\omega(\phi_r^{-1}(\bc))+n^{\eps_T},
\end{align}
with high probability.
To see~\eqref{eq:few activated weight on vertex testers}, note that for every $\xX\in\sX_{\sqcup I}$, the probability that $\xX$ contains an activated vertex is at most $|I|\mu$.
Hence, an application of Theorem~\ref{thm:McDiarmid} and a union bound yield~\eqref{eq:few activated weight on vertex testers} with probability, say, at least $1-\eul^{-n^{\eps}}$.

\medskip
Next, we consider tuples $\xX\in\sX_{\sqcup {I}}$ that have not been embedded onto $\bc$ by $\phi_r$.
To that end, we fix $m\in[|I|]$ and $J\in\binom{I}{m}$, and a pair of disjoint sets $\JX,\JV\subseteq J$.
We aim to control the $\omega$-weight on the tuples in $E_{\phi_r}(\bc,J,\JX,\JV)$ as defined in~\eqref{eq:def leftover g-edge}.
Analogously as in Claim~\ref{claim:leftover for g}, we can employ \ind{r}\ref{edge weight} and~\ref{treffer auf g} for the edge testers that we defined for $(\omega,\bc)$ in~\eqref{eq:def initial weight vertex tester}. 
Proceeding as in Claim~\ref{claim:leftover for g}\COMMENT{In particular, we have to sum over all possible patterns.} yields that\looseness=-1
\begin{align}\label{eq:weight leftover for g}
\omega\left(E_{\phi_r}(\bc,J,\JX,\JV)\right)
\leq
\left(\IND\{\JXV=\emptyset \}+2\eps_T\right)
\prod_{i\in  J} d_i \frac{\omega(\sX_{\sqcup {I}})}{n^{|I|-|J|}}+ n^{2\eps_T}.
\end{align}
We can now proceed similarly as in the proof of Claim~\ref{claim:prob to use G_B edge} in Step~\ref{step:prob to use G_B edge}; that is, we consider different cases and sets of tuples $\xX\in \sX_{\sqcup {I}}$ that could potentially be embedded onto $\bc$, and derive an upper bound on the expected $\omega$-weight used in such a case. 
In the end, a union bound over all these cases gives us an upper bound on the total expected $\omega$-weight of tuples that are embedded onto $\bc$ during the completion.
Therefore, we consider different choices for $\JX,\JV\subseteq J$.

\medskip
For $\xX=\{x_i\}_{i\in I} \in E_{\phi_r}(\bc,J,\JX=\emptyset,\JV=\emptyset)\sm E_{\phi_r}(\bc,J,\JX=\emptyset,\JV\neq\emptyset)$, in order that $\{x_i\}_{i\in J}$ can be mapped onto $\{c_i\}_{i\in J}$ during the RPP, it must hold that the vertices $\{x_i\}_{i\in J}$ and $\{c_i\}_{i\in J}$ become activated because $\JV=\emptyset$. 
Since $\phi_r(\{x_i\}_{i\in J})\neq\{c_i\}_{i\in J}$, this happens with probability at most $\mu^{m+1}$.
Further, by~\ref{item:random packing}, the activated vertices $\{x_i\}_{i\in J}$ are mapped onto $\{c_i\}_{i\in J}$ with probability at most $2^m\prod_{i\in J}(\mu d_in)^{-1}$. 
Altogether, this implies that the expected weight of edges in $E_{\phi_r}(\bc,J,\JX=\emptyset,\JV=\emptyset)\sm E_{\phi_r}(\bc,J,\JX=\emptyset,\JV\neq\emptyset)$ that are mapped onto $\bc$ during the RPP is at most
\begin{align}
\nonumber&\mu^{m+1}2^{m}\prod_{i\in J}(\mu d_i n)^{-1}
\omega\left(E_{\phi_r}\left(
\bc, J,\JX=\emptyset,\JV=\emptyset
\right)\right)
\\\label{eq:prob c is used v unembeddedby only activated vertices}
\stackrel{\eqref{eq:weight leftover for g}}{\leq}
&\mu^{m+1}2^{m}\prod_{i\in J}(\mu d_i n)^{-1}
\left(2\prod_{i\in J} d_i \frac{\omega(\sX_{\sqcup {I}})}{n^{|I|-|J|}}+n^{2\eps_T}\right)
\leq
\mu^{4/5}\frac{\omega(\sX_{\sqcup {I}})}{n^{|I|}}+n^{3\eps_T}.
\end{align}

\medskip
Next, we consider the set $E_{\phi_r}(\bc,J,\JX,\JV)$ for disjoint but fixed $\JX,\JV\subseteq J$ such that $\JXV=\JX\cup\JV\neq\emptyset$.
The vertices $\{x_i\}_{i\in J}$ are mapped onto $\{c_i\}_{i\in J}$ with probability at most $2^m\prod_{i\in J}(\mu d_in)^{-1}$ by~\ref{item:random packing}.
Hence, 
the expected weight of edges in $E_{\phi_r}(\bc,J,\JX,\JV)$ that are mapped onto $\bc$ during the RPP is at most 
\begin{align}\label{eq:prob c is used v unembedded}
2^{m}\prod_{i\in J}(\mu d_i n)^{-1}
\omega\left(E_{\phi_r}\left(
\bc, J,\JX,\JV
\right)\right)
\stackrel{\eqref{eq:weight leftover for g}}{\leq}
\eps_T^{1/2}\frac{\omega(\sX_{\sqcup {I}})}{n^{|I|}}+n^{3\eps_T}.
\end{align}

\medskip
Altogether, a union bound over all $m\in[|I|]$, $J\in\binom{I}{m}$, and all sets $\JX, \JV\subseteq J$ that we considered to be fixed in~\eqref{eq:prob c is used v unembeddedby only activated vertices} and \eqref{eq:prob c is used v unembedded} together with an application of Lemma~\ref{lem:convenient freedman} yields that 
\begin{align*}
\omega\left(
\left\{
\xX\in \sX_{\sqcup I}\colon \phi_r(\xX)\neq\bc, \phi^{|\cH|}(\xX)=\bc
\right\}
\right)
\leq \mu^{1/2}\frac{\omega(\sX_{\sqcup {I}})}{n^{|I|}}+n^{4\eps_T}
\end{align*}
with probability at least, say, $1-\eul^{-n^{2\eps_T}}$.\COMMENT{
Let $X=\omega\left(
\left\{
\xX\in \sX_{\sqcup I}\colon \phi_r(\xX)\neq\bc, \phi^{|\cH|}(\xX)=\bc
\right\}
\right)$ and for each tuple $\xX\in \sX_{\sqcup I}$ with $\phi_r(\xX)\neq\bc$ let $X_\xX$ be the random variable taking value $\omega(\xX)\leq\alpha^{-1}$ if $\phi^{|\cH|}(\xX)=\bc$ and $0$ otherwise.   
We have $X=\sum_\xX X_\xX$ and by~\eqref{eq:prob c is used v unembeddedby only activated vertices} and \eqref{eq:prob c is used v unembedded}, 
we have that
$$\sum_\xX \mathbb{E}'[X_\xX]\leq \mu^{3/4}\frac{\omega(\sX_{\sqcup {I}})}{n^{|I|}}+n^{7\eps_T/2}.$$
Hence, we can apply Lemma~\ref{lem:convenient freedman}, which establishes the displayed equation with probability at least
$$1-2\exp\left(-
\frac{\mu^{3/4}\frac{\omega(\sX_{\sqcup {I}})}{n^{|I|}}+n^{7\eps_T/2}}{6\alpha^{-1}}
\right)
\geq
1-\exp\left(
- n^{2\eps_T}
\right).
$$
}

Combining this with~\eqref{eq:few activated weight on vertex testers} yields
\begin{align*}
\omega((\phi^{|\cH|})^{-1}(\bc))
&
\stackrel{\hphantom{\text{\ind{r}\ref{vertex tester}}}}{=}
(1\pm\mu^{1/2})\omega(\phi_r^{-1}(\bc))
\pm  \mu^{1/2}\frac{\omega(\sX_{\sqcup {I}})}{n^{|I|}}
\pm 2n^{4\eps_T}
\\
&\stackrel{\text{\ind{r}\ref{vertex tester}}}{=}
(1\pm\alpha)\frac{\omega(\sX_{\sqcup {I}})}{n^{|I|}}\pm n^{\alpha}.
\end{align*}
This establishes conclusion~\ref{item: vertex tester1} of Lemma~\ref{lem:main matching} and completes the proof of Claim~\ref{claim:random packing}.
\endclaimproof

\begin{substep}\label{step:finishing the completion}
Finishing the completion
\end{substep}

As the RPP outputs $\phi^{|\cH|}$ with positive probability by Claim~\ref{claim:random packing}, we obtain a packing $\phi:=\phi^{|\cH|}$ of $\cH$ into $G$ which clearly satisfies conclusions~\ref{item:partition1} and~\ref{item: vertex tester1} of Lemma~\ref{lem:main matching}.
Note that by the construction of $\phi^{|\cH|}$, we have that $\phi|_{X_i^H}=\phi_r|_{X_i^H\sm Y_i^H}\cup \phi^H|_{Y_i^H}$ for all $H\in\cH$, $i\in[r]$. 
Since $|Y_i^H|=(1\pm\eps_T^{1/2})\mu n$, we therefore merely modified $\phi_r$ to obtain $\phi$ and thus,~\ind{r}\ref{set tester} easily implies conclusion~\ref{item:set testers1} of Lemma~\ref{lem:main matching}.\COMMENT{
To establish~\ref{item:set testers1}, let $(W,Y_1,\ldots,Y_\ell)\in \cW_{set}$ with $W\subseteq V_i$ for $i\in[r]$ and let $H_{m_1},\ldots, H_{m_\ell}\in\cH$ be such that $Y_j\subseteq X_i^{H_{m_j}}$ for each $j\in[\ell]$.
We have that $$|W\cap\bigcap_{j\in [\ell]}\phi(Y_j\cap Y_i^{H_{m_j}})|\leq 2\mu n,$$
which together with~\ind{r}\ref{set tester} yields conclusion~\ref{item:set testers1}.}
This completes the proof of Lemma~\ref{lem:main matching}.
\endproof

\section{Proof of the main results}\label{sec:mainproofs}

In this section we prove Theorem~\ref{thm:main_new} and Theorem~\ref{thm:normalwithweights}.

\lateproof{Theorem~\ref{thm:main_new}}
\noindent
Our general approach is as follows.
Given a blow-up instance $(\cH,G,R,\cX,\cV)$, we refine the vertex partition~$\cX$ of the graphs in~$\cH$ using Lemma~\ref{lem:refining partitions} and we randomly refine the vertex partition~$\cV$ of~$G$ accordingly. 
Afterwards we can apply Lemma~\ref{lem:main matching} to obtain the required packing of $\cH$ into $G$.

Suppose {$1/n\ll\eps\ll1/t\ll\beta\ll \alpha,1/k$} for a new parameter $\beta$.
For each $(W,Y_1,\ldots,Y_m)\in\cW_{set}$ with $W\subseteq V_i$, $i\in[r]$, and each $\ell\in[m]$, let $\omega_{Y_{\ell}}\colon \bigcup_{H\in\cH}X_i^H\to\{0,1\}$ be such that $\omega_{Y_{\ell}}(x)=\IND{\{x\in Y_{\ell} \}}$, and let $\cW_Y$ be the set containing all those weight functions.\COMMENT{$\cW_{Y}$ is used to control the size of the refined $Y_{\ell,j}$.}
Further, for all $\rR\in E(R)$, let $\omega_\rR\colon\sX_{\sqcup {\rR}}\to\{0,1\}$ be defined by $\omega_\rR(\xX):=\IND\{\xX\in E(\cH) \}$, and let $\cW_\cH$ be the set containing all those weight functions.
We apply Lemma~\ref{lem:refining partitions} to $\cH$ with weight functions $\{\omega\colon (\omega,\bc)\in\cW_{ver} \}\cup \cW_{Y}\cup \cW_\cH$.
This yields a refined partition $\cX'=(X_{i,j}^H)_{H\in\cH,i\in[r],j\in[\beta^{-1}]}$ of $\cH$ such that for all $H\in\cH$ and $i\in[r]$, the partitions $(X_{i,j}^H)_{j\in[\beta^{-1}]}$ of $X_i^H$ satisfy the conclusions~\ref{item:split1}--\ref{item:split weights} of Lemma~\ref{lem:refining partitions}.

Let $R'$ be the $k$-graph with vertex set $[r]\times[\beta^{-1}]$ and edge set $$\left\{ \{(i_\ell,j_\ell) \}_{\ell\in[k]}\colon \{i_\ell\}_{\ell\in[k]} \in E(R), j_\ell\in[\beta^{-1}] \right\}.$$
Note that $\Delta(R')\leq \alpha^{-1}\beta^{-(k-1)}$ because $\Delta(R)\leq\alpha^{-1}$.
Let $n':=\beta n$.

Employing conclusion~\ref{item:split weights} of Lemma~\ref{lem:refining partitions} for the weight functions in $\cW_\cH$ implies for all $\{(i_\ell,j_\ell) \}_{\ell\in[k]}\in E(R')$ with $\rR:=\{i_\ell\}_{\ell\in[k]}\in E(R)$ that
\begin{align*}
\sum_{H\in \cH} e_H(X_{i_1,j_1}^H,\ldots,X_{i_k,j_k}^H) 
&\leq (1+\eps)\beta^k \omega_\rR(\sX_{\sqcup {\rR}})+ n^{1+\eps}
\\
&= (1+\eps)\beta^k e_\cH(\sX_{\sqcup {\rR}})+ n^{1+\eps}
\leq (1-\alpha/2)dn'^k,
\end{align*}
because by assumption, $e_\cH(\sX_{\sqcup {\rR}})\leq(1-\alpha)dn^k$ and $d\geq n^{-\eps}$.

Further, note that Lemma~\ref{lem:refining partitions}\ref{item:split1} implies for each $H\in\cH$ that $H[X_{\cup \rR}^H]$ is a matching if $\rR\in E(R')$ and empty if $\rR\in\binom{[r]\times[\beta^{-1}]}{k}\sm E(R')$. 

According to the refinement $\cX'$ of $\cX$, we claim that there exists a refined partition $\cV'=(V_{i,j})_{i\in[r],j\in[\beta^{-1}]}$ of $\cV$, where $(V_{i,j})_{j\in[\beta^{-1}]}$ is a partition of~$V_i$ for every $i\in[r]$ such that
\begin{enumerate}[label=(\alph*)]
\item\label{item:refined W} $|W\cap V_{i,j}|=\beta|W|\pm\beta^{3/2}n$ for all $(W,Y_1,\ldots,Y_m)\in \cW_{set}$ and $j\in[\beta^{-1}]$ with $W\subseteq V_i$;

\item\label{item:refined blowup instance} $\sB':=(\cH,G,R',\cX',\cV')$ is an $(\eps^{1/2},t,d)$-typical, $\beta^{-k}$-bounded blow-up instance of size $(n',k,\beta^{-1}r)$ with $n'=\beta n$.
\end{enumerate}
{The existence of such a partition $\cV'$ can be seen by a probabilistic argument.
For all $i\in[r]$ and $j\in[\beta^{-1}]$, we take disjoint subset $V_{i,j}$ of $V_i$ of size exactly $|X_{i,j}^H|$ uniformly and independently at random. 
We analyse the probability that~\ref{item:refined W} or~\ref{item:refined blowup instance} are not satisfied. 
To that end, we consider the slightly different random experiment where we assign for all $i\in[r]$, every vertex in $V_i$ uniformly and independently at random to some $V_{i,j}$ for $j\in[\beta^{-1}]$. 
For this experiment and a fixed $i\in[r]$, we consider the bad events that~\ref{item:refined W} or~\ref{item:refined blowup instance} are not satisfied; in such a case, we say the experiment \emph{fails} (in step $i$).
Standard properties of the multinomial distribution yield that $|V_{i,j}|=|X_{i,j}^H|$ for all $j\in[\beta^{-1}]$, $H\in\cH$ with probability at least $\Omega(n^{-\beta^{-1}})$. 
Hence, together with Theorem~\ref{thm:McDiarmid} and a union bound, this yields that the original experiment fails in step $i$ with probability, say, at most~$\eul^{-n^{1/2}}$.
Since in fact we take $V_{i,j}$ of size exactly $|X_{i,j}^H|$,
this altogether implies the existence of a refined partition $\cV'$ of $\cV$ satisfying~\ref{item:refined W} and~\ref{item:refined blowup instance} with positive probability.}\COMMENT{For the $k$-multinomial coefficient, we obtain with Stirling's formula that
$$\binom{n}{\frac{n}{k},\ldots,\frac{n}{k}}=\frac{n!}{\left(\frac{n}{k}\right)!^k}\sim \frac{\sqrt{2\pi n}\left(\frac{n}{e}\right)^n}{\left(\sqrt{\frac{2\pi n}{k}}\left(\frac{n}{ke}\right)^{n/k}\right)^k}=\frac{k^{k/2}k^n}{\left(2\pi n \right)^{(k-1)/2}}.$$
This implies that 
$$\prob{\forall j\in[\Delta^2]\colon\left|\tau_i(j)\right|=\frac{n}{\Delta^2}}=\Omega(n^{-\Delta^2}),$$
for every $i\in[r]$. Thus,
$$\prob{\forall i\in[r], j\in[\Delta^2]\colon\left|\tau_i(j)\right|=\frac{n}{\Delta^2}}=\Omega(n^{-r\Delta^2}).$$
}

We show how to adapt the vertex and set testers from the original blow-up instance to the blow-up instance $\sB'$.
For each $(W,Y_1,\ldots,Y_m)\in \cW_{set}$ and distinct $H_1,\ldots,H_m\in\cH$ such that $W\subseteq V_i$ for some $i\in[r]$ and $Y_{\ell}\subseteq X_i^{H_{\ell}}$ for all $\ell\in[m]$, 
we define $(W_j,Y_{1,j},\ldots,Y_{m,j})$ by setting $W_j:=W\cap V_{i,j}$ and $Y_{\ell,j}:=Y_{\ell}\cap X_{i,j}^{H_{\ell}}$ for all $j\in[\beta^{-1}], \ell\in[m]$. 
By~\ref{item:refined W}, we conclude that $|W_j|=\beta |W|\pm\beta^{3/2}n$. 
Employing conclusion~\ref{item:split weights one H, |I|=1} of Lemma~\ref{lem:refining partitions} for the weight function $\omega_{Y_{\ell}}\in\cW_Y$, we have that 
\begin{align}\label{eq:size Y_l,j}
|Y_{\ell,j}|=\omega_{Y_{\ell}}(X_{i,j}^{H_{\ell}})=\beta \omega_{Y_{\ell}}(X_i^{H_{\ell}}) \pm \beta^{3/2}n =\beta |Y_{\ell}|\pm \beta^{3/2}n.
\end{align}
Let $\cW_{set}':=\{(W_j,Y_{1,j},\ldots,Y_{m,j})\colon 
j\in[\beta^{-1}],(W,Y_1,\ldots,Y_m)\in \cW_{set} \}$.
For each $(\omega,\bc)\in\cW_{ver}$ with centres $\bc=\{c_i\}_{i\in I}$  in $I\subseteq[r]$ and multiset $\{j_i\}_{i\in I}\subseteq[\beta^{-1}]$ such that $c_i\in V_{i,j_i}$ for each $i\in I$, let 
 $\omega':=\omega|_{\bigcup_{H\in \cH,i\in I}X_{i,j_i}^H}$ and $\cW_{ver}':=\{(\omega',\bc)\colon (\omega,\bc)\in \cW_{ver} \}$.

Hence, we can apply Lemma~\ref{lem:main matching} to $\sB'$ with set testers $\cW_{set}'$ and vertex testers~$\cW_{ver}'$ as follows:\COMMENT{Note that $n'^{2\log n'}=(\beta n)^{2\log(\beta n)}
\geq(\beta n)^{3\log n/2}\geq\beta^{-2}n^{\log n}\geq |\cW_{set}'|,|\cW_{ver}'|.$}

\begin{center}
\begin{tabular}{r||c|c|c|c|c|c}
parameter & $n'$  & $\eps^{1/2}$ & $t$ & {$\beta^k$} & $d$ &  $r\beta^{-1}$
\\ \hline
plays the role of & $n$  & $\eps$  & $t$ & $\alpha$ & $d$  & $r$ 
\vphantom{\Big(}
\end{tabular}
\end{center} 
This yields a packing of $\cH$ into $G$ such that 
\begin{enumerate}[label=(\Roman*)]
\item\label{item:I} $\phi(X_{i,j}^H)= V_{i,j}$ for all $i\in[r]$,  $j\in[\beta^{-1}]$, $H\in\cH$;
\item\label{item:II} $|W_j\cap \bigcap_{\ell\in [m]}\phi(Y_{\ell,j})|= |W_j||Y_{1,j}|\cdots |Y_{m,j}|/n'^m \pm {\beta^k} n'$ for all $(W_j,Y_{1,j},\ldots ,Y_{\ell,j})\in \cW_{set}'$;
\item \label{item:III} 
$\omega'( \phi^{-1}(\bc))
)
=(1\pm\beta^k)
\omega'(
\bigcup_{H\in\cH}(\bigsqcup_{i\in I}X_{i,j_i}^H)
)
/n'^{|I|}
\pm n'^{\beta^k}$ for all $(\omega',\bc)\in \cW_{ver}'$ with centres $\bc=\{c_i\}_{i\in I}$ in $I\subseteq[r]$ and multiset $\{j_i\}_{i\in I}\subseteq[\beta^{-1}]$ such that $c_i\in V_{i,j_i}$ for each $i\in I$.
\end{enumerate}

\medskip
For $(W,Y_1,\ldots,Y_m)\in\cW_{set}$, we conclude that
\begin{align*}
\textstyle\left|W\cap \bigcap_{\ell\in [m]}\phi(Y_{\ell})\right| 
&\stackrel{\hphantom{\rm \ref{item:II},\ref{item:refined W},\eqref{eq:size Y_l,j}}}{=} \sum_{j\in[\beta^{-1}]}\textstyle\left|W_j\cap \bigcap_{\ell\in [m]}\phi(Y_{\ell,j})\right|\\
&\stackrel{\rm \ref{item:II},\ref{item:refined W},\eqref{eq:size Y_l,j}}{=} \sum_{j\in[\beta^{-1}]}\left(\frac{\beta^{m+1}\left(|W||Y_1|\cdots|Y_m|\pm\beta^{1/3}n^{m+1}\right)}{(\beta n)^m}\pm{\beta^k} n'\right)\\
&\stackrel{\hphantom{\rm \ref{item:II},\ref{item:refined W},\eqref{eq:size Y_l,j}}}{=}|W||Y_1|\cdots|Y_m|/n^m \pm\alpha n.
\end{align*}
\COMMENT{We have $$|W_j||Y_{1,j}|\cdots |Y_{m,j}|=\beta^{m+1}\Big(|W||Y_1|\cdots|Y_m|\pm\beta^{1/3}n^{m+1}\Big)$$
because $|W_j|=\beta |W|\pm\beta^{3/2}n$ by~\ref{item:refined W} and $|Y_{\ell,j}|=\beta |Y_{\ell}|\pm\beta^{3/2}n$ by~\eqref{eq:size Y_l,j}. Hence,~\ref{item:II} implies that 
\begin{align*}
\sum_{j\in[\beta^{-1}]}\left|W_j\cap \bigcap_{\ell\in [m]}\phi(Y_{\ell,j})\right|
&=\sum_{j\in[\beta^{-1}]}\left(\frac{\beta^{m+1}\left(|W||Y_1|\cdots|Y_m|\pm\beta^{1/3}n^{m+1}\right)}{(\beta n)^m}\pm{\beta^2} n'\right)
\\&=\frac{|W||Y_1|\cdots|Y_\ell|}{n^m}\pm\beta^{1/3}n\pm {\beta^2} n
\\&=|W||Y_1|\cdots|Y_\ell|/n^m \pm\alpha n.
\end{align*}
}
This establishes Theorem~\ref{thm:main_new}\ref{item:set testers}.

In order to establish Theorem~\ref{thm:main_new}\ref{item: vertex tester}, we fix $(\omega,\bc)\in\cW_{ver}$  with centres $\bc=\{c_i\}_{i\in I}$ for $I\subseteq[r]$ and multiset $\{j_i\}_{i\in I}\subseteq[\beta^{-1}]$ such that $c_i\in V_{i,j_i}$ for each $i\in I$, and we fix the corresponding tuple $(\omega',\bc)\in \cW_{ver}'$.
We conclude that 
\begin{align*}
\omega( \phi^{-1}(\bc))
&\stackrel{\rm\ref{item:I}}{=}
~\omega'( \phi^{-1}(\bc))
\stackrel{\rm \ref{item:III}}{=}
(1\pm\beta^k)
\frac{\omega'\left(
\bigcup_{H\in\cH}\left(\bigsqcup_{i\in I}X_{i,j_i}^H\right)
\right)
}
{n'^{|I|}}
\pm n'^{\beta^k}
\\
&\stackrel{\hphantom{\rm\ref{item:split weights}}}{=}
(1\pm\beta^k)
\frac{\omega\left(
\bigcup_{H\in\cH}\left(\bigsqcup_{i\in I}X_{i,j_i}^H\right)
\right)
}
{n'^{|I|}}
\pm n'^{\beta^k}
\\&\stackrel{{\rm\ref{item:split weights}}}{=}
(1\pm\beta^k)
\frac{(1\pm\eps)\beta^{|I|} \omega(\sX_{\sqcup {I}})
\pm n^{1+\eps}
}
{(\beta n)^{|I|}}
\pm n'^{\beta^k}
=(1\pm\alpha) \omega(\sX_{\sqcup {I}})/n^{|I|}
\pm n^\alpha,
\end{align*}
where we employed conclusion~\ref{item:split weights} of Lemma~\ref{lem:refining partitions} in the penultimate equation.
This establishes Theorem~\ref{thm:main_new}\ref{item: vertex tester} and completes the proof.
\endproof

We now proceed to the proof of Theorem~\ref{thm:normalwithweights}.
The highlevel strategy is similar as in the proof of Theorem~\ref{thm:main_new}.
Additionally, we group the hypergraphs in $\cH$ into $P=\plog n$ many collections of hypergraphs $\cH_1,\ldots,\cH_P$ and accordingly we partition the edge set of the host graph $G$ into $G_1,\ldots,G_P$ subgraphs. 
Afterwards, we partition the vertex sets of the graphs in $\cH$ via Lemma~\ref{lem:refining partitions r=1} and randomly partition the vertex set of each $G_p$ for $p\in [P]$ accordingly. 
Then, we can iteratively apply Lemma~\ref{lem:main matching} for each $p\in[P]$ to map $\cH_p$ into $G_p$. 
Note that this yields a packing of $\cH$ into~$G$, and
considering $P=\plog n$ partitions enables us to establish conclusions~\ref{item:set testers2} and~\ref{item: vertex tester2} of Theorem~\ref{thm:normalwithweights}.

\lateproof{Theorem~\ref{thm:normalwithweights}}
We set $P:=\log^{t}n$ and suppose $1/n\ll\eps\ll1/t\ll\beta\ll \alpha,1/k$ for  a new parameter $\beta$.

First, we group the graphs in $\cH$ into $P$ collections $\cH_1,\ldots,\cH_P$ with roughly equally many edges. 
That is, we claim that there exists a partition of $\cH$ into $P$ collections of graphs $\cH_1,\ldots,\cH_P$ such that each $H\in\cH$ belongs to exactly one $\cH_p$ for $p\in[P]$, and for each $p\in[P]$, we have
\begin{align}\label{eq:edges sliced Hs}
e(\cH_p)\leq (1+\eps)P^{-1}e(\cH)+n^{1+\eps},
\end{align}
and for every $\omega\in\cW_{ver}$ with $\omega(V(\cH))\geq n^{1+\eps}$ and each $p\in[P]$, we have
\begin{align}\label{eq:sliced weight}
\omega(V(\cH_p))
=
(1\pm\eps)P^{-1}\omega(V(\cH)).
\end{align}
The existence of $\cH_1,\ldots,\cH_P$ can be easily seen by assigning every graph $H\in\cH$ to one collection~$\cH_p$ for $p\in[P]$ uniformly and independently at random.

We now aim to apply Lemma~\ref{lem:refining partitions r=1} to each $\cH_p$.
Let $p\in[P]$ be fixed.
For each $(W,Y_1,\ldots,Y_m)\in\cW_{set}$ and each $\ell\in[m]$ such that $Y_\ell\subseteq V(H_\ell)$ for $H_\ell\in\cH_p$, let $\omega_{Y_{\ell}}\colon V(H_\ell)\to\{0,1\}$ be such that $\omega_{Y_{\ell}}(x)=\IND{\{x\in Y_{\ell} \}}$, and let $\cW_Y$ be the set containing all those weight functions.\COMMENT{$\cW_{Y}$ is used to control the size of the refined $Y_{\ell,j}$.}
Further, let $\omega_{edges}\colon\bigcup_{H\in \cH} \binom{V(H)}{k}_\prec\to\{0,1\}$ be defined by $\omega_{edges}((x_1,\ldots,x_k)):=\IND\{\{x_1,\ldots,x_k\}\in E(\cH_p) \}$, and let~$\cW_{\cH_p}$ be the set containing all those weight functions.
We apply Lemma~\ref{lem:refining partitions r=1} to~$\cH_p$ with weight functions $\{\omega|_{\cup V(\cH_p)}\colon (\omega,\bc)\in\cW_{ver} \}\cup \cW_{Y}\cup \cW_{\cH_p}$.
This yields a partition $\cX_p=(X_{j}^{H})_{H\in\cH_p,j\in[\beta^{-1}]}$ of $\cH_p$ such that for all $H\in\cH_p$, the partitions $(X_{j}^H)_{j\in[\beta^{-1}]}$ of $V(H)$ satisfy the conclusions~\ref{item:split1}--\ref{item:split weights} of Lemma~\ref{lem:refining partitions r=1}.

For each $p\in [P]$, let $R_p$ be the $k$-graph with vertex set $[\beta^{-1}]$ and $\rR\in\binom{[\beta^{-1}]}{k}$ is an edge in~$R_p$ if~$H[X_{\sqcup\rR}^H]$ is non-empty for some $H\in \cH_p$.
Clearly, $\Delta(R_p)\leq \beta^{-k}$.
Let $n':=\beta n$.

Employing conclusion~\ref{item:split' weights} of Lemma~\ref{lem:refining partitions r=1} for the weight functions in $\cW_{\cH_p}$ yields for all $\rR=\{i_1,\ldots,i_k\}\in E(R_p)$ that
\begin{align}
\nonumber
\sum_{H\in \cH_p}
e_H(X_{\sqcup\rR}^H)
&= \omega_{edges}(
\textstyle \bigcup_{H\in\cH_p}(X_{i_1}^H\times\ldots \times X_{i_k}^H)
)
\leq 
(1+\beta^{1/2})\beta^k\omega_{edges}(V(\cH_p))+n^{1+\eps}
\\&=
(1+\beta^{1/2})\beta^kk! e(\cH_p)+n^{1+\eps}
\stackrel{\eqref{eq:edges sliced Hs}}{\leq}
(1-\alpha/2)P^{-1}dn'^k
\end{align}
because by assumption, $e(\cH)\leq (1-\alpha)e(G)$. 

\medskip
Now, we want to prepare $G$ accordingly to $\cH_1,\ldots,\cH_P$ and their partitions.
To that end, we first partition $G$ into $P$ edge-disjoint spanning subgraphs $G_1,\ldots,G_P$ such that $G_p$ is $(\eps^{1/2},t,P^{-1}d)$-typical for every $p\in[P]$. 
The existence of $G_1,\ldots,G_P$ can be seen by assigning every edge in $G$ to one subgraph $G_p$ for $p\in[P]$ uniformly and independently at random.

Further, we claim that there exist  partitions $\cV_p=(V_j)_{j\in[\beta^{-1}]}$ of $V(G_p)$ according to the partition~$\cX_p$ of~$\cH_p$, such that for every $p\in[P]$ and $\cV_p=(V_j)_{j\in[\beta^{-1}]}$, we have
\begin{enumerate}[label=(\alph*)]
\item\label{item:refined W r=1} $|W\cap V_{j}|=\beta|W|\pm\beta^{3/2}n$ for all $(W,Y_1,\ldots,Y_m)\in \cW_{set}$ and $j\in[\beta^{-1}]$,

\item\label{item:refined blowup instance r=1} $\sB_p:=(\cH_p,G_p,R_p,\cX_p,\cV_p)$ is an $(\eps^{1/2},t,P^{-1}d)$-typical, $\beta^{-k}$-bounded blow-up instance of size $(n',k,\beta^{-1})$ with $n'=\beta n$,
\end{enumerate}

as well as

\begin{enumerate}[label=(\alph*)]
\setcounter{enumi}{2}
\item\label{item:centres collide} 
$\sum_{p\in[P]\colon\text{centres collide}} \omega(V(\cH_p))\leq \beta^{1/2}\omega(V(\cH))$ for all $(\omega,\bc)\in\cW_{ver}$ {with $\omega(V(\cH))\geq n^{1+\eps}$}, where we say that the \emph{centres $\bc$ collide (with respect to $\cV_p$)} if $|\bc\cap V_j|\geq 2$ for some $V_j\in\cV_p$. 
\end{enumerate}
The existence of such partitions $\cV_p$ can be seen by assigning for every $p\in[P]$, every vertex to some~$V_j$ for $j\in[\beta^{-1}]$ uniformly and independently at random.
Theorem~\ref{thm:McDiarmid} and a union bound establish~\ref{item:refined W r=1} with probability, say, at least $1-\eul^{-n^{1/2}}$.
For~\ref{item:refined blowup instance r=1}, note that standard properties of the multinomial distribution yield for $p\in[P]$ that $|V_j|=|X_j^H|$ for all $j\in[\beta^{-1}]$ and $H\in\cH$ with probability at least $\Omega(n^{-\beta^{-1}})$. 
For~\ref{item:centres collide}, note that for $p\in[P]$, the probability that the centres~$\bc$ collide with respect to $\cV_p$ is at most $k^2\beta$.
By~\eqref{eq:sliced weight}, we therefore expect that at most 
$\sum_{p\in[P]}k^2\beta\omega(V(\cH_p))
\leq
2k^2\beta \omega(V(\cH))$
weight of $\omega$ collides in~\ref{item:centres collide}. 
Since $P$ grows sufficiently fast in terms of $n$, we can establish concentration. {That is, Theorem~\ref{thm:McDiarmid} and a union bound that yield~\ref{item:centres collide} with probability, say, at least $1-n^{-\log n}$.
Hence, a final union bound yields} the existence of these partitions $\cV_p$ for $p\in[P]$ satisfying~\ref{item:refined W r=1}--\ref{item:centres collide} with positive probability.

Next, we iteratively apply Lemma~\ref{lem:main matching} to $\sB_p$ for $p\in[P]$ which yields a packing $\phi_p$ of $\cH_p$ into~$G_p$.
Let us first explain how we adapt the vertex and set testers from the original blow-up instance to the blow-up instance $\sB_p$.

For $p\in[P]$, we define 
\begin{align*}
\cW_{ver}(p):=
\Big\{
(\omega_p,\bc)\colon&
(\omega,\bc)\in\cW_{ver}, \text{ centres $\bc=\{c_i\}_{i\in I}$ do not collide with respect to $\cV_p$,}
\\&\omega_p:=\omega|_{\bigcup_{H\in \cH_p,i\in I} X_{j_i}^H}
\Big\},
\end{align*}
where we define $\{j_i\}_{i\in I}\subseteq[\beta^{-1}]$ as the indices such that $c_i\in V_{j_i}$ for centres $\bc=\{c_i\}_{i\in I}$ that do not collide with respect to $\cV_p=(V_j)_{j\in[\beta^{-1}]}$. 

For all $j\in[\beta^{-1}]$, $(W,Y_1,\ldots,Y_m)\in\cW_{set}$ and distinct $H_1,\ldots,H_m\in\cH$ such that $Y_\ell\subseteq V(H_\ell)$, we define $Y_{\ell,j}:=Y_\ell\cap X_j^{H_\ell}$  for each $\ell\in[m]$.
For $j\in[\beta^{-1}], p\in[P]$, we define
\begin{align*}
\cW_{set}(j,p):=
\Big\{
\big(
W_j(p), \{Y_{\ell,j}\}_{\ell\in[m]\colon Y_{\ell,j} \subseteq V(\cH_p)}
\big)
\colon
(W,Y_1,\ldots,Y_m)\in\cW_{set},~Y_{\cup [m]}\cap V(\cH_p)\neq\emptyset
\Big\},
\end{align*}
where we define $W_j(p)$ for $(W,Y_1,\ldots,Y_m)\in\cW_{set}$ recursively by $$W_j(p):=W_j(p-1) \cap \bigcap_{\ell\in[m]\colon Y_{\ell,j}\subseteq V(\cH_{p-1})}\phi_{p-1}(Y_{\ell,j})$$ with $\phi_0$ being the empty function and thus, $W_j(1)=W_j(0):=W\cap V_j$.
By~\ref{item:refined W r=1}, we have that $|W\cap V_j|=\beta|W|\pm\beta^{3/2}n$.
By employing conclusion~\ref{item:split' weights one H, |I|=1} of Lemma~\ref{lem:refining partitions r=1} for the weight function $\omega_{Y_\ell}\in\cW_Y$, we have that 
\begin{align}\label{eq:size Y_l,j r=1}
|Y_{\ell,j}|=\omega_{Y_{\ell}}(X_{j}^{H_{\ell}})
\stackrel{\ref{item:split' weights one H, |I|=1}}{=}
\beta \omega_{Y_{\ell}}(V({H_{\ell}})) \pm \beta^{3/2}n =\beta |Y_{\ell}|\pm \beta^{3/2}n.
\end{align}

\medskip
Hence, we iteratively apply Lemma~\ref{lem:main matching} for every $p\in[P]$ to $\sB_p$ with set testers $\bigcup_{j\in[\beta^{-1}]}\cW_{set}(j,p)$ and vertex testers $\cW_{ver}(p)$ as follows:
\begin{center}
\begin{tabular}{r||c|c|c|c|c|c}
parameter & $n'$  & $\eps^{1/2}$ & $t$ & {$\beta^k$} & $P^{-1}d$ &  $\beta^{-1}$
\\ \hline
plays the role of & $n$  & $\eps$  & $t$ & $\alpha$ & $d$  &  $r$ 
\vphantom{\Big(}
\end{tabular}
\end{center} 
For every $p\in[P]$, this yields a packing $\phi_p$ of $\cH_p$ into $G_p$ such that 
\begin{enumerate}[label=(\Roman*)]
\item\label{item:I r=1} $\phi_p(X_{j}^H)= V_{j}$ for all $j\in[\beta^{-1}]$, $H\in\cH_p$;
\item\label{item:II r=1} $|W_j(p)
\cap \bigcap_{\ell\in[m]\colon Y_{\ell,j} \subseteq V(\cH_{p})}\phi_{p}(Y_{\ell,j})|
=
|W_j(p)|\prod_{\ell\in[m]\colon Y_{\ell,j} \subseteq V(\cH_{p})} n'^{-1}|Y_{\ell,j}|
\pm \beta^kn'
$
for all \linebreak$(
W_j(p), \{Y_{\ell,j}\}_{\ell\in[m]\colon Y_\ell \subseteq V(\cH_p)}
)\in \cW_{set}(j,p)$, $j\in[\beta^{-1}]$;
\item \label{item:III r=1} 
$\omega_p( \phi_p^{-1}(\bc))
)
=(1\pm\beta^k)
\omega_p(
\bigcup_{H\in\cH_p}(\bigsqcup_{i\in I}X_{j_i}^H)
)
/n'^{|I|}
\pm n'^{\beta^k}$ for all $(\omega_p,\bc)\in \cW_{ver}(p)$ with centres $\bc=\{c_i\}_{i\in I}$ that do not collide with respect to $\cV_p$ and $\{j_i\}_{i\in I}\subseteq[\beta^{-1}]$ such that $c_i\in V_{j_i}$ for each $i\in I$.
\end{enumerate}
Let $\phi:=\bigcup_{p\in[P]}\phi_p$ and note that $\phi$ is a packing of $\cH$ into $G$.

\medskip
For $(W,Y_1,\ldots,Y_m)\in\cW_{set}$, we conclude that
\begin{align*}
\textstyle\left|W\cap \bigcap_{\ell\in [m]}\phi(Y_{\ell})\right| 
&\stackrel{\hphantom{\rm \ref{item:refined W},\eqref{eq:size Y_l,j r=1}}}{=}
\sum_{j\in[\beta^{-1}]}\textstyle\left|W_j(P)\cap \bigcap_{\ell\in [m]\colon Y_{\ell,j}\subseteq V(\cH_P)}\phi_P(Y_{\ell,j})\right| 
\\
&\stackrel{\rm ~~\ref{item:II r=1}~~~}{=}
 \sum_{j\in[\beta^{-1}]}
 \left(|W_j(0)||Y_{1,j}|\cdots|Y_{m,j}|/n'^m \pm m\beta^kn'\right)
\\
&\stackrel{\rm\ref{item:refined W r=1},\eqref{eq:size Y_l,j r=1}}{=} 
|W||Y_1|\cdots|Y_m|/n^m \pm\alpha n.
\end{align*}
This establishes Theorem~\ref{thm:normalwithweights}\ref{item:set testers2}.

\medskip
In order to establish Theorem~\ref{thm:normalwithweights}\ref{item: vertex tester2}, we fix $(\omega,\bc)\in\cW_{ver}$ with centres $\bc=\{c_i\}_{i\in I}$.\COMMENT{In fact, we have to assume for \ref{item:centres collide} and~\eqref{eq:sliced weight} that $\omega(V(\cH))\geq n^{1+\eps}$ but in the other case it is plain to establish Theorem~\ref{thm:normalwithweights}\ref{item: vertex tester2}.}
Recall that $\supp(\omega)\subseteq E(\cH)$ for $|I|\geq 2$, and thus, if the centres $\bc$ collide with respect to one of the partitions~$\cV_p$, then $\omega(\phi_p^{-1}(\bc))=0$.
Therefore, we can consider the corresponding tuples $(\omega_p,\bc)\in\cW_{ver}(p)$ for $p\in[P]$ and conclude that
\begin{align*}
\allowdisplaybreaks
\omega(\phi^{-1}(\bc))
&
\stackrel{\hphantom{\rm\ref{item:III r=1}}}{=}
\sum_{p\in[P]}\omega(\phi_p^{-1}(\bc))
=\sum_{p\in[P]\colon\text{no collision}}\omega_p(\phi_p^{-1}(\bc))
\\&\allowdisplaybreaks
\stackrel{\rm\ref{item:III r=1}}{=}
\sum_{p\in[P]\colon\text{no collision}}
\left(
(1\pm\beta^k)
\frac{\omega_p(
\bigcup_{H\in\cH_p}(\bigsqcup_{i\in I}X_{j_i}^H)
)}{n'^{|I|}}
\pm n'^{\beta^k}
\right)
\\&\allowdisplaybreaks
\stackrel{\rm\ref{item:split' weights}}{=}
\sum_{p\in[P]\colon\text{no collision}}
\left(
(1\pm\beta^k)
\frac{(1\pm\beta^{1/2})\beta^{|I|}\omega(V(\cH_p))}
{n'^{|I|}}
\pm n'^{\beta^k}
\right)
\\&
\stackrel{~\rm\ref{item:centres collide}~}{=}
(1\pm\alpha)\frac{\omega(V(\cH))}{n^{|I|}}\pm n^\alpha.
\end{align*}
This establishes Theorem~\ref{thm:normalwithweights}\ref{item: vertex tester2} and completes the proof.
\endproof

\vspace{-.1cm}
\section{Concluding remarks}\label{sec:conclusion}

In this paper we have proved a functional tool that allows for approximate decompositions of (sparse) (multipartite) quasirandom $k$-graphs into any collection of spanning bounded degree $k$-graphs (with slightly fewer edges).
Our main result extends the strong decomposition result due to Kim, K\"uhn, Osthus and Tyomkyn into two directions, both from the dense to the sparse regime and from graphs to hypergraphs.
It thereby answers Question~\ref{questionHyper} of Keevash, and Kim, K\"uhn, Osthus and Tyomkyn
and Question~\ref{question} of Kim, K\"uhn, Osthus and Tyomkyn.

We conclude with a selection of immediate applications of Theorem~\ref{thm:simple}
including an asymptotic solution to a hypergraph Oberwolfach problem asked by Glock, K\"uhn and Osthus~\cite{GKO:ta}.
We first consider several natural hypergraph analogues of graph decomposition questions and then turn to problems concerning simplicial complexes. Along the way we state a few conjectures and open problems. We also believe that our main results will have further significantly more complex applications in the future.\looseness=-2

\subsection{Applications to hypergraph decompositions}\label{sec:hypergraphs}

As we pointed out in the introduction,
the Conjectures~\ref{conj:R}--\ref{conj:OW} intrigued mathematicians for decades.
In the following we propose several conjectures of similar spirit for $k$-graphs.
Our main results imply approximate versions thereof.

Recall that the Oberwolfach problem in~\ref{conj:OW} asks for a decomposition of $K_{n}$ into $(n-1)/2$ copies of a graph on $n$ vertices that is the disjoint union of cycles.
There are many definitions for cycles in $k$-graphs and tight cycles are among the most well studied cycles.
A $k$-graph is a tight cycle if its vertex set can be cyclically ordered and the edge set consists of all $k$-sets that appear consecutively in this ordering.
We refer to the number of vertices in a tight cycle as its length.
One potential version of a hypergraph Oberwolfach problem has recently been asked by Glock, K\"uhn and Osthus.

\begin{conj}[Hypergraph Oberwolfach problem; Glock, K\"uhn and Osthus~\cite{GKO:ta}]\label{conj:GKO}
Let $k\geq 3$ and suppose $n$ is sufficiently large in terms of $k$ and $k$ divides $\binom{n-1}{k-1}$.
Suppose $F$ is a $k$-graph on~$n$ vertices that is the disjoint union of tight cycles each of length at least $2k-1$.
Then there is a decomposition of $K_n^{(k)}$ into copies of $F$.
\end{conj}

Clearly, Theorem~\ref{thm:simple} yields an approximate solution of Conjecture~\ref{conj:GKO}.

We think that an even stronger result is true.

\begin{conj}[Hypergraph Oberwolfach problem]\label{conj:wOW}
Let $k\geq 3$ and suppose $n$ is sufficiently large in terms of $k$ and $k$ divides $\binom{n-1}{k-1}$.
Suppose $F$ is a $k$-graph on $n$ vertices that is the disjoint union of tight cycles each of length at least $k+2$.
Then there is a decomposition of $K_n^{(k)}$ into copies of $F$.
\end{conj}

Observe that Conjecture~\ref{conj:GKO} includes the natural generalization of Walecki's theorem to hypergraphs, namely decompositions into Hamilton cycles.
This has already been conjectured by Bailey and Stevens~\cite{BS:10} (and when $n$ and $k$ are coprime by Baranyai~\cite{Bar:79} and independently by Katona) and there are a few results that provide approximate decompositions of quasirandom graphs into Hamilton cycles (of various types); see for example~\cite{BF:12,FK:12,FKL:12}.

Whenever we allow cycles of length $k+1$ and the cycle factor consists (essentially) of cycles of length $k+1$, we suspect that there are more divisibility obstructions present. Hence we pose the following problem.\looseness=-2 

\begin{problem}[Hypergraph Oberwolfach problem]
Let $k\geq 3$ and suppose $n$ is sufficiently large in terms of $k$ and $k$ divides $\binom{n-1}{k-1}$.
Which disjoint unions of tight cycles whose length add up to $n$ admit a decomposition of $K_n^{(k)}$?
\end{problem}

It immediately follows from Theorem~\ref{thm:simple} that the hypergraph Oberwolfach problems are approximately true in a sense that $K_n^{(k)}$ contains $(1-o(1))\binom{n-1}{k-1}/k$ disjoint copies of~$F$ (for any choice of~$F$ as above);
in fact, we can we take any collection of $(1-o(1))\binom{n-1}{k-1}/k$ cycle factors.

\medskip

Similarly as for cycles, there is more then one notion for trees in $k$-graphs.
Let us stick to the following recursive definition of tree to which we refer as a $k$-tree.
A single edge is a $k$-tree.
A $k$-tree with $\ell$ edges can be constructed from a $k$-tree with $\ell-1$ edges $T$ by adding a vertex $v$ and an edge that contains $v$ and a $(k-1)$-set that is contained in an edge of $T$.
For this definition, we propose the following generalization of Ringel's conjecture.

\begin{conj}\label{conj:HyperRingel}
Let $k,n\in \bN\sm \{1\}$.
Suppose $T$ is a $k$-tree with $n$ edges.
Then $K_{kn + k - 1}^{(k)}$ admits a decomposition into copies of $T$.
\end{conj}

Observe that similarly as for Ringel's conjecture, the order of the complete graph needs to be at least $kn+k-1$ if we allow $T$ to be any tree with $n$ edges as the natural generalization of a star shows.
It is an easy exercise to show this conjecture for stars.

There is a conjecture related to Ringel's conjecture for bipartite graphs due to Graham and H\"aggkvist
 stating that $K_{n,n}$ can be decomposed into $n$ copies of any tree with $n$ edges.
We propose here the following strengthening.

\begin{conj}\label{conj:k-partite}
Suppose $k,n\in \bN\sm\{1\}$ and $T$ is a $k$-tree with $n$ edges.
Then there is a decomposition of the complete balanced $k$-partite graph on $kn$ vertices.
\end{conj}

The tree packing conjecture has arguably the least obvious strengthening to $k$-graphs
and there may be more than one.
We propose the following one.

\begin{conj}\label{conj:hypertree packing}
Suppose $k,n\in \bN\sm\{1\}$.
Let $\cT$ be a family of $k$-trees such that~$\cT$ contains $\binom{n-i-1}{k-2}$ trees with $i$ edges for $i\in [n-k+1]$.
Then $K_n^{(k)}$ admits a decomposition into $\cT$.
\end{conj}

It follows directly from Theorem~\ref{thm:simple} that Conjectures~\ref{conj:HyperRingel} and~\ref{conj:hypertree packing} are approximately true when restricted to bounded degree trees (and similarly an approximate version of Conjecture~\ref{conj:k-partite} follows from our main theorem, see Theorem~\ref{thm:multi} in Section~\ref{sec:main results}).

\subsection{Applications to simplicial complexes}\label{sec:simplicial complexes}

Generalizing long-studied and nowadays classical combinatorial questions to higher dimensions appears to be a challenging but insightful theme.
There are several results considering $k$-dimensional permutations
and it was Linial and Meshulam~\cite{LM:06} who introduced a random model for simplicial complexes whose probability measure is the same as those of a binomial random $(d+1)$-graph; simply add all $d'$-faces for $d'\leq d-1$ with points in~$[n]$ and add every potential $d$-face independently with probability $p$.
Recently, Linial and Peled investigated and determined the threshold in $Y_d(n,p)$ for the emergence of a what may considered as an analogoue of a giant component~\cite{LP:16}.

With this topological viewpoint of treating a $k$-graph as a simplicial complex,
a cycle in a graph is simply an object homomorphic to $\bS^1$
and hence a $2$-dimensional Hamilton cycle a collection of $2$-faces containing every vertex and which is homomorphic to~$\bS^2$.
In~\cite{LT:19}, Luria and Tessler determined the threshold in $Y_2(n,p)$ for the appearance of such a Hamilton cycle (another suitable term may be spanning triangulation of the sphere).

An analogue of Dirac's theorem was proved by Georgakopoulos, Haslegrave, Narayanan and Montgomery;
to be precise, when every pair of vertices is contained in at least $n/3+o(n)$ edges/2-faces of a $3$-graph $G$, then there is a spanning triangulation of the sphere in $G$.
This bound remains the same when we replace $\bS^2$ by any other compact surface without boundary~\cite{GHMN:18}.

Instead of only asking for a single triangulation of some surface, we can of course also investigate decompositions into (spanning) triangulations of surfaces.
Our results imply that every quasirandom simplicial complex (this in particular includes almost all graphs in $Y_d(n,p)$) can even be almost decomposed into any list of triangulations of any kind of manifolds provided the density does not vanish too quickly with $n$ and as long as every vertex is contained in at most a bounded number of $d$-faces.
The triangulations may even be chosen in advance.
Hence a precise statement for 2-complexes is as follows.

\begin{cor}\label{cor:simplicial_complex}
For all $\alpha>0$,
there exist $n_0,h\in \bN$ and $\eps>0$ such that the following holds for all $n\geq n_0$.
Suppose $G$ is an $(\eps,h,d)$-typical $3$-graph on $n$ vertices with $d\geq n^{-\eps}$ and 
$H_1,\ldots,H_\ell$ are spanning triangulations of $\bS^2$ where every vertex is contained in at most $\alpha^{-1}$ $2$-faces and $\ell\leq (1-\alpha)dn^2/12$.
Then $G$ contains edge-disjoint copies of $H_1,\ldots,H_\ell$ such that every $2$-face of $G$ is contained in at most one $H_i$.
\end{cor}

We omit statements for higher dimensions as they follow in the obvious way from Theorem~\ref{thm:simple}.
We wonder whether there is always an actual decomposition (subject to certain divisibility conditions).
This might be easier than a decomposition into tight Hamilton cycles as the structure of tight Hamilton cycles seems to be more restrictive.

\bibliographystyle{amsplain_v2.0customized}
\bibliography{References_approx}

\providecommand{\bysame}{\leavevmode\hbox to3em{\hrulefill}\thinspace}
\providecommand{\MR}{\relax\ifhmode\unskip\space\fi MR }
\providecommand{\MRhref}[2]{%
  \href{http://www.ams.org/mathscinet-getitem?mr=#1}{#2}
}
\providecommand{\href}[2]{#2}
\begin{thebibliography}{10}

\bibitem{ABCT:19}
P.~Allen, J.~B\"ottcher, D.~Clemens, and A.~Taraz, \emph{Perfectly packing
  graphs with bounded degeneracy and many leaves}, arXiv:1906.11558 (2019).

\bibitem{ABHKP:16}
P.~Allen, J.~B\"ottcher, H.~H\`an, Y.~Kohayakawa, and Y.~Person, \emph{Blow-up
  lemmas for sparse graphs}, arXiv:1612.00622 (2016).

\bibitem{ABHP:17}
P.~Allen, J.~B\"ottcher, J.~Hladk\'y, and D.~Piguet, \emph{Packing degenerate
  graphs}, Adv. Math.~\textbf{354} (2019), 106739.

\bibitem{ABSS:20}
P.~Allen, J.~B{\"o}ttcher, J.~Skokan, and M.~Stein, \emph{Regularity
  inheritance in pseudorandom graphs}, Random Structures Algorithms~\textbf{56}
  (2020), 306--338.

\bibitem{ARR:98}
N.~Alon, V.~R{\"o}dl, and A.~Ruci\'nski, \emph{Perfect matchings in
  $\epsilon$-regular graphs}, Electron. J. Combin.~\textbf{5} (1998), 1--4.

\bibitem{BS:10}
R.~F.~Bailey and B.~Stevens, \emph{Hamiltonian decompositions of complete
  {$k$}-uniform hypergraphs}, Discrete Math.~\textbf{310} (2010), 3088--3095.

\bibitem{BF:12}
D.~Bal and A.~Frieze, \emph{Packing tight {H}amilton cycles in uniform
  hypergraphs}, SIAM J. Discrete Math.~\textbf{26} (2012), 435--451.

\bibitem{Bar:79}
Z.~Baranyai, \emph{The edge-coloring of complete hypergraphs. {I}}, J. Combin.
  Theory Ser. B~\textbf{26} (1979), 276--294.

\bibitem{BHPT:16}
J.~B\"{o}ttcher, J.~Hladk\'{y}, D.~Piguet, and A.~Taraz, \emph{An approximate
  version of the tree packing conjecture}, Israel J. Math.~\textbf{211} (2016),
  391--446.

\bibitem{CG:02}
F.~Chung and R.~Graham, \emph{Sparse quasi-random graphs},
  Combinatorica~\textbf{22} (2002), 217--244.

\bibitem{ehard:phd}
S.~Ehard, \emph{Embeddings and decompositions of graphs and hypergraphs}, PhD
  thesis (2020).

\bibitem{EGJ:19a}
S.~Ehard, S.~Glock, and F.~Joos, \emph{Pseudorandom hypergraph matchings},
  Combin. Probab. Comput.~\textbf{29} (2020), 868--885.

\bibitem{EJ:19}
S.~Ehard and F.~Joos, \emph{A short proof of the blow-up lemma for approximate
  decompositions}, arXiv:2001.03506 (2020).

\bibitem{FLM:17}
A.~Ferber, C.~Lee, and F.~Mousset, \emph{Packing spanning graphs from separable
  families}, Israel J. Math.~\textbf{219} (2017), 959--982.

\bibitem{freedman:75}
D.~A.~Freedman, \emph{On tail probabilities for martingales}, Ann.
  Probab.~\textbf{3} (1975), 100--118.

\bibitem{FK:12}
A.~Frieze and M.~Krivelevich, \emph{Packing {H}amilton cycles in random and
  pseudo-random hypergraphs}, Random Structures Algorithms~\textbf{41} (2012),
  1--22.

\bibitem{FKL:12}
A.~Frieze, M.~Krivelevich, and P.-S.~Loh, \emph{Packing tight {H}amilton cycles
  in 3-uniform hypergraphs}, Random Structures Algorithms~\textbf{40} (2012),
  269--300.

\bibitem{GHMN:18}
A.~Georgakopoulos, J.~Haslegrave, R.~Montgomery, and B.~Narayanan,
  \emph{Spanning surfaces in $3$-graphs}, arXiv:1808.06864 (2018).

\bibitem{GJKKO:18}
S.~Glock, F.~Joos, J.~Kim, D.~K\"uhn, and D.~Osthus, \emph{Resolution of the
  {O}berwolfach problem}, arXiv:1806.04644 (2018).

\bibitem{GKLO:ta}
S.~Glock, D.~K\"uhn, A.~Lo, and D.~Osthus, \emph{The existence of designs via
  iterative absorption: hypergraph {$F$}-designs for arbitrary~{$F$}}, Mem.
  Amer. Math. Soc. (to appear).

\bibitem{GKO:ta}
S.~Glock, D.~K\"uhn, and D.~Osthus, \emph{Extremal aspects of graph and
  hypergraph decomposition problems}, Surveys in Combinatorics, London
  Mathematical Society Lecture Note Series, to appear.

\bibitem{HHM:20}
H.~H{\`a}n, J.~Han, and P.~Morris, \emph{Factors and loose {H}amilton cycles in
  sparse pseudo-random hypergraphs}, Proceedings of the Fourteenth Annual
  ACM-SIAM Symposium on Discrete Algorithms, SIAM, 2020, pp.~702--717.

\bibitem{JKKO:19}
F.~Joos, J.~Kim, D.~K\"uhn, and D.~Osthus, \emph{Optimal packings of bounded
  degree trees}, J. Eur. Math. Soc.~\textbf{21} (2019), 3573--3647.

\bibitem{JK:21}
F.~Joos and M.~K\"uhn, \emph{Fractional cycle decompositions in hypergraphs},
  arXiv:2101.05526 (2021).

\bibitem{keevash:11}
P.~Keevash, \emph{A hypergraph blow-up lemma}, Random Structures
  Algorithms~\textbf{39} (2011), 275--376.

\bibitem{keevash:14}
\bysame, \emph{The existence of designs}, arXiv:1401.3665 (2014).

\bibitem{keevash:18b}
\bysame, \emph{The existence of designs~{II}}, arXiv:1802.05900 (2018).

\bibitem{keevash:18c}
\bysame, \emph{Hypergraph matchings and designs}, Proceedings of the
  {I}nternational {C}ongress of {M}athematicians---{R}io de {J}aneiro 2018.
  {V}ol. {IV}. {I}nvited lectures, World Sci. Publ., Hackensack, NJ, 2018,
  pp.~3113--3135.

\bibitem{KS:20a}
P.~Keevash and K.~Staden, \emph{The generalised {O}berwolfach problem},
  arXiv:2004.09937 (2020).

\bibitem{KS:20b}
\bysame, \emph{Ringel's tree packing conjecture in quasirandom graphs},
  arXiv:2004.09947 (2020).

\bibitem{KKOT:19}
J.~Kim, D.~K\"uhn, D.~Osthus, and M.~Tyomkyn, \emph{A blow-up lemma for
  approximate decompositions}, Trans. Amer. Math. Soc.~\textbf{371} (2019),
  4655--4742.

\bibitem{KSS:97}
J.~Koml\'os, G.~N.~S\'ark\"ozy, and E.~Szemer\'edi, \emph{Blow-up lemma},
  Combinatorica~\textbf{17} (1997), 109--123.

\bibitem{LM:06}
N.~Linial and R.~Meshulam, \emph{Homological connectivity of random
  2-complexes}, Combinatorica~\textbf{26} (2006), 475--487.

\bibitem{LP:16}
N.~Linial and Y.~Peled, \emph{On the phase transition in random simplicial
  complexes}, Ann. of Math.~\textbf{184} (2016), 745--773.

\bibitem{LT:19}
Z.~Luria and R.~J.~Tessler, \emph{A sharp threshold for spanning 2-spheres in
  random 2-complexes}, Proc. Lond. Math. Soc.~\textbf{119} (2019), 733--780.

\bibitem{mcdiarmid:89}
C.~McDiarmid, \emph{On the method of bounded differences}, Surveys in
  combinatorics, 1989 ({N}orwich, 1989), London Math. Soc. Lecture Note Ser.
  141, Cambridge Univ. Press, 1989, pp.~148--188.

\bibitem{MRS:16}
S.~Messuti, V.~R\"{o}dl, and M.~Schacht, \emph{Packing minor-closed families of
  graphs into complete graphs}, J. Combin. Theory Ser. B~\textbf{119} (2016),
  245--265.

\bibitem{MPS:20}
R.~Montgomery, A.~Pokrovskiy, and B.~Sudakov, \emph{A proof of {R}ingel's
  {C}onjecture}, arXiv:2001.02665 (2020).

\bibitem{RR:99}
V.~R\"odl and A.~Ruci\'nski, \emph{Perfect matchings in {$\epsilon$}-regular
  graphs and the blow-up lemma}, Combinatorica~\textbf{19} (1999), 437--452.

\bibitem{thomason:87b}
A.~Thomason, \emph{Pseudo-random graphs}, Proceedings of Random Graphs,
  Pozna{\'n} 1985, Annals of Discrete Math.~33, North Holland, 1987,
  pp.~307--331.

\bibitem{thomason:87}
\bysame, \emph{Random graphs, strongly regular graphs and pseudo-random
  graphs}, Surveys in Combinatorics,  123, London Mathematical Society Lecture
  Note Series, 1987, pp.~173--195.

\bibitem{wilson:72a}
R.~M.~Wilson, \emph{An existence theory for pairwise balanced designs {I}.
  {C}omposition theorems and morphisms}, J. Combin. Theory Ser.~A~\textbf{13}
  (1972), 220--245.

\bibitem{wilson:72b}
\bysame, \emph{An existence theory for pairwise balanced designs {II}. {T}he
  structure of {PBD}-closed sets and the existence conjectures}, J. Combin.
  Theory Ser.~A~\textbf{13} (1972), 246--273.

\bibitem{wilson:75}
\bysame, \emph{An existence theory for pairwise balanced designs {III}. {P}roof
  of the existence conjectures}, J. Combin. Theory Ser.~A~\textbf{18} (1975),
  71--79.

\end{thebibliography}
\end{document}